\documentclass{article}
\usepackage{times}
\usepackage{amsmath,amssymb}
\usepackage{mathrsfs}
\usepackage{mathabx}\changenotsign 
\usepackage{dsfont}
\usepackage{bbm}
\usepackage{scalerel}
\usepackage{graphicx}
\usepackage{algorithm}
\usepackage{adjustbox}
\usepackage{array}
\usepackage[noend]{algpseudocode}
\usepackage{algcompatible}
\usepackage[T1]{fontenc}
\usepackage{nicefrac}
\usepackage{booktabs} 
\usepackage{amsfonts}
\usepackage{etex}
\usepackage{natbib}
\usepackage{prettyref}
\usepackage{subfiles}
\usepackage{subcaption} 
\usepackage[letterpaper, left=1in, right=1in, top=1in, bottom=1in]{geometry}
\PassOptionsToPackage{hypertexnames=false}{hyperref}  
\usepackage{parskip}
\usepackage[dvipsnames]{xcolor}
\usepackage[colorlinks=true, linkcolor=blue!70!black, citecolor=blue!70!black,urlcolor=black,breaklinks=true]{hyperref}
\hypersetup{linktoc=all}

\usepackage{microtype}
\usepackage{hhline}

\makeatletter
\newcommand{\neutralize}[1]{\expandafter\let\csname c@#1\endcsname\count@}
\makeatother

\usepackage{algorithm}

\usepackage{natbib}
\bibpunct{(}{)}{;}{a}{,}{,}

\usepackage{amsthm}
\usepackage{mathtools}
\usepackage{amsmath}
\usepackage{bbm}
\usepackage{amsfonts}
\usepackage{amssymb}

\usepackage{xpatch}

\usepackage{thmtools}
\usepackage{thm-restate}
\declaretheorem[name=Theorem,parent=section]{theorem}
\declaretheorem[name=Lemma,parent=section]{lemma}
\declaretheorem[name=Assumption, parent=section]{assumption}
\declaretheorem[name=Condition, parent=section]{condition}
\declaretheorem[name=Example,style=definition]{example}
\declaretheorem[name=Remark, style=definition]{remark}
\declaretheorem[name=Proposition, parent=section]{proposition}

\makeatletter
  \renewenvironment{proof}[1][Proof]%
  {%
   \par\noindent{\bfseries\upshape {#1.}\ }%
  }%
  {\qed \\ \newline}
  \makeatother

\newtheorem{definition}{Definition}[section]

\xpatchcmd{\proof}{\itshape}{\normalfont\proofnameformat}{}{}

\usepackage[nameinlink,capitalize]{cleveref}

\renewcommand{\eqref}[1]{\texorpdfstring{\hyperref[#1]{Eq.~(\ref*{#1})}}{Eq.~(\ref*{#1})}}

\crefformat{equation}{#2(#1)#3}
\Crefformat{equation}{#2(#1)#3}

\Crefformat{figure}{#2Figure #1#3}
\Crefname{assumption}{Assumption}{Assumptions}
\Crefformat{assumption}{#2Assumption #1#3}
\Crefname{subsubsection}{Section}{Sections}
\crefformat{subsubsection}{#2Section #1#3}
\Crefformat{subsubsection}{#2Section #1#3}

\usepackage{crossreftools}
\usepackage{xparse}

\ExplSyntaxOn
\DeclareDocumentCommand{\XDeclarePairedDelimiter}{mm}
 {
  \__egreg_delimiter_clear_keys: 
  \keys_set:nn { egreg/delimiters } { #2 }
  \use:x 
   {
    \exp_not:n {\NewDocumentCommand{#1}{sO{}m} }
     {
      \exp_not:n { \IfBooleanTF{##1} }
       {
        \exp_not:N \egreg_paired_delimiter_expand:nnnn
         { \exp_not:V \l_egreg_delimiter_left_tl }
         { \exp_not:V \l_egreg_delimiter_right_tl }
         { \exp_not:n { ##3 } }
         { \exp_not:V \l_egreg_delimiter_subscript_tl }
       }
       {
        \exp_not:N \egreg_paired_delimiter_fixed:nnnnn 
         { \exp_not:n { ##2 } }
         { \exp_not:V \l_egreg_delimiter_left_tl }
         { \exp_not:V \l_egreg_delimiter_right_tl }
         { \exp_not:n { ##3 } }
         { \exp_not:V \l_egreg_delimiter_subscript_tl }
       }
     }
   }
 }

\keys_define:nn { egreg/delimiters }
 {
  left      .tl_set:N = \l_egreg_delimiter_left_tl,
  right     .tl_set:N = \l_egreg_delimiter_right_tl,
  subscript .tl_set:N = \l_egreg_delimiter_subscript_tl,
 }

\cs_new_protected:Npn \__egreg_delimiter_clear_keys:
 {
  \keys_set:nn { egreg/delimiters } { left=.,right=.,subscript={} }
 }

\cs_new_protected:Npn \egreg_paired_delimiter_expand:nnnn #1 #2 #3 #4
 {
  \mathopen{}
  \mathclose\c_group_begin_token
   \left#1
   #3
   \group_insert_after:N \c_group_end_token
   \right#2
   \tl_if_empty:nF {#4} { \c_math_subscript_token {#4} }
 }
\cs_new_protected:Npn \egreg_paired_delimiter_fixed:nnnnn #1 #2 #3 #4 #5
 {
  \mathopen{#1#2}#4\mathclose{#1#3}
  \tl_if_empty:nF {#5} { \c_math_subscript_token {#5} }
 }
\ExplSyntaxOff

\XDeclarePairedDelimiter{\supnorm}{
  left=\lVert,
  right=\rVert,
  subscript=\infty
  }

\usepackage{mathtools}
\usepackage{amsmath}
\usepackage{amsxtra}
\usepackage{enumitem}      
\usepackage{mathtools}
\usepackage{bbm}
\usepackage{amsfonts}
\usepackage{amssymb}
\let\vec\undefined
\usepackage{MnSymbol} 

\usepackage{xpatch}
\newcommand\blfootnote[1]{%
  \begingroup
  \renewcommand\thefootnote{}\footnote{#1}%
  \addtocounter{footnote}{-1}%
  \endgroup
}

\newtheorem{claim}{Claim}
{\newtheorem{corollary}{Corollary}} 

\usepackage{prettyref}
\newcommand{\pref}[1]{\prettyref{#1}}

\newcommand{\savehyperref}[2]{\texorpdfstring{\hyperref[#1]{#2}}{#2}}
\newrefformat{eq}{\savehyperref{#1}{\textup{(\ref*{#1})}}}
\newrefformat{eqn}{\savehyperref{#1}{Equation~\ref*{#1}}}
\newrefformat{con}{\savehyperref{#1}{Conjecture~\ref*{#1}}}
\newrefformat{lem}{\savehyperref{#1}{Lemma~\ref*{#1}}}
\newrefformat{def}{\savehyperref{#1}{Definition~\ref*{#1}}}
\newrefformat{line}{\savehyperref{#1}{line~\ref*{#1}}}
\newrefformat{thm}{\savehyperref{#1}{Theorem~\ref*{#1}}}
\newrefformat{corr}{\savehyperref{#1}{Corollary~\ref*{#1}}}
\newrefformat{sec}{\savehyperref{#1}{Section~\ref*{#1}}}
\newrefformat{app}{\savehyperref{#1}{Appendix~\ref*{#1}}}
\newrefformat{ass}{\savehyperref{#1}{Assumption~\ref*{#1}}}
\newrefformat{ex}{\savehyperref{#1}{Example~\ref*{#1}}}
\newrefformat{fig}{\savehyperref{#1}{Figure~\ref*{#1}}}
\newrefformat{alg}{\savehyperref{#1}{Algorithm~\ref*{#1}}}
\newrefformat{rem}{\savehyperref{#1}{Remark~\ref*{#1}}}
\newrefformat{conj}{\savehyperref{#1}{Conjecture~\ref*{#1}}}
\newrefformat{prop}{\savehyperref{#1}{Proposition~\ref*{#1}}}
\newrefformat{proto}{\savehyperref{#1}{Protocol~\ref*{#1}}}
\newrefformat{prob}{\savehyperref{#1}{Problem~\ref*{#1}}}
\newrefformat{claim}{\savehyperref{#1}{Claim~\ref*{#1}}}
\newrefformat{subsec}{\savehyperref{#1}{Subsection~\ref*{#1}}}
\newrefformat{subsubsec}{\savehyperref{#1}{Subsection~\ref*{#1}}}
\newrefformat{sec}{\savehyperref{#1}{Section~\ref*{#1}}}
\newrefformat{ineq}{\savehyperref{#1}{Inequality~\ref*{#1}}}

\newcommand\numberthis{\addtocounter{equation}{1}\tag{\theequation}}
\allowdisplaybreaks

\DeclarePairedDelimiter{\abs}{\lvert}{\rvert} 
\DeclarePairedDelimiter{\brk}{[}{]}
\DeclarePairedDelimiter{\crl}{\{}{\}}
\DeclarePairedDelimiter{\prn}{(}{)}
\DeclarePairedDelimiter{\nrm}{\|}{\|}
\DeclarePairedDelimiter{\tri}{\langle}{\rangle}

\def\ddefloop#1{\ifx\ddefloop#1\else\ddef{#1}\expandafter\ddefloop\fi}
\def\ddef#1{\expandafter\def\csname bb#1\endcsname{\ensuremath{\mathbb{#1}}}}
\ddefloop ABCDEFGHIJKLMNOPQRSTUVWXYZ\ddefloop
\def\ddefloop#1{\ifx\ddefloop#1\else\ddef{#1}\expandafter\ddefloop\fi}
\def\ddef#1{\expandafter\def\csname b#1\endcsname{\ensuremath{\mathbf{#1}}}}
\ddefloop ABCDEFGHIJKLMNOPQRSTUVWXYZ\ddefloop
\def\ddef#1{\expandafter\def\csname c#1\endcsname{\ensuremath{\mathcal{#1}}}}
\ddefloop ABCDEFGHIJKLMNOPQRSTUVWXYZ\ddefloop
\def\ddef#1{\expandafter\def\csname h#1\endcsname{\ensuremath{\widehat{#1}}}}
\ddefloop ABCDEFGHIJKLMNOPQRSTUVWXYZabcdefghijklmnopqrsuvwxyz\ddefloop    
\def\ddef#1{\expandafter\def\csname hc#1\endcsname{\ensuremath{\widehat{\mathcal{#1}}}}}
\ddefloop ABCDEFGHIJKLMNOPQRSTUVWXYZ\ddefloop
\def\ddef#1{\expandafter\def\csname t#1\endcsname{\ensuremath{\widetilde{#1}}}}
\ddefloop ABCDEFGHIJKLMNOPQRSTUVWXYZ\ddefloop
\def\ddef#1{\expandafter\def\csname tc#1\endcsname{\ensuremath{\widetilde{\mathcal{#1}}}}}
\ddefloop ABCDEFGHIJKLMNOPQRSTUVWXYZ\ddefloop

\newcommand{\ball}{\mathbb{B}}

\usepackage{tikz}

\newcommand{\grad}{\nabla}

\renewcommand{\epsilon}{\varepsilon}

\newcommand{\DERIV}{\mathrm{d}}

\newcommand{\OPNORM}{\mathrm{op}{}}

\newcommand{\vecW}{\pmb{w}}
\newcommand{\vecX}{\pmb{x}}
\newcommand{\vecY}{\pmb{y}}
\newcommand{\vecZ}{\pmb{z}}
\newcommand{\vecV}{\pmb{v}}
\newcommand{\vecB}{\pmb{B}}
\newcommand{\vecU}{\pmb{u}}
\newcommand{\vecE}{\pmb{e}}
\newcommand{\vecP}{\pmb{p}}

\newcommand{\vecA}{\pmb{a}}
\newcommand{\vecb}{\pmb{b}}
\newcommand{\vecDelta}{\pmb{\delta}}
\newcommand{\vecOrigin}{\vec{\pmb{0}}}
\newcommand{\vecEps}{\pmb{\varepsilon}}
\newcommand{\vecAlpha}{\pmb{\alpha}}
\newcommand{\vecXi}{\pmb{\xi}}
\newcommand{\vecZeta}{\pmb{\zeta}}
\newcommand{\vecBeta}{\pmb{\beta}}
\newcommand{\vecPsi}{\pmb{\psi}}
\newcommand{\vecLambda}{\pmb{\Lambda}}
\newcommand{\vechatPsi}{\pmb{\hat{\psi}}}
\newcommand{\vecXiS}{\vecXi_{\vecU}}
\newcommand{\vecXiSperp}{\vecXi_{\vecV}}
\newcommand{\vecYtilde}{\tilde{\vecY}}
\newcommand{\vecXbar}{\overline{\vecX}}

\newcommand{\lambdamin}{\lambda_{\textsc{min}}}
\newcommand{\matI}{\pmb{I}}
\newcommand{\matD}{\pmb{D}}
\newcommand{\hatmatD}{\hat{\pmb{D}}}
\newcommand{\matA}{\pmb{A}}
\newcommand{\matH}{\pmb{H}}
\newcommand{\matM}{\pmb{M}}
\newcommand{\matU}{\pmb{U}}
\newcommand{\matLambda}{\pmb{\Lambda}}
\newcommand{\matHS}{\pmb{H}_{\cS}}
\newcommand{\matHSperp}{\pmb{H}_{\cS^{\perp}}}
\newcommand{\projS}{\pmb{\cP}_{\cS}}
\newcommand{\projSperp}{\pmb{\cP}_{\cS^{\perp}}}

\newcommand{\projSgradF}{\grad_{\vecU} F}
\newcommand{\projSperpgradF}{\grad_{\vecV} F}

\DeclareMathOperator{\polylog}{polylog}

\newcommand{\cmax}{c_{\text{max}}}

\newcommand{\Tthres}{\normalfont t_{\text{thres}}}
\newcommand{\Toracle}{\normalfont t_{\text{oracle}}}
\newcommand{\Gthres}{\normalfont g_{\text{thres}}}
\newcommand{\Fthres}{\normalfont f_{\text{thres}}}
\newcommand{\Tnoise}{\normalfont t_{\text{noise}}}
\newcommand{\cKexit}{\mathcal{K}_{\text{exit}}}

\usepackage[utf8]{inputenc} 
\usepackage[T1]{fontenc}    
\usepackage{hyperref}       
\usepackage{url}            
\usepackage{booktabs}       
\usepackage{amsfonts}       
\usepackage{nicefrac}       
\usepackage{microtype}      
\usepackage{xcolor}         

\title{Efficiently Escaping Saddle Points under Generalized Smoothness via Self-Bounding Regularity}

\author{%
Daniel Yiming Cao$^{\dag}$ \qquad August Y. Chen$^{\dag}$ \qquad Karthik Sridharan$^{\dag}$ \qquad Benjamin Tang$^{\dag}$
\vspace{10pt} 
\\
\small{$^\dag$Cornell University} 
}

\begin{document}

\maketitle

\begin{abstract}
We study the optimization of non-convex functions that are not necessarily smooth (gradient and/or Hessian are Lipschitz) using first order methods. Smoothness is a restrictive assumption in machine learning in both theory and practice, motivating significant recent work on finding first order stationary points of functions satisfying generalizations of smoothness with first order methods.
We develop a novel framework that lets us systematically study the convergence of a large class of first-order optimization algorithms (which we call decrease procedures) under generalizations of smoothness.
We instantiate our framework to analyze the convergence of first order optimization algorithms to first and \textit{second} order stationary points under generalizations of smoothness. As a consequence, we establish the first convergence guarantees for first order methods to second order stationary points under generalizations of smoothness. We demonstrate that several canonical examples fall under our framework, and highlight practical implications.
\blfootnote{ $^\star$Alphabetical ordering. \null\;\; 
{\scriptsize~~\texttt{Emails:}~\{dyc33, ayc74, ks999, bt283\}@cornell.edu}} 
\end{abstract}

\section{Introduction}\label{sec:intro}
A widely studied problem in machine learning (ML) and optimization is finding a First Order Stationary Point (FOSP) of a generic function $F$ with domain $\mathbb{R}^d$, defined as follows:
\[ \text{Given a tolerance $\epsilon>0$, find $\vecW$ such that } \nrm*{\grad F(\vecW)} \le \epsilon\numberthis\label{eq:fospproblem}.\]
The methods of choice in theory and practice for this task are Gradient Descent (GD), Stochastic Gradient Descent (SGD), and variants thereof. Under the additional assumption of (second-order) \textit{smoothness} on $F$, i.e. that the gradient $\grad F$ is Lipschitz with parameter $L>0$, this task is well-understood. In several settings -- such as with access to exact gradients, stochastic gradients, Hessian-Vector Products, and the exact Hessian -- we have matching upper and lower bounds. The literature on this problem is extensive; for a subset see e.g.  \citet{ghadimi2013stochastic, johnson2013accelerating, fang2018spider, fang2019sharp, foster2019complexity, arjevani2020second, carmon2020lower, carmon2021lower}.

However, for many non-convex functions $F$, FOSPs are uninformative. A significant and difficult problem established in the literature for over a decade -- which carries strong theoretical and practical implications in optimization for machine learning -- is establishing efficient rates for finding a Second Order Stationary Point (SOSP). In many non-convex optimization problems such as Phase Retrieval and Matrix Square Root \citep{ge2015escaping, jin2017escape, ge2017no, sun2018geometric}, SOSPs are global minima. Finding a SOSP is defined as follows:
\[ \text{Given a tolerance $\epsilon>0$, find $\vecW$ such that }\nrm*{\grad F(\vecW)} \le \epsilon, \grad^2 F(\vecW) \succeq -\sqrt{\epsilon}\matI,\numberthis\label{eq:sospproblem}\]
where $\succeq$ denotes the PSD order, $\matI$ is the $d \times d$ identity matrix, and $\grad^2 F(\vecW)$ is the Hessian of $F$.\footnote{\text{There are several definitions of a SOSP; see \pref{rem:sospdefnosmooth} for why we use this definition here.}}

Under the additional \textit{Hessian Lipschitz} assumption, that the operator norm of the Hessian $\grad^2 F$ in addition to the gradient $\grad F$ is Lipschitz, this task is also well-understood. Under these regularity assumptions, finding SOSPs is classical under exact oracle access to the full Hessian $\grad^2 F$. Decades ago, it was shown that cubic regularization and trust region methods succeed \citep{nesterov2006cubic, conn2000trust}, with a matching lower bound in \citet{arjevani2020second}. Motivated by the success of non-convex optimization in ML via first order methods, solving this problem \pref{eq:sospproblem} with first order methods has seen much recent study \citep{ge2015escaping, jin2017escape, fang2019sharp, arjevani2020second, jin2021nonconvex}. We have matching upper and lower bounds in several cases, such as for SGD which is perhaps most relevant to ML \citep{fang2019sharp, arjevani2020second}.

However, in many optimization problems in ML, the gradient and Hessian of the loss function is not Lipschitz. This was observed empirically through extensive experiments of \citet{zhang2019gradient} on LSTMs and of \citet{crawshaw2022robustness} on transformers. We provide theoretical examples in \pref{subsec:examplesmainbody}. As such, a line of work began in \citet{zhang2019gradient} on studying finding FOSPs under weaker regularity assumptions, see e.g. \citep{zhang2020improved, jin2021smoothness, crawshaw2022robustness, reisizadeh2023variance, li2023convergence, wang2024provable, hong2024convergence, gaash2025convergence, yu2025convergence}. The regularity assumption generally made is $(L_0, L_1)$-smoothness: $\nrm*{\grad^2 F(\vecW)}_{\OPNORM} \le L_0+L_1 \nrm*{\grad F(\vecW)}$ for all $\vecW\in\mathbb{R}^d$ for some $L_0, L_1 \ge 0$. This allows for arbitrarily polynomial growth rates of $F$ in $\nrm*{\vecW}$.
The guarantees in \citet{zhang2019gradient} and follow-up works generally hold for adaptive methods, presented as theoretical justification for gradient clipping.

The authors of \citet{li2023convex}, under a milder regularity assumption than \citet{zhang2019gradient}, studied finding FOSPs via \textit{fixed-step-size} GD and SGD rather than adaptive methods. In particular, \textit{\citet{li2023convex} demonstrated clipping is not necessary for $(L_0, L_1)$-smooth functions.} Related works extended this analysis to Nesterov's Accelerated Gradient Descent \citep{li2023convergence, hong2024convergence}. \citet{xie2024trust} studied finding SOSPs under $(L_0, L_1)$-smoothness and a similar assumption that for all $\vecW$, in a small neighborhood of $\vecW$, the Hessian of $F$ is Lipschitz with parameter $M_0+M_1 \nrm*{\grad F(\vecW)}$. However, their algorithm is \textit{second-order} and requires the \textit{full} Hessian, analogous to classical work \citep{nesterov2006cubic, conn2000trust}. This contrasts with recent developments of finding SOSPs using first order methods when $F$ has Lipschitz gradient and Hessian, which are more pertinent to ML where first-order algorithms are the only tractable method \citep{ge2015escaping, jin2017escape, fang2019sharp, arjevani2020second, jin2021nonconvex}.

\subsection{Our Contributions}
In this work, we develop a novel framework to study \textbf{non-asymptotic} guarantees finding FOSPs and SOSPs via first-order methods, for functions whose gradient and/or Hessian are not Lipschitz. Central to our work is the following regularity assumption:
\begin{assumption}[Second-Order Self-Bounding Regularity]\label{ass:selfbounding}
$F$ is twice differentiable, and there exists a non-decreasing function $\rho_1: \bbR_{\ge 0} \mapsto \bbR_{\ge 0}$ such that $\nrm*{\grad^2 F(\vecW)}_{\OPNORM} \le \rho_1\prn*{F(\vecW)} \text{ for all }\vecW\in\mathbb{R}^d$.
\end{assumption}
This assumption implies the relevant Hessian operator norm is upper bounded by a function of the function value. It was also made in \citet{priorpaper} for the different task of studying global convergence of GD/SGD, where it was shown that \pref{ass:selfbounding} holds for many canonical non-convex optimization problems. Some quantitative control of the Hessian is necessary for non-asymptotic guarantees of finding FOSPs \citep{kornowski2024hardness}. In \pref{ex:breaklit}, we show these prior assumptions are not satisfied by a natural univariate function. We show in \pref{prop:comparetolit} that \pref{ass:selfbounding} \textbf{generalizes} $(L_0,L_1)$-smoothness and its extension from \citet{li2023convex}, and that
\[ (L_0, L_1)\text{-smoothness }(\nrm*{\grad^2 F}_{\OPNORM} \le L_0+L_1\nrm*{\grad F})\implies\text{ \pref{ass:selfbounding} with }\rho_{1}(x)=\tfrac32L_{0}+4L_{1}^{2}x.\]
For finding SOSPs, we impose the following additional regularity assumption:
\begin{assumption}[Third-Order Self-Bounding Regularity]\label{ass:thirdorderselfbounding}
$F$ satisfies \pref{ass:selfbounding}, and either:
\begin{itemize}
    \item $F$ is three-times differentiable everywhere, and for some non-decreasing function $\rho_2:\mathbb{R}_{\ge 0} \rightarrow \mathbb{R}_{\ge 0}$, $\nrm*{\grad^3 F(\vecW)}_{\OPNORM} \le \rho_2\prn*{F(\vecW)}$ for all $\vecW\in\mathbb{R}^d$.
    \item Or for some constant $\delta>0$ and some non-decreasing function $\rho_2:\mathbb{R}_{\ge 0} \rightarrow \mathbb{R}_{\ge 0}$, for all $\vecW, \vecW' \in\mathbb{R}^d$ with $\nrm*{\vecW-\vecW'} \le \delta$, we have $\nrm*{\grad^2 F(\vecW) - \grad^2 F(\vecW')}_{\OPNORM} \le \rho_2\prn*{F(\vecW)} \nrm*{\vecW-\vecW'}$.
\end{itemize}
\end{assumption}
\pref{ass:thirdorderselfbounding} naturally extends \pref{ass:selfbounding}, and generalizes the Hessian Lipschitz assumption ubiquitous in the literature on \textit{non-asymptotic} rates for finding SOSPs. (We note that the works \citet{lee2016gradient, lee2019first} established \textit{asymptotic} guarantees for GD finding SOSPs without the Hessian Lipschitz assumption, and note their proof strategy uses Lipschitzness of the gradient in a crucial way.)
In \pref{subsec:examplesmainbody}, we show several canonical non-convex losses with non-Lipschitz gradient and Hessian satisfy \pref{ass:thirdorderselfbounding}. \pref{ass:thirdorderselfbounding} covers several growth rates of interest (e.g. univariate self-concordant functions satisfying \pref{ass:selfbounding}). It also subsumes that of \citet{xie2024trust}, which to our knowledge is the only other result on finding SOSPs under generalized smoothness (but uses the full Hessian). Under the assumptions of \citet{xie2024trust}, an explicit, simple form for $\rho_2(\cdot)$ can be found. We detail all of this in \pref{ex:thirdordersmoothnessjustification}. 

Furthermore, \pref{ass:thirdorderselfbounding} encompasses several examples of Distributionally Robust Optimization (DRO) problems. \citet{xie2024trust} very interestingly demonstrates that under mild assumptions, the objective of DRO satisfies their Assumption 3, see Theorem 3 therein. Assumption 3 of \citet{xie2024trust} is subsumed by \pref{ass:thirdorderselfbounding} as per our \pref{ex:thirdordersmoothnessjustification}. Thus our results apply to DRO. DRO is a general optimization problem that has significant applications in fairness in machine learning and in learning under distribution shifts; see \citet{xie2024trust} for more discussion.

We now introduce the following standard definition, which, when combined with \pref{ass:selfbounding} and \pref{ass:thirdorderselfbounding}, forms the core of our argument, as we explain in \pref{subsec:highlevelideaframework}.
\begin{definition}\label{def:regularpoint}
For a function $F$ and threshold $\alpha$, the $\alpha$-sublevel set of $F$ is $\cL_{F,\alpha} = \{\vecW: F(\vecW) \le \alpha\}$.
\end{definition}
Now, our contributions are as follows:
\begin{enumerate}
\item \textbf{We develop a novel, systematic framework} detailed in \pref{sec:mainidea} and \pref{thm:generalframework} to study the convergence of first order methods to FOSPs and SOSPs under \pref{ass:selfbounding} and \pref{ass:thirdorderselfbounding} respectively. The core idea is in \pref{subsec:highlevelideaframework}. \textbf{Our framework lets us systematically analyze existing practical, and widely used first-order optimization algorithms in the challenging generalized smooth setting.}

\item \textbf{Main Results, non-asymptotic convergence to SOSPs:} Under \pref{ass:thirdorderselfbounding}, we establish efficient rates for first-order optimization algorithms finding SOSPs. See \pref{thm:escapesecondordergd} for Perturbed GD \citep{jin2017escape} and \pref{thm:sgdsecondorder} for Restarted SGD \citep{fang2019sharp}. The dependence on $\epsilon, d$ matches that in the smooth setting, and in particular is polylogarithmic in $d$. This is particularly pertinent for ML applications, where the ambient dimension is so large that the second-order methods of \citet{xie2024trust} are not feasible. 

\item Non-asymptotic convergence to FOSPs: Under \pref{ass:selfbounding}, we establish efficient rates for GD, Adaptive GD, and SGD finding FOSPs. See \pref{thm:gdfirstorder}, \pref{thm:adaptivegdguarantee}, and \pref{thm:sgdfirstorder} respectively. The dependence on $\epsilon, d$ again matches that in the smooth setting.

\item We provide examples and practical implications in \pref{subsec:examplesmainbody}. Our examples are direct corollaries of \pref{thm:escapesecondordergd}, \pref{thm:sgdsecondorder}. They show variants of GD/SGD globally optimize non-convex `strict-saddle' losses from ML with non-Lipschitz gradient and Hessian.
\end{enumerate}

\textbf{Notation:} 
$\ball(\vecP,R)$ denotes the Euclidean $l_2$ ball centered at $\vecP\in\mathbb{R}^d$ with radius $R \ge 0$, with boundary. 
By shifting, we assume WLOG that $F$ attains a minimum value of $0$. We follow the convention that $F$ is smooth, specifically $L$-smooth, if $\nrm*{\grad^2 F} \le L$ holds globally. 
We always let $\vecW_0$ denote the initialization of a given algorithm (which is clear from context) unless stated otherwise.

\section{Main Idea}\label{sec:mainidea}
\subsection{High Level Idea}\label{subsec:highlevelideaframework}
One classic analysis of GD on smooth functions to converge to a FOSP goes by establishing decrease per iterate, via the so-called `Descent Lemma' \citep{bubeck2015convex}. For $L$-smooth functions, setting the step size $\eta = \frac1{L}$ in GD,
\begin{align*}
F(\vecW_{t+1}) &\le F(\vecW_t) - \eta \prn*{1-\frac12 L\eta} \nrm*{\grad F(\vecW_t)}^2 = F(\vecW_t)-\frac1{2L} \nrm*{\grad F(\vecW_t)}^2.\numberthis\label{eq:descentlemma}
\end{align*}
Such an analysis fails if $F$ is not $L$-smooth. Following the above recipe under \pref{ass:selfbounding}, as such a bound $L<\infty$ need not exist, one must set $\eta=0$ and does not obtain any convergence rate. 

\textbf{Core Insight 1: } The first simple but powerful insight in our work is that many optimization algorithms such as GD decrease the function value at each iterate (with high probability) when $\eta$ is appropriately chosen as a function of the smoothness (Hessian operator norm) at the \textit{current} iterate.

Specifically, consider iterates of GD initialized at some $\vecW_0$. For step size $\eta$ small enough in terms of $\nrm*{\grad^2 F (\vecW_0)}$, the next iterate $\vecW_1$ of GD is sufficiently `local' (see \pref{corr:gradcontrol}). This lets us upper bound $\nrm*{\grad^2 F}$ along the segment $\overline{\vecW_0 \vecW_1}$ by an increasing function $L_1(F(\vecW_0))$ of $F(\vecW_0)$ (see \pref{lem:boundlocalsmoothnesseasy}). Thus, for appropriate $\eta$ in terms of $F(\vecW_0)$, we obtain $F(\vecW_1) \le F(\vecW_0)$, and so $\vecW_1$ lies in the $F(\vecW_0)$-sublevel set $\cL_{F,F(\vecW_0)}$.

\textbf{Core Insight 2:} Crucially, we can `chain together' this decrease. 
By \pref{ass:selfbounding}, the aforementioned argument goes through at any $\vecW$ in the $F(\vecW_0)$-sublevel set $\cL_{F,F(\vecW_0)}$ -- in particular, at $\vecW_1$. Consequently, this \textit{same} step size $\eta$ is small enough to ensure $F(\vecW_2) \le F(\vecW_1) \le F(\vecW_0)$, and so forth through all the iterates of GD.
Moreover, this argument yields a convergence rate. As each iterate is in $\cL_{F, F(\vecW_0)}$, if the gradient norm is at least $\epsilon$ at each iterate, we obtain decrease of at least $\tfrac{\epsilon^2}{2L_1(F(\vecW_0))}$ per iterate analogously to \pref{eq:descentlemma}. Too many iterations contradicts that $F$ is lower bounded by 0 (recall Notation), so we must reach an iterate $\vecW_t$ which is a FOSP within $\tfrac{2L_1(F(\vecW_0))F(\vecW_0)}{\epsilon^2}$ iterates. 

\textbf{Generalizing the argument:} This idea is powerful enough to readily analyze SGD and variants of GD/SGD which find SOSPs. Rather than a single iterate where decrease need not hold, we consider a sequence of consecutive $\Tthres$ iterates. We show with high probability, the last iterate in this sequence decreases function value for $\vecW \in \cL_{F,F(\vecW_0)}$. To do so, recall the analyses of first-order optimization algorithms often establish decrease by considering `local' behavior. Locally around $\vecW \in \cL_{F,F(\vecW_0)}$, \pref{ass:selfbounding} and \pref{ass:thirdorderselfbounding} give enough control over the relevant derivatives to do so.

Then the above argument still goes through, with a fixed step size defined in terms of $F(\vecW_0)$. We group the iterates of the algorithm into `blocks' of length $\Tthres$, and establish $F(\vecW_{\Tthres}) \le F(\vecW_0)$ and so forth (rather than establishing $F(\vecW_2) \le F(\vecW_1) \le F(\vecW_0)$ for consecutive iterates). 

\subsection{The Formal Framework}\label{subsec:formalframework}
Consider a set of interest $\cS$, e.g. FOSPs or SOSPs with tolerance $\epsilon$.
We begin by presenting a simpler version of our formal framework. Consider a deterministic update procedure $\cA:\mathbb{R}^d\rightarrow \mathbb{R}^d$, where the output of $\cA$ denotes the future iterate of the algorithm. For example, $\cA(\vecW) = \vecW - \eta \grad F(\vecW)$ for GD. 
Following \pref{subsec:highlevelideaframework}, we consider algorithms that decrease function value in the $F(\vecW_0)$-sublevel set $\cL_{F, F(\vecW_0)}$ if they have not reached $\cS$. The following definition formalizes this property:
\begin{definition}[Special case of Decrease Procedure in \pref{def:highprobdecreasealg}]
Consider a set of interest $\cS$, a decrease threshold $\Delta > 0$, a point $\vecU_0$, and a deterministic procedure $\cA$ to compute the next iteration. 
We say $\cA$ forms a $(\cS, \Toracle(\vecU_0), \Delta(\vecU_0), \vecU_0)$-decrease procedure if  computing $\cA(\vecU_0)$ takes at most $\Toracle(\vecU_0)$ oracle calls, and one of the following holds:
\[ 1)\text{ } F(\cA(\vecU_0)) < F(\vecU_0) - \Delta(\vecU_0),\hspace{0.7cm}\text{ or }\hspace{0.7cm}2)\text{ } \cA(\vecU_0) \cap \cS \neq \{ \}. \]
\end{definition}
Here 1) means that the subsequent iterate has smaller function value, and 2) means that the rule of output $\cA_2$ outputs a sequence of candidate vectors, one of which is in $\cS$. 

Then, \pref{thm:generalframework} states that if $\cA$ is a decrease procedure for all $\vecU_0$ in $\cL_{F, F(\vecW_0)}$, we can bound the number of oracle calls for $\cA$ to output a candidate vector in $\cS$, e.g. for GD to output a FOSP. We prove it arguing as in \pref{subsec:highlevelideaframework}, `chaining together' the decrease per iterate in $\cL_{F, F(\vecW_0)}$. Then as $F$ is lower bounded, 1) in \pref{def:highprobdecreasealg} cannot occur too often, so 2) must occur at some point.

We now generalize this to randomized procedures $\cA$ which can output several candidate vectors. 
\paragraph{Framework in full generality. } Consider an update procedure $\cA:\mathbb{R}^d\rightarrow \mathbb{R}^d \times \bigcup_{n=0}^{\infty} (\mathbb{R}^d)^n$ (possibly randomized).
We now consider a map $\cA=(\cA_1, \cA_2)$, $\cA:\mathbb{R}^d \rightarrow \mathbb{R}^d \times \bigcup_{n=0}^{\infty} (\mathbb{R}^d)^n$ defined as follows:
\[ \text{For all }\vecU\in\mathbb{R}^d, \cA(\vecU) = \prn*{\vecP_1, \vecP_2} \text{ for }\vecP_1 \in \mathbb{R}^d, \vecP_2 \in \bigcup_{n=0}^{\infty} (\mathbb{R}^d)^n, \text{ and define }\cA_1(\vecU) := \vecP_1, \cA_2(\vecU) := \vecP_2.\]
Intuitively, $\cA_1$ computes a future iterate $\cA_1(\vecU)$. $\cA_2$ outputs a sequence of candidate vectors in $\mathbb{R}^d$, among which we hope one lies in $\cS$ (e.g. different candidate models in statistical learning).

However, the output of $\cA_1$ need not correspond to the `next iterate' in the traditional sense. For SGD, $\cA_1$ does \textit{not} output the next iterate of SGD, but rather the iterate produced by SGD after $K_0 > 1$ steps. 
This is necessary to guarantee decrease; a single step of SGD need not decrease the value of $F$, but with high probability and large enough $K_0$, a consecutive `block' of $K_0$ iterates will. We will lay this out concretely next in \pref{subsec:exinframework}.

\begin{remark}\label{rem:whymultipleoutput}
Often $\cA_2$ will output a single vector in $\mathbb{R}^d$, which we hope lies in $\cS$, but this is not always the case. Consider guarantees for GD or SGD, which upper bound $\frac1T \sum_{t=1}^T \nrm*{\grad F(\vecW_t)}^2 \le \epsilon^2$ or $\frac1T \sum_{t=1}^T \nrm*{\grad F(\vecW_t)} \le \epsilon$. This only ensures a single $\vecW_t \in \cS, 1 \le t \le T$ where $\cS$ is the set of FOSPs to tolerance $\epsilon$ (e.g. \citet{zhang2019gradient}, \citet{jin2021smoothness}, \citet{li2023convergence}, \citet{xie2024trust} and many others). Consequently $(\vecW_1, \ldots, \vecW_T)$ is our sequence of candidate vectors, and the guarantee obtained is that $\vecW_t \in \cS$ for some $1 \le t \le T$. We thus allow for $\cA_2$ to output multiple candidate vectors.
\end{remark}

The following definition formalizes a common property of optimization algorithms we study:
\begin{definition}[Decrease Procedure]\label{def:highprobdecreasealg}
Consider a set of interest $\cS$, a confidence parameter $\delta > 0$, a decrease threshold $\Delta > 0$, a point $\vecU_0$, and a procedure $\cA$ to compute the next iteration. 
We say $\cA$ forms a $(\cS, \Toracle(\vecU_0), \Delta(\vecU_0), \delta(\vecU_0), \vecU_0)$-decrease procedure if with probability at least $1-\delta(\vecU_0)$ over the randomness in $\cA$ to compute $\cA(\vecU_0)$ from $\vecU_0$, computing $\cA(\vecU_0)$ takes at most $\Toracle(\vecU_0)$ oracle calls, and one of the following holds:
\[ 1)\text{ } F(\cA_1(\vecU_0)) < F(\vecU_0) - \Delta(\vecU_0),\hspace{0.7cm}\text{ or }\hspace{0.7cm}2)\text{ } \cA_2(\vecU_0) \cap \cS \neq \{ \}. \]
\end{definition}
Here 1) means that the subsequent iterate has smaller function value, and 2) means that the rule of output $\cA_2$ outputs a sequence of candidate vectors, one of which is in $\cS$. $\cA$ forms a $(\cS, \Toracle(\vecU_0), \Delta(\vecU_0), \delta(\vecU_0), \vecU_0)$-decrease procedure if 1) or 2) occurs with high probability. 

\textbf{Informal Theorem:} For analogous reasons as before, we will establish that if $\cA$ is a decrease procedure for all $\vecU_0$ in $\cL_{F,F(\vecW_0)}$, we can bound the number of oracle calls for $\cA_2$ to output a candidate vector lying in $\cS$. Formally, this is \pref{thm:generalframework}. 

\subsection{Examples Subsumed by Framework}\label{subsec:exinframework}
We demonstrate that a host of first-order optimization algorithms are covered in our framework, and highlight the general recipe for using our framework.
\paragraph{GD:} Starting from $\vecU$, the next iterate of GD with step size $\eta>0$ is $\vecU - \eta \grad F(\vecU)$.
\begin{enumerate}
    \item For $\epsilon>0$, let $\cS =\{\vecW:\nrm*{\grad F(\vecW)} \le \epsilon\}$, the set of FOSPs.
    \item For all $\vecU_0 \in \mathbb{R}^d$, let $\cA(\vecU_0) = \prn*{\vecU_0 - \eta \grad F(\vecU_0), \vecU_0}$. Hence, $\cA_1(\vecU_0) = \vecU_0 - \eta \grad F(\vecU_0)$, $\cA_2(\vecU_0)=\vecU_0$, and $\Toracle(\vecU_0)=1$.
    \item \textbf{In \pref{claim:gdfirstorderdecreaseprocedure}, we establish that if $F$ is satisfies \pref{ass:selfbounding}, then $\cA$ is a decrease procedure for all $\vecU_0 \in \cL_{F, F(\vecW_0)}$, for suitable $\eta$ depending on $F(\vecW_0)$.}
Our result for GD, \pref{thm:gdfirstorder}, subsequently follows by our general framework \pref{thm:generalframework}.
\end{enumerate}

\paragraph{Adaptive GD:}
Starting from $\vecU$, the next iterate of Adaptive GD is $\vecU - \eta_{\vecU} \grad F(\vecU)$,
where $\eta_{\vecU}>0$ is an adaptive step size that depends on $\vecU$. 
\begin{enumerate}
    \item For $\epsilon>0$, let $\cS =\{\vecW:\nrm*{\grad F(\vecW)} \le \epsilon\}$, the set of FOSPs.
    \item For all $\vecU_0 \in \mathbb{R}^d$, let $\cA(\vecU_0) = \prn*{ \vecU_0 - \eta_{\vecU_0} \grad F(\vecU_0), \vecU_0}$. Hence, $\cA_1(\vecU_0) = \vecU_0 - \eta \grad F(\vecU_0)$, $\cA_2(\vecU_0)=\vecU_0$, and $\Toracle(\vecU_0)=1$.
    \item \textbf{In \pref{claim:adaptivegdfirstorderdecreaseprocedure}, we establish that if $F$ is satisfies \pref{ass:selfbounding}, then $\cA$ is a decrease procedure for all $\vecU_0 \in \cL_{F, F(\vecW_0)}$, for suitable $\eta_{\vecU}$ depending on $F(\vecW_0)$ and $\nrm*{\grad F(\vecU)}$.} Our result for Adaptive GD, \pref{thm:adaptivegdguarantee}, then follows by \pref{thm:generalframework}.
\end{enumerate}
However, for SGD and other randomized algorithms involving randomness, 1) in \pref{def:highprobdecreasealg} does not hold deterministically. This is where the generality in our framework is powerful. For SGD, by concentration inequalities we show that 1) is true \textit{with high probability over a long enough `block' of subsequent iterates}, as long as none of the iterates in the block have small gradient. We then define $\cA$ so that $\cA_1$ outputs the composition of the SGD steps in the block, and $\cA_2$ outputs all the iterates of the block. The resulting guarantee is that one of the points among all the blocks lies in $\cS$.
\paragraph{SGD:} Starting from $\vecU$, letting $\grad f(\vecU; \vecZeta)$ be a stochastic gradient oracle where $\vecZeta$ is a minibatch sample, the next iterate of SGD is $\vecU - \eta \grad f(\vecU;\vecZeta)$ where $\eta>0$ is the step size.
\begin{enumerate}
    \item For $\epsilon>0$, let $\cS =\{\vecW:\nrm*{\grad F(\vecW)} \le \epsilon\}$, the set of FOSPs.
    \item Consider any $K_0 \ge 1$. For all $\vecU_0 \in \mathbb{R}^d$, let $\vecP_0=\vecU_0$, and define a sequence $(\vecP_i)_{0 \le i \le K_0}$ via $\vecP_i = \vecP_{i-1} - \eta \grad f(\vecP_{i-1};\vecZeta_{i})$, where the $\vecZeta_{i}$ are i.i.d. minibatch samples. Note this sequence can be equivalently defined by repeatedly composing the function $\vecU \rightarrow \vecU - \eta \grad f(\vecU; \vecZeta)$.
We then define $\cA(\vecU_0) = \prn*{ \vecP_{K_0}, (\vecP_i)_{0 \le i \le K_0-1}}$, hence $\cA_1(\vecU_0)=\vecP_{K_0}$, $\cA_2(\vecU_0)=(\vecP_i)_{0 \le i \le K_0-1}$. Note all the $\vecP_i$ are a function of $\vecU_0$ and the randomness in the stochastic gradient oracle $\grad f(\cdot;\cdot)$. We let $\Toracle(\vecU_0)=K_0$, which need not equal 1. This procedure is clearly SGD, with its iterates divided into blocks of length $K_0$.

\item \textbf{In \pref{claim:sgdfirstorderdecreasesequence}, we establish that if $F$ is satisfies \pref{ass:selfbounding} and $\grad f(\cdot;\cdot)$ satisfies \pref{ass:noiseregularity}, then $\cA$ is a decrease procedure for all $\vecU_0 \in \cL_{F, F(\vecW_0)}$ for suitable algorithm parameters.} Our result for SGD, \pref{thm:sgdfirstorder}, then follows by \pref{thm:generalframework}.
\end{enumerate}
\paragraph{SOSP-finding algorithms:} We now study finding SOSPs using first order methods under our regularity assumptions. We analyze two algorithms to achieve this under exact and stochastic gradients, respectively Perturbed GD (\pref{alg:perturbedgd}, \citet{jin2017escape}) and Restarted SGD (\pref{alg:restartedsgd}, \citet{fang2019sharp}). We remark that our framework likely subsumes many other algorithms. 

\paragraph{Perturbed GD:} This algorithm, formally written in \pref{alg:perturbedgd}, \pref{sec:perturbedgdproofs}, is as follows. At $\vecU$,
\begin{itemize}
    \item If $\nrm*{\grad F(\vecU)} > \Gthres$ for some appropriate $\Gthres$, the algorithm simply runs a step of GD. 
    \item Else, \pref{alg:perturbedgd} adds uniform noise from a ball with particular radius and runs GD for $\Tthres$ iterations for suitably chosen $\Tthres$, yielding $\vecU'$. We check if $F(\vecU') - F(\vecU) \le -\Fthres$ for some appropriate $\Fthres$. If decrease does not occur, we return $\vecU$; if decrease occurred, we go back to the If/Else with $\vecU'$ in place of $\vecU$. 
\end{itemize}
Notice now that the oracle complexity $\Toracle$, probability $\delta$, and amount of decrease $\Delta$ depend on the location $\vecU$. Our framework readily subsumes this example as follows. 
\begin{enumerate}
\item For $\epsilon>0$, let $\cS =\{\vecW:\nrm*{\grad F(\vecW)} \le \epsilon, \grad^2 F(\vecW) \succeq -\sqrt{\epsilon} \matI\}$, the set of SOSPs.
\item For all $\vecU_0 \in \mathbb{R}^d$, if $\nrm*{\grad F(\vecU_0)} > \Gthres$, we let 
\[ \cA(\vecU_0) = \prn*{\vecU_0 - \eta \grad F(\vecU_0), \vecU_0}, \text{ hence }\cA_1(\vecU_0)=\vecU_0 - \eta \grad F(\vecU_0), \cA_2(\vecU_0) = \vecU_0.\]

Otherwise if $\nrm*{\grad F(\vecU_0)} \le \Gthres$, we let $\vecP_0 = \vecU_0 + \vecXi$ where $\vecXi$ is uniform from $\ball(\vecOrigin, r)$, and define a sequence $(\vecP_i)_{0 \le i \le \Tthres}$ via $\vecP_i = \vecP_{i-1} - \eta \grad F(\vecP_{i-1})$.
We then define
\[ \cA(\vecU_0) = \prn*{\vecP_{\Tthres}, \vecU_0}, \text{ hence }\cA_1(\vecU_0)=\vecP_{\Tthres}, \cA_2(\vecU_0) = \vecU_0.\]
Thus
\[ \Toracle(\vecU_0) = \begin{cases} \Tthres &: \nrm*{\grad F(\vecU_0)} \le \Gthres \\ 1 &:\nrm*{\grad F(\vecU_0)} >\Gthres. \end{cases}\]
This is identical to \pref{alg:perturbedgd}, and highlights why $\Toracle, \delta, \Delta$ need to depend on $\vecU_0$. 
\item \textbf{In \pref{claim:perturbedgddecreasesequence}, we establish that if $F$ satisfies \pref{ass:thirdorderselfbounding}, then $\cA$ is a decrease procedure for all $\vecU_0 \in \cL_{F, F(\vecW_0)}$ for suitable algorithm parameters.} Our result for Perturbed GD, \pref{thm:escapesecondordergd}, then follows by \pref{thm:generalframework}.
\end{enumerate}

\paragraph{Restarted SGD:} This algorithm, formally written in \pref{alg:restartedsgd}, \pref{sec:restartedsgdproofs}, works as follows. Take $B = \tilde{\Theta}(\epsilon^{0.5})$, $K_0 = \Tilde{\Theta}(\epsilon^{-2})$. Consider an anchor point $\vecU$, first taken to be the initialization $\vecW_0$. The algorithm runs SGD until its iterates first escape the ball $\ball(\vecU, B)$, tracking at most $K_0$ iterations.
\begin{itemize}
    \item If an escape occurs within $K_0$ iterations, letting $\vecU'$ be the first iterate that escaped $\ball(\vecU, B)$, the algorithm sets $\vecU'$ to be the anchor point and runs the same procedure.
    \item If these $K_0$ iterates do not escape within $K_0$ iterations, return their \textit{average}.
\end{itemize}
We cover Restarted SGD in our framework as follows.
\begin{enumerate}
    \item For $\epsilon>0$, let $\cS =\{\vecW:\nrm*{\grad F(\vecW)} \le \epsilon, \grad^2 F(\vecW) \succeq -\sqrt{\epsilon} \matI\}$, the set of SOSPs.
    \item For all $\vecU_0 \in \mathbb{R}^d$, let $\vecP_0=\vecU_0$. We define a sequence $(\vecP_i)_{0 \le i \le K_0}$ via $\vecP_i = \vecP_{i-1} - \eta (\grad f(\vecP_{i-1};\vecZeta_{i})+\tilde{\sigma} \Lambda^{i})$, where $\grad f(\cdot;\cdot)$ is our stochastic gradient oracle, the $\vecZeta_{i}$ are i.i.d. minibatch samples, the $\Lambda^{i} \sim \ball(\vecOrigin,1)$ are i.i.d., and $\tilde{\sigma}$ is a parameter governing the noise level. Note this sequence can be equivalently defined by repeatedly composing the function $\vecU \rightarrow \vecU - \eta (\grad f(\vecU; \vecZeta)+\tilde{\sigma}\Lambda)$.
If it exists, let $i, 1 \le i \le K_0$ be the minimal index such that $\nrm*{\vecP_i - \vecP_0} > B$. Otherwise let $i=K_0$. In either case, we define
\[ \cA(\vecU_0) = \prn*{\vecP_i, \frac1{i} \sum_{t=0}^{i-1} \vecP_t}, \text{ hence }\cA_1(\vecU_0)=\vecP_i, \cA_2(\vecU_0) = \frac1{i} \sum_{t=0}^{i-1} \vecP_t. \]
We let $\Toracle(\vecU_0)=K_0$.\footnote{Defining $i$ as above, note that we can compute $\cA(\vecU_0)$ using $i$ rather than $K_0$ oracle calls, but this change does not affect runtime beyond constant factors.} This is clearly identical to \pref{alg:restartedsgd}. 
\item \textbf{In \pref{claim:restartedsgddecreasesequence}, we establish that if $F$ satisfies \pref{ass:thirdorderselfbounding} and $\grad f(\cdot;\cdot)$ satisfies \pref{ass:noiseregularity} and \pref{ass:gradoraclecontrol}, then $\cA$ is a decrease procedure for all $\vecU_0 \in \cL_{F, F(\vecW_0)}$ for suitable algorithm parameters.} Our result for Restarted GD, \pref{thm:sgdsecondorder}, then follows by \pref{thm:generalframework}.
\end{enumerate}

\begin{theorem}[General Framework]\label{thm:generalframework}
Consider a given initialization $\vecW_0$ of $\cA$ and a desired set $\cS$. Define a sequence $(\vecW_t)_{t \ge 0}$ recursively by $\vecW_{t+1}=\cA_1(\vecW_t)$. Suppose that for all $\vecU_0 \in \cL_{F,F(\vecW_0)}$, $\cA$ forms a $(\cS, \Toracle(\vecU_0), \Delta(\vecU_0), \delta(\vecU_0), \vecU_0)$-decrease procedure. 
Define $\overline{\Delta} = \inf_{\vecU \in \cL_{F, F(\vecW_0)} }\frac{\Delta(\vecU)}{\Toracle(\vecU)}$.
Then with probability at least
\[1-\sup_{\vecU \in \cL_{F, F(\vecW_0)} }\delta(\vecU) \cdot \sup_{\vecU \in \cL_{F, F(\vecW_0)} }\crl*{\frac{F(\vecW_0)}{\Delta(\vecU)}},\text{ upon making } N = \frac{F(\vecW_0)}{\overline{\Delta}}+\sup_{\vecU \in \cL_{F, F(\vecW_0)} }\Toracle(\vecU)\]
oracle calls, there exists $\vecW_t \in (\vecW_t)_{t \ge 0}$ such that $\cA_2(\vecW_t) \cap \cS \neq \{ \}$. I.e. for some $\vecW_t$, $\cA_2(\vecW_t)$ will output a sequence of candidate vectors, one of which is in $\cS$. Furthermore, if the output of $\cA_2$ has length at most $S$, then the number of candidate vectors outputted is at most $S \cdot \sup_{\vecU \in \cL_{F, F(\vecW_0)} }\crl*{\frac{F(\vecW_0)}{\Delta(\vecU)}}$.
\end{theorem}
Our full proof is in \pref{sec:generalframeworkpf}.\footnote{The extra second term in the sum defining $N$ occurs as $\Toracle, \Delta, \delta$ have $\vecU_0$-dependence.} The proof formalizes the main idea from \pref{subsec:highlevelideaframework}, by `chaining together' the decrease per iterate in $\cL_{F, F(\vecW_0)}$. Then as $F$ is lower bounded, 1) in \pref{def:highprobdecreasealg} cannot occur too many times, so 2) must occur at some point.
\begin{remark}\label{rem:frameworklocalnecessity}
To verify $\cA$ is a decrease procedure in $\cL_{F, F(\vecW_0)}$, we can systematically port over analyses in the literature. As discussed in \pref{subsec:highlevelideaframework}, $\vecU_0$ being in $\cL_{F, F(\vecW_0)}$ allows us to show the algorithm is `local', crucially giving us quantitative control over the relevant derivatives. 
We view this as a core \textit{strength} of our work; our framework allows us to \textit{systematically} extend results from the smooth setting to generalizations of smoothness. 
\end{remark}

\section{Convergence Results}\label{sec:convergenceresults}
Here we systematically obtain our convergence results for the algorithms listed in \pref{subsec:exinframework}, by formally showing that they are decrease procedures. \textbf{Our main results are \pref{thm:escapesecondordergd}, \pref{thm:sgdsecondorder}: that under \pref{ass:thirdorderselfbounding}, variants of GD/SGD can find SOSPs.} We note our dependence on $\epsilon, d$ for \pref{thm:gdfirstorder}, \pref{thm:adaptivegdguarantee}, \pref{thm:sgdfirstorder}, and \pref{thm:sgdsecondorder} match lower bounds for smooth functions \citep{carmon2020lower, carmon2021lower, arjevani2020second}, and hence are optimal in this setting too.\footnote{Dependence on $\epsilon$ in \pref{thm:sgdfirstorder} and on $\epsilon, d$ in \pref{thm:sgdsecondorder} are tight up to $\log$ factors.} We present examples and implications of our results in \pref{subsec:examplesmainbody}.

\begin{remark}[Dependence on Initialization]\label{rem:stepsizeoninitialization}
In our results, the step size $\eta$ here depends only on $\rho_1(F(\vecW_0))$, a fixed value depending only on initialization. 
Moreover, the expressions on $\eta$ depending on $\rho_1(F(\vecW_0))$ in our results and proofs to follow are only an \textbf{upper bound} for working step sizes. We do not need to know these exact values. Therefore, all that is needed is an upper bound on fixed quantities such as $\rho_1(F(\vecW_0))$; hence a working step size $\eta$ for our algorithms in practice and theory can be found using cross validation or binary search. 

Letting $\eta(\vecW_0)$ be an upper bound on the step size $\eta$ needed to guarantee convergence, we note by searching over $\log\prn*{\eta(\vecW_0)}$ with binary search, we will find an $\eta$ with a constant factor 2 of $\eta(\vecW_0)$. This $\log$ factor will be logarithmic in $\epsilon, d$, and will only change the claimed iteration complexity by a universal constant factor. The latter is because the amount of decrease in the definition of Decrease Procedure will in turn only change by a universal constant multiple.
\end{remark}

\begin{remark}[On Adaptivity]\label{rem:adaptivityremark}
Our results hold for non-adaptive versions of GD/SGD and their variants. That said, one can interpret cross validation or binary search over $\eta$ 
as adaptive algorithms in their own right. As mentioned above, it is relatively straightforward to obtain analogous results to our current ones for cross validation or binary search. In the learning from data setting, one can make the cross validation result formal using classic techniques.
\end{remark}

\subsection{Gradient Descent}\label{subsec:fospgd}
\begin{theorem}[GD for FOSP]\label{thm:gdfirstorder}
Suppose $F$ satisfies \pref{ass:selfbounding}. Run GD initialized at $\vecW_0$, with step size $\eta = \frac1{L_1(\vecW_0)}$ where $L_1(\vecW_0)$ is defined in \pref{eq:L1def}. Then letting 
\[ T=\frac{2 F(\vecW_0) L_1(\vecW_0)}{\epsilon^2}, \text{ within $T+1$ oracle calls to $\grad F(\cdot)$,} \]
GD will output $T$ candidate vectors $(\vecP_1, \ldots, \vecP_T)$, one of which satisfies $\nrm*{\grad F(\vecP_t)} \le \epsilon$.
\end{theorem}
We prove \pref{thm:gdfirstorder} here to show our strategy's simplicity. The following Lemmas, proved in \pref{subsec:technicalhelperresults}, help show GD is `local' for $\vecW \in \cL_{F,F(\vecW_0)}$.
\begin{corollary}\label{corr:gradcontrol}
For $F$ satisfying \pref{ass:selfbounding}, we have $\nrm*{\grad F(\vecW)} \leq \rho_0\prn*{F(\vecW)}$, where $\rho_0:\mathbb{R}_{\ge 0}\rightarrow\mathbb{R}_{\ge 0}$ is a non-decreasing function given by $\rho_0(x) = \rho_1(x)\sqrt{2\theta(x)}$, where $\theta(x) = \int_0^x \frac1{\rho_1(v)} \DERIV v$.
\end{corollary}
\begin{lemma}\label{lem:boundfuncvalueradius}
Under \pref{ass:selfbounding}, for $\vecX, \vecY$ with $\nrm*{\vecY-\vecX} \le \frac{1}{\rho_0(F(\vecX)+1)}$, $F(\vecY) - F(\vecX) \le 1$.
\end{lemma}
Combining the above with \pref{ass:selfbounding} immediately gives:
\begin{lemma}\label{lem:boundlocalsmoothnesseasy}
Suppose $F$ satisfies \pref{ass:selfbounding}. Defining $\rho_0$ as in \pref{corr:gradcontrol}, let
\[  L_1(\vecW_0) = \max\crl*{1, \rho_0(F(\vecW_0)+1), \rho_0(F(\vecW_0))\rho_0(F(\vecW_0)+1), \rho_1(F(\vecW_0)+1)}.\numberthis\label{eq:L1def} \]
Then for all $\vecW \in \cL_{F,F(\vecW_0)}$, $\nrm*{\grad^2 F(\vecU)}_{\OPNORM} \le L_1(\vecW_0)$ for all $\vecU \in \ball\prn*{\vecW, \rho_0(F(\vecW_0)+1)^{-1}}$.
\end{lemma}

\begin{proof}[Proof of \pref{thm:gdfirstorder}]
Use \pref{thm:generalframework} with $\cS = \{\vecW:\nrm*{\grad F(\vecW)} \le \epsilon\}$, defining $\cA$ as in \pref{subsec:exinframework}. Upon applying \pref{thm:generalframework}, the following Claim directly proves \pref{thm:gdfirstorder}:
\begin{claim}\label{claim:gdfirstorderdecreaseprocedure}
For any $\vecU_0$ in $\cL_{F,F(\vecW_0)}$, $\cA$ is a $(\cS, 1, \tfrac{\epsilon^2}{2 L_1(\vecW_0)}, 0, \vecU_0)$-decrease procedure.
\end{claim}
To prove \pref{claim:gdfirstorderdecreaseprocedure}, note for $\vecU_0 \in \cS$, by definition of $\cA_2$ that $\cA_2(\vecU_0) = (\vecU_0) \in \cS$. Now if $\vecU_0 \not\in \cS$ (i.e. $\nrm*{\grad F(\vecU_0)} > \epsilon$), consider $\vecU_1 = \cA_1(\vecU_0) = \vecU_0 - \eta \grad F(\vecU_0)$. 
By \pref{corr:gradcontrol} and as $F(\vecU_0) \le F(\vecW_0)$, $\nrm*{\grad F(\vecU_0)} \le \rho_0(F(\vecU_0)) \le \rho_0(F(\vecW_0))$, so by choice of $\eta$,
\[ \nrm*{\vecU_1- \vecU_0} = \eta \nrm*{\grad F(\vecU_0)} \le \eta \rho_0(F(\vecW_0)) \le \rho_0(F(\vecW_0)+1)^{-1}.\]
By \pref{lem:boundlocalsmoothnesseasy}, for all $\vecP$ in the line segment $\overline{\vecU_0 \vecU_1}$, $\nrm*{\grad^2 F(\vecP)}_{\OPNORM} \le L_1(\vecW_0)$.
By \pref{lem:secondordersmoothnessineq}, which only depends on the smoothness constant in the segment between the two iterates (see \pref{subsec:backgroundmath}),
\[ F(\vecU_1) \le F(\vecU_0) - \eta \nrm*{\grad F(\vecU_0)}^2 + \tfrac{L_1(\vecW_0) \eta^2}2 \cdot \nrm*{\grad F(\vecU_0)}^2 < F(\vecU_0) - \tfrac{\epsilon^2}{2L_1(\vecW_0)},\]
as $\nrm*{\grad F(\vecU_0)} > \epsilon$ and by our choice of $\eta$. This proves \pref{claim:gdfirstorderdecreaseprocedure}, completing the proof. 
\end{proof}
Note it is critical here that $\vecU_0$ is in the $F(\vecW_0)$-sublevel set.
Also, to satisfy \pref{corr:gradcontrol}, $\rho_0(x)$ just needs to be a non-decreasing pointwise upper bound of $\rho_1(x) \sqrt{2\theta(x)}$. For example when $F$ is $(L_0,L_1)$-smooth, we show in \pref{prop:L0L1rho0ex} that we can take $\rho_0(x) = 2L_0^{1/2} x^{1/2} + 5L_1^2 L_0^{-1/2} x^{3/2}$. 

\subsection{Adaptive Gradient Descent}
Our proof and framework readily adapt to Adaptive GD, as discussed \pref{subsec:exinframework}. It is even easier as Adaptive GD is automatically `local' via gradient clipping. Our proof is in \pref{subsec:AdaptiveGDFOSP}.
\begin{theorem}[GD for FOSP]\label{thm:adaptivegdguarantee}
Suppose $F$ satisfies \pref{ass:selfbounding}. Run Adaptive GD initialized at $\vecW_0$, with adaptive step size $\eta_{\vecW_t} = \min\crl*{\frac1{L'_1(\vecW_0)}, \frac1{\rho_0(F(\vecW_0)+1)\nrm*{\grad F(\vecW_t)}}}$ where $L'_1(\vecW_0) = \rho_1\prn*{F(\vecW_0)+1}$. 
Let $T=\frac{2 F(\vecW_0)}{\min\crl*{\frac{L_1'(\vecW_0)}{\rho_0(F(\vecW_0)+1)^2}, \frac{\epsilon^2}{L'_1(\vecW_0)}}}$. Within $T+1$ oracle calls to $\grad F(\cdot)$, Adaptive GD will output $T$ candidate vectors $(\vecP_1, \ldots, \vecP_T)$, one of which satisfies $\nrm*{\grad F(\vecP_t)} \le \epsilon$.
\end{theorem}

\subsection{Stochastic Gradient Descent}
We make the following assumption on the stochastic gradient oracle:
\begin{assumption}\label{ass:noiseregularity}
The stochastic gradient oracle $\grad f(\cdot;\cdot)$ is unbiased (i.e. $\mathbb{E}_{\vecZeta}\brk*{\grad f(\cdot;\vecZeta)} = \grad F(\cdot)$), and for a non-decreasing function $\sigma:\bbR^+ \mapsto \bbR^+$ and all $\vecW$, $\vecZeta$, $\nrm*{\grad f(\vecW; \vecZeta) - \grad F(\vecW)}^2 \le \sigma(F(\vecW))^2$.
\end{assumption}

In many problems of interest in ML, noise scales with function value \citep{wojtowytsch2023stochastic, wojtowytsch2024stochastic}; \pref{ass:noiseregularity} captures this setting. Note we do not assume a global bound on $\nrm*{\grad F}$ or $F$, thus noise is \textit{unbounded}. We show in \pref{rem:sgdfirstordersubgaussian} that one can extend \pref{thm:sgdfirstorder} to when $\nrm*{\grad f(\vecW; \vecZeta) - \grad F(\vecW)}$ is sub-Gaussian with parameter $\sigma(F(\vecW))$ with a longer technical argument. We also note that bounding $L_2$ gradient error in terms of function value has been studied -- denoted by the expected smoothness assumption -- in \citet{gower2019sgd, gower2021stochastic}.

\begin{theorem}[SGD for FOSP]\label{thm:sgdfirstorder}
Suppose $F$ satisfies \pref{ass:selfbounding} and that the stochastic gradient oracle $\grad f(\cdot;\cdot)$ satisfies \pref{ass:noiseregularity}. For any $\delta \in (0,1)$, run SGD initialized at $\vecW_0$, for a given fixed step size $\eta \le \tilde{O}(\epsilon^2)$ depending on $\epsilon$, $\delta$, and $F(\vecW_0)$. Then with probability at least $1-\delta$, within
\[ T=\tilde{O}\prn*{\frac{1}{\epsilon^4} \cdot \polylog(1/\epsilon, 1/\delta)}\text{ oracle calls to }\grad f(\cdot;\cdot),\]
SGD will output $T$ candidate vectors $\vecW$, one of which satisfies $\nrm*{\grad F(\vecW)} \le \epsilon$.
\end{theorem}
Here $\tilde{O}(\cdot)$ hides additional $F(\vecW_0)$-dependence. Our full proof is in \pref{subsec:SGDFOSP}. As discussed in \pref{subsec:exinframework}, the idea is similar to the proof of \pref{thm:gdfirstorder}, except we now establish high-probability decrease over blocks of consecutive iterates using concentration inequalities.

\subsection{Perturbed Gradient Descent}\label{subsec:perturbedgdresults}
\begin{theorem}[Perturbed GD for SOSP]\label{thm:escapesecondordergd}
Suppose $F$ satisfies \pref{ass:thirdorderselfbounding}. For any $\delta\in(0,1)$, run Perturbed GD (\pref{alg:perturbedgd}, from \citet{jin2017escape}) initialized at $\vecW_0$, with appropriate step size $\eta$ and other parameters depending on $\epsilon, \delta, d$, and $F(\vecW_0)$. Then with probability at least $1-\delta$, within 
\[ T=O\prn*{\frac{1}{\epsilon^2} \log^4\prn*{\frac{d}{\epsilon \delta}}}\text{ oracle calls to }\grad F(\cdot),\]
Perturbed GD outputs $T$ candidates $\vecW$, one of which satisfies $\nrm*{\grad F(\vecW)} \le \epsilon, \grad^2 F(\vecW) \succeq -\sqrt{\epsilon} \matI$.
\end{theorem}
\begin{remark}\label{rem:sospdefnosmooth}
Here we find $\vecW$ with $\grad^2 F(\vecW) \succeq -\sqrt{\epsilon} \matI$, which is most sensible without Lipschitz Hessian. 
\end{remark}
For Perturbed GD here in \pref{subsec:perturbedgdresults}, asymptotic notation hides universal constants and dependence on $F(\vecW_0)$.
The full proof is in \pref{sec:perturbedgdproofs}; here we give the main ideas. 
Define $\cA, \Toracle(\vecU_0), \cS$ as in \pref{subsec:exinframework} for Perturbed GD. Consider $\Gthres = \tilde{\Theta}\prn*{\epsilon}$, $\Fthres=\tilde{\Theta}\prn*{\epsilon^{1.5}}$ defined in \pref{alg:perturbedgd}. Let
\[ \Delta(\vecU_0) = \begin{cases}\Fthres &:\nrm*{\grad F(\vecU_0)} \le \Gthres \\ \frac{\eta}2 \cdot \Gthres^2 &: \nrm*{\grad F(\vecU_0)} > \Gthres. \end{cases}\]
The central Claim is as follows, from which \pref{thm:escapesecondordergd} follows directly via \pref{thm:generalframework}:
\begin{claim}\label{claim:perturbedgddecreasesequence}
For all $\vecU_0 \in \cL_{F,F(\vecW_0)}$, $\cA$ is a $(\cS, \Toracle(\vecU_0), \Delta(\vecU_0), \frac{d L_1(\vecW_0)}{\sqrt{\epsilon}} e^{-\chi}, \vecU_0)$-decrease procedure, where $\chi = \Theta\prn*{\log\prn*{\frac{d}{\epsilon^{2.5} \delta}}}$ and $L_1(\vecW_0)$ is defined in \pref{eq:L1def}.
\end{claim}
Perturbed GD is a decrease procedure only in $\cL_{F, F(\vecW_0)}$ where we have quantitative control on $F$ and its derivatives -- \textit{using our framework is crucial}. 
To prove \pref{claim:perturbedgddecreasesequence}, we note the analysis of Perturbed GD in \citet{jin2017escape} only considers `local' points close to the current iterate the algorithm. Thus we can apply similar analysis, using \pref{lem:boundfuncvalueradius}, \pref{lem:boundlocalsmoothnesseasy}, and the similar \pref{lem:formalizingperturbationopnrom} to give enough control over the derivatives of $F$ between these `local' points close to $\vecU_0 \in \cL_{F,F(\vecW_0)}$.

\subsection{Restarted Stochastic Gradient Descent}\label{subsec:mainsgdsecondoder}
In addition to \pref{ass:noiseregularity}, we will make the following mild assumption on the error of the stochastic gradient oracle, a relaxation of Assumption 1 of \citet{fang2019sharp}. 
\begin{assumption}\label{ass:gradoraclecontrol}
For every $\vecW, \vecZeta$, $\nrm*{\grad^2 f(\vecW;\vecZeta)}_{\OPNORM} \le \rho_3\prn*{\nrm*{\grad f(\vecW; \vecZeta)}, F(\vecW)}$, where $\rho_3(\cdot, \cdot):\mathbb{R}_{\ge 0}\times \mathbb{R}_{\ge 0} \rightarrow \mathbb{R}_{\ge 0}$ is non-decreasing in both arguments.
\end{assumption}
Note if $f(\cdot;\vecZeta)$ satisfies the regularity assumptions of \citet{zhang2019gradient} or \citet{li2023convex} for every $\vecZeta$, then \pref{ass:gradoraclecontrol} is satisfied. However, \pref{ass:gradoraclecontrol} goes well beyond these assumptions, allowing for the operator norm of $\grad^2 f(\cdot;\vecZeta)$ to also diverge in $F(\vecW)$.\footnote{While the above assumes that $f(\cdot;\vecZeta)$ is twice differentiable, it can be easily phrased in terms of $\grad f(\cdot;\vecZeta)$.}
\begin{theorem}[Restarted SGD for SOSP]\label{thm:sgdsecondorder}
Suppose $F$ satisfies \pref{ass:thirdorderselfbounding} and $\grad f(\cdot;\cdot)$ satisfies \pref{ass:noiseregularity} and \pref{ass:gradoraclecontrol}. For any $\delta\in(0,1)$, run Restarted SGD (\pref{alg:restartedsgd}, the same algorithm from \citet{fang2019sharp}) initialized at $\vecW_0$, with appropriate step size $\eta$ and other parameters depending on $\epsilon$, $\delta$, $d$, and $F(\vecW_0)$. 
Then with probability at least $1-\delta$, upon making 
\[ T=\tilde{O}\prn*{\frac1{\epsilon^{3.5}}}\text{ oracle calls to }\grad f(\cdot;\cdot), \]
Restarted SGD outputs $T$ candidates $\vecW$, one of which satisfies $\nrm*{\grad F(\vecW)} \le \epsilon, \grad^2 F(\vecW)\succeq -\sqrt{\epsilon} \matI$.
\end{theorem}
Here $\tilde{O}(\cdot)$ only hides constant factors, $F(\vecW_0)$-dependent constants, and logarithmic factors in $d, 1/\epsilon, 1/\delta$.
We specify the exact parameters and detail the proof in \pref{sec:restartedsgdproofs}. The proof follows our framework instantiated for Restarted GD as in \pref{subsec:exinframework}. The crux again is establishing that the algorithm is a decrease procedure in the $F(\vecW_0)$-sublevel set, done in \pref{claim:restartedsgddecreasesequence}.

\subsection{Examples}\label{subsec:examplesmainbody}
Several interesting problems in ML and optimization, such as Phase Retrieval and Matrix PCA, can be globally optimized by finding a SOSP (but not a FOSP), and satisfy \pref{ass:thirdorderselfbounding}. See \pref{sec:examplesproofs} for these verifications. Thus \pref{thm:escapesecondordergd} and \pref{thm:sgdsecondorder} immediately imply we can solve the following problems, with no customized analysis required. 
\paragraph{Phase Retrieval:} We reconstruct a hidden vector $\vecW^* \in \mathbb{R}^d$ with $\nrm*{\vecW^*} = 1$ using phaseless observations $\cS = \{(\vecA_j, y_j)\}$ where $y_j = \langle \vecA_j, \vecW^* \rangle^2$, $\vecA_j \sim \mathcal{N}(\vecOrigin, \matI_d)$. The population loss is $F_{\text{pr}}(\vecW) = \mathbb{E}_{\vecA\sim \mathcal{N}(\vecOrigin, \matI_d)}\brk*{\prn*{\tri*{\vecA, \vecW}^2 - \tri*{\vecA,\vecW^*}^2}^2}$.

\paragraph{Matrix PCA:} Given a $d \times d$ symmetric positive definite (PD) matrix $\matM$, we aim to find $\vecW\in\mathbb{R}^d$ (the first principal component) minimizing $F_{\text{pca}}(\vecW) = \frac12 \nrm*{\vecW \vecW^{\top} - \matM}_F^2$.

\subsection{Practical Implications and Simulations} 
Our results show under generalizations of smoothness, unlike with Lipschitz gradient/Hessian, the larger the loss is at initialization (larger $F(\vecW_0)$) and larger self-bounding functions $\rho_1(\cdot)$ shrink the `window' for choosing a working $\eta$. Specifically, with larger loss at initialization, the smaller the largest working step size is, in contrast to optimizing smooth functions. \textit{This implies in practice, for losses with non-Lipschitz gradient/Hessian, one should tune $\eta$ based on suboptimality at initialization.} 

In \pref{sec:experiments}, we validate this finding through simulations with GD and SGD on several natural smooth and generalized smooth functions, namely $F(\vecW) = \nrm*{\matA \vecW}^p$ for $p=2,3,4,5,6$. Our simulations show the above theoretical conclusions match behavior in practice, validating the practical implications of our theoretical results on which step sizes successfully optimize generalized smooth functions.

\section{Conclusion}\label{sec:conclusion}
We present a systematic framework to analyze the convergence of first order methods to FOSPs and SOSPs under generalizations of smoothness, extending key results in finding SOSPs via first-order methods to this setting. Our work \textit{elucidates fundamental behavior of first-order optimization algorithms}, showing that `chaining together high-probability decrease' enables their success under generalizations of smoothness. Our framework applies for many other algorithms (e.g. Langevin Dynamics) and sets of interest $\cS$ (e.g. higher order stationary points, or minima with good generalization properties). 
It can also inform the design of new optimization algorithms, by designing procedures which are decrease procedures. 
These promising directions are left for future research.

\section{Acknowledgments}\label{sec:acknowledgements}
We thank Dylan J. Foster and Ayush Sekhari for discussions, and Anthony Bao, Fan Chen, and Albert Gong for useful suggestions on the presentation of our manuscript.

\newpage
\bibliography{sources.bib}



\appendix
\onecolumn

\newpage
\paragraph{Additional Notation:} For a matrix $\matM$, $\lambda_{\min}(\matM)$ denotes its minimum eigenvalue, and $\lambda_r(\matM)$ denotes its $r$-th largest eigenvalue. Thus $\lambda_1(\matM) \ge \lambda_2(\matM) \ge \ldots$. We denote the $k \times k$ identity matrix by $\matI_k$. We use $\ball^k(\vecP,R)$ to denote the full $k$-dimensional $l_2$-ball centered at $\vecP\in\mathbb{R}^k$ with radius $R$, including the boundary. When $k$ is not specified explicitly, $\ball(\vecP, R)$ refers to the $l_2$-ball in $\mathbb{R}^d$, following Notation. All logarithms in the following are the natural logarithm. For an event $\cS$, $1_{\cS}$ denotes the indicator function. In the following, the norm $\nrm*{\cdot}$ of matrices and higher-order tensors refers to the operator norm unless otherwise stated. The norm $\nrm*{\cdot}$ of vectors refers to $l_2$-Euclidean norm.

\tableofcontents

\section{Technical Preliminaries}\label{sec:preliminaries}
\subsection{Helpful Background Lemmas}\label{subsec:backgroundmath}
We will use the following classical inequalities from optimization to show we still have some notion of control if we have local bounds on the relevant derivatives.
\begin{lemma}\label{lem:secondordersmoothnessineq}
Suppose $F$ is twice differentiable, and for all $\vecU \in \overline{\vecX \vecY}$ (the line segment) we have $\nrm*{\grad^2 F(\vecU)}_{\OPNORM} \leq L$. Then, we have
\[
F(\vecY) \leq F(\vecX) + \tri*{\grad F(\vecX), \vecY-\vecX} + \frac{L}{2} \nrm*{
\vecY-\vecX}^2.
\]
\end{lemma}
\begin{proof}
This follows by the proof of Lemma 3.4 in \citet{bubeck2015convex}. In particular, one can readily verify that $\vecX+t(\vecY-\vecX) \in \overline{\vecX \vecY}$ for all $t \in [0,1]$. Hence for all $t \in [0,1]$ and $\vecU$ in the line segment between $\vecX$ and $\vecX+t(\vecY-\vecX)$, $\nrm*{\grad^2 F(\vecU)}_{\OPNORM} \le L$. Thus,
\begin{align*}
\abs*{F(\vecY)-F(\vecX)-\tri*{\grad F(\vecX),\vecY-\vecX}} &= \abs*{\int_0^1 \tri*{\grad F(\vecX+t(\vecY-\vecX)), \vecY-\vecX}\DERIV t - \tri*{\grad F(\vecX), \vecY-\vecX}}  \\
&= \abs*{\int_0^1 \tri*{\grad F(\vecX+t(\vecY-\vecX))-\grad F(\vecX), \vecY-\vecX}\DERIV t }  \\
&\le \abs*{\int_0^1 L t\nrm*{\vecY-\vecX}^2 \DERIV t} = \frac{L}2 \nrm*{\vecY-\vecX}^2.
\end{align*}
This gives the desired result.
\end{proof}
Analogously, one can show the following by considering the local second-order approximation around $\vecX$.
\begin{lemma}\label{lem:thirdordersmoothnessineq}
Suppose $F$ is twice differentiable, and for all $\vecU \in \overline{\vecX \vecY}$ (again the line segment), we have 
\[ \nrm*{\grad^2 F(\vecU) - \grad^2 F(\vecX)}_{\OPNORM} \leq L \nrm*{\vecU-\vecX}. \]
Then,
\[
F(\vecY) \leq F(\vecX) + \tri*{\grad F(\vecX), \vecY - \vecX} + \frac12 (\vecY - \vecX)^\top \grad^2 F(\vecX) (\vecY - \vecX) + \frac{L}{6} \nrm*{\vecY - \vecX}^3.
\]
\end{lemma}
\begin{proof}
Similarly to the proof of \pref{lem:secondordersmoothnessineq}, we show this via the proof of Lemma 1 in \citet{nesterov2006cubic}.
Analogously as in the proof of \pref{lem:secondordersmoothnessineq}, one can readily verify that for any $\vecY' \in \overline{\vecX \vecY}$, $\vecX+t(\vecY'-\vecX) \in \overline{\vecX \vecY}$ holds for all $t \in [0,1]$. Hence for all $t \in [0,1]$, applying the condition of this Lemma,
\[ \nrm*{\grad^2 F(\vecX+t(\vecY'-\vecX)) - \grad^2 F(\vecX)}_{\OPNORM} \le L t\nrm*{\vecY'-\vecX}.\]
Thus for any $\vecY' \in \overline{\vecX \vecY}$, by Cauchy-Schwartz and the above, we obtain
\begin{align*}
\nrm*{\grad F(\vecY')-\grad F(\vecX)-\tri*{\grad^2 F(\vecX),\vecY'-\vecX}} &= \nrm*{\int_0^1 \tri*{\grad^2 F(\vecX+t(\vecY'-\vecX)), \vecY'-\vecX}\DERIV t - \tri*{\grad^2 F(\vecX), \vecY'-\vecX}}  \\
&= \nrm*{\int_0^1 \tri*{\grad^2 F(\vecX+t(\vecY'-\vecX))-\grad^2 F(\vecX), \vecY'-\vecX}\DERIV t }  \\
&\le \abs*{\int_0^1 L t\nrm*{\vecY'-\vecX}^2 \DERIV t} = \frac{L}2 \nrm*{\vecY'-\vecX}^2.
\end{align*}
Applying the above relation for $\vecY'=\vecX+t(\vecY-\vecX)$ which is in $\overline{\vecX \vecY}$ for all $t \in [0,1]$, we obtain
\begin{align*}
&\abs*{F(\vecY) - F(\vecX) - \tri*{\grad F(\vecX), \vecY-\vecX} - \frac12 \tri*{\grad^2 F(\vecX)(\vecY-\vecX), \vecY-\vecX}} \\
&= \abs*{\int_0^1 \tri*{\grad F(\vecX+t(\vecY-\vecX))- \grad F(\vecX)- t \grad^2 F(\vecX)(\vecY-\vecX), \vecY-\vecX} \DERIV t} \\
&= \abs*{\int_0^1 \tri*{\grad F(\vecX+t(\vecY-\vecX))- \grad F(\vecX)- \grad^2 F(\vecX) \cdot t(\vecY-\vecX), \vecY-\vecX} \DERIV t} \\
&\le \int_0^1 \nrm*{\vecY-\vecX} \cdot \frac{L}2 \nrm*{t(\vecY-\vecX)}^2 \DERIV t = \frac{L}6 \nrm*{\vecY-\vecX}^3.
\end{align*}
This gives the desired result.
\end{proof}

We will also use the following Lemmas.
\begin{lemma}\label{lem:upperboundouterproductopnorm}
For vectors $\vecA$, $\vecb$, the matrix operator norm $\nrm*{\vecA\vecb^{\top}}_{\OPNORM} \le \nrm*{\vecA}\nrm*{\vecb}$.
\end{lemma}
\begin{proof}
Consider any unit vector $\vecX$. By Cauchy-Schwartz and associativity, we have
\[ \vecX^{\top} \prn*{\vecA\vecb^{\top} }\vecX \le \tri*{\vecX, \vecA} \tri*{\vecX, \vecb} \le \nrm*{\vecX}^2 \nrm*{\vecA}\nrm*{\vecb}=\nrm*{\vecA}\nrm*{\vecb}.\]
The conclusion follows by definition of operator norm.
\end{proof}
\begin{lemma}\label{lem:comparefuncstechnical}
Consider any non-negative, continuous function $g(x)$ such that $\lim_{x \rightarrow \infty} g(x) = \infty$ and such that $g(x)>0$ on $[1,\infty)$. Then on $[1,\infty)$, $g(x)$ can be lower bounded by a strictly positive, infinitely differentiable, strictly increasing function $\tilde{g}(x)$, where $\tilde{g}$ has domain $[1,\infty)$.
\end{lemma}
\begin{proof}
We will explicitly construct such a $\tilde{g}$ in terms of $g$. First, since $\lim_{x \rightarrow \infty} g(x) = \infty$, for all $i\ge 1$, there exists $t_i \in [1,\infty)$ such that $g(x) \ge i+1$ for all $x \ge t_i$. We furthermore can clearly assume $2 \le t_1 < t_2 < \cdots$, by increasing each $t_N$ if necessary. Also let $t_0 = 1$. Thus $\bigcup_{i \ge 0} [t_i, t_{i+1})$ forms a disjoint union of $[1,\infty)$.

Now, let $c=\min\prn*{1, \inf_{x \in [1,t_1]} g(x)}>0$; the strict inequality here holds as $t_1<\infty$ and as $g$ is continuous. Define a sequence $\{b_i\}_{i \ge 0}$ by $b_0 = c/2, b_1=c$, and $b_i = i$ for all $i \ge 2$. Thus $b_0 < b_1 < \cdots$. Furthermore, this construction of $\{b_i\}_{i \ge 0}$ implies for all $i\ge 0$, we have $g(x) \ge b_{i+1}$ for all $x \in [t_i, t_{i+1}]$.

Now construct $\tilde{g}(x)$ as follows. For all $i \ge 0$, we let $\tilde{g}(x)$ equal a function $h_i(x)$ defined on $[t_i, t_{i+1}]$ such that $h_i(t_i)=b_i$, $h_i(t_{i+1})=b_{i+1}$, where we define $h_i$ as follows. We first define $h:[0,1]\rightarrow [0,1]$ such that $h$ is infinitely differentiable, $h(0)=0$, $h(1)=1$, $h^{(n)}(0)=h^{(n)}(1)=0$ for all $n \ge 1$ where $h^{(n)}$ denotes the $n$-th derivative, and $h'(x) > 0$ for all $x \in (0,1)$. To this end we use a construction from \citet{chen2024optimization}: let 
\[ h(x) = \frac{e^{-\frac1{x^2}}}{e^{-\frac1{x^2}}+e^{-\frac1{1-x^2}}}\text{ on }(0,1),\]
and extend $h$ to $[0,1]$ by $h(0)=0, h(1)=1$. We justify these claims about $h$ shortly below. Now we let 
\[ h_i(x) = (b_{i+1}-b_i)\cdot h\prn*{\frac{x-t_i}{t_{i+1}-t_i}}+b_i \text{ for all }i \ge 0.\]
We now check $h$ satisfies the claimed properties. 
\begin{itemize}
    \item In \citet{chen2024optimization}, it is argued that $h$ maps to $[0,1]$, $h(0)=0$, $h(1)=1$, and that $h$ is infinitely differentiable. It is also argued in \citet{chen2024optimization}, Lemma 11.5, that $h'(x)$ (which is called $\tilde{p}(x)$ there) is non-negative on $[0,1]$. 

\item Next, we check $h^{(n)}(0)=h^{(n)}(1)=0$ for all $n \ge 1$. Via a straightforward induction outlined in \citet{chen2024optimization}, one can check that $\prn*{e^{-\frac1{x^2}}}^{(n)}=0$, $\prn*{e^{-\frac1{1-x^2}}}^{(n)}=0$ for all $n \ge 1$ (following the standard convention in analysis that $0 \cdot \infty=0$, see e.g. \citet{folland1999real}). Now let $f(x) = e^{-\frac1{x^2}}$, $g(x) = e^{-\frac1{x^2}}+e^{-\frac1{1-x^2}}$, thus $h=f/g$. Consequently $f^{(n)}(0)=0$, $f^{(n)}(1)=0$, $g^{(n)}(0)=0$, $g^{(n)}(1)=0$ for all $n \ge 1$. As $g>0$ always holds in $[0,1]$ as shown in \citet{chen2024optimization} and can be easily checked, we have $f=gh$. A straightforward induction gives $f^{(n)} = \sum_{k=0}^n \binom{n}k g^{(k)} h^{(n-k)}$ where $\binom{n}k$ is the binomial coefficient. We thus obtain $gh^{(n)} = f^{(n)} - \sum_{k=0}^{n-1} \binom{n}k g^{(k)} h^{(n-k)}$. For any $n \ge 1$, taking $x=0, 1$ in this expression for $h(x)$ and noting at least one of $k, n-k \ge 1$ for $0 \le k \le n-1$ implies $g(0)h^{(n)}(0)=g(1)h^{(n)}(1)=0$. Recalling $g(x)>0$ on $[0,1]$ proves $h^{(n)}(0)=h^{(n)}(1)=0$ for $n \ge 1$, as requested.

\item Finally, we check that $h'(x)>0$ for all $x \in (0,1)$. Consider any $x \in (0,1)$. By a calculation in Lemma 11.5, \citet{chen2024optimization}, we have $h'(x) > 0$ if and only if $q(x) = \frac2{x^3}\prn*{e^{-\frac1{x^2}}+e^{-\frac1{1-x^2}}} + e^{\frac{-1}{x^2}} \cdot \frac{2}{x^3} + e^{-\frac1{1-x^2}} \cdot \frac{-2x}{(1-x^2)^2} > 0$. If $x \in [\frac{\sqrt{2}}2, 1)$, directly following the proof of Lemma 11.5 in \citet{chen2024optimization} establishes that $q(x)>0$. Otherwise if $x \in (0,\frac{\sqrt{2}}2)$, note the strict inequality $\frac1{x^3} > \frac{x}{(1-x^2)^2}$, which in turn implies $q(x)>0$.
\end{itemize}
By the above properties of $h$, it follows from the Chain Rule that for all $i \ge 0$, $h_i$ satisfies the following properties:
\begin{itemize}
    \item $h_i(t_i) = b_i$, $h_i(t_{i+1})=b_{i+1}$, and $h_i(x) \in [b_i, b_{i+1}]$ for all $x \in [t_i, t_{i+1}]$.
    \item $h_i$ is infinitely differentiable.
    \item $h_i'(x) > 0$ for $x \in (t_i, t_{i+1})$, and for all $x \in [t_i, t_{i+1}], h_i'(x) \ge 0$.
    \item For all $n \ge 1$, $h_i^{(n)}(t_i)=h_i^{(n)}(t_{i+1})=0$, where again $h_i^{(n)}$ denotes the $n$-th derivative.
\end{itemize}
Finally, we check that $\tilde{g}$ has the desired properties:
\begin{itemize}
\item $\tilde{g}$ is well-defined: This follows because for all $i \ge 1$, we have $h_i(t_i) = h_{i-1}(t_i) = b_i$.
\item $\tilde{g}$ is strictly positive: This follows because $h_i(x) \in [b_i, b_{i+1}] \subseteq (0,\infty)$ for all $x \in [t_i, t_{i+1}]$.
\item $\tilde{g}$ is continuous, and moreover is infinitely differentiable: Continuity of $\tilde{g}$ follows because each $h_i$ is infinitely differentiable, and hence continuous, combined with the fact that for all $i \ge 1$, we have $h_i(t_i) = h_{i-1}(t_i) = b_i$. Infinite differentiability of $\tilde{g}$ follows because each $h_i$ is infinitely differentiable, and because for all $n \ge 1$ and all $i \ge 0$, $h_i^{(n)}(t_i)=h_i^{(n)}(t_{i+1})=0$.
\item $\tilde{g}(x) \le g(x)$ always holds for $x \in [1,\infty)$: Recall for all $i \ge 0$, we have $g(x) \ge b_{i+1}$ for all $x\in [t_i, t_{i+1}]$. Since we have $\tilde{g}(x)=h_i(x) \le b_{i+1}$ for all $x\in [t_i, t_{i+1}]$, it follows that for all $x\in [t_i, t_{i+1}]$, $\tilde{g}(x) \le g(x)$. The result follows upon recalling that $\bigcup_{i \ge 0} [t_i, t_{i+1})$ forms a disjoint union of $[1,\infty)$.
\item $\tilde{g}$ is strictly increasing: Consider any $x_1 < x_2, x_1, x_2 \in [1,\infty)$. Since $x_1<x_2$, and recalling that $\bigcup_{i \ge 0} [t_i, t_{i+1})$ forms a disjoint union of $[1,\infty)$, it follows that for some $j \ge 0$, $(x_1, x_2) \cap (t_j, t_{j+1}) \neq \emptyset$. This intersection is open, and therefore contains some open interval $(a,b) \subseteq (t_j, t_{j+1})$. Let $c' = \inf_{x \in [\frac{2a+b}3, \frac{a+2b}3]} h_j'(x) > 0$, where the strict inequality follows as $[\frac{2a+b}3, \frac{a+2b}3]  \subseteq (t_j, t_{j+1})$, and by continuity of $h'_j$ on the compact $[\frac{2a+b}3, \frac{a+2b}3]$. Since we have $h'_i(x) \ge 0$ for all $x \in [t_i, t_{i+1}]$ for all $i \ge 0$, we obtain
\[ \tilde{g}(x_2) \ge 0+c' \cdot \frac{b-a}3 + \tilde{g}(x_1) > \tilde{g}(x_1).\]
This proves that $\tilde{g}$ is strictly increasing as claimed.
\end{itemize}
Thus, we have constructed a function $\tilde{g}$ that satisfies the requested properties.
\end{proof}

\subsection{Comparison of Assumptions with Literature}\label{subsec:detailsforcomparewithlit}
Here, we establish that our regularity conditions are more general than those of literature.
\begin{proposition}\label{prop:comparetolit}
If $\nrm*{\grad^2 F(\vecW)} \le l(\grad F(\vecW))$ for non-decreasing, differentiable sub-quadratic $l$ (where sub-quadratic means that $\lim_{x\rightarrow\infty} \frac{l(x)}{x^2}=0$), then our \pref{ass:selfbounding} is satisfied for some non-decreasing $\rho_1(x)$. In this generality, $\rho_1(x)$ depends on $l(x)$, and can be found explicitly from the construction from \pref{lem:comparefuncstechnical}. 

Furthermore, suppose $F$ is $(L_0,L_1)$-smooth, that $\nrm*{\grad^2 F(\vecW)} \le L_0+L_1 \nrm*{\grad F(\vecW)}$ for $L_0,L_1 \ge 0$. Then \pref{ass:selfbounding} is satisfied with $\rho_1(x)=\tfrac32L_{0}+4L_{1}^{2}x$.
\end{proposition}
\begin{proof}
Essentially this follows from Lemma 3.5, \citet{li2023convex}, where it is shown that these assumptions of \citet{zhang2019gradient}, \citet{li2023convex} imply an upper bound on $\nrm*{\grad F(\vecW)}$ in terms of an increasing function of $F(\vecW)$; combining with the assumptions of \citet{zhang2019gradient, li2023convex} implies that $\nrm*{\grad^2 F(\vecW)}$ is upper bounded in terms of an increasing function of $F(\vecW)$. 

\paragraph{Proof for general $l$:} Consider any $\vecW\in\mathbb{R}^d$. By Lemma 3.5 of \citet{li2023convex}, 
\begin{align*}
\nrm*{\grad F(\vecW)}^2 &\leq 2 \ell \prn*{2 \nrm*{\grad F(\vecW)}} \cdot F(\vecW).
\end{align*}
This implies 
\[\frac{4\nrm*{\grad F(\vecW)}^2}{\ell \prn*{2 \nrm*{\grad F(\vecW)}}} \le 8F(\vecW).\]
Let $2\nrm*{\grad F(\vecW)} = t$. Consider when $t \ge 2$. Then the left hand side equals $\frac{t^2}{l(t)}$. Note that WLOG, we can add 1 to $l(\cdot)$ so that $l(t) \ge 1$ for $t \ge 1$. Thus $\frac{t^2}{l(t)}$ is continuous on $[1,\infty)$, and furthermore is positive on this interval. Now note $\lim_{x\rightarrow\infty}\frac{x^2}{l(x)}=\infty$ by the condition (including after adding 1 WLOG), and thus by \pref{lem:comparefuncstechnical}, $\frac{x^2}{l(x)}$ is lower bounded by some strictly increasing function $\tilde{g}(x)$ on $[2,\infty)$. Therefore, $\tilde{g}$ is invertible and so we have 
\[ \tilde{g}(2\nrm*{\grad F(\vecW)}) \le \frac{4\nrm*{\grad F(\vecW)}^2}{\ell \prn*{2 \nrm*{\grad F(\vecW)}}} \le 8F(\vecW) \implies \nrm*{\grad F(\vecW)} \le \frac12 \tilde{g}^{-1}(8F(\vecW)).\]
Then by the assumptions of \citet{li2023convergence}, it holds that
\[ \nrm*{\grad^2 F(\vecW)} \le l\prn*{\frac12 \tilde{g}^{-1}(8F(\vecW))}.\]
Else when $t < 2$, we have $\nrm*{\grad F(\vecW)} \le 1$, and by the assumptions of \citet{li2023convergence}, we have $\nrm*{\grad^2 F(\vecW)} \le l(1)$.

Thus the assumptions of \citet{li2023convergence} imply that the following always holds:
\[ \nrm*{\grad^2 F(\vecW)} \le l\prn*{\frac12 \tilde{g}^{-1}(8F(\vecW))}+l(1).\]
We thus can take $\rho_1(x) = l\prn*{\frac12 \tilde{g}^{-1}(8x)}+l(1)$, which is clearly non-negative. It remains to check that $l\prn*{\frac12 \tilde{g}^{-1}(8x)}$ is non-decreasing. As $l$ is non-decreasing, as compositions of non-decreasing functions are non-decreasing, it remains to check that $\frac12 \tilde{g}^{-1}(8x)$ is non-decreasing. Since $\tilde{g}$ is non-decreasing, $\tilde{g}^{-1}$ is non-decreasing as well, and this completes the proof.

\paragraph{Proof for $(L_0,L_1)$-smoothness:} First, when $L_1 = 0$ the result is immediate, so from here on out suppose $L_1 > 0$. By Lemma~3.5 from \cite{li2023convex} we have for all $\vecW\in\mathbb{R}^d$,
\begin{align*}
\nrm*{\grad F(\vecW)}^2 &\leq 2 \ell \prn*{2 \nrm*{\grad F(\vecW)}} \cdot F(\vecW),
\end{align*}
where $\ell(x) = L_0 + L_1(x)$ for $L_0, L_1 \geq 0$. We thus obtain:
\begin{align*}
    \nrm*{\grad F(\vecW)}^2 
    &\leq 2 \prn*{L_0 + 2L_1 \nrm*{\grad F(\vecW)}} \cdot F(\vecW) \\
    &= 2L_0 F(\vecW) + 4L_1 \nrm*{\grad F(\vecW)} F(\vecW). 
\end{align*}
Rewriting this inequality, we get
\begin{align*}
    \nrm*{\grad F(\vecW)}^2 
    - 4L_1 \nrm*{\grad F(\vecW)} F(\vecW) 
    - 2L_0 F(\vecW) 
    &\leq 0.
\end{align*}
Consider the quadratic $x^2-4L_1 F(\vecW) \cdot x - 2L_0 F(\vecW)$. The coefficient on the quadratic term is positive, and the quadratic is non-negative when $x=\nrm*{\grad F(\vecW)}$. Thus $\nrm*{\grad F(\vecW)}$ must be no larger than the largest root of $x^2-4L_1 F(\vecW) \cdot x - 2L_0 F(\vecW)$, and we obtain
\begin{align*}
\nrm*{\grad F(\vecW)} 
&\leq \frac{1}{2} \prn*{ 4L_1 F(\vecW) + \sqrt{16L_1^2 F(\vecW)^2 + 8L_0 F(\vecW)} } \\
&\leq 2L_{1}F(\vecW)+\sqrt{(2L_{1}F(\vecW))^{2}+2L_{0}F(\vecW)} \numberthis\label{eq:grad-bound}
\end{align*}
If $F(\vecW)=0$, the above immediately implies $\nrm*{\grad F(\vecW)} = 0$. Otherwise, recall by shifting (in Notation) that $F(\vecW) \ge 0$ always holds, so suppose $F(\vecW)>0$. Recall also from earlier that it suffices to show the result for $L_1>0$. Applying the inequality
\(\sqrt{a^{2}+b}\le a+\tfrac{b}{2a}\), valid for all $a > 0, b \ge 0$ with
\(a=2L_{1}F(\vecW)>0,\;b=2L_{0}F(\vecW) \ge 0\), we obtain
\begin{align*}
\sqrt{(2L_{1}F(\vecW))^{2}+2L_{0}F(\vecW)}\;\le\;2L_{1}F(\vecW) + \frac{L_{0}}{2L_{1}}.
\end{align*}
Substituting into \pref{eq:grad-bound} gives that for all $\vecW$ with $F(\vecW)>0$, we have
\begin{align}
\nrm*{\grad F(\vecW)} \leq \frac{L_{0}}{2L_{1}}+4L_{1}F(\vecW).\numberthis\label{eq:g_bound} 
\end{align}
By the argument earlier, if $F(\vecW)=0$, the above bound \pref{eq:g_bound} holds too. Thus \pref{eq:g_bound} holds for all $\vecW\in\mathbb{R}^d$. Now inserting \pref{eq:g_bound} into the definition of $(L_0,L_1)$-smoothness gives
\begin{align*}
  \nrm*{\grad^{2}F(\vecW)}
  \;\le\;
  L_{0}+L_{1}\prn*{\tfrac{L_{0}}{2L_{1}}+4L_{1}F(\vecW)}
   = \tfrac32L_{0}+4L_{1}^{2}F(\vecW).
\end{align*}
Hence \pref{ass:selfbounding} is satisfied with the increasing function $\rho_{1}(x)=\frac32L_{0}+4L_{1}^{2}x$.
\end{proof}

\begin{proposition}\label{prop:L0L1rho0ex}
When $F$ is $(L_0,L_1)$-smooth, letting $\rho_0(x) =2L_0^{1/2} x^{1/2} + \frac{5L_1^2}{L_0^{1/2}} x^{3/2}$, we have $\nrm*{\grad F(\vecW)} \le \rho_0\prn*{F(\vecW)}$.
\end{proposition}
\begin{proof}
By \pref{prop:comparetolit}, we can take $\rho_{1}(x)=\frac32L_{0}+4L_{1}^{2}x$ in this case. As noted in \pref{subsec:fospgd}, we need to show that $2L_0^{1/2} x^{1/2} + \frac{5L_1^2}{L_0^{1/2}} x^{3/2}$ is a pointwise upper bound on 
\[ \rho_1(x) \sqrt{2\theta(x)} \text{ where }\theta(x) = \int_0^x \frac1{\rho_1(v)} \DERIV v.\]
To this end note for each $x \ge 0$ that $\theta(x) \le x \cdot \frac{1}{\frac32 L_0}=\frac{2}{3L_0} x$, thus for each $x \ge 0$,
\[ \rho_1(x) \sqrt{2\theta(x)} \le \prn*{\frac32L_{0}+4L_{1}^{2}x} \sqrt{\frac{4}{3L_0} x} \le 2L_0^{1/2} x^{1/2} + \frac{5L_1^2}{L_0^{1/2}} x^{3/2}.  \]
This completes the proof.
\end{proof}

\begin{example}\label{ex:breaklit}
We now provide a natural example of a univariate function that satisfies our regularity assumptions but does not necessarily satisfy those of \citet{li2023convergence} for non-convex optimization. Namely, consider the univariate function:
\[
F(x) = 1 - \log(\cos(1 + x)), 0 \le x < \frac{\pi}2-1.
\]
The argument here is in radians. The first derivative is:
\[
F'(x) = \tan(1 + x).
\]
The second derivative is:
\[
F''(x) = \sec^2(1 + x).
\]
Thus as $\tan^2(\theta)+1=\sec^2(\theta)$, $F$ satisfies the ODE: 
\[ F''(x) = F'(x)^2 + 1.\numberthis\label{eq:odecounterex}\]
Suppose that $F$ satisfied the conditions of \citet{li2023convergence} for non-convex optimization on the relevant domain, thus for all $0 \le x < \frac{\pi}2-1$, we would have
\[ F''(x) \leq \ell(F'(x)),\] 
for some sub-quadratic $l(\cdot)$.

Then by \pref{eq:odecounterex} and noting $F'(x) > 0$ on the domain, we obtain for all $0 \le x < \frac{\pi}2-1$
\[ 1 \le 1+\frac1{F'(x)^2} = \frac{F'(x)^2+1}{F'(x)^2} = \frac{F''(x)}{F'(x)^2} \leq \frac{ \ell(F'(x))}{F'(x)^2}. \]
As $l$ is subquadratic, there exists $x' < \infty$ such that $l(x)/x^2 < 1$ for all $x > x'$. Noting $F'(x) \rightarrow \infty$ for $x \rightarrow \frac{\pi}2-1$ yields a contradiction. 

Consequently $F$ does not satisfy the conditions of \citet{li2023convergence} for non-convex optimization. However, we show that $F$ satisfies \pref{ass:selfbounding}.
Rewriting \( F''(x) \) in terms of \( F(x) \), note that:
\[
\cos(1 + x) = e^{1 - F(x)},
\]
and thus:
\[
F''(x) = \sec^2(1 + x) = \frac{1}{\cos^2(1 + x)} = e^{2(F(x) - 1)}.
\]
Hence we can define the increasing, non-negative function
\[
\rho_1(t) = e^{2(t - 1)},
\]
which satisfies:
\[
F''(x) \leq \rho_1(F(x)).
\]
Thus $F$ satisfies \pref{ass:selfbounding} (in the relevant domain).
\end{example}
We now discuss \pref{ass:thirdorderselfbounding}.
\begin{example}\label{ex:thirdordersmoothnessjustification}
First, we show that \pref{ass:thirdorderselfbounding} captures several univariate functions of interest. Notice also if $F(\vecW)$ is a sum of functions satisfying \pref{ass:thirdorderselfbounding}, Triangle Inequality implies that $F(\vecW)$ also satisfies \pref{ass:thirdorderselfbounding}. 
\begin{itemize}
    \item Polynomials: Consider whenever $F(x)$ is a linear combination of monomials $x^p$ for $p \ge 1$, combined with a constant term. We claim $F(x)$ satisfies \pref{ass:thirdorderselfbounding}. By linearity of derivative and Triangle Inequality, it suffices to prove this whenever $F(x)=x^p$ for $p \ge 1$ as the constant term vanishes, and then add up all the non-decreasing, non-negative functions on the right hand side to form $\rho_1$ and $\rho_2$. To this end note $F''(x) = p(p-1)x^{p-2}$, thus
    \[ |F''(x)| = p(p-1)x^{p-2} \le p(p-1)(x^p + 1) = p(p-1)(F(x)+1).\]
    Similarly, $F'''(x)=p(p-1)(p-2)x^{p-3}$, thus
    \[ |F'''(x)|=p(p-1)(p-2)x^{p-3} \le p(p-1)(p-2)(1+F(x)).\]
    Noting $p(p-1)(1+t)$ and $p(p-1)(p-2)(1+t)$ are non-decreasing and non-negative for $t \ge 0$, combined with our earlier remarks that it suffices to prove this result when $F(x)=x^p$, completes the proof. 
    \item Single-exponential functions: Consider when $F(x)=a^{x} = e^{x\ln a}$ for $a>1$. Then $F''(x)=(\ln a)^2 e^{x \ln a}$, $F'''(x)=(\ln a)^3 e^{x \ln a}$, and so we can take $\rho_1(t)=(\ln a)^2 t, \rho_2(t) = (\ln a)^3 t$.
    \item Doubly-exponential functions: Consider when $F(x)=a^{b^x} = e^{\ln a e^{x\ln b}}$ for $a,b>1$. Thus 
    \[ F'(x) = e^{\ln a e^{x\ln b}} \cdot \ln a e^{x \ln b} \cdot \ln b = \ln a \ln b F(x) e^{x \ln b}.\]
    It follows that 
    \[ F''(x) = \ln a \ln b \prn*{F'(x) e^{x \ln b} + \ln b F(x) e^{x \ln b}} = (\ln a) (\ln b)^2 F(x) (e^{2x \ln b} \ln a + e^{x \ln b}).\]
    This then implies
    \begin{align*}
    F'''(x) &=  (\ln a) (\ln b)^2 F(x) (e^{2x \ln b} 2\ln a \ln b + e^{x \ln b} \ln b) \\
    &\hspace{1in}+ (\ln a) (\ln b)^2 (e^{2x \ln b} \ln a + e^{x \ln b}) \ln a \ln b F(x) e^{x \ln b} \\
    &= (\ln a) (\ln b)^3 F(x) \prn*{2e^{2 x \ln b} \ln a + e^{x \ln b} + e^{3x \ln b} (\ln a)^2 + e^{2x \ln b} \ln a}.
    \end{align*}
    Notice
    \[ e^{x \ln b} \le e^{\ln a e^{x \ln b}} - 1 < F(x),\]
    therefore we have 
    \begin{align*}
    F''(x) &\le (\ln a)(\ln b)^2 F(x) \prn*{F(x)^2 \ln a + F(x)},\\
    F'''(x) &\le (\ln a) (\ln b)^3 F(x) \prn*{F(x)^3 (\ln a)^2 + 3F(x)^2 \ln a + F(x)}.
    \end{align*}
    We thus can take
    \begin{align*}
    \rho_1(t) = (\ln a)(\ln b)^2 t\prn*{t^2 \ln a + t},\\
    \rho_2(t) &= (\ln a) (\ln b)^3 t(t^3 (\ln a)^2 +3t^2 \ln a + t),
    \end{align*}
    which are clearly non-negative and non-decreasing on $[0,\infty)$.
    \item Next we highlight the natural example of any self-concordant function $F:\mathbb{R}\rightarrow\mathbb{R}$. Thus 
    \[ |F'''(x)| \le 2F''(x)^{3/2} \le 2|F''(x)|^{3/2}.\]
    Suppose $F$ satisfies \pref{ass:selfbounding}. Then there exists a non-negative, non-decreasing $\rho_1$ such that $|F''(x)| \le \rho_1(F(x))$. Thus,
    \[ |F'''(x)| \le 2\rho_1(F(x))^{3/2}.\]
    Since $\rho_1$ is non-negative and non-decreasing, $\rho_2(t) := 2\rho_1(t)^{3/2}$ is as well, and thus \pref{ass:thirdorderselfbounding} is satisfied.
\end{itemize}
Next, we show that the regularity assumptions Assumptions 1 and 3 of \citet{xie2024trust}, which they need for their guarantees finding SOSPs, are less general than \pref{ass:thirdorderselfbounding} when $F$ is twice-differentiable. To do so we show they imply \pref{ass:thirdorderselfbounding}, and are hence subsumed by \pref{ass:thirdorderselfbounding}. 

When $F$ is twice-differentiable, their Assumption 1 implies $(L_0, L_1)$-smoothness. As shown in \pref{prop:L0L1rho0ex}, this means that
\[ \nrm*{\grad F(\vecW)} \le \rho_0\prn*{F(\vecW)}\text{ where }\rho_0(x) = 2L_0^{1/2} x^{1/2} + \frac{5L_1^2}{L_0^{1/2}} x^{3/2}.\]
Their Assumption 3 implies for $M_0, M_1 \ge 0$ and some $\delta>0$ that for all $\vecW,\vecW'$ with $\nrm*{\vecW-\vecW'} \le \delta$,
\[ \nrm*{\grad^2 F(\vecW) - \grad^2 F(\vecW')}_{\OPNORM} \le \nrm*{\vecW-\vecW'}\prn*{M_0+M_1\nrm*{\grad F(\vecW)}} . \]
Combining this with the earlier display gives for all $\vecW,\vecW'$ with $\nrm*{\vecW-\vecW'} \le \delta$,
\[ \nrm*{\grad^2 F(\vecW) - \grad^2 F(\vecW')}_{\OPNORM} \le \nrm*{\vecW-\vecW'}\prn*{M_0+M_1 \rho_0(F(\vecW))},\]
where $\rho_0(x) = 2L_0^{1/2} x^{1/2} + \frac{5L_1^2}{L_0^{1/2}} x^{3/2}$. 
We thus see that $F$ satisfies \pref{ass:thirdorderselfbounding} with the non-decreasing, non-negative function $\rho_2(x)=M_0+M_1 \prn*{2L_0^{1/2} x^{1/2} + \frac{5L_1^2}{L_0^{1/2}} x^{3/2}}$, where the latter two properties are evident as $\rho_0(\cdot)$ is non-decreasing and non-negative.
\end{example}

\subsection{Proofs of Technical Results}\label{subsec:technicalhelperresults}
Now, we prove general results used throughout our work. 
We prove \pref{corr:gradcontrol}, which gives us control over the gradient:
\begin{proof}[Proof of \pref{corr:gradcontrol}]
Applying Lemma 11, \citet{priorpaper} with $\Phi$ in place of $F$, we obtain
\[ \nrm*{\grad F(\vecW)} \le \rho(F(\vecW))\sqrt{2\theta(F(\vecW))} = \rho_0(F(\vecW)),\]
where $\theta(\cdot)$ is defined as in the statement of \pref{corr:gradcontrol}. To prove $\rho_0(x)$ is increasing, simply note $\theta$ and thus $\sqrt{\theta}$ are clearly increasing, and are both non-negative. $\rho_1$ is non-decreasing and non-negative as well, thus $\rho_0$ is non-decreasing and non-negative.
\end{proof}
We also prove the central \pref{lem:boundfuncvalueradius}, which is very important to our results: it lets us control the change in function value under our regularity assumptions. We first state the following Lemma from \citet{li2023convex}, a generalization of Gronwall's Inequality:
\begin{lemma}[Lemma A.3, \citet{li2023convex}]\label{lem:helperlemgronwall}
Let $\alpha:[a,b] \rightarrow[0,\infty)$ and $\beta:[0,\infty)\rightarrow[0,\infty)$ be two continuous functions. Suppose $\alpha'(t) \le \beta(\alpha(t))$ almost everywhere over $(a,b)$. Let $\phi(u) = \int_0^u \frac1{\beta(v)}\DERIV v$. Then for all all $t \in [a,b]$,
\[ \phi(\alpha(t)) \le \phi(\alpha(a))-a+t.\]
\end{lemma}
This allows us to prove \pref{lem:boundfuncvalueradius}, which is an extension of Lemma A.4, \citet{li2023convex}:
\begin{proof}[Proof of \pref{lem:boundfuncvalueradius}]
The proof is essentially identical to the proof of Lemma A.4, \citet{li2023convex}. Let $\vecZ(t) = (1-t)\vecX+t\vecY$, $\alpha(t)=F(\vecZ(t))$. Then for all $t \in (0,1)$, we obtain
\begin{align*}
\alpha'(t) &= \lim_{s\rightarrow t}\frac{\alpha(s)-\alpha(t)}{s-t} \\
&\le \lim_{s\rightarrow t}\frac{\abs*{F(\vecZ(s))-F(\vecZ(t))}}{s-t} \\
&= \abs*{\lim_{s\rightarrow t}\frac{F(\vecZ(s))-F(\vecZ(t))}{s-t}} \\
&= \abs*{\frac{\DERIV}{\DERIV t}F(\vecZ(t))} \\
&= \abs*{\grad F(\vecZ(t))^{\top} (\vecY-\vecX)} \\
&\le \rho_0(F(\vecZ(t)))\nrm*{\vecY - \vecX},
\end{align*}
the last step using $\nrm*{\grad F(\vecW)} \le \rho_0(F(\vecW))$. Let $\beta(x) = \nrm*{\vecY-\vecX} \rho_0(x)$ and let $\phi(u) = \int_0^u \frac1{\beta(v)} \DERIV v$. Thus, $\alpha'(t) \le \beta(\alpha(t))$ almost everywhere. Applying \pref{lem:helperlemgronwall} gives 
\[ \phi(F(\vecY)) = \phi(\alpha(1)) \le \phi(\alpha(0))+1 = \phi(F(\vecX))+1. \]
Let $\psi(u) = \nrm*{\vecY-\vecX} \phi(u) = \int_0^u \frac1{\rho_0(v)} \DERIV v$, which is clearly strictly increasing. Consequently we obtain from the above and assumption on $\vecY$ that
\begin{align*}
\psi(F(\vecY)) &\le \psi(F(\vecX))+\nrm*{\vecY-\vecX} \\
&\le \psi(F(\vecX)) + \frac{1}{\rho_0(F(\vecX)+1)} \\
&\le \int_0^{F(\vecX)} \frac1{\rho_0(v)}\DERIV v + \int_{F(\vecX)}^{F(\vecX)+1} \frac1{\rho_0(v)} \DERIV v \\
&= \int_0^{F(\vecX)+1}\frac1{\rho_0(v)} \DERIV v = \psi(F(\vecX)+1).
\end{align*}
Since $\psi$ is strictly increasing, taking inverses implies
\[ F(\vecY) \le F(\vecX)+1,\]
as desired.
\end{proof}
We also introduce the following Lemma, which lets us exploit \pref{ass:thirdorderselfbounding} to control the Lipschitz constant of the Hessian of $F$.
\begin{lemma}\label{lem:thirdorderlipschitzlem}
Suppose $F$ satisfies \pref{ass:thirdorderselfbounding}. Suppose $\vecX,\vecY \in \mathbb{R}^d$ are such that $\nrm*{\vecY-\vecX} \le r$ for some $r>0$. Then 
\[ \nrm*{\grad^2 F(\vecX) - \grad^2 F(\vecY)}_{\OPNORM} \le \nrm*{\vecX - \vecY} \cdot \sup_{\vecU \in \overline{\vecX \vecY}} \rho_2\prn*{F(\vecU)}.\]
In particular, we have 
\[ \nrm*{\grad^2 F(\vecX) - \grad^2 F(\vecY)}_{\OPNORM} \le \nrm*{\vecX - \vecY} \cdot \sup_{\vecU \in \ball(\vecY, r)} \rho_2\prn*{F(\vecU)}.\]
\end{lemma}
\begin{proof}
Consider $\delta>0$, either from \pref{ass:thirdorderselfbounding} if the second case of \pref{ass:thirdorderselfbounding} holds, and otherwise set to some arbitrary positive real. Similar to the proof of \pref{lem:boundfuncvalueradius}, divide the line segment between $\vecX, \vecY$ into $N = \frac{\nrm*{\vecX-\vecY}}{\delta}$ equally spaced segments of length $\delta$ between points $\vecX_i$, where we define $\vecX_0=\vecX, \vecX_1, \ldots, \vecX_{N-1}, \vecX_N = \vecY$. Thus $\nrm*{\vecX-\vecY}=N\delta$.

Suppose for all $\vecU \in \overline{\vecX\vecY}$ we have $\nrm*{\grad^3 F(\vecU)}_{\OPNORM} \le L$. Consider any $\vecX', \vecY'$ in the line segment $\overline{\vecX \vecY}$. Applying this for $\vecX' + t(\vecY'-\vecX')$ for $t \in [0,1]$, which always lies in the line segment $\overline{\vecX \vecY}$, we obtain
\[ \nrm*{\grad^2 F(\vecY') - \grad^2 F(\vecX')}_{\OPNORM} \le \nrm*{\int_0^1 \tri*{\grad^3 F(\vecX'+t(\vecY'-\vecX')), \vecY'-\vecX'} \DERIV t} \le L\nrm*{\vecY'-\vecX'}. \]
Consequently irrespective of which case of \pref{ass:thirdorderselfbounding} holds, because $\nrm*{\vecX_i - \vecX_{i-1}} \le \delta$, we have for each $i, 1 \le i \le N$ that 
\[ \nrm*{\grad^2 F(\vecX_i) - \grad^2 F(\vecX_{i-1})}_{\OPNORM} \le \nrm*{\vecX_i - \vecX_{i-1}} \sup_{\vecU \in \overline{\vecX \vecY}} \rho_2\prn*{F(\vecU)}.\]
Now Triangle Inequality gives
\begin{align*}
\nrm*{\grad^2 F(\vecX) - \grad^2 F(\vecY)}_{\OPNORM} &\le \sum_{i=1}^N \nrm*{\grad^2 F(\vecX_i) - \grad^2 F(\vecX_{i-1})}_{\OPNORM} \\
&\le \sum_{i=1}^N \nrm*{\vecX_i - \vecX_{i-1}} \sup_{\vecU \in \overline{\vecX \vecY}} \rho_2\prn*{F(\vecU)} \\
&\le N \delta \cdot \sup_{\vecU \in \overline{\vecX \vecY}} \rho_2\prn*{F(\vecU)} \\
&= \nrm*{\vecX-\vecY} \sup_{\vecU \in \overline{\vecX \vecY}} \rho_2\prn*{F(\vecU)},
\end{align*}
as desired.
\end{proof}
We will also generalize the proof of \pref{thm:gdfirstorder} to show that GD, when initialized in the $F(\vecW_0)$-sublevel set $\cL_{F,F(\vecW_0)}$ with appropriate step size defined in terms of $F(\vecW_0)$, never increases function value.
\begin{lemma}\label{lem:gdnoincrease}
Consider any $\vecW_0\in\mathbb{R}^d$, and consider iterates $\{\vecU_t\}_{t \ge 0}$ of GD initialized at any $\vecU_0 \in \cL_{F,F(\vecW_0)}$, the $F(\vecW_0)$-sublevel set. If the step size $\eta$ of GD is at most $\frac1{L_1(\vecW_0)}$ where $L_1(\cdot)$ is defined as per \pref{eq:L1def}, then $F(\vecU_t) \le F(\vecU_0)$ for all $t \ge 0$.
\end{lemma}
\begin{proof}
It suffices to prove this for $t=1$; a simple inductive argument then establishes this for all $t \ge 0$. We have $\vecU_1 = \vecU_0 - \eta \grad F(\vecU_0)$. 
By \pref{corr:gradcontrol} and because $\vecU_0 \in \cL_{F,F(\vecW_0)}$, $\nrm*{\grad F(\vecU_0)} \le \rho_0(F(\vecU_0)) \le \rho_0(F(\vecW_0))$. Thus by choice of $\eta$ and definition of $L_1(\vecW_0)$,
\[ \nrm*{\vecU_1- \vecU_0} = \eta \nrm*{\grad F(\vecU_0)} \le \eta \rho_0(F(\vecW_0)) \le \frac1{\rho_0(F(\vecW_0)+1)}.\]
By \pref{lem:boundlocalsmoothnesseasy}, because $\vecU_0 \in \cL_{F,F(\vecW_0)}$, for all $\vecP$ in the line segment $\overline{\vecU_0 \vecU_1}$we have $\nrm*{\grad^2 F(\vecP)}_{\OPNORM} \le L_1(\vecW_0)$.
By \pref{lem:secondordersmoothnessineq}, it follows that
\begin{align*}
F(\vecU_1) &\le F(\vecU_0) - \eta \nrm*{\grad F(\vecU_0)}^2 + \frac{L_1(\vecW_0) \eta^2}2 \cdot \nrm*{\grad F(\vecU_0)}^2 \\
&\le F(\vecU_0) + \nrm*{\grad F(\vecU_0)}^2 \cdot \prn*{-\eta + \frac{L_1(\vecW_0)\eta^2}2}.
\end{align*}
Noting $-\eta + \frac{L_1(\vecW_0)\eta^2}2 \le 0$ for $\eta \in \brk*{0, \frac{2}{L_1(\vecW_0)}}$, the conclusion follows.
\end{proof}

\section{Proof of Framework}\label{sec:generalframeworkpf}
\begin{proof}[Proof of \pref{thm:generalframework}]
For convenience, for all $n \ge 0$, define $p_n := 1-n\cdot\sup_{\vecU \in \cL_{F, F(\vecW_0)} }\delta(\vecU)$.
Also let $T = \sup_{\vecU \in \cL_{F, F(\vecW_0)} }\crl*{\frac{F(\vecW_0)}{\Delta(\vecU)}}$.
\begin{lemma}\label{lem:frameworkinduction}
For any $n \ge 0$, let $\cE_n$ be the event that the sequence of iterates $(\vecW_t)_{0 \le t \le n-1}$ satisfies either:
\begin{enumerate}
    \item The event $\cE_{n,1}$: For all $0 \le t \le n-1$, $F(\cA_1(\vecW_{t})) < F(\vecW_t) - \Delta(\vecW_t)$.
    \item The event $\cE_{n,2}$: There exists $\vecW_t \in (\vecW_t)_{0 \le t \le n-1}$ such that $\cA_2(\vecW_t) \cap \cS \neq \{ \}$, and for all $\vecW_s$ with $0 \le s<t$, we have $F(\cA_1(\vecW_{s})) < F(\vecW_s) - \Delta(\vecW_s)$. 
\end{enumerate}
That is, $\cE_n = \cE_{n,1} \cup \cE_{n,2}$. Then over the randomness in $\cA$, we have $\mathbb{P}(\cE_n) \ge p_n$ for all $n \ge 0$.
\end{lemma}
\begin{proof}
We proceed by induction on $n$. The base case $n=0$ is vacuously evident, and the case $n=1$ follows immediately by the definition of a decrease procedure from \pref{def:highprobdecreasealg} and hypotheses of \pref{thm:generalframework}.

For the inductive step, suppose \pref{lem:frameworkinduction} is true for some $n \ge 1$; we show it is for $n+1$. By the inductive hypothesis, we know that $\mathbb{P}(\cE_n) \ge p_n$. We aim to show $\mathbb{P}(\cE_n) \ge p_{n+1}$. If $p_n \le 0$ there is nothing to prove, so suppose now that $n$ is such that $p_n > 0$.
\begin{enumerate}
    \item Let $p = \mathbb{P}\prn*{\cE_{n,2}|\cE_n}$. Note $\cE_{n,2} \subseteq \cE_{n+1,2} \subseteq \cE_{n+1}$. 
    
    \item Let $\cB := \cE_{n,1} \cap \cE_{n,2}^c$. Thus, if $\cB$ occurs, then all the $(\vecW_t)_{0 \le t \le n-1}$ are such that $F(\cA_1(\vecW_t)) < F(\vecW_t) - \Delta(\vecW_t)$, but $\cE_{n,2}$ did not occur. Note $\cE_n$ is the disjoint union $\cE_{n, 2} \sqcup \cB$, so $\mathbb{P}(\cB | \cE_n)=1-p$. 
    
    Under $\cB$, we know $\vecW_{n} = \cA(\vecW_{n-1})$ is such that $F(\vecW_{n}) \le F(\vecW_0)$. Hence $\vecW_{n} \in \cL_{F,F(\vecW_0)}$. Therefore, conditioned on $\cB$, by the hypotheses of \pref{thm:generalframework} we have with probability at least $p_0$ that either $F(\cA_1(\vecW_n)) < F(\vecW_n) - \Delta(\vecW_n)$ or $\cA_2(\vecW_n) \cap \cS \neq \{\}$. 
    
    Let $\cC$ be the event that $F(\cA_1(\vecW_n)) < F(\vecW_n) - \Delta(\vecW_n)$ occurs. Let $\cD$ be the event that $\cA_2(\vecW_n) \cap \cS \neq \{\}$ occurs but $\cC$ does not occur. 
    Recall that $\vecW_n \in \cL_{F,F(\vecW_0)}$ conditioned on $\cB$. Furthermore recall that $\cA(\vecW_n)$ is only a function of $\vecW_n$, and none of the $(\vecW_t)_{0 \le t \le n-1}$. 
    Thus the definition of decrease procedure, \pref{def:highprobdecreasealg}, implies that 
    \[ \mathbb{P}\prn*{\cC \sqcup \cD | \cB} \ge p_0.\]
    Now since $\mathbb{P}(\cB) = \mathbb{P}(\cB | \cE_n) \mathbb{P}(\cE_n) \ge (1-p)p_n > 0$, Bayes' Rule implies
    \begin{align*}
    \mathbb{P}\prn*{(\cB \cap \cC) \sqcup (\cB \cap \cD) | \cB} &= \frac{\mathbb{P}(\cB \cap ((\cB \cap \cC) \sqcup (\cB \cap \cD)))}{\mathbb{P}(\cB)} \\
    &= \frac{\mathbb{P}(\cB \cap (\cC \sqcup \cD))}{\mathbb{P}(\cB)} =\mathbb{P}\prn*{\cC \sqcup \cD | \cB} \ge p_0.
    \end{align*}
    Note $\cB \cap \cC$ implies that $\cE_{n+1,1}$ occurs, since under $\cB \cap \cC$ we have $F(\cA_1(\vecW_t)) < F(\vecW_t) - \Delta(\vecW_t)$ for all $0 \le t \le n$. 
    Similarly, $\cB \cap \cD$ implies that $\cE_{n+1,2}$ occurs, since under $\cB \cap \cD$ we have $F(\cA_1(\vecW_t)) < F(\vecW_t) - \Delta(\vecW_t)$ for $0 \le t \le n-1$ and $\cA_2(\vecW_n) \cap \cS \neq \{\}$.
\end{enumerate}
Thus recalling $\cE_{n,2}, \cB$ are disjoint, we see that $\cE_{n+1}$ contains the following disjoint union of events:
\[ \cE_{n+1} \supseteq \cE_{n,2} \sqcup (\cB \cap \cC) \sqcup (\cB \cap \cD).\]
The above observations imply via Bayes' Rule that
\begin{align*}
\mathbb{P}(\cE_{n+1}) &\ge \mathbb{P}\prn*{\cE_{n,2} \sqcup (\cB \cap \cC) \sqcup (\cB \cap \cD)} \\
&= \mathbb{P}\prn*{\cE_{n,2}} + \mathbb{P}\prn*{(\cB \cap \cC) \sqcup (\cB \cap \cD)} \\
&= \mathbb{P}\prn*{\cE_{n,2} | \cE_n} \mathbb{P}\prn*{\cE_n} + \mathbb{P}\prn*{(\cB \cap \cC) \sqcup (\cB \cap \cD) | \cB} \mathbb{P}(\cB | \cE_n) \mathbb{P}(\cE_n) \\
&= \mathbb{P}(\cE_n)\prn*{ p +  \mathbb{P}\prn*{(\cB \cap \cC) \sqcup (\cB \cap \cD) | \cB} \cdot (1-p)} \\
&\ge p_n \prn*{p + p_0 (1-p)}\\
&\ge p_n (p_0 p + p_0 (1-p)) = p_n p_0 \ge p_{n+1}.
\end{align*}
Here we used that $\mathbb{P}(\cE_n) \ge p_n$, $p_n p_0 \ge p_{n+1}$ which follows immediately from the definition of $p_n$, $p_0 \le 1$, and simple manipulations. The inductive step, and hence the proof, is thus complete.
\end{proof}
Using \pref{lem:frameworkinduction} now readily proves the following:
\begin{claim}\label{claim:intuitiveclaim}
Let $\cE$ be the event that there exists $\vecW_t$ with $\vecW_t \in (\vecW_t)_{0 \le t \le T-1}$ such that $\cA_2(\vecW_t) \cap \cS \neq \{ \}$, and for all $\vecW_s$ with $0 \le s<t$, we have $F(\cA_1(\vecW_{s})) < F(\vecW_s) - \Delta(\vecW_s)$. Then $\mathbb{P}(\cE) \ge p_{T}$.
\end{claim}
\begin{proof}[Proof of \pref{claim:intuitiveclaim}]
Apply \pref{lem:frameworkinduction} with $n=T$. Following the notation from there, we have that the event $\cE_{T} = \cE_{T, 1} \sqcup \cE_{T,2}$ has probability at least $p_{T}$. 

Suppose that $\cE_{T,1}$ occurs. Note $\cE_{T,1}$ implies that $\vecW_t \in \cL_{F,F(\vecW_0)}$ for all $0 \le t \le T$. Therefore 
\[ \Delta(\vecW_t) \ge \inf_{\vecU\in\cL_{F,F(\vecW_0)}}\Delta(\vecU)\text{ for all }0 \le t \le T. \numberthis\label{eq:frameworkpflowerbounddecrease}\]
Moreover, telescoping the direct implication of $\cE_{T,1}$ gives that
\[ F(\vecW_T) < F(\vecW_0) - \sum_{t=0}^{T-1} \Delta(\vecW_t).\numberthis\label{eq:frameworkpfdecreasetelescope} \]
Combining \pref{eq:frameworkpflowerbounddecrease} and \pref{eq:frameworkpfdecreasetelescope} and recalling that we shifted WLOG so $F$ has minimum value 0 (see Notation) gives
\[ T \inf_{\vecU\in \cL_{F, F(\vecW_0)}} \Delta(\vecU) \le \sum_{t=0}^{T-1} \Delta(\vecW_t) < F(\vecW_0) - F(\vecW_T) \le F(\vecW_0).\]
This contradicts our choice of $T$. 

Thus $\cE_{T,1}$ cannot occur, and so $\cE_{T,2}$ must occur, i.e. $\cE_{T} = \cE_{T,2}$. Note $\cE_{T,2}$ is exactly the event $\cE$. Thus 
\[ \mathbb{P}(\cE) = \mathbb{P}(\cE_{T,2}) = \mathbb{P}(\cE_{T}) \ge p_{T},\]
as desired.
\end{proof}
Conditioning on the event $\cE$ from \pref{claim:intuitiveclaim}, by \pref{claim:intuitiveclaim}, we immediately recover the desired guarantee on the output, probability, and number of candidate vectors stated in \pref{thm:generalframework}. The only part remaining to prove \pref{thm:generalframework} is to establish the bound $N = \frac{F(\vecW_0)}{\overline{\Delta}}+\sup_{\vecU \in \cL_{F, F(\vecW_0)} }\Toracle(\vecU)$ on the number of oracle calls. 

To this end, condition on $\cE$ from \pref{claim:intuitiveclaim} in all of the following, and follow the notation from there, in particular the definition of $\vecW_t$. Directly, we obtain that the number of oracle calls is at most $\sum_{i=0}^t \Toracle(\vecW_i)$ (the last term $\Toracle(\vecW_t)$ in the sum appears since computing $\cA(\vecW_t)$ and $\cA(\vecW_t)$ takes at most $\Toracle(\vecW_t)$ oracle calls). We now upper bound this sum.

As we are conditioning on $\cE$ and since we assumed WLOG by shifting that $F$ has minimum value 0, we have
\[ F(\vecW_{i+1}) - F(\vecW_{i}) < -\Delta(\vecW_{i}) < 0 \text{ for all } 0 \le i \le t-1 \implies \sum_{i=0}^{t-1} \Delta(\vecW_i) < F(\vecW_0)-F(\vecW_t) \le F(\vecW_0). \numberthis\label{eq:frameworkpfdecreasesumfinish}\]
The above also implies $F(\vecW_i) \le F(\vecW_0)$, i.e. $\vecW_i \in \cL_{F,F(\vecW_0)}$, for all $0 \le i \le t$. Therefore, $\Toracle(\vecW_{i}) \le \sup_{\vecU\in\cL_{F, F(\vecW_0)}} \Toracle(\vecU)$ for all $0 \le i \le t$. Thus \pref{eq:frameworkpfdecreasesumfinish} gives
\[ \frac{F(\vecW_0)}{\sum_{i=0}^{t-1} \Toracle(\vecW_{i})} > \frac{\sum_{i=0}^{t-1} \Delta(\vecW_i)}{\sum_{i=0}^{t-1} \Toracle(\vecW_{i})} \ge \min_{0 \le i \le t-1} \frac{\Delta(\vecW_i)}{\Toracle(\vecW_i)} \ge \overline{\Delta},\]
where the last inequality uses the elementary inequality $\frac{\sum_{i=1}^{k'} a_i}{\sum_{i=1}^{k'} b_i} \ge \min_i \frac{a_i}{b_i}$ for $a_i \ge 0, b_i > 0$, that $\vecW_i \in \cL_{F,F(\vecW_0)}$ for all $0 \le i \le t-1$, and the definition of $\overline{\Delta}$.
Rearranging and recalling $\Toracle(\vecW_t) \le \sup_{\vecU\in\cL_{F, F(\vecW_0)}} \Toracle(\vecU)$ as justified above, we obtain 
\[ \sum_{i=0}^{t} \Toracle(\vecW_{i}) \le  \sup_{\vecU\in \cL_{F, F(\vecW_0)}} \Toracle(\vecU) + \sum_{i=0}^{t-1} \Toracle(\vecW_{i})  \le  \sup_{\vecU\in \cL_{F, F(\vecW_0)}} \Toracle(\vecU) +\frac{F(\vecW_0)}{\overline{\Delta}}.\]
This yields the desired conclusion on oracle complexity, completing the proof.
\end{proof}

\section{First Order Convergence Proofs}\label{sec:firstorderproofs}
\subsection{Proofs for Adaptive GD}\label{subsec:AdaptiveGDFOSP}
\begin{proof}
As with the proof of \pref{thm:gdfirstorder}, we use \pref{thm:generalframework}. We again have $\cS = \{\vecW:\nrm*{\grad F(\vecW)} \le \epsilon\}$, and recall the choice of $\eta$ from \pref{thm:adaptivegdguarantee}. Now we let $\cA(\vecU_0) = \prn*{\vecU_0 - \eta_{\vecU_0} \grad F(\vecU_0), \vecU_0}$. Thus$\cA_1(\vecU_0) = \vecU_0 - \eta \grad F(\vecU_0)$, $\cA_2(\vecU_0)=\vecU_0$, and $\Toracle(\vecU_0)=1$.

\begin{claim}\label{claim:adaptivegdfirstorderdecreaseprocedure}
For any $\vecU_0$ in the $F(\vecW_0)$-sublevel set $\cL_{F,F(\vecW_0)}$, $\cA$ is a $(\cS, 1, \min\crl*{\frac{L_1'(\vecW_0)}{2\rho_0(F(\vecW_0)+1)^2}, \frac{\epsilon^2}{2 L'_1(\vecW_0)}}, 0, \vecU_0)$-decrease procedure.
\end{claim}
To show this, analogously to the proof of \pref{thm:gdfirstorder}, for any $\vecU_0 \not \in \cS$ in the $F(\vecW_0)$-sublevel set $\cL_{F,F(\vecW_0)}$, we will show that the function will deterministically decrease by strictly greater than $\min\crl*{\frac{L_1'(\vecW_0)}{\rho_0(F(\vecW_0)+1)^2}, \frac{\epsilon^2}{2 L'_1(\vecW_0)}}$ at the next iterate. By definition of $\cA_2$, exactly as with the proof of \pref{thm:gdfirstorder}, we conclude via \pref{thm:generalframework} upon showing \pref{claim:adaptivegdfirstorderdecreaseprocedure}. 

To show \pref{claim:adaptivegdfirstorderdecreaseprocedure}, by choice of step size, we have $\eta_{\vecU_0} \nrm*{\grad F(\vecU_0)} \le \frac1{\rho_0(F(\vecW_0)+1)}$. Thus 
\[ \nrm*{\vecU_1-\vecU_0} \le \frac1{\rho_0(F(\vecW_0)+1)} \le \frac1{\rho_0(F(\vecU_0)+1)}. \]
Now combining \pref{lem:boundfuncvalueradius} with \pref{ass:selfbounding}, and because $\vecU_0 \in \cL_{F,F(\vecW_0)}$, we see for all $\vecP \in \overline{\vecU_0 \vecU_1}$, $\nrm*{\grad^2 F(\vecP)}_{\OPNORM} \le L'_1(\vecW_0)$ where $L'_1(\vecW_0)$ is defined as in the statement of \pref{thm:adaptivegdguarantee}. We thus obtain by \pref{lem:secondordersmoothnessineq},
\[ F(\vecU_1) \le F(\vecU_0) - \eta \nrm*{\grad F(\vecU_0)}^2 + \frac{L'_1(\vecW_0) \eta^2}2 \cdot \nrm*{\grad F(\vecU_0)}^2.\numberthis\label{eq:adaptivegddescent}\]
Recall that $\vecU_0\not\in\cS$, so $\nrm*{\grad F(\vecU_0)} > \epsilon$. We break into cases:
\begin{enumerate}
    \item If $\nrm*{\grad F(\vecU_0)} > \frac{L'_1(\vecW_0)}{\rho_0(F(\vecW_0)+1)}$, then $\eta_{\vecU_0} = \frac1{\rho_0(F(\vecW_0)+1)\nrm*{\grad F(\vecU_0)}}$. In this case, substituting into \pref{eq:adaptivegddescent} gives
    \begin{align*}
    F(\vecU_1) &\le F(\vecU_0) - \eta \nrm*{\grad F(\vecU_0)}^2 + \frac{L'_1(\vecW_0) \eta^2}2 \cdot \nrm*{\grad F(\vecU_0)}^2 \\
    &= F(\vecU_0) - \frac1{\rho_0(F(\vecW_0)+1)} \nrm*{\grad F(\vecU_0)} + \frac{L'_1(\vecW_0)}{2\rho_0(F(\vecW_0)+1)^2} \\
    &< F(\vecU_0) - \frac12 \cdot \frac{L'_1(\vecW_0)}{\rho_0(F(\vecW_0)+1)^2}.
    \end{align*}
    \item Else if $\nrm*{\grad F(\vecU_0)} \le L'_1(\vecW_0)$, then $\eta_{\vecU_0} = \frac1{L_1'(\vecW_0)}$. In this case, substituting into \pref{eq:adaptivegddescent} gives
    \begin{align*}
    F(\vecU_1) &\le F(\vecU_0) - \eta \nrm*{\grad F(\vecU_0)}^2 + \frac{L'_1(\vecW_0) \eta^2}2 \cdot \nrm*{\grad F(\vecU_0)}^2 \\
    &\le F(\vecU_0) - \frac{\nrm*{\grad F(\vecU_0)}^2}{2L'_1(\vecW_0)} < F(\vecU_0) - \frac{\epsilon^2}{2L'_1(\vecW_0)},
    \end{align*}
    where we used that $\nrm*{\grad F(\vecU_0)} > \epsilon$.
\end{enumerate}
In either case, for $\nrm*{\grad F(\vecU_0)} > \epsilon$ we have that
\begin{align*}
F(\vecU_1) &< F(\vecU_0)-\min\crl*{\frac{L'_1(\vecW_0)}{2\rho_0(F(\vecW_0)+1)^2}, \frac{\epsilon^2}{2 L'_1(\vecW_0)}}.
\end{align*}
This proves \pref{claim:adaptivegdfirstorderdecreaseprocedure}. By our framework \pref{thm:generalframework}, the proof is complete.
\end{proof}

\subsection{Proofs for SGD for FOSPs}\label{subsec:SGDFOSP}
Here, we prove \pref{thm:sgdfirstorder}. We first introduce technical preliminaries, which will also be used in \pref{sec:restartedsgdproofs}.
\begin{theorem}[Vector-Valued Azuma-Hoeffding, Theorem 3.5 in \citet{pinelis1994optimum}]\label{thm:subgaussianazumahoeffding}
Let $\vecEps_1, \ldots, \vecEps_K \in \mathbb{R}^d$ be such that for all $k$, $\mathbb{E}\brk*{\vecEps_k | \mathfrak{F}^{k-1}}=0$, $\nrm*{\vecEps_k}^2 \le \sigma_k^2$. Then for any $\lambda>0$,
\[ \mathbb{P}\prn*{\nrm*{\sum_{k=1}^K \vecEps_k} \ge \lambda} \le 4\exp\prn*{-\frac{\lambda^2}{4\sum_{k=1}^K \sigma_k^2}}.\]
\end{theorem}
Note the bound here is dimension free, so this result does not follow directly from standard Azuma-Hoeffding. Such a result can also be found in \citet{kallenberg1991some, zhang2005learning, fang2019sharp}.

\begin{theorem}[Data-Dependent Concentration Inequality, Lemma 3 in \citet{rakhlin2012making}]\label{thm:subgaussiandatadependentconcentration}
Let $\epsilon_1, \ldots, \epsilon_K \in \mathbb{R}$ be such that for all $k$, $\mathbb{E}\brk*{\epsilon_k | \mathfrak{F}^{k-1}}=0$, $\mathbb{E}\brk*{\epsilon_k^2 | \mathfrak{F}^{k-1}} \le \sigma_k^2$. Furthermore suppose that $\mathbb{P}\prn*{\epsilon_k \le b | \mathfrak{F}^{k-1}}=1$. Letting $V_K = \sum_{k=1}^K \sigma_k^2$, for any $\delta  < 1/e$, $K \ge 4$, we have
\[ \mathbb{P}\prn*{\sum_{k=1}^K \epsilon_k  > 2\max\crl*{2\sqrt{V_k}, b\sqrt{\log(1/\delta)}} \sqrt{\log(1/\delta)}} \le \delta \log(K).\]
\end{theorem}
Such a result is also presented in \citet{zhang2005learning, bartlett2008high, fang2019sharp}.

We will first prove \pref{thm:sgdfirstorder} in the case where $\nrm*{\grad f(\vecW;\vecZeta) - \grad F(\vecW)}$ is bounded by $\sigma(F(\vecW))$. As noted in \citet{fang2019sharp}, these same inequalities hold when the martingale difference is not bounded or almost-surely bounded but rather the norms are sub-Gaussian with parameter $\sigma_k$. Thus after the proof, we remark how to straightforwardly generalize \pref{thm:sgdfirstorder} to the case when $\nrm*{\grad f(\vecW;\vecZeta) - \grad F(\vecW)}$ is sub-Gaussian with parameter $\sigma(F(\vecW))$ in \pref{rem:sgdfirstordersubgaussian}.

Now, we prove \pref{thm:sgdfirstorder}.
\begin{proof}
We use our framework \pref{thm:generalframework} with $\cS=\{\vecW:\nrm*{\grad F(\vecW)} \le \epsilon\}$. Recall as per the discussion of SGD in our framework in \pref{subsec:exinframework}, we let $\vecP_0=\vecU_0$, and define a sequence $(\vecP_i)_{0 \le i \le K_0}$ via
\[ \vecP_i = \vecP_{i-1} - \eta \grad f(\vecP_{i-1};\vecZeta_{i}),\]
where the $\vecZeta_{i}$ are minibatch samples i.i.d. across different $i$. Note this sequence can be equivalently defined by repeated compositions of the function $\vecU \rightarrow \vecU - \eta \grad f(\vecU; \vecZeta)$.

We now let $\cA(\vecU_0) = \prn*{ \vecP_{K_0}, (\vecP_i)_{0 \le i \le K_0-1}}$, hence $\cA_1(\vecU_0)=\vecP_{K_0}$, $\cA_2(\vecU_0)=(\vecP_i)_{0 \le i \le K_0-1}$. Thus $\Toracle(\vecU_0)=K_0$. Also note the noise $\vecXi_t$ used defining $\cA$ are independent across different $t$. 

For appropriate $\eta=\tilde{\Theta}(\epsilon^2)$, $K_0=\tilde{\Theta}(\epsilon^{-2})$ depending only on $\epsilon,\delta,F(\vecW_0)$ and polylogarithmically in $1/\delta$, which we define below, we establish the following \pref{claim:sgdfirstorderdecreasesequence}: 
\begin{claim}\label{claim:sgdfirstorderdecreasesequence}
For any $\vecU_0$ in the $F(\vecW_0)$-sublevel set $\cL_{F,F(\vecW_0)}$, $\cA$ is a $(\cS, K_0, \frac{\eta K_0 \epsilon^2}{4}, p, \vecU_0)$-decrease procedure, where $p=\frac{\delta \eta K_0 \epsilon^2}{4(F(\vecW_0)+1)}$.
\end{claim}
Then using \pref{thm:generalframework}, we then directly conclude \pref{thm:sgdfirstorder}.

To show \pref{claim:sgdfirstorderdecreasesequence}, consider any $\vecU_0$ in the $F(\vecW_0)$-sublevel set but not in $\cS$. 
Following the notation from above, consider a `block' of $K_0$ consecutive iterates of SGD starting at $\vecP_0 = \vecU_0$. We establish that with probability at least $1-p$, if none of the iterates $\{\vecP_0=\vecU_0, \ldots, \vecP_{K_0-1}\}$ lie in $\cS$, then $F(\vecP_{K_0}) < F(\vecP_0)-\Delta$ where $\Delta = \frac{\eta K_0 \epsilon^2}4$. Then recalling the definitions of $\cA_2$, we immediately conclude \pref{claim:sgdfirstorderdecreasesequence}.

\paragraph{Definitions and Parameters:} For convenience, define
\begin{align*}
L_0(\vecW_0) &= \rho_0(F(\vecW_0)+1),\\ L_1(\vecW_0)&=\rho_1(F(\vecW_0)+1),\\
\sigma_1(\vecW_0) &= \sigma(F(\vecW_0)+1), \\
B(\vecW_0)&=\sigma_1(\vecW_0)^2+\frac18 \sigma_1(\vecW_0) L_0(\vecW_0).
\end{align*}
Also define
\[ \vecXi_{t+1} = \grad f(\vecP_t;\vecZeta_{t+1}) - \grad F(\vecP_t),\]
where $\vecZeta_{t+1}$ denotes the i.i.d. minibatch samples. Note by \pref{ass:noiseregularity} that $\mathbb{E}\brk*{\vecXi_{t+1}}=0$, where expectation is with respect to $\vecZeta_{t+1}$.

In particular, we choose these parameters as follows: 
\begin{align*}
\tilde{\eta} &= \frac{\epsilon^2}{\tilde{L}(\vecW_0) \log(1/\epsilon)^6 \log(1/\delta)^6} \\
K_0 &= \frac{C(\vecW_0)}{\epsilon^2} \log(1/\tilde{\eta})^2 \log(1/\delta)^2 \log(1/\epsilon)^2, \\
\eta &= \frac1{\max\crl*{1, \rho_0(F(\vecW_0)+1)}} \cdot \tilde{\eta},
\end{align*}
where 
\begin{align*}
C(\vecW_0) &= 128 B(\vecW_0) \lor 64(F(\vecW_0)+1)^2,\\
\tilde{L}'(\vecW_0) &= 8L_1(\vecW_0)(L_0(\vecW_0)^2+\sigma_1(\vecW_0)^2) \lor 2L_0(\vecW_0) \lor 4\sigma_1(\vecW_0),\\
\tilde{L}(\vecW_0) &= \tilde{L}'(\vecW_0)^2 C(\vecW_0)^2 \lor (3\sqrt{2} \log(\tilde{L}(\vecW_0)))^8 \lor (3\sqrt{2})^8.
\end{align*}
\begin{remark}
Note that $C, \tilde{L}', \tilde{L}$ depend only polynomially in terms of the self-bounding functions $\rho_0, \rho_1$, $\sigma$, and $F(\vecW_0)$.
\end{remark}
Note we can assume WLOG that $\epsilon$ and the desired probability $\delta$ are at most some small enough \textit{universal} constants in $(0,1)$; by doing so, the result does not change up to universal constant, and hence is identical under the $O(\cdot)$. Consequently we may assume WLOG that $\tilde{\eta}$ and $\eta$ are at most some small enough universal constant in $(0,1)$ and that $K_0 \ge 4$.
\begin{claim}\label{claim:sgdproofclaimparams}
For $\epsilon,\delta$ small enough \textit{universal} constants, the above choice of parameters satisfies the following properties:
\begin{align*}
&\max\crl*{1, \rho_0(F(\vecW_0)+1)}\eta \\
&= \tilde{\eta} \le \min\crl*{\frac{\epsilon^2}{8L_1(\vecW_0)(L_0(\vecW_0)^2+\sigma_1(\vecW_0)^2)}, \frac1{2K_0 L_0(\vecW_0)}, \frac1{4\sigma_1(\vecW_0)\sqrt{K_0 \log(4K_0/p)}}},\numberthis\label{eq:etaparamchocesgdfirstorder}\\
K_0 \epsilon^2 &\ge 128 B(\vecW_0) \log\prn*{\frac{2 \log K_0}{p}}.\numberthis\label{eq:K0paramchocesgdfirstorder}
\end{align*}
\end{claim}
For the sake of brevity, we prove \pref{claim:sgdproofclaimparams} after the our main proof. Checking this is a matter of elementary, albeit tedious, univariate inequalities.

Again, our plan is to apply \pref{thm:generalframework} by showing decrease with high probability for a block of $K_0$ iterates starting at $\vecP_0$. 

\paragraph{Notation:} Let $\mathfrak{F}^t$ denote the filtration of all information up through $\vecP_t$, but \textit{not} including the minibatch sample $\vecZeta_{t+1}$. Let $\cK$ be a stopping time denoting the first $t$ such that $\vecP_t \not\in \ball\prn*{\vecP_0, \frac1{\rho_0(F(\vecW_0)+1)}}$, i.e. the escape time of the iterates beginning at $\vecP_0$ from $\ball\prn*{\vecP_0, \frac1{\rho_0(F(\vecW_0)+1)}}=\ball\prn*{\vecU_0, \frac1{\rho_0(F(\vecW_0)+1)}}$.

We first detail two high probability events we will condition on for the remainder of the proof:
\begin{itemize}
\item By Vector-Valued Azuma Hoeffding \pref{thm:subgaussianazumahoeffding}, for a given $1 \le t \le K_0$ we have with probability at least $1-\frac{p}{2K_0}$,
\[ \nrm*{\eta \sum_{k=1}^t \vecXi_k} \le 2\eta \sqrt{\log(48K_0/p) \sum_{k=1}^t \sigma\prn*{F(\vecP_{k-1})}^2} = 2\eta \sqrt{\log(4K_0/p) \sum_{k=0}^{t-1} \sigma\prn*{F(\vecP_{k-1})}^2}.\]
This follows since each $\mathbb{E}\brk*{\vecXi_k | \mathfrak{F}^{k-1}} = 0$ as the stochastic gradient oracle is unbiased, and as $\nrm*{\vecXi_k} \le \sigma\prn*{F(\vecP_{k-1})}$ by \pref{ass:noiseregularity}. 

Thus by Union Bound, with probability at least $1-p/2$, we have for all $1 \le t \le K_0$ that
\[ \nrm*{\eta \sum_{k=1}^t \vecXi_k} \le2\eta \sqrt{\log(4K_0/p) \sum_{k=0}^{t-1} \sigma\prn*{F(\vecP_{k})}^2}.\numberthis\label{eq:noiseupperboundsgdfirstorder}\]
Denote this event by $\cE_1$, so $\mathbb{P}(\cE_1) \ge 1-p/2$.
\item We define a stochastic process with the following trick to derive uniform bounds. Define the following sequence of real numbers:
\[ Y_{t} := -\eta \tri*{\grad F(\vecP_t), \vecXi_{t+1}} 1_{t < \cK}.\]
Notice $1_{t<\cK}$ is $\mathfrak{F}^{t}$-measurable, as $\{t<\cK\}$ holds if and only if $\vecP_1, \ldots, \vecP_t \in \ball\prn*{\vecP_0, \frac1{\rho_0(F(\vecW_0)+1)}}$.

Clearly $\grad F(\vecP_t)$ is also $\mathfrak{F}^t$-measurable. Thus as the stochastic gradient oracle is unbiased (i.e. $\mathbb{E}\brk*{\vecXi_{t+1} | \mathfrak{F}^t}=0$),
\[ \mathbb{E}\brk*{Y_t} = \mathbb{E}\brk*{\tri*{\grad F(\vecP_t), \vecXi_{t+1}} 1_{t < \cK} | \mathfrak{F}^{t}} = 0. \]
For $t \ge \cK$ we have $Y_t \equiv 0$. For $t < \cK$, we have $\vecP_t \in \ball\prn*{\vecP_0, \frac1{\rho_0(F(\vecW_0)+1)}}$. Consequently by \pref{lem:boundfuncvalueradius} and \pref{corr:gradcontrol} we have
\begin{align*}
\abs*{Y_t} &\le \eta \abs*{\tri*{\grad F(\vecP_t), \vecXi_{t+1}}} 
\le \eta \nrm*{\grad F(\vecP_t)} \nrm*{\vecXi_{t+1}} 
\le \eta \rho_0(F(\vecW_0)+1) \nrm*{\vecXi_{t+1}}.
\end{align*}
Moreover by \pref{ass:noiseregularity} and \pref{lem:boundfuncvalueradius}, 
\[ \nrm*{\vecXi_{t+1}} \le \sigma(F(\vecP_t)) \le \sigma(F(\vecW_0)+1) = \sigma_1(\vecW_0). \]
In particular, recall that $\vecXi_{t+1}$ is the difference between the gradient oracle and actual gradient at $\vecP_t$.

By the above arguments, both of the following inequalities hold deterministically:
\begin{align*}
\abs*{Y_t} &\le \eta \nrm*{\grad F(\vecP_t)}\sigma_1(\vecW_0),\\
\abs*{Y_t} &\le \eta \rho_0(F(\vecW_0)+1) \sigma_1(\vecW_0) = \eta L_0(\vecW_0) \sigma_1(\vecW_0). 
\end{align*} 
We now apply both of these bounds in Data-Dependent Concentration Inequality, \pref{thm:subgaussiandatadependentconcentration} (whose conditions hold because we can assume $\delta, \epsilon$ are at most given universal constants, so $K_0 \ge 4, 2\log K_0/p > e$). Consequently we obtain with probability at least $1-\frac{p}2$ that
\begin{align*}
- \eta \sum_{t=0}^{K_0-1}\tri*{\grad F(\vecP_t), \vecXi_{t+1}} 1_{t < \cK} &\le 2\eta L_0(\vecW_0) \sigma_1(\vecW_0)\log \prn*{\frac{2\log K_0}{p}} \bigvee\\
&4\sqrt{\eta^2 \sigma_1(\vecW_0)^2 \sum_{t=0}^{K_0-1} \nrm*{\grad F(\vecP_t)}^2} \sqrt{\log \prn*{\frac{2\log K_0}{p}}}.\numberthis\label{eq:firstorderupperboundsgdfirstorder}
\end{align*}
Denote this event by $\cE_2$, so $\mathbb{P}(\cE_2) \ge 1-p/2$.
\end{itemize}
For the rest of this proof, we condition on $\cE_1 \cap \cE_2$. By the above, $\cE_1 \cap \cE$ occurs with probability at least $1-p$. Denote $\cE=\cE_1 \cap \cE_2$.

A-priori, these bounds are not particularly useful, especially in our more challenging setting under \pref{ass:gradoraclecontrol} where noise can depend on function value. However conditioned on $\cE$, we prove that SGD is sufficiently `local', in particular that $\nrm*{\vecP_t - \vecU_0} \le 1$ for all $t, 1 \le t \le K_0$. This will then give us control over function value via \pref{lem:boundfuncvalueradius}, which then allow us to make use of these bounds in a more standard way.
\begin{lemma}\label{lem:firstordersgdlocal}
Conditioned on $\cE_1$ (and hence conditioned on $\cE$), for all $t, 1 \le t \le K_0$, we have 
\[ \nrm*{\vecP_t - \vecP_0}=\nrm*{\vecP_t - \vecU_0} \le \frac1{\rho_0(F(\vecW_0)+1)}.\]
\end{lemma}
\begin{proof}
We go by induction on $t$. Notice after $t$ iterates, 
\[ \vecP_t = \vecW_0 - \eta \sum_{k=0}^{t-1} \grad F(\vecP_k) - \eta \sum_{k=1}^{t} \vecXi_k.\]
For the base case $t=1$, we have from \pref{corr:gradcontrol} that $\nrm*{\grad F(\vecW_0)} \le \rho_0(F(\vecW_0)) \le L_0(\vecW_0)$. From the definition of the high-probability event $\cE_1$ and properties of $\eta$ from \pref{claim:sgdproofclaimparams}, and as $\sigma_1(\vecW_0) \ge \sigma(\vecW_0)$), it follows that
\[ \nrm*{\eta \vecXi_1} \le 2\eta \sigma(F(\vecW_0)) \sqrt{K_0 \log\prn*{4K_0/p}}\le \frac1{2\rho_0(F(\vecW_0)+1)}. \] 
Consequently by properties of $\eta$ from \pref{claim:sgdproofclaimparams}, 
\[ \nrm*{\vecP_1-\vecP_0} \le \nrm*{\eta\grad F(\vecW_0)}+\nrm*{\eta \vecXi_0} \le \frac1{\rho_0(F(\vecW_0)+1)}.\]
This finishes the proof of the base case.

Now suppose \pref{lem:firstordersgdlocal} holds for all $1 \le k \le t-1$; we will show it for $t$. From \pref{lem:boundfuncvalueradius}, for all $k \le t-1$, we have 
\[ \nrm*{\grad F(\vecP_{k})} \le \rho_0(F(\vecW_0)+1) \le L_0(\vecW_0). \]
Thus for each $k$, we have
\[ \sigma(F(\vecP_k)) \le \sigma\prn*{F(\vecW_0)+1} = \sigma_1(\vecW_0).\] 
Thus conditioned on $\cE_1$ we obtain
\begin{align*}
\nrm*{\vecP_t-\vecP_0} &\le \nrm*{\eta \sum_{k=0}^{t-1} \grad F(\vecP_k)} + \nrm*{\eta \sum_{k=1}^{t}\vecXi_k} \\
&\le \eta K_0 L_0(\vecW_0) + 2\eta \sqrt{\log\prn*{4K_0/p}\sum_{k=0}^{K_0-1} \sigma_1(\vecW_0)^2} \\
&= \eta K_0 L_0(\vecW_0) + 2\eta \sigma_1(\vecW_0) \sqrt{K_0 \log\prn*{4K_0/p}}\\
&\le \frac1{2\rho_0(F(\vecW_0)+1)} + \frac1{2\rho_0(F(\vecW_0)+1)}=\frac1{\rho_0(F(\vecW_0)+1)}. 
\end{align*}
Here we used the choice of $\eta$ from \pref{claim:sgdproofclaimparams} and the upper bound \pref{eq:noiseupperboundsgdfirstorder} on $\nrm*{\eta \sum_{k=1}^{t}\vecXi_k}$ implied by $\cE_1$. This completes the induction.
\end{proof}
Now that we know the iterates of SGD are `sufficiently local' for $K_0$ iterations via \pref{lem:firstordersgdlocal}, the finish is straightforward. Condition on $\cE$ for the rest of the proof. Consider any $0 \le t \le K_0-1$. $\cE$ implies for all $\vecP \in \overline{\vecP_{t-1}\vecP_t}$, writing $\vecP = \theta \vecP_{t-1}+(1-\theta)\vecP_t$ for $\theta \in [0,1]$, that we have
\[ \nrm*{\vecP - \vecP_0} \le \theta\nrm*{\vecP_{t-1}-\vecP_0}+(1-\theta)\nrm*{\vecP_{t}-\vecP_0} \le \prn*{1-\theta+\theta} \cdot \frac1{\rho_0(F(\vecW_0)+1)}=\frac1{\rho_0(F(\vecW_0)+1)}.\]
Consequently $F(\vecP) \le \rho_0(F(\vecW_0)+1)$, so the above combined with \pref{ass:selfbounding} gives
\[ \nrm*{\grad^2 F(\vecP)} \le L_1(\vecW_0).\numberthis\label{eq:sgdfirstorderhessianbound}\]
We also obtain from \pref{lem:firstordersgdlocal} together with \pref{corr:gradcontrol} and \pref{ass:noiseregularity} that for all $0 \le t \le K_0$,
\begin{align*}
\nrm*{\vecXi_t} &\le \sigma\prn*{F(\vecW_0)+1} = \sigma_1(\vecW_0), \\
\nrm*{\grad F(\vecP_t)} &\le \rho_0\prn*{F(\vecW_0)+1} =L_0(\vecW_0).\numberthis\label{eq:sgdfirstordernoisebound}
\end{align*} 
Now by \pref{lem:secondordersmoothnessineq} and \pref{eq:sgdfirstorderhessianbound}, 
\begin{align*}
F(\vecP_{t+1}) &\le F(\vecP_t) - \eta \tri*{\grad F(\vecP_t), \grad f(\vecP_t; \vecZeta_{t+1})} + \frac{\eta^2 L_1(\vecW_0)}2 \nrm*{\grad f(\vecP_t; \vecZeta_{t+1})}^2\\
&\le F(\vecP_t) - \eta \nrm*{\grad F(\vecP_t)}^2 - \eta \tri*{\grad F(\vecP_t), \vecXi_{t+1})} + \eta^2 L_1(\vecW_0) \prn*{\nrm*{\grad F(\vecP_t)}^2 + \nrm*{\vecXi_{t+1}}^2}.
\end{align*}
The last step uses the definition of $\vecXi_{t+1}$ and Young's Inequality.

Summing and telescoping the above for $0 \le t \le K_0-1$, and applying \pref{eq:sgdfirstordernoisebound}, gives
\begin{align*}
F(\vecP_{K_0}) &\le F(\vecP_0) - \eta \sum_{t=0}^{K_0-1} \nrm*{\grad F(\vecP_t)}^2 - \eta \sum_{t=0}^{K_0-1}\tri*{\grad F(\vecP_t), \vecXi_{t+1}} \\
&\hspace{1in}+ \eta^2 K_0 L_0(\vecW_0)^2 L_1(\vecW_0) + \eta^2 K_0 \sigma^2_1(\vecW_0) L_1(\vecW_0).\numberthis\label{eq:sgdfirstordertelescopeineq}
\end{align*}
Now, conditioned on $\cE$, we upper bound 
\[ - \eta \sum_{t=0}^{K_0-1}\tri*{\grad F(\vecP_t), \vecXi_{t+1}} \]
using \pref{eq:firstorderupperboundsgdfirstorder}. Under $\cE$, by \pref{lem:firstordersgdlocal} and \pref{lem:boundfuncvalueradius}, we have $\vecP_t \in \ball\prn*{\vecP_0, \frac1{\rho_0(F(\vecW_0)+1)}}$ for all $1 \le t \le K_0$, which implies that $t < \cK$ for all $1 \le t \le K_0$. Therefore
\[ - \eta \sum_{t=0}^{K_0-1}\tri*{\grad F(\vecP_t), \vecXi_{t+1}}=- \eta \sum_{t=0}^{K_0-1}\tri*{\grad F(\vecP_t), \vecXi_{t+1}}1_{t<\cK}.\]
Now AM-GM gives
\begin{align*}
&4\sqrt{\eta^2 \sigma_1(\vecW_0)^2\sum_{t=0}^{K_0-1} \nrm*{\grad F(\vecP_t)}^2} \sqrt{\log \prn*{\frac{2\log K_0}{p}}} \\
&\le 2\eta \prn*{\frac14 \sum_{t=0}^{K_0-1} \nrm*{\grad F(\vecP_t)}^2 + 8\sigma_1(\vecW_0)^2 \log \prn*{\frac{2\log K_0}{p}}}.
\end{align*}
Combining with \pref{eq:firstorderupperboundsgdfirstorder}, we obtain
\begin{align*}
- \eta \sum_{t=0}^{K_0-1}\tri*{\grad F(\vecP_t), \vecXi_{t+1}} &= - \eta \sum_{t=0}^{K_0-1}\tri*{\grad F(\vecP_t), \vecXi_{t+1}}1_{t<\cK} \\
&\le \frac{\eta}2\sum_{t=0}^{K_0-1} \nrm*{\grad F(\vecP_t)}^2 +16\eta  B(\vecW_0)\log \prn*{\frac{2\log K_0}{p}}.
\end{align*}
Combining with \pref{eq:sgdfirstordertelescopeineq} gives
\begin{align*}
F(\vecP_{K_0}) &\le F(\vecP_0) - \frac{\eta}2 \sum_{t=0}^{K_0-1} \nrm*{\grad F(\vecP_t)}^2 + 16\eta  B(\vecW_0)\log \prn*{\frac{2\log K_0}{p}}+ \eta^2 K_0 L_0(\vecW_0)^2 L_1(\vecW_0) \\
&\hspace{1in}+ \eta^2 K_0 \sigma^2_1(\vecW_0) L_1(\vecW_0).\numberthis\label{eq:firstordersgdprefinal}
\end{align*}
Suppose that $\nrm*{\grad F(\vecP_t)} > \epsilon$ for all $0 \le t \le K_0-1$. Then the above gives
\begin{align*}
F(\vecP_{K_0}) &< F(\vecP_0) - \frac{\eta K_0 \epsilon^2}2 + 16\eta  B(\vecW_0)\log \prn*{\frac{2\log K_0}{p}}\\
&\hspace{1in}+ \eta^2 K_0 L_0(\vecW_0)^2 L_1(\vecW_0)+ \eta^2 K_0 \sigma^2_1(\vecW_0) L_1(\vecW_0). 
\end{align*}
To make use of this bound, by our choice of $\eta$, \pref{claim:sgdproofclaimparams} implies that
\[ \eta^2 K_0 L_0(\vecW_0)^2 L_1(\vecW_0)+ \eta^2 K_0 \sigma^2_1(\vecW_0) L_1(\vecW_0) \le \frac{\eta K_0 \epsilon^2}8.\]
By choice of $K_0$, \pref{claim:sgdproofclaimparams} implies that
\[ 16 \eta B(\vecW_0) \log \prn*{\frac{2\log K_0}{p}} \le \frac{\eta K_0 \epsilon^2}{8}. \]
The above was all conditioned on $\cE$, which occurred with probability at least $1-p$. Thus by \pref{eq:firstordersgdprefinal}, we obtain that with this same probability which is at least $1-p$, if none of $\vecP_0, \ldots, \vecP_{K_0-1}$ have gradient norm larger than $\epsilon$, we have
\[ F(\vecP_{K_0}) < F(\vecP_0) - \frac{\eta K_0 \epsilon^2}4=F(\vecU_0) - \frac{\eta K_0 \epsilon^2}4.\]
This establishes that $\cA$ is a $(\cS, K_0+1, \frac{\eta K_0 \epsilon^2}{4}, p, \vecU_0)$-decrease procedure. Following our initial observations, we conclude via \pref{thm:generalframework}.
\end{proof}
Now we prove \pref{claim:sgdproofclaimparams}.
\begin{proof}[Proof of \pref{claim:sgdproofclaimparams}]
We first prove \pref{eq:K0paramchocesgdfirstorder}. Recall we chose 
\[ K_0 = \frac{C(\vecW_0)}{\epsilon^2} \log(1/\tilde{\eta})^2 \log(1/\delta)^2 \log(1/\epsilon)^2.\]
Furthermore recall $p = \frac{\delta \tilde{\eta} K_0 \epsilon^2}{4(F(\vecW_0)+1)}$. Thus, \pref{eq:K0paramchocesgdfirstorder} holds if and only if 
\[ C(\vecW_0) \log(1/\tilde{\eta})^2 \log(1/\delta)^2 \log(1/\epsilon)^2 \ge 128 B(\vecW_0) \log\prn*{\frac{8\log K_0 \cdot (F(\vecW_0)+1)}{\delta \tilde{\eta} K_0 \epsilon^2}}. \]
As $C(\vecW_0) \ge 128 B(\vecW_0) \lor 64(F(\vecW_0)+1)^2$, again using the expression for $K_0$, it suffices to prove
\[ \log(1/\tilde{\eta})^2 \log(1/\delta)^2 \log(1/\epsilon)^2 \ge \log\prn*{\frac{\log K_0 }{C(\vecW_0)^{1/2} \delta \tilde{\eta} \log(1/\tilde{\eta})^2 \log(1/\delta)^2}}.\]
As $\log(1/\delta)$, $\log(1/\tilde{\eta})$ are both larger than 1, it suffices to prove 
\begin{align*}
&\log(1/\tilde{\eta})^2 \log(1/\delta)^2 \log(1/\epsilon)^2 \\
&\ge \log(1/\tilde{\eta})+\log(1/\delta) \\
&+\log\prn*{\frac{\log C(\vecW_0) + \log(1/\epsilon^2) + 2\log \log (1/\tilde{\eta})+ 2\log \log (1/\delta)+ 2\log \log (1/\epsilon)}{C(\vecW_0)^{1/2}}}.
\end{align*}
Since $C(\vecW_0) \ge 64$, it satisfies $\log C(\vecW_0) < C(\vecW_0)^{1/2}$, so it suffices to prove 
\begin{align*}
&\log(1/\tilde{\eta})^2 \log(1/\delta)^2 \log(1/\epsilon)^2 \\
&\ge \log(1/\tilde{\eta})+\log(1/\delta) \\
&+ \log\prn*{1+2\log(1/\epsilon) + 2\log \log (1/\tilde{\eta})+ 2\log \log (1/\delta)+ 2\log \log (1/\epsilon)}.
\end{align*}
By comparing `degrees', we conclude recalling we can assume WLOG that $\delta,\epsilon,\tilde{\eta}$ are smaller than some universal constant.

Now we prove \pref{eq:etaparamchocesgdfirstorder}. We will prove that 
\[ \tilde{\eta} \le \frac1{ \tilde{L}'(\vecW_0) K_0 \sqrt{\log(4K_0/p)}}.\numberthis\label{eq:newetaparamchocesgdfirstorder}\]
After proving \pref{eq:newetaparamchocesgdfirstorder}, recalling our choice of $K_0 > 1/\epsilon^2$ directly implies \pref{eq:etaparamchocesgdfirstorder}. To show \pref{eq:newetaparamchocesgdfirstorder}, equivalently, we want to show 
\[ \tilde{\eta} \log(1/\tilde{\eta})^2 \sqrt{\log(4K_0/p)} \le \frac{\epsilon^2}{\tilde{L}'(\vecW_0) C(\vecW_0) \log(1/\delta)^2 \log(1/\epsilon)^2}.\]
Recalling the definition of $p$, this holds if and only if 
\[ \tilde{\eta} \log(1/\tilde{\eta})^2 \sqrt{\log\prn*{\frac{16(F(\vecW_0)+1)}{\delta \tilde{\eta} \epsilon^2}}}\le\frac{\epsilon^2}{\tilde{L}'(\vecW_0) C(\vecW_0) \log(1/\delta)^2 \log(1/\epsilon)^2}.\]
Now we explicitly recall our expression for $\tilde{\eta} = \frac{\epsilon^2}{\tilde{L}(\vecW_0) \log(1/\epsilon)^6 \log(1/\delta)^6}$. Plugging this in and recalling $\tilde{L}(\vecW_0) \ge \tilde{L}'(\vecW_0)^2 C(\vecW_0)^2$, it suffices to prove 
\begin{align*}
&\frac1{\tilde{L}(\vecW_0)^{1/2} \log(1/\epsilon)^6 \log(1/\delta)^6} \log\prn*{\frac{\tilde{L}(\vecW_0) \log(1/\epsilon)^6 \log(1/\delta)^6}{\epsilon^2}}^2 \\
&\cdot \sqrt{\log\prn*{\frac{16(F(\vecW_0)+1) \tilde{L}(\vecW_0) \log(1/\epsilon)^6 \log(1/\delta)^6}{\delta \epsilon^4}}} \\
&\le \frac{1}{\log(1/\delta)^2 \log(1/\epsilon)^2}.
\end{align*}
Thus it suffices to prove:
\begin{align*}
&\frac{18}{\tilde{L}(\vecW_0)^{1/2}} \log\prn*{\frac{\tilde{L}(\vecW_0) \log(1/\epsilon) \log(1/\delta)}{\epsilon}}^2 \sqrt{\log\prn*{\frac{16(F(\vecW_0)+1) \tilde{L}(\vecW_0) \log(1/\epsilon) \log(1/\delta)}{\delta \epsilon}}} \\
&\le \log(1/\delta)^4 \log(1/\epsilon)^4.
\end{align*}
Recall $\tilde{L}(\vecW_0)^{1/8} \ge 3\sqrt{2} \log(\tilde{L}(\vecW_0)) \lor 3\sqrt{2}$ and so 
\begin{align*}
&\frac{3\sqrt{2}}{\tilde{L}(\vecW_0)^{1/4}}\log\prn*{\frac{\tilde{L}(\vecW_0) \log(1/\epsilon) \log(1/\delta)}{\epsilon}} \\
&\le \frac{3\sqrt{2}}{\tilde{L}(\vecW_0)^{1/4}} \prn*{\log(1/\epsilon)+\log \log (1/\epsilon)+\log \log (1/\delta) + \log \tilde{L}(\vecW_0)} \\
&\le \frac1{\tilde{L}(\vecW_0)^{1/8}}\prn*{1+\log(1/\epsilon)+\log \log (1/\epsilon)+\log \log (1/\delta)}.
\end{align*}
Thus it suffices to show 
\begin{align*}
&\frac1{\tilde{L}(\vecW_0)^{1/4}}\prn*{1+\log(1/\epsilon)+\log \log (1/\epsilon)+\log \log (1/\delta)}^2 \\
&\cdot \sqrt{\log\prn*{\frac{16(F(\vecW_0)+1) \tilde{L}(\vecW_0) \log(1/\epsilon) \log(1/\delta)}{\delta \epsilon}}} \\
&\le \log(1/\delta)^4 \log(1/\epsilon)^4.
\end{align*}
To this end recall $\tilde{L}(\vecW_0)^{1/8} \ge \log (16(F(\vecW_0)+1) \tilde{L}(\vecW_0))$, thus
\begin{align*}
&\frac1{\tilde{L}(\vecW_0)^{1/8}} \log\prn*{\frac{16(F(\vecW_0)+1) \tilde{L}(\vecW_0) \log(1/\epsilon) \log(1/\delta)}{\delta \epsilon}} \\
&=\frac1{\tilde{L}(\vecW_0)^{1/8}} \prn*{\log (16(F(\vecW_0)+1) \tilde{L}(\vecW_0)) + \log(1/\delta)+\log(1/\epsilon) + \log \log (1/\delta))+ \log \log (1/\epsilon)} \\
&\le 1+\log(1/\delta)+\log(1/\epsilon) + \log \log (1/\delta))+ \log \log (1/\epsilon).
\end{align*}
Therefore it suffices to show 
\begin{align*}
&\prn*{1+\log(1/\epsilon)+\log \log (1/\epsilon)+\log \log (1/\delta)}^2 \\
&\cdot \prn*{1+\log(1/\delta)+\log(1/\epsilon) + \log \log (1/\delta))+ \log \log (1/\epsilon)}^{1/2} \\
&\le \log(1/\delta)^4 \log(1/\epsilon)^4.
\end{align*}
Evidently the above holds for small enough universal constants $\delta,\epsilon$ (compare `degrees'), so we conclude the proof.
\end{proof}
\begin{remark}\label{rem:sgdfirstordersubgaussian}
We also discuss how to extend this result to when the $\nrm*{\vecXi_t}$ has sub-Gaussianity parameter $\sigma(F(\vecP_t))$. The extension is straightforward. Again, we aim to prove \pref{claim:sgdfirstorderdecreasesequence}. For the rest of this remark, follow the notation from the proof for SGD above. Besides applying \pref{thm:subgaussianazumahoeffding}, \pref{thm:subgaussiandatadependentconcentration} when the relevant random variables are sub-Gaussian, which still hold true as mentioned in \citet{fang2019sharp}, the only other time we used that $\nrm*{\vecXi_t} \le \sigma(F(\vecP_t))$ holds deterministically is to derive \pref{eq:sgdfirstordertelescopeineq}. 

We apply \pref{thm:subgaussianazumahoeffding}, \pref{thm:subgaussiandatadependentconcentration} identically to the proof earlier. This time, we have for $t < \cK$ that $\vecXi_{t+1}$ is sub-Gaussian with parameter $\sigma_1(\vecW_0)$, thanks to the same trick of multiplying with $1_{t<\cK}$ when applying \pref{thm:subgaussiandatadependentconcentration}.

The only change is as follows: in the definition $\cE$, add in the intersection the event $\cE_3$ that for all $1 \le t \le K_0$, $\nrm*{\vecXi_t}^2 \le \sigma(F(\vecP_t))^2 \log(K_0/p)$, where $p$ is defined the same as before. We control the probability of $\cE_3$ via the following Lemma:
\begin{lemma}[Equivalent of Lemma 12, \citet{priorpaper}]\label{lem:sgdsubgaussianlem}
With probability at least $1-p$, we have for all $1 \le t \le K_0$,
\[ \nrm*{\vecXi_t}^2 \le \sigma(F(\vecP_t))^2 \log(K_0/p).\]
\end{lemma}
\begin{proof}
By \pref{ass:noiseregularity}, with probability $1-\frac{p}{K_0}$, we have 
\[\frac{\nrm*{\vecXi_t}^2}{\sigma(F(\vecP_t))^2} \le \log(K_0/p). \]
A Union Bound finishes the proof.
\end{proof}
Now we condition on $\cE=\cE_1 \cap \cE_2 \cap \cE_3$, which has probability at least $1-2p$ by combining our earlier argument with \pref{lem:sgdsubgaussianlem}. Note this only changes the resulting guarantee by a universal constant. We still have \pref{lem:firstordersgdlocal}, which does not require an upper bound on \textit{each} $\nrm*{\vecXi_t}$ in its proof but simply uses concentration from event $\cE_1$. 

Thus, conditioned on $\cE$, we still have $F(\vecP_t) \le F(\vecW_0)+1$ by \pref{lem:firstordersgdlocal}, \pref{lem:boundfuncvalueradius}, and as $\vecU_0 \in \cL_{F,F(\vecW_0)}$. Now conditioned on $\cE$, by  \pref{lem:sgdsubgaussianlem}, we still have the following upper bound for all $1 \le t \le K_0$:
\[ \nrm*{\vecXi_t}^2 \le \sigma(F(\vecW_0)+1)^2 \log(K_0 / p) = \sigma_1(\vecW_0)^2 \log(K_0/p).\]
Therefore conditioned on $\cE$, we can still derive a bound analogous to \pref{eq:sgdfirstordertelescopeineq}. This resulting bound changes by only a $\log(K_0/p)$ factor (from \pref{lem:sgdsubgaussianlem}, see the above display); moreover recall $K_0, p$ depend polynomially in $\delta, 1/\epsilon$. By adjusting $\eta$ smaller by a $\polylog(K_0/p)$ factor, the same proof as above goes through, up to changing quantities by polylogarithmic factors.
\end{remark}

\section{Perturbed GD finding Second Order Stationary Points}\label{sec:perturbedgdproofs}
\subsection{Proof using the Framework}
Here we prove \pref{thm:escapesecondordergd}. We instantiate \pref{alg:perturbedgd} formally here. The parameters of \pref{alg:perturbedgd} will depend on $L_1(\vecW_0), L_2(\vecW_0)$, which are defined in \pref{eq:L1def}, \pref{eq:L2defperturbedgd} respectively, and depend only on $\rho_1, \rho_2, F(\vecW_0)$. Given a desired success probability $1-\delta$ for $\delta>0$, a tolerance $\epsilon>0$, and $F(\vecW_0), L_1(\vecW_0), L_2(\vecW_0)$, the algorithm's other parameters are defined in terms of as follows:
\begin{enumerate}
    \item $c \le \cmax$ is a universal constant, where $\cmax$ is a universal constant defined in \pref{lem:escapesaddlepointgdfinal}.
    \item $\tilde{\epsilon}= \frac{\epsilon}{L_2(\vecW_0)}$.
    \item $\chi \leftarrow 4\max\crl*{\log\prn*{\frac{2d L_1(\vecW_0)^2 F(\vecW_0) }{c^2\tilde{\epsilon}^{2.5} \delta}}, 5}$.
    \item $\eta \leftarrow \frac{c}{L_1(\vecW_0)}$. 
    \item $r \leftarrow \frac{\sqrt{c}\tilde{\epsilon}}{\chi^2 L_1(\vecW_0)}$.
    \item $\Gthres \leftarrow \frac{\sqrt{c}}{\chi^2}\tilde{\epsilon}$.
    \item $\Fthres\leftarrow \frac{c}{\chi^3}\sqrt{\frac{\tilde{\epsilon}^3}{L_2(\vecW_0)}}$.
    \item $\Tthres \leftarrow \frac{\chi}{c^2} \frac{L_1(\vecW_0)}{\sqrt{L_2(\vecW_0)\tilde{\epsilon}}}$
\end{enumerate}

\begin{algorithm}[tb]
   \caption{Perturbed Gradient Descent, modified from \citet{jin2017escape}. }
   \label{alg:perturbedgd}
\begin{algorithmic}
   \STATE $\tilde{\epsilon}= \frac{\epsilon}{L_2(\vecW_0)}$, $\chi \leftarrow 4\max\crl*{\log\prn*{\frac{2d L_1(\vecW_0)^2 F(\vecW_0)}{c\tilde{\epsilon}^{2.5} \delta}}, 5}$, $\eta \leftarrow \frac{c}{L_1(\vecW_0)}$, $r \leftarrow \frac{\sqrt{c}\tilde{\epsilon}}{\chi^2 L_1(\vecW_0)}$, $\Gthres \leftarrow \frac{\sqrt{c}}{\chi^2}\tilde{\epsilon}$, $\Fthres\leftarrow \frac{c}{\chi^3}\sqrt{\frac{\tilde{\epsilon}^3}{L_2(\vecW_0)}}$, $\Tthres \leftarrow \frac{\chi}{c^2} \frac{L_1(\vecW_0)}{\sqrt{L_2(\vecW_0)\tilde{\epsilon}}}$. Here $c$ refers to a small enough universal constant upper bounded by $\cmax$ in \pref{lem:escapesaddlepointgdfinal}.  
   \STATE 
   \WHILE{True}
   \IF{$\nrm*{\grad F(\vecW_t)} \le \Gthres$}
   \STATE $\tilde{\vecW}_t \leftarrow \vecW_t$, $\Tnoise \leftarrow t$
   \STATE $\vecW_t \leftarrow \tilde{\vecW}_t + \vecXi_t$, $\vecXi_t$ uniform from $\ball(\vecOrigin, r)$
   \STATE $s\leftarrow 0$
        \WHILE{$s<\Tthres$}
        \STATE $\vecW_{t+1}=\vecW_t-\eta \grad F(\vecW_t)$, $s \leftarrow s+1$, $t \leftarrow t+1$
        \ENDWHILE
        \IF{$F(\vecW_{t})-F(\tilde{\vecW}_{\Tnoise}) > -\Fthres$}
        \STATE Return $\tilde{\vecW}_{\Tnoise}$
        \ENDIF
   \ELSE{}
   \STATE $\vecW_{t+1}=\vecW_t-\eta \grad F(\vecW_t)$, $t \leftarrow t+1$
   \ENDIF
   \ENDWHILE
\end{algorithmic}
\end{algorithm}

\begin{proof}[Proof of \pref{thm:escapesecondordergd} given \pref{lem:escapesaddlepointgdfinal}]
We will first prove the following Lemma, which will define $L_2(\vecW_0)$ and explain its significance.
\begin{lemma}\label{lem:formalizingperturbationopnrom}
Define $L_1(\vecW_0)$ as in \pref{eq:L1def}, and define
\begin{align*}
L_2(\vecW_0) &= \max\crl*{1,L_1(\vecW_0), \rho_2\prn*{F(\vecW_0)+1}}.\numberthis\label{eq:L2defperturbedgd}
\end{align*}
Then we have the following:
\begin{enumerate}
    \item Suppose $\vecU$ is such that $\nrm*{\vecU-\tilde{\vecW}} \le \frac1{\rho_0(F(\vecW_0)+1)}$, where $\tilde{\vecW} \in \cL_{F, F(\vecW_0)}$, the $F(\vecW_0)$-sublevel set. Then under \pref{ass:selfbounding} (and in particular under \pref{ass:thirdorderselfbounding}),
    \[ \nrm*{\grad^2 F(\vecU)}_{\OPNORM} \le L_1(\vecW_0).\]
    \item Suppose that $\vecU_1, \vecU_2$ are such that $\nrm*{\vecU_1-\tilde{\vecW}}, \nrm*{\vecU_2-\tilde{\vecW}} \le \frac1{\rho_0(F(\vecW_0)+1)}$, where $\tilde{\vecW} \in \cL_{F, F(\vecW_0)}$. Then 
    \[\nrm*{\grad^2 F(\vecU_1) - \grad^2 F(\vecU_2)}_{\OPNORM} \le L_2(\vecW_0)\nrm*{\vecU_1 - \vecU_2}.\]
\end{enumerate}
\end{lemma}
\begin{remark}
Note $L_1(\vecW_0), L_2(\vecW_0) \ge 1$, and $L_2(\vecW_0) \ge L_1(\vecW_0)$.
\end{remark}
\begin{proof}[Proof of \pref{lem:formalizingperturbationopnrom}]
Recall by \pref{corr:gradcontrol} that $\nrm*{\grad F(\vecW)} \le \rho_0(F(\vecW))$. Now by \pref{lem:boundfuncvalueradius} and as $\tilde{\vecW} \in \cL_{F, F(\vecW_0)}$, for any $\vecU'$ with $\nrm*{\vecU' - \tilde{\vecW}} \le \frac1{\rho_0(F(\vecW_0)+1)} \le \frac1{\rho_0(F(\tilde{\vecW})+1)}$, we have $F(\vecU') \le F(\tilde{\vecW}) + 1$.
The first part now directly follows by \pref{ass:selfbounding}. 

The second part now follows by noting the line segment $\overline{\vecU_1 \vecU_2}$ is contained in $\ball\prn*{\tilde{\vecW}, \frac1{\rho_0(F(\vecW_0)+1)}}$ via Triangle Inequality, recalling $\tilde{\vecW} \in \cL_{F, F(\vecW_0)}$, and then applying \pref{lem:thirdorderlipschitzlem} and \pref{lem:boundfuncvalueradius}. 
\end{proof}
We now prove \pref{thm:escapesecondordergd} by instantiating our framework. 

Define $\tilde{\epsilon}=\frac{\epsilon}{L_2(\vecW_0)}$ as we did earlier, and note $L_2(\vecW_0) \ge 1$. It suffices to show for $\tilde{\epsilon} \le 1$, that with probability at least $1-\delta$, we will return $\vecW$ such that $\nrm*{\grad F(\vecW)} \le \tilde{\epsilon}$, $\grad^2 F(\vecW) \succeq -\sqrt{L_2(\vecW_0) \tilde{\epsilon}} \matI$ in $T = O\prn*{\frac{L_1(\vecW_0) \max\crl*{F(\vecW_0),1} \chi^4}{\tilde{\epsilon}^2}}=O\prn*{\frac{L_1(\vecW_0) L_2(\vecW_0)^2 \max\crl*{F(\vecW_0),1} \chi^4}{\epsilon^2}}$ oracle calls.\footnote{The $\max\crl*{1,F(\vecW_0)}$ is a proof artifact.}

Now let the set of interest
\[ \cS = \{\vecW:\nrm*{\grad F(\vecW)} \le \Gthres, \grad^2 F(\vecW) \succeq -\sqrt{L_2(\vecW_0) \tilde{\epsilon}} \matI\}.\]
Note $\Gthres \le \tilde{\epsilon}$, so $\vecW \in \cS$ immediately implies $\nrm*{\grad F(\vecW)} \le \tilde{\epsilon}$, $\grad^2 F(\vecW) \succeq -\sqrt{L_2(\vecW_0) \tilde{\epsilon}} \matI$. Also note it suffices to show the result for all $\tilde{\epsilon} \le \frac1{100L_2(\vecW_0)}$; otherwise for larger $\tilde{\epsilon}$ we can just apply the result for $\tilde{\epsilon} = \frac1{100L_2(\vecW_0)}$. Thus as $L_2(\vecW_0) \ge 1$, we can assume $\tilde{\epsilon} \le 1$. Clearly, we also can assume WLOG that $\Tthres \ge 1$.

As in \pref{subsec:exinframework}, we make the following definitions for \pref{alg:perturbedgd}. For all $\vecU_0 \in \mathbb{R}^d$, if $\nrm*{\grad F(\vecU_0)} > \Gthres$, we let 
\[ \cA(\vecU_0) = \prn*{\vecU_0 - \eta \grad F(\vecU_0), \vecU_0}, \text{ hence }\cA_1(\vecU_0)=\vecU_0 - \eta \grad F(\vecU_0), \cA_2(\vecU_0) = \vecU_0.\]
Otherwise if $\nrm*{\grad F(\vecU_0)} \le \Gthres$, we let $\vecP_0 = \vecU_0 + \vecXi$ where $\vecXi$ is uniform from $\ball(\vecOrigin, r)$, and define a sequence $(\vecP_i)_{0 \le i \le \Tthres}$ via
\[ \vecP_i = \vecP_{i-1} - \eta \grad F(\vecP_{i-1}).\]
When then take 
\[ \cA(\vecU_0) = \prn*{\vecP_{\Tthres}, \vecU_0}, \text{ hence }\cA_1(\vecU_0)=\vecP_{\Tthres}, \cA_2(\vecU_0) = \vecU_0.\]
We then have 
\[ \Toracle(\vecU_0) = \begin{cases} \Tthres &: \nrm*{\grad F(\vecU_0)} \le \Gthres \\ 1 &:\nrm*{\grad F(\vecU_0)} >\Gthres. \end{cases}\]
We also define
\[  \Delta(\vecU_0) = \begin{cases}  \Fthres &: \nrm*{\grad F(\vecU_0)} \le \Gthres \\ \frac{\eta}2 \cdot \Gthres^2 &:\nrm*{\grad F(\vecU_0)} >\Gthres. \end{cases}\]
We now establish the crucial \pref{claim:perturbedgddecreasesequence}: for all $\vecU_0 \in \cL_{F,F(\vecW_0)}$, $\cA$ is a $(\cS, \Toracle(\vecU_0), \Delta(\vecU_0), \frac{d L_1(\vecW_0)}{\sqrt{L_2(\vecW_0)\tilde{\epsilon}}} e^{-\chi}, \vecU_0)$-decrease procedure. (Recall $\tilde{\epsilon} = \frac{\epsilon}{L_2(\vecW_0)}$.)

To do this, we use the following crucial Lemma ensuring high-probability decrease around saddle points \textit{in the $F(\vecW_0)$-sublevel set}:
\begin{lemma}[Equivalent of Lemma 13, \citet{jin2017escape}]\label{lem:escapesaddlepointgdfinal}
There exists a universal constant $\cmax \le 1$ such that the following occurs. Suppose we start with a $\tilde{\vecW} \in \cL_{F, F(\vecW_0)}$, that is in the $F(\vecW_0)$-sublevel set, satisfying the following conditions:
\[
\nrm*{\grad F(\tilde{\vecW})} \leq \Gthres \quad \text{and} \quad \lambda_{\min}(\nabla^2 F(\tilde{\vecW})) \le -\sqrt{L_2(\vecW_0)\tilde{\epsilon}}.
\]
Now let $\vecP_0 = \tilde{\vecW} + \vecZeta$, where $\vecZeta$ is sampled uniformly from $\ball(\vecOrigin, r)$ where $r$ is defined in \pref{lem:escapesaddlepoints}, and let $\{\vecP_t\}$ be the iterates of gradient descent starting from $\vecP_0$. Then when the step size $\eta \leq \frac{\cmax}{L_1(\vecW_0)}$, with probability at least $1 - \frac{d L_1(\vecW_0)}{\sqrt{L_2(\vecW_0)\tilde{\epsilon}}} e^{-\chi}$, we have:
\[
F(\vecP_{\Tthres}) - F(\tilde{\vecW}) < -\Fthres.
\]
\end{lemma}
The variables in the above are defined in \pref{alg:perturbedgd}. As noted earlier, because we work in the generalized smooth setting, the details require significant care compared to the proof of Lemma 13 in \citet{jin2017escape}. 

With \pref{lem:escapesaddlepointgdfinal}, we have the ingredients to prove \pref{thm:escapesecondordergd}. First we establish \pref{claim:perturbedgddecreasesequence}.
\begin{proof}[Proof of \pref{claim:perturbedgddecreasesequence}]
We prove this by breaking into the following cases:
\begin{itemize}
    \item Suppose $\nrm*{\grad F(\vecU_0)} > \Gthres$. Then $\vecU_1 = \cA_1(\vecU_0) = \vecU_0-\eta \grad F(\vecU_0)$.

Our condition on $\eta$ implies that 
\[ \eta \le \frac1{L_1(\vecW_0)} \le \frac1{\rho_0(F(\vecW_0)) \rho_0(F(\vecW_0)+1)}.\] 
As $\vecU_0  \in \cL_{F, F(\vecW_0)}$, we have by \pref{corr:gradcontrol},
\[ \nrm*{\vecU_1-\vecU_0} = \eta \nrm*{\grad F(\vecU_0)} \le \eta \rho_0(F(\vecU_0)) \le \eta \rho_0(F(\vecW_0)) \le \frac1{\rho_0(F(\vecW_0)+1)}. \]
Consequently, by \pref{lem:boundfuncvalueradius}, 
\[ F(\vecP) \le F(\vecU_0)+1\le F(\vecW_0)+1 \text{ for  all } \vecP \in \overline{\vecU_0 \vecU_1}.\]
Now by \pref{lem:secondordersmoothnessineq} and \pref{ass:selfbounding},
\begin{align*}
F(\vecU_1) &\le F(\vecU_0) - \eta \nrm*{\grad F(\vecU_0)}^2 + \frac{L_1(\vecW_0) \eta^2}2  \nrm*{\grad F(\vecU_0)}^2 \\
&\le F(\vecU_0) - \frac{\eta}2\nrm*{\grad F(\vecU_0)}^2\\
&< F(\vecU_0) - \frac{\eta}2 \cdot \Gthres^2 = F(\vecU_0) - \Delta(\vecU_0).
\end{align*}

    \item Else suppose $\nrm*{\grad F(\vecU_0)} \le \Gthres$. Then $\vecU_0$ is perturbed, and we consider the sequence of the next $\Tthres$ iterates $\vecP_0=\vecU_0+\vecXi, \vecP_1, \ldots, \vecP_{\Tthres}$.
    
Consider the event $\cE$ from  \pref{lem:escapesaddlepointgdfinal}, which occurs with probability at least $1-\frac{dL_1(\vecW_0)}{\sqrt{L_2(\vecW_0)\tilde{\epsilon}}} e^{-\chi}$. 
Under $\cE$, for such $\vecU_0$, we have:
\begin{itemize}
    \item Either
    \[ F(\vecP_{\Tthres}) - F(\vecU_0) < -\Fthres,\]
    that is 
    \[ F(\vecU_1) = F(\vecP_{\Tthres}) < F(\vecU_0) - \Fthres.\]
    \item Or
    \[ \lambdamin\prn*{\grad^2 F(\vecU_0)} \ge -\sqrt{\tilde{\epsilon} L_2(\vecW_0)}, \text{ hence }\vecU_0 \in \cS.\]
\end{itemize}
\end{itemize}
In all cases, by definition of $\cA$, we conclude that $\cA$ is a $(\cS, \Toracle(\vecU_0), \Delta(\vecU_0), \frac{d L_1(\vecW_0)}{\sqrt{L_2(\vecW_0)\tilde{\epsilon}}} e^{-\chi}, \vecU_0)$ decrease procedure for $\vecU_0 \in \cL_{F,F(\vecW_0)}$.
\end{proof}
Consider these two cases, and recall the definition of $\overline{\Delta}$ from \pref{thm:generalframework}. Using the definition of $\eta, \Gthres, \Fthres$, we obtain for $c$ a small enough universal constant,
\begin{align*}
\overline{\Delta} &\ge \frac12 \min \crl*{\frac{c^2 \tilde{\epsilon}^2}{2L_1(\vecW_0) \chi^4}, \frac{c^3 \tilde{\epsilon}^2}{\chi^4 L_1(\vecW_0)}} \\
&\ge \frac{c^3 \tilde{\epsilon}^2}{\chi^4 L_1(\vecW_0)}.
\end{align*}
Combining with \pref{thm:generalframework}, and note $\Toracle(\vecU_0) \le \Tthres\le \frac{\max\crl*{1,F(\vecW_0)}}{\overline{\Delta}}$ for $\tilde{\epsilon} \le 1$. We thus obtain the desired oracle complexity of $O\prn*{\frac{L_1(\vecW_0) \max\crl*{F(\vecW_0),1} \chi^4}{\tilde{\epsilon}^2}}=O\prn*{\frac{L_1(\vecW_0) L_2(\vecW_0)^2 \max\crl*{F(\vecW_0),1} \chi^4}{\epsilon^2}}$ to obtain an iterate in $\cS$.\footnote{Note $\Tthres$ generally does not decrease with $F(\vecW_0)$, and this is why the $\max\crl*{1,F(\vecW_0)}$ comes in.}

We finally show the desired probability of success. Through \pref{thm:generalframework}, since $\chi \ge 18$ and by definition of $\chi$, we can verify that the probability of failure is at most
\begin{align*}
&\frac{dL_1(\vecW_0)}{\sqrt{L_2(\vecW_0)\tilde{\epsilon}}} e^{-\chi} \cdot \sup_{\vecU \in \cL_{F, F(\vecW_0)}}\crl*{\frac{F(\vecW_0)}{\Delta(\vecW)}}\\
&\le\frac{dL_1(\vecW_0)}{\sqrt{L_2(\vecW_0)\tilde{\epsilon}}} e^{-\chi} \cdot \frac{F(\vecW_0)}{\frac{c^2 \tilde{\epsilon}^2}{2\chi^4 L_1(\vecW_0)\sqrt{L_2(\vecW_0)}} } \\
&\le \chi^4 e^{-\chi} \frac{2d L_1(\vecW_0)^2 F(\vecW_0)}{c^2 \tilde{\epsilon}^{2.5}} \\
&\le e^{-\chi/4} \cdot \frac{2 F(\vecW_0) d L^2_1(\vecW_0)}{c\tilde{\epsilon}^{2.5}}\\
&\le \delta.
\end{align*}
This completes the proof, assuming \pref{lem:escapesaddlepointgdfinal}.
\end{proof}

\subsection{Proving the key Lemma}
We now prove \pref{lem:escapesaddlepointgdfinal} to complete the proof. The rest of the proof is similar to that of \citet{jin2017escape}, but hinges crucially on the fact that the analysis in \citet{jin2017escape} is `local'.

Consider any $\gamma > 0$, and define the `units' in a similar way as \citet{jin2017escape}, but now in terms of $L_1(\vecW_0), L_2(\vecW_0) > 0$ defined earlier. First let the new `condition number' be $\kappa = \kappa(\vecW_0) := \frac{L_1(\vecW_0)}{\gamma}$ (note this is \textit{not} the real condition number, but rather is the `effective condition number' of $\grad^2 F$ in $\cL_{F,F(\vecW_0)}$). Now define the following positive reals:
\begin{align*}
\cF_1 &= \eta L_1(\vecW_0) \frac{\gamma^3}{L_2(\vecW_0)^2} \log^{-3}\prn*{\frac{d \kappa}{\delta}},\\
\cF_2 &= \frac{\log\prn*{\frac{d \kappa}{\delta}}}{\eta \gamma},\\
\cG &= \sqrt{\eta L_1(\vecW_0) } \frac{\gamma^2}{L_2(\vecW_0)}\log^{-2}\prn*{\frac{d \kappa}{\delta}},\\
\cL &= \sqrt{\eta L_1(\vecW_0)} \frac{\gamma}{L_2(\vecW_0)} \log^{-1}\prn*{\frac{d \kappa}{\delta}}.
\end{align*}
Our goal is to prove the following. 
\begin{lemma}[equivalent of Lemma 14 in \citet{jin2017escape}]\label{lem:escapesaddlepoints}
There exists a universal constant $\cmax$ such that the following holds. For any $F$ satisfying the conditions of \pref{thm:escapesecondordergd}, for any $\delta \in \left(0, \frac{d\kappa}{e}\right]$, suppose we start with a point $\tilde{\vecW} \in \cL_{F, F(\vecW_0)}$ satisfying the following conditions for some $\gamma>0$, where $\cG$ is defined as above:
\[
\|\nabla F(\tilde{\vecW})\| \leq \cG \quad \text{and} \quad \lambda_{\min}(\nabla^2 F(\tilde{\vecW})) \leq -\gamma.
\]
Let $\vecP_0 = \tilde{\vecW} + \vecZeta$, where $\vecZeta$ is sampled from the uniform distribution over a ball with radius $\frac{\cL}{\kappa \cdot \log\left(\frac{d\kappa}{\delta}\right)} := r$ and where $\cL$ is defined as above. Let $\{\vecP_t\}$ be the iterates of gradient descent starting from $\vecP_0$. Then, when the step size $\eta \leq \frac{\cmax}{L_1(\vecW_0)}$, with probability at least $1 - \delta$, we have the following for any $T \geq \frac{1}{\cmax}\cF_2$:
\[
F(\vecP_T) - F(\tilde{\vecW}) < -\cF_1.
\]
\end{lemma}
Plugging in $\gamma = \sqrt{L_2(\vecW_0) \tilde{\epsilon}}$, $\eta = \frac{\cmax}{L_1(\vecW_0)}$, $\delta=\frac{d L_1(\vecW_0)}{\sqrt{L_2(\vecW_0)\tilde{\epsilon}}} e^{-\chi}$ into the above expressions for $\cF_1, \cF_2, \cG, \cL$, using $c \le \cmax$, and directly applying \pref{lem:escapesaddlepoints}, we immediately obtain \pref{lem:escapesaddlepointgdfinal}. The rest of \pref{sec:perturbedgdproofs} is thus devoted to proving \pref{lem:escapesaddlepoints}.

\begin{remark}\label{rem:perturbedgdremarksimplify}
Note it suffices to prove \pref{lem:escapesaddlepoints} for $\delta$ and $\gamma$ smaller than universal constants, as the result \pref{thm:escapesecondordergd} will remain identical under the $O(\cdot)$. Thus we can assume WLOG that $\log(d\kappa/\delta)$ is larger than some universal constant, and that $\gamma \le \frac1{60}$. Also notice by our choice of step size $\eta \le \frac{\cmax}{L_1(\vecW_0)}$ and the assumption $\gamma \le \frac{1}{60}$, for $c\le \cmax \le \frac1{12100}$ we obtain
\[ \kappa \ge 1, r \le 1.\]
This in turn implies
\begin{align*}
\cG &\le \cL,\\
\cF_2 &\ge 40,\\
\cL &\le \sqrt{\eta L_1(\vecW_0)} \cdot  \frac{\gamma}{L_2(\vecW_0)}\cdot  \log^{-1}\prn*{\frac{d \kappa}{\delta}}\\
&\le \frac1{6600} \cdot \min\crl*{1, \frac1{\rho_0(F(\vecW_0)+1)},\frac1{\rho_0(F(\vecW_0)) \rho_0(F(\vecW_0)+1)}},
\end{align*}
where the second line uses that 
\[ L_2(\vecW_0) \ge L_1(\vecW_0)\ge \max\crl*{1,\rho_0(F(\vecW_0)+1),\rho_0(F(\vecW_0)) \rho_0(F(\vecW_0)+1)}. \]
As \textit{these assumptions come with no loss of generality}, we make these assumptions for the rest of the proof.
\end{remark}

To show \pref{lem:escapesaddlepoints}, again as in \citet{jin2017escape}, we prove that the width of the stuck region is not too large.
\begin{lemma}[equivalent of Lemma 15 in \citet{jin2017escape}]\label{lem:widthnotlarge}
There exists a universal constant $\cmax$ such that the following occurs. For any $\delta \in \left(0, \frac{d\kappa}{e}\right]$, let $F$ and $\tilde{\vecW}$ satisfy the conditions in \pref{lem:escapesaddlepoints}. Without loss of generality, by rotational symmetry, let $\vecE_1$ be the minimum eigenvector of $\nabla^2 F(\tilde{\vecW})$. Consider two gradient descent sequences $\{\vecU_t\}$ and $\{\vecX_t\}$ with initial points $\vecU_0, \vecX_0$ satisfying (again, denote the radius $r = \frac{\cL}{\kappa \cdot \log\left(\frac{d\kappa}{\delta}\right)}$):
\[
\nrm*{\vecU_0 - \tilde{\vecW}} \leq r, \quad \vecX_0 = \vecU_0 \pm \mu \cdot r \cdot \vecE_1, \quad \mu \in \brk*{\frac{\delta}{2\sqrt{d}}, 1}.
\]
Then for any step size $\eta \leq \frac{\cmax}{L_1(\vecW_0)}$, and any $T \geq \frac{1}{\cmax}\cF_2$, we have:
\[
\min\{F(\vecU_T) - F(\vecU_0), F(\vecX_T) - F(\vecX_0)\} \le -2.5\cF_1.
\]
\end{lemma}
Now, we prove \pref{lem:escapesaddlepoints} given \pref{lem:widthnotlarge}.
\begin{proof}[Proof of \pref{lem:escapesaddlepoints} given \pref{lem:widthnotlarge}]
Recall as per \pref{rem:perturbedgdremarksimplify} that
\[ \nrm*{\vecP_0-\tilde{\vecW}} \le r \le \cL \le \frac1{\rho_0(F(\vecW_0)+1)}. \]
Also recall $\tilde{\vecW}  \in \cL_{F, F(\vecW_0)}$. Thus by \pref{lem:formalizingperturbationopnrom} we obtain for all $\vecU \in \overline{\vecP_0 \tilde{\vecW}}$ that
\begin{align*}
\nrm*{\grad^2 F(\vecU)}_{\OPNORM} &\le L_1(\vecW_0).
\end{align*}
Therefore by \pref{lem:secondordersmoothnessineq},
\[ F(\vecP_0) \le F(\tilde{\vecW}) +\nrm*{\grad F(\tilde{\vecW})} r + \frac{L_1(\vecW_0)}2 r^2 \le  F(\tilde{\vecW}) +\cG r + \frac{L_1(\vecW_0)}2 r^2=F(\tilde{\vecW}) + \cF_1,\]
where we can readily verify from \pref{rem:perturbedgdremarksimplify} that $\cG r + \frac{L_1(\vecW_0)}2 r^2 \le \cF_1$. 

Now let the stuck region be the set of points $\vecP_0$ in $\ball(\tilde{\vecW}, r)$ such that 
\[ F(\vecP_T) - F(\vecP_0) \ge -2.5 \cF_1.\]
Define the unstuck points by the complement of the stuck points. 

We upper bound the volume of the stuck region as done in \citet{jin2017escape}; this step does not use gradient and Hessian Lispchitzness. Let $1_{\text{Stuck Region}}(\cdot)$ be the indicator function of the stuck region. Write all $\vecW\in\mathbb{R}^d$ as $\vecW = (\vecW^{(1)}, \vecW^{(-1)})$, where $\vecW^{(1)}$ is the component of $\vecW$ along $\vecE_1$ direction and $\vecW^{(-1)}$ is the component of $\vecW$ along the orthogonal complement of $\vecE_1$. By \pref{lem:widthnotlarge}, for any $\vecW \in \ball(\tilde{\vecW}, r)$, 
\begin{align*}
1_{\text{Stuck region}}(\vecW) \DERIV\vecW &= 1_{\text{Stuck region}}(\vecW) \DERIV\vecW^{(-1)} \int_{\tilde{\vecW} - \sqrt{r^2 - \nrm*{\tilde{\vecW}^{(-1)} - \vecW^{(-1)}}^2}}^{\tilde{\vecW} + \sqrt{r^2 - \nrm*{\tilde{\vecW}^{(-1)} - \vecW^{(-1)}}^2}} \DERIV\vecW^{(1)} \\
&\le \DERIV\vecW^{(-1)} \cdot 2 \cdot \frac{\delta}{2\sqrt{d}} r.
\end{align*}
Using this, we have:
\begin{align*}
\text{Volume}(\text{Stuck region}) &= \int_{\ball^d(\tilde{\vecW},r)} 1_{\text{Stuck region}}(\vecW) \DERIV\vecW \\
&= \int_{\ball^{d-1}(\tilde{\vecW},r)} 1_{\text{Stuck region}}(\vecW) \DERIV\vecW^{(-1)} \int_{\tilde{\vecW} - \sqrt{r^2 - \nrm*{\tilde{\vecW}^{(-1)} - \vecW^{(-1)}}^2}}^{\tilde{\vecW} + \sqrt{r^2 - \nrm*{\tilde{\vecW}^{(-1)} - \vecW^{(-1)}}^2}} \DERIV\vecW^{(1)} \\
&\le \int_{\ball^{d-1}(\tilde{\vecW},r)} \DERIV\vecW^{(-1)} \cdot 2 \cdot \frac{\delta}{2\sqrt{d}} r. \\
&= \text{Volume}(\ball^{d-1}(\vecOrigin,r)) \cdot \frac{\delta r}{\sqrt{d}}.
\end{align*}
Then letting $\Gamma(\cdot)$ denote the Gamma function, we have the following ratio:
\begin{align*}
\frac{\text{Volume}(\text{Stuck region})}{\text{Volume}(\ball(\tilde{\vecW}, r))} &\le \frac{\delta r}{\sqrt{d}} \cdot \frac{\text{Volume}(\ball^{d-1}(\vecOrigin,r))}{\text{Volume}(\ball^{d}(\vecOrigin,r))} \\
&= \frac{\delta}{\sqrt{\pi d}} \cdot \frac{\Gamma\left(\frac{d}{2} + 1\right)}{\Gamma\left(\frac{d}{2} + \frac{1}{2}\right)} \\
&\leq \frac{\delta}{\sqrt{\pi d}} \cdot \sqrt{\frac{d}{2} + \frac{1}{2}} \le \delta.
\end{align*}
Here we use the following property of the Gamma function: for $x \ge 0$, $\frac{\Gamma(x+1)}{\Gamma(\frac{x}{2} + \frac{1}{2}} \leq \sqrt{{x + \frac{1}{2}}}$. 

This directly implies that with probability at least $1-\delta$, $\vecP_0$ is an unstuck point. 
Consequently with probability at least $1-\delta$, for any $T \ge \frac1{\cmax}\cF_2$, we have 
\[ F(\vecP_T) - F(\tilde{\vecW}) = F(\vecP_T) - F(\vecP_0) + F(\vecP_0) - F(\tilde{\vecW}) \le -2.5\cF_1+\cF_1=-1.5 \cF_1 < -\cF_1.\]
This proves \pref{lem:escapesaddlepoints}.
\end{proof}
Now we prove \pref{lem:widthnotlarge}, which we do with an analogous strategy as \citet{jin2017escape} by coupling two gradient descent sequences. We have the following two Lemmas, analogous to Lemmas 16, 17 in \citet{jin2017escape}. Again, the reason why they hold in our setting under generalized smoothness is because they all concern `local' behavior around points in the sublevel set of $F(\vecW_0)$. Consequently \pref{lem:boundfuncvalueradius} and \pref{ass:thirdorderselfbounding} ensure we have the required `local' smoothness properties.

Again define $\matH, \tilde{F}_{\vecY}(\vecX)$ analogously to page 20, \citet{jin2017escape}, as follows:
\[ \matH := \grad^2 F(\tilde{\vecW}), \tilde{F}_{\vecY}(\vecX) := F(\vecY)+\tri*{\grad F(\vecY), \vecX - \vecY} + \frac12 (\vecX - \vecY)^\top \matH  (\vecX - \vecY).\numberthis\label{eq:approxfunction}\]
That is, $\tilde{F}_{\vecY}$ is a quadratic approximation of $F$, Taylor expanded about $\tilde{\vecW}$.

The aforementioned Lemmas are as follows:
\begin{lemma}[equivalent of Lemma 16 in \citet{jin2017escape}]\label{lem:technicallemma1}
Letting $\hat{c}=11$, there exists a universal constant $\cmax\le \frac1{12100}$ such that following holds. For any $\delta \in (0, \frac{d\kappa}{e}]$, consider $F, \tilde{\vecW}, r$ as in \pref{lem:escapesaddlepoints}. For any $\vecU_0$ with $\nrm*{\vecU_0 - \tilde{\vecW}} \le 2r = \frac{2\cL}{\kappa \cdot \log\left(\frac{d\kappa}{\delta}\right)}$, define
\[ T = \min \crl*{\inf_t \crl*{t \mid \tilde{F}_{\vecU_0}(\vecU_t)-F(\vecU_0) \leq -3\cF_1}, \hat{c} \cF_2}.\]
Then for any $\eta \le \frac{\cmax}{L(\vecW_0)}$, we have for all $t<T$ that $\nrm*{\vecU_t - \tilde{\vecW}} \le 150\cL \hat{c}$.
\end{lemma}
\begin{lemma}[equivalent of Lemma 17 in \citet{jin2017escape}]\label{lem:technicallemma2}
Letting $\hat{c}=11$, there exists a universal constant $\cmax\le \frac1{12100}$ such that the following holds. For any $\delta \in \left(0, \frac{d\kappa}{e}\right]$, consider $F, \tilde{\vecW}, r$ as in \pref{lem:escapesaddlepoints}, and sequences $\{\vecU_t\}$, $\{\vecX_t\}$ satisfying the conditions in \pref{lem:widthnotlarge}. Define: 
\[
T = \min \left\{
\inf_t 
\left\{ t \mid \tilde{F}_{\vecX_0}(\vecX_t) - F(\vecX_0) \leq -3\cF_1 \right\}, \; \hat{c}\cF_2
\right\}.
\]
Then, for any $\eta \leq \frac{\cmax}{L_1(\vecW_0)}$, if $\nrm*{\vecU_t - \tilde{\vecW}} \leq 150 \cL\hat{c}$ for all $t < T$, we will have $T < \hat{c}\cF_2$. Equivalently, this means that 
\[
\inf_t 
\left\{ t : \tilde{F}_{\vecX_0}(\vecX_t) - F(\vecX_0) \leq -3\cF_1 \right\} < \hat{c}\cF_2,
\]  
i.e. that we escaped the saddle point.
\end{lemma}
\begin{proof}[Proof of \pref{lem:widthnotlarge} given \pref{lem:technicallemma1}, \pref{lem:technicallemma2}]
Choosing $\cmax$ to be the minimum of the $\cmax$ from \pref{lem:technicallemma1}, \pref{lem:technicallemma2}, we can ensure both Lemmas hold. Clearly this preserves that $\cmax \le \frac1{12100}$.

Define
\[ T^{\star} = \hat{c} \cF_2, T' = \inf \crl*{t: \tilde{F}_{\vecU_0}(\vecU_t) - F(\vecU_0) \le -3\cF_1}.\]
We break into cases on $T'$ versus $T^{\star}$:
\begin{itemize}
\item $T' \le T^{\star}$: By \pref{lem:technicallemma1}, $\nrm*{\vecU_{T'-1} - \tilde{\vecW}} \le 150\cL \hat{c}$. Since $\cL \le \frac1{6600} \cdot \frac1{\rho_0(F(\vecW_0)+1)}$ from \pref{rem:perturbedgdremarksimplify} and $\hat{c}=11$, this yields 
\[ \nrm*{\vecU_{T'-1} - \tilde{\vecW}} \le 150 \cL\hat{c}\le \frac14 \cdot \frac1{\rho_0(F(\vecW_0)+1)}.\]
Thus because $\tilde{\vecW} \in \cL_{F,F(\vecW_0)}$, by \pref{lem:formalizingperturbationopnrom}, we have
\[ \nrm*{\grad^2 F(\vecU)} \le L_1(\vecW_0) \text{ for all } \vecU \in \overline{\vecU_{T'-1} \tilde{\vecW}}.\]
Thus, recalling $\cG \le \cL$ from \pref{rem:perturbedgdremarksimplify}, we obtain
\begin{align*}
\nrm*{\grad F(\vecU_{T'-1})} &\le \nrm*{\grad F(\tilde{\vecW})} + L_1(\vecW_0) \nrm*{\vecU_{T'-1} - \tilde{\vecW}} \\
&\le \cG + 150 \hat{c} L_1(\vecW_0) \cL \le \cL + 150 \hat{c} L_1(\vecW_0) \cL.
\end{align*}
Therefore, as $\eta L_1(\vecW_0) \le \cmax \le 1$,
\begin{align*}
\nrm*{\vecU_{T'}-\tilde{\vecW}} &\le \nrm*{\vecU_{T'-1}-\tilde{\vecW}} + \eta \nrm*{\grad F(\vecU_{T'-1})} \\
&\le 150 \cL \hat{c} + \cL + 150 \hat{c} \cdot \eta L_1(\vecW_0) \cL \le (300 \hat{c}+1) \cL \numberthis\label{eq:perturbgdlocal1}
\end{align*}
Recalling $\kappa, \log\prn*{\frac{d\kappa}{\delta}} \ge 1$, the conditions of \pref{lem:widthnotlarge} give 
\[ \nrm*{\vecU_0 - \tilde{\vecW}} \le r \le \cL. \numberthis\label{eq:perturbgdlocal2}\]
Combining \pref{eq:perturbgdlocal1}, \pref{eq:perturbgdlocal2} and applying Triangle Inequality gives
\[ \nrm*{\vecU_{T'} - \vecU_0} \le (300\hat{c}+2)\cL.\numberthis\label{eq:perturbgdlocal3}\]
Also by \pref{eq:perturbgdlocal2}, we have $\nrm*{\vecU_0-\tilde{\vecW}} \le\cL \le \frac1{\rho_0(F(\vecW_0)+1)}$. Thus as $\tilde{\vecW} \in \cL_{F, F(\vecW_0)}$, by \pref{lem:formalizingperturbationopnrom} we obtain
\[ \nrm*{\grad^2 F(\vecU_0)} \le L_1(\vecW_0).\numberthis\label{eq:perturbgdopboundez}\]
Moreover, by Triangle Inequality we obtain that for any $\vecU \in \overline{\vecU_{0}\vecU_{T'}}$, we have 
\[ \nrm*{\vecU - \tilde{\vecW}} \le (300\hat{c}+2) \cL = 3302 \cL \le \frac1{\rho_0(F(\vecW_0)+1)}.\] 
As $\tilde{\vecW} \in \cL_{F, F(\vecW_0)}$, \pref{lem:formalizingperturbationopnrom} implies for all such $\vecU_1, \vecU_2 \in \overline{\vecU_{0}\vecU_{T'}}$ that
\[ \nrm*{\grad^2 F(\vecU_1) - \grad^2 F(\vecU_2)}_{\OPNORM} \le \nrm*{\vecU_1 - \vecU_2} L_2(\vecW_0).\]
Now applying \pref{lem:thirdordersmoothnessineq}, and by choosing $\eta = \frac{c}{L(\vecW_0)}$ for a small enough universal constant $c$, we obtain:
\begin{align*}
&F(\vecU_{T'}) - F(\vecU_0) \\
&\le \grad F(\vecU_0)^\top (\vecU_{T'} - \vecU_0) + \frac12 (\vecU_{T'} - \vecU_0)^\top \grad^2 F(\vecU_0) (\vecU_{T'} - \vecU_0) + \frac{L_2(\vecW_0)}{6} \nrm*{\vecU_{T'} - \vecU_0}^3 \\
&\le \tilde{F}_{\vecU_0}(\vecU_{T'})-F(\vecU_0) + \frac{L_2(\vecW_0)}2 \nrm*{\vecU_{T'} - \vecU_0}^2 \nrm*{\vecU_0 - \tilde{\vecW}}+ \frac{L_2(\vecW_0)}{6} \nrm*{\vecU_{T'} - \vecU_0}^3 \\
&\le -3\cF_1 + O(L_1(\vecW_0) \cL^3) \\
&= -3\cF_1 + O(\sqrt{\eta L_1(\vecW_0)} \cF_1) \le -2.5 \cF_1.
\end{align*}
Here we used \pref{eq:perturbgdopboundez}, \pref{eq:perturbgdlocal2}, \pref{eq:perturbgdlocal3}, and that $\cL \le 1$ as per \pref{rem:perturbedgdremarksimplify}.
In the above, $O(\cdot)$ only hides universal constants as $\hat{c}=11$ is a universal constant, and so these final inequalities can be made to hold by choosing $\cmax$ a sufficiently small universal constant.

Since $\tilde{\vecW}  \in \cL_{F, F(\vecW_0)}$ and $\eta \le \frac2{L_1(\vecW_0)}$, \pref{lem:gdnoincrease} shows that gradient descent will not increase value (this is essentially the same as several steps the proof of \pref{thm:gdfirstorder}, combined with induction). Thus for all $T \ge T'$ and hence for all $T \ge \frac1{\cmax} \cF_2 \ge \hat{c} \cF_2 \ge T'$ along this gradient descent trajectory, we have 
\[ F(\vecU_T) - F(\vecU_0) \le F(\vecU_{T'}) - F(\vecU_0) \le -2.5 \cF_1.\]
\item $T' > T^{\star}$: In this case, by \pref{lem:technicallemma1}, we know $\nrm*{\vecU_t - \tilde{\vecW}} \leq 150\cL \hat{c}$ for all $t < T^\star=\hat{c}\cF_2$. 

Define 
\[
T'' = \inf_t \left\{t \mid \tilde{F}_{\vecX_0}(\vecX_t) - F(\vecX_0) \leq -3\cF_1\right\}.
\]
Since $\nrm*{\vecU_t - \tilde{\vecW}} \leq 150\cL \hat{c}$ for all $t < T^\star=\hat{c}\cF_2$, it follows that $\nrm*{\vecU_t - \tilde{\vecW}} \leq 150\cL \hat{c}$ for all $t<\min\crl*{T'', T^{\star}}$. Thus by \pref{lem:technicallemma2}, we have that $\min\crl*{T'', T^{\star}} < T^{\star}$, and so $T'' < T^\star$. Applying the same argument as in the first case to the $\{\vecX_t\}$, we have that for all $T \geq \frac{1}{c_{\text{max}}}\cF_2$ that
\[ F(\vecX_T) - F(\vecX_0) \le -2.5 \cF_1.\]
\end{itemize}
This proves \pref{lem:widthnotlarge}.
\end{proof}
\begin{remark}
Note that $\tilde{\vecW} \in \cL_{F, F(\vecW_0)}$ is central to this argument, unlike the Lipschitz gradient and Hessian case from \citet{jin2017escape}.
\end{remark}

\subsection{Proof of Escaping Saddles Lemmas}
Now we prove \pref{lem:technicallemma1}, \pref{lem:technicallemma2}.
\begin{proof}[Proof of \pref{lem:technicallemma1}]
We follow the proof of Lemma 16, \citet{jin2017escape}. Again, we aim to show that if the function value does not decrease, then all the iterates must remain constrained in a small ball. This is done by analyzing the dynamics of the iterates and decomposing the $d$-dimensional space into two subspaces: a subspace $S$, which is the span of the negative enough eigenvectors of the Hessian, and its orthogonal complement.

The main difference now is that now we cannot directly control relevant operator norms with global Lipschitz properties of the gradient and Hessian. However, it turns out that the proof of this Lemma will follow induction on the iterate $\vecU_t$, and consequently we will obtain that all of the prior iterates $\vecU_{t'}$ for $t'<t$ are close enough to $\tilde{\vecW}$. By a similar argument as in \pref{lem:escapesaddlepoints}, since $\tilde{\vecW} \in \cL_{F, F(\vecW_0)}$, this lets us upper bound the gradient of these points. By the Gradient Descent update rule, this in turn implies the current iterate is also close to $\tilde{\vecW}$, and thus we obtain bounds on the relevant derivatives in terms of $L_1(\vecW_0)$, $L_2(\vecW_0)$ for all points in the convex hull of the relevant iterates.

We begin the argument. Analogously to \citet{jin2017escape}, since $\delta \in \left(0, \frac{d\kappa}{e}\right]$, we always have $\log\left(\frac{d\kappa}{\delta}\right) \geq 1$. By the gradient descent update function, we have
\[
\vecU_{t+1} = \vecU_t - \eta \grad F(\vecU_t).
\]
This can be expanded as:
\[
\vecU_{t+1} = \vecU_t - \eta \nabla F(\vecU_0) - \eta \prn*{\int_0^1 \nabla^2 F(\theta(\vecU_t - \vecU_0) + \vecU_0)\DERIV\theta} (\vecU_t - \vecU_0).
\]
Recall the definition $\matH = \grad^2 F(\tilde{\vecW})$.
Let $\Delta_t$ be defined as:
\[
\Delta_t := \int_0^1 \nabla^2 F(\theta(\vecU_t - \vecU_0) + \vecU_0)\DERIV\theta - \matH.
\]
Substituting, we obtain:
\[
\vecU_{t+1} = (\matI - \eta \matH - \eta \Delta_t)(\vecU_t - \vecU_0) - \eta \nabla F(\vecU_0) + \vecU_0.
\]
Note we do not immediately have an upper bound on the operator norm of $\Delta_t$. In particular this is because $t$ could diverge (logarithmically) in the dimension, only being upper bounded by $\cF_2$.

We now compute the projections of $\vecU_t-\vecU_0$ in different eigenspaces of $\matH$. Define $\cS$ as the subspace spanned by all eigenvectors of $\matH$ whose eigenvalues are less than $-\frac{\gamma}{\hat{c} \log\left(\frac{d\kappa}{\delta}\right)}$. Let $\cS^c$ denote the subspace of the remaining eigenvectors. Let $\vecAlpha_t$ and $\vecBeta_t$ denote the projections of $\vecU_t-\vecU_0$ onto $\cS$ and $\cS^c$ respectively, i.e., $\vecAlpha_t = \cP_{\cS} (\vecU_t-\vecU_0)$, and $\vecBeta_t = \cP_{\cS^c} (\vecU_t-\vecU_0)$. 

We can decompose the update equations for $\vecU_{t+1}$ into:
\[
\vecAlpha_{t+1} = (\matI - \eta \matH)\vecAlpha_t - \eta \cP_{\cS} \Delta_t (\vecU_t-\vecU_0) - \eta \cP_{\cS} \grad F(\vecU_0),
\]
\[
\vecBeta_{t+1} = (\matI - \eta \matH)\vecBeta_t - \eta \cP_{\cS^c} \Delta_t (\vecU_t-\vecU_0) - \eta \cP_{\cS^c} \grad F(\vecU_0). 
\]
By the definition of $T$, we know for all $t < T$:
\begin{align*}
-3\cF_1 < \tilde{F}_{\vecU_0}(\vecU_t) - F(\vecU_0) &= \nabla F(\vecU_0)^\top (\vecU_t-\vecU_0) - \frac{1}{2} (\vecU_t-\vecU_0)^\top \matH (\vecU_t-\vecU_0) \\
&\leq \nabla F(\vecU_0)^\top (\vecU_t-\vecU_0) - \frac{\gamma}{2} \frac{\nrm*{\vecAlpha_t}^2}{\hat{c} \log\left(\frac{d\kappa}{\delta}\right)} + \frac{1}{2} \vecBeta_t^\top \matH \vecBeta_t.
\end{align*}
Evidently we have $\nrm*{\vecU_t-\vecU_0}^2 = \nrm*{\vecAlpha_t}^2 + \nrm*{\vecBeta_t}^2$, and thus the above rearranges to
\[ \nrm*{\vecU_t-\vecU_0}^2 \leq \frac{2\hat{c} \log\left(\frac{d\kappa}{\delta}\right)}{\gamma} \prn*{3\cF_1 + \nabla F(\vecU_0)^\top (\vecU_t-\vecU_0) + \frac{1}{2} \vecBeta_t^\top \matH \vecBeta_t } + \nrm*{\vecBeta_t}^2.\numberthis\label{eq:perturbedgdlem1utu0initbound} \]
Now we control $\nrm*{\grad F(\vecU_0)}$. We use the fact that $\tilde{\vecW} \in \cL_{F, F(\vecW_0)}$ to give us the necessary control over this quantity. Similar ideas were used in the proof of \pref{lem:widthnotlarge}, and will continue to be used in the rest of the proofs of \pref{lem:technicallemma1}, \pref{lem:technicallemma2}. In particular, recall as per \pref{rem:perturbedgdremarksimplify} that
\[ \nrm*{\vecU_0-\tilde{\vecW}} \le 2r \le 2 \cL \le \frac1{\rho_0(F(\vecW_0)+1)}. \]
Thus by \pref{lem:formalizingperturbationopnrom}, as $\tilde{\vecW} \in cL_{F,F(\vecW_0)}$, we obtain
\[ \nrm*{\grad^2 F(\vecU)} \le L_1(\vecW_0) \text{ for all }\vecU \in \overline{\vecU_0 \tilde{\vecW}}.\]
Consequently,
\[ \nrm*{\grad F(\vecU_0)- \grad F(\tilde{\vecW})} \le L_1(\vecW_0) \nrm*{\vecU_0-\tilde{\vecW}} \le 2r L_1(\vecW_0) = 2\cG, \]
which implies 
\[ \nrm*{\grad F(\vecU_0)} \le \nrm*{\grad F(\tilde{\vecW})} +2\cG = 3\cG.\numberthis\label{eq:perturbedgdlem1fu0gradbound}\]
This gives us an analogous bound on $\nrm*{\grad F(\vecU_0)}$ as in the proof of Lemma 16, \citet{jin2017escape}. Substituting this bound on $\nrm*{\grad F(\vecU_0}$ into \pref{eq:perturbedgdlem1utu0initbound}, we obtain
\[
\nrm*{\vecU_t-\vecU_0}^2 \leq 14 \max \crl*{ \frac{\cG \hat{c} \log\left(\frac{d\kappa}{\delta}\right)}{\gamma} \nrm*{\vecU_t-\vecU_0}, \frac{\cF_1 \hat{c} \log\left(\frac{d\kappa}{\delta}\right)}{\gamma}, \frac{\vecBeta_t^\top \matH \vecBeta_t \hat{c} \log\left(\frac{d\kappa}{\delta}\right)}{\gamma}, \nrm*{\vecBeta_t}^2 }.
\]
In turn this implies
\[
\nrm*{\vecU_t-\vecU_0} \leq 14  \max \crl*{ \frac{\cG \hat{c} \log\left(\frac{d\kappa}{\delta}\right)}{\gamma}, \sqrt{\frac{\cF_1 \hat{c} \log\left(\frac{d\kappa}{\delta}\right)}{\gamma}}, \sqrt{\frac{\vecBeta_t^\top \matH \vecBeta_t \hat{c} \log\left(\frac{d\kappa}{\delta}\right)}{\gamma}}, \nrm*{\vecBeta_t} }.\numberthis\label{eq:perturbedgdlem1utu0bound}
\]

\paragraph{The key induction:} Now, we induct on $t$ to prove 
\[ \nrm*{\vecU_t-\vecU_0} \le 148\cL \hat{c} \text{ for all }t<T.\numberthis\label{eq:goalforlem1}\]
Clearly this implies \pref{lem:technicallemma1}, upon recalling $\nrm*{\vecU_0-\tilde{\vecW}} \le 2r = 2\cL \le \hat{c}\cL$ by our choice $\hat{c}=11$.

The base case $t=0$ is evident.

Now for the inductive step, suppose \pref{eq:goalforlem1} is true for all $\tau \le t$ such that $t+1 < T$. We show it is true for $t+1$.

Due to the above bound \pref{eq:perturbedgdlem1utu0bound}, it suffices to upper bound $\nrm*{\vecBeta_{t+1}}, \vecBeta_{t+1}^\top \matH \vecBeta_{t+1}$. We note as in the proof of Lemma 16 of \citet{jin2017escape} that letting 
\[ \vecDelta_t := \cP_{\cS^c}\prn*{\Delta_t (\vecU_t-\vecU_0) + \grad F(\vecU_0)},\]
we have by the Triangle Inequality and properties of projections that
\[ \nrm*{\vecDelta_{t}} \le \nrm*{\Delta_t}_{\OPNORM} \nrm*{\vecU_t-\vecU_0}+\nrm*{\grad F(\vecU_0)}.\numberthis\label{eq:perturbedgdlem1control1}\]
Furthermore, we have by definition of the update rule for $\vecBeta_{t+1}$ that
\[ \vecBeta_{t+1} = (\matI - \eta \matH)\vecBeta_t + \eta \vecDelta_t.\numberthis\label{eq:perturbedgdlem1betarecurseexact}\]
Thus,
\[ \nrm*{\vecBeta_{t+1}} \le \nrm*{(\matI - \eta \matH)\vecBeta_t} + \eta \vecDelta_t \le \nrm*{\vecBeta_t} + \eta \nrm*{\matH \vecBeta_t} + \eta \vecDelta_t.\numberthis\label{eq:perturbedgdlem1betarecurse}\]
Now, consider any $\tau, 0 \le \tau \le t$. We upper bound $\nrm*{\Delta_{\tau}}_{\OPNORM}$. Rewrite 
\[ \Delta_{\tau} = \int_0^1 \prn*{\grad^2 F(\theta(\vecU_{\tau}-\vecU_0)+\vecU_0) -\grad^2 F(\vecU_0)} \DERIV \theta + \grad^2 F(\vecU_0) - \grad^2 F(\tilde{\vecW}).\]
Clearly, as per \pref{rem:perturbedgdremarksimplify},
\[ \nrm*{\vecU_0 - \tilde{\vecW}} \le 2r \le 2\cL \le \frac1{\rho_0(F(\vecW_0)+1)}.\]
Recalling $\tilde{\vecW} \in \cL_{F,F(\vecW_0)}$ and applying \pref{lem:formalizingperturbationopnrom} gives 
\[ \nrm*{\grad^2 F(\vecU_0) - \grad^2 F(\tilde{\vecW})}_{\OPNORM} \le L_2(\vecW_0) \nrm*{\vecU_0 - \tilde{\vecW}}.\numberthis\label{eq:perturbedgdlem1hessianlipschitz1}\]
Moreover by inductive hypothesis, we know that $\nrm*{\vecU_{\tau}-\vecU_0} \le148 \cL \hat{c}$. Consequently as $\hat{c}=11 \ge 1$ and following \pref{rem:perturbedgdremarksimplify}, for all $\theta \in [0,1]$, we have
\[ \nrm*{\prn*{\theta (\vecU_{\tau} - \vecU_0)+\vecU_0}-\tilde{\vecW}} \le 2\cL+148\hat{c}\cL \le \frac1{\rho_0(F(\vecW_0)+1)}. \]
Since $\tilde{\vecW} \in \cL_{F, F(\vecW_0)}$, it follows by \pref{lem:formalizingperturbationopnrom} that
\[ \nrm*{\grad^2 F(\theta(\vecU_{\tau}-\vecU_0)+\vecU_0) -\grad^2 F(\vecU_0)}_{\OPNORM} \le L_2(\vecW_0) \nrm*{\vecU_{\tau} - \vecU_0} \text{ for all }\theta \in [0,1].\numberthis\label{eq:perturbedgdlem1hessianlipschitz2} \]
Hence by Triangle Inequality, from \pref{eq:perturbedgdlem1hessianlipschitz1} and \pref{eq:perturbedgdlem1hessianlipschitz2}, we have
\[ \nrm*{\Delta_t}_{\OPNORM} \le L_2(\vecW_0)\prn*{\nrm*{\vecU_{\tau}-\vecU_0} + \nrm*{\vecU_0-\tilde{\vecW}}} \le L_2(\vecW_0)(148\cL \hat{c}+\nrm*{\vecU_0-\tilde{\vecW}}).\numberthis\label{eq:perturbedgdlem1control2}\]
Proceeding from here is now exactly the same as in \citet{jin2017escape}. We detail the argument for completeness.

Combining \pref{eq:perturbedgdlem1control1}, \pref{eq:perturbedgdlem1control2}, \pref{eq:perturbedgdlem1fu0gradbound} and applying the inductive hypothesis and the condition of \pref{lem:escapesaddlepoints} that $\nrm*{\vecU_0 - \tilde{\vecW}} \le 2r$, gives
\begin{align*}
\nrm*{\vecDelta_{\tau}} &\leq L_2(\vecW_0)(148\cL \hat{c}+\nrm*{\vecU_0-\tilde{\vecW}})\nrm*{\vecU_{\tau}-\vecU_0} + \nrm*{\nabla F(\vecU_0)}\\
&\leq L_2(\vecW_0)\cdot 148\hat{c}\left(148\hat{c} + \frac{2}{\kappa \cdot \log\left(\frac{d\kappa}{\delta}\right)}\right) \cL^2 + 3\cG.
\end{align*}
Plugging in the choice of $\cL$, and choosing a small enough constant $c_{\max} \leq \prn*{\frac{1}{2 \cdot 148\hat{c}(148\hat{c}+2)}}^2$ and choosing step size $\eta < \frac{c_{\max}}{L_1(\vecW_0)}$, gives for any $0 \le \tau \le t$:
\[
\nrm*{\vecDelta_{\tau}} \leq \crl*{148\hat{c}\left(148\hat{c} + \frac{2}{\kappa \cdot \log\left(\frac{d\kappa}{\delta}\right)}\right)\sqrt{\eta L_1(\vecW_0)} + 3}\cG \leq 3.5\cG.\numberthis\label{eq:perturbedgdlem1deltabound}
\]
We now bound $\nrm*{\vecBeta_{t+1}}, \vecBeta_{t+1}^\top \matH \vecBeta_{t+1}$, which combining with \pref{eq:perturbedgdlem1utu0bound} finishes the induction and thus the proof.
\begin{itemize}
    \item In order to bound $\nrm*{\vecBeta_{t+1}}$, combining \pref{eq:perturbedgdlem1betarecurse} with \pref{eq:perturbedgdlem1deltabound} and recalling the definition of $\cS$ and $\vecBeta_t$ gives:
\[
\nrm*{\vecBeta_{t+1}} \leq \left(1 + \frac{\eta\gamma}{\hat{c} \log\left(\frac{d\kappa}{\delta}\right)}\right)\nrm*{\vecBeta_t} + 3.5\eta \cG.
\]
Since $\nrm*{\vecBeta_0} = 0$ and $t + 1 \leq T$, by applying the above relation recursively, we have:
\[
\nrm*{\vecBeta_{t+1}} \leq \sum_{\tau=0}^T 3.5\left(1 + \frac{\eta\gamma}{\hat{c} \log\left(\frac{d\kappa}{\delta}\right)}\right)^\tau \eta \cG \leq 3.5 \cdot 3 \cdot T \eta \cG \leq 10.5\cL \hat{c}. \numberthis\label{eq:perturbedgdlem1betat}
\]
In the above we used $T \leq \hat{c}\cF$, which also implies $\left(1 + \frac{\eta\gamma}{\hat{c} \log\left(\frac{d\kappa}{\delta}\right)}\right)^T \le \left(1 + \frac{\eta\gamma}{\hat{c} \log\left(\frac{d\kappa}{\delta}\right)}\right)^{\hat{c}\cF} \leq 3$ (one can find an easy upper bound on $\cF$ based on its definition and check using $L_2(\vecW_0) \ge L_1(\vecW_0) \ge 1$ that this is the case).
\item Now for bounding $\vecBeta_{t+1}^\top \matH \vecBeta_{t+1}$, notice we can also write the update equation \pref{eq:perturbedgdlem1betarecurseexact} for $\vecBeta_t$ as:
\[
\vecBeta_t = \eta\sum_{\tau=0}^{t-1} (\matI - \eta \matH)^\tau \vecDelta_{t-1-\tau}.
\]
As $\matH$ is symmetric this gives:
\[
\vecBeta_{t+1}^\top \matH \vecBeta_{t+1} = \eta^2 \sum_{\tau_1=0}^t \sum_{\tau_2=0}^t \vecDelta^{\top}_{t-1-\tau_1} (\matI - \eta \matH)^{\tau_1} \matH (\matI - \eta \matH)^{\tau_2} \vecDelta_{t-1-\tau_2}.
\]
Thus we have:
\[
\vecBeta_{t+1}^\top \matH\vecBeta_{t+1} \leq \eta^2 \sum_{\tau_1=0}^t \sum_{\tau_2=0}^t \nrm*{\vecDelta_{t-1-\tau_1}} \| (\matI - \eta \matH)^{\tau_1} \matH (\matI - \eta \matH)^{\tau_2} \| \nrm*{\vecDelta_{t-1-\tau_2}}.
\]
Since for $0 \le \tau_1, \tau_2 \le t$ we have $\nrm*{\vecDelta_{t-1-\tau_1}}, \nrm*{\vecDelta_{t-1-\tau_2}} \leq 3.5\cG$ as argued earlier, we have:
\[
\vecBeta_{t+1}^\top \matH \vecBeta_{t+1} \leq 3.5^2\eta^2 \cG^2 \sum_{\tau_1=0}^t \sum_{\tau_2=0}^t \nrm*{ (\matI - \eta \matH)^{\tau_1} \matH (\matI - \eta \matH)^{\tau_2} }.
\]
Let the eigenvalues of $\matH$ be $\{\lambda_i\}$. Thus for any $\tau_1, \tau_2 \geq 0$, the eigenvalues of $(\matI - \eta \matH)^{\tau_1} \matH(\matI - \eta \matH)^{\tau_2}$ are $\{\lambda_i (1 - \eta \lambda_i)^{\tau_1 + \tau_2}\}$. We now detail a calculation from \citet{jin2017escape}. Letting $g_t(\lambda) := \lambda(1 - \eta \lambda)^t$ and setting its derivative to zero yields
\[
\nabla g_t(\lambda) = (1 - \eta \lambda)^t - t \eta \lambda (1 - \eta \lambda)^{t-1} = 0.
\]
It is easy to check that $\lambda_t^\star = \frac{1}{(1 + t)\eta}$ is the unique maximizer, and $g_t(\lambda)$ is monotonically increasing in $(-\infty, \lambda_t^\star]$.

This gives:
\[
\nrm*{ (\matI - \eta \matH)^{\tau_1} \matH (\matI - \eta \matH)^{\tau_2} } = \max_i \lambda_i (1 - \eta \lambda_i)^{\tau_1 + \tau_2} \leq \hat{\lambda}(1 - \eta \hat{\lambda})^{\tau_1 + \tau_2} \leq \frac{1}{(1 + \tau_1 + \tau_2)\eta},
\]
where $\hat{\lambda} = \min\{\ell, \lambda_{\tau_1 + \tau_2}^\star\}$. Therefore, we have:
\[
\vecBeta_{t+1}^\top \matH \vecBeta_{t+1} \leq 3.5^2\eta \cG^2 \sum_{\tau_1=0}^t \sum_{\tau_2=0}^t \frac{1}{1 + \tau_1 + \tau_2}.
\]
To bound the sum note:
\[
\sum_{\tau_1=0}^t \sum_{\tau_2=0}^t \frac{1}{1 + \tau_1 + \tau_2} = \sum_{\tau=0}^{2t} \min\{1 + \tau, 2t + 1 - \tau\} \cdot \frac{1}{1 + \tau} \leq 2t + 1 < 2T.
\]
Thus:
\[
\vecBeta_{t+1}^\top \matH \vecBeta_{t+1} \leq 2 \cdot 3.5^2\eta T \cG^2 \leq \frac{3.5^2 \cL^2 \gamma \hat{c}}{\log\left(\frac{d\kappa}{\delta}\right)}. \numberthis\label{eq:perturbedgdlem1quadformbound}
\]
\end{itemize}
Finally, substituting the previous upper bounds \pref{eq:perturbedgdlem1betat}, \pref{eq:perturbedgdlem1quadformbound} for $\nrm*{\vecBeta_t}$, $\vecBeta_{t+1}^\top \matH \vecBeta_{t+1}$ into our prior display \pref{eq:perturbedgdlem1utu0bound} for $\nrm*{\vecU_t-\vecU_0}$, we obtain:
\[
\nrm*{\vecU_t-\vecU_0} \leq 14  \max \crl*{ \frac{\cG \hat{c} \log\left(\frac{d\kappa}{\delta}\right)}{\gamma}, \sqrt{\frac{\cF_1 \hat{c} \log\left(\frac{d\kappa}{\delta}\right)}{\gamma}}, \sqrt{\frac{\vecBeta_t^\top \matH \vecBeta_t \hat{c} \log\left(\frac{d\kappa}{\delta}\right)}{\gamma}}, \nrm*{\vecBeta_t} } \leq 148 \cL \hat{c}.
\]
This finishes the induction, and hence the proof of the Lemma.
\end{proof}

\begin{proof}[Proof of \pref{lem:technicallemma2}]
Again, we aim to show that if all iterates from $\vecU_0$ are contained in a small ball, then the iterates from $\vecX_0$ decrease function value. As with the proof of \pref{lem:technicallemma1}, the proof combines the proof idea of Lemma 17, \citet{jin2017escape} with the self-bounding framework. This time it goes through even easier, because the required new bounds that we need from the relevant iterates being `local' hold not due to induction, but rather from a direct application of \pref{lem:technicallemma1}.

Define $\vecV_t = \vecX_t - \vecU_t$. By the assumptions of this Lemma we have that $\vecV_0 = \pm \mu \left[\frac{\cL}{\kappa \cdot \log\left(\frac{d \kappa}{\delta}\right)}\right] \vecE_1$ where $\mu \in \left[\frac{\delta}{2\sqrt{d}}, 1\right]$. Consequently 
\[ \frac{\delta}{2\sqrt{d}} \cdot r \le \nrm*{\vecV_0} \le r.\numberthis\label{eq:perturbedgdlem2v0bound}\]
Recall the definition 
\[ \matH=\grad^2 F(\tilde{\vecW}) \]
as per \pref{eq:approxfunction}. Also define
\[
\Delta_t' := \int_0^1 \nabla^2 F(\vecU_t + \theta \vecV_t)\DERIV \theta - \matH.
\]
Exactly as in the proof of Lemma 17, \citet{jin2017escape}, by directly writing the update equations, we have
\begin{align*}
\vecU_{t+1} + \vecV_{t+1} = \vecX_{t+1} &= \vecX_t - \eta \grad F(\vecX_t) \\
&= \vecU_t + \vecV_t - \eta \grad F(\vecU_t+\vecV_t) \\
&= \vecU_t + \vecV_t - \eta \grad F(\vecU_t) - \eta \prn*{\int_0^1 \grad^2 F(\vecU_t+\theta \vecV_t) \DERIV \theta} \vecV_t \\
&= \vecU_t + \vecV_t - \eta \grad F(\vecU_t) - \eta \prn*{\matH + \Delta_t'} \vecV_t \\
&= \vecU_t - \eta \grad F(\vecU_t) + (\matI - \eta \matH - \eta \Delta'_t) \vecV_t.
\end{align*}
Hence as $\vecU_{t+1} =  \vecU_t - \eta \grad F(\vecU_t)$, we obtain
\[ \vecV_{t+1} = (\matI-\eta \matH - \eta \Delta'_t) \vecV_t.\numberthis\label{eq:perturbedgdlem2vtrecurse}\]
The difference from the proof of Lemma 17, \citet{jin2017escape} is now that we do not immediately have an upper bound on $\nrm*{\Delta'_t}_{\OPNORM}$ without global Lipschitzness of the gradient and Hessian. However, similarly as in the proof of \pref{lem:technicallemma1}, we can obtain such a bound using the self-bounding framework, since the point $\tilde{\vecW}$ in question is in the $F(\vecW_0)$-sublevel set $\cL_{F, F(\vecW_0)}$.

Note by hypothesis on $\vecU_0$ from \pref{lem:widthnotlarge} and as $\nrm*{\vecV_0} \le r$ by \pref{eq:perturbedgdlem2v0bound},
\[ \nrm*{\vecX_0-\tilde{\vecW}} \le \nrm*{\vecU_0-\tilde{\vecW}}+\nrm*{\vecV_0} \le r+r=2r.\]
Applying \pref{lem:technicallemma1} directly to the $\{\vecX_t\}$ implies that 
\[ \nrm*{\vecX_t - \tilde{\vecW}} \le 150 \cL \hat{c} \text{ for all } t<T.\]
By assumption of this Lemma, we have 
\[ \nrm*{\vecU_t - \tilde{\vecW}} \le 150 \cL \hat{c}\text{ for all } t<T.\]
Triangle Inequality thus gives
\[ \nrm*{\vecV_t} \le 300 \cL \hat{c}, \nrm*{\vecU_t - \vecU_0} \le 300\cL \hat{c}\text{ for all } t<T.\]
Therefore for all $0 \le \theta \le 1$,
\[ \vecU_t + \theta \vecV_t \in \ball(\tilde{\vecW}, 600\cL \hat{c}).\]
Note as per \pref{rem:perturbedgdremarksimplify},
\[ 600 \cL \hat{c} = 6600\cL \le \frac1{\rho_0(F(\vecW_0)+1)}. \]
As $\tilde{\vecW} \in \cL_{F, F(\vecW_0)}$, it follows from \pref{lem:formalizingperturbationopnrom} that 
\[ \nrm*{\grad^2 F(\vecU_t + \theta \vecV_t) - \grad^2 F(\vecU_t)}_{\OPNORM} \le L_2(\vecW_0) \cdot \theta \vecV_t\text{ for all }\theta \in [0,1].\numberthis\label{eq:technicallemma2smoothness1}\]
Similarly, by the above bound 
\[ \nrm*{\vecU_t-\tilde{\vecW}} \le 150 \cL \hat{c} \le \frac1{\rho_0(F(\vecW_0)+1)}\]
and as $\tilde{\vecW} \in \cL_{F, F(\vecW_0)}$, \pref{lem:formalizingperturbationopnrom} proves that
\[ \nrm*{\grad^2 F(\vecU_t) - \grad^2 F(\tilde{\vecW})}_{\OPNORM} \le L_2(\vecW_0) \nrm*{\vecU_t - \tilde{\vecW}}.\numberthis\label{eq:technicallemma2smoothness2}\]
Now, rewrite
\[ \Delta'_t = \int_0^1 \prn*{\grad^2 F(\vecU_t + \theta \vecV_t) - \grad^2 F(\vecU_t)} \DERIV \theta + \grad^2 F(\vecU_t) - \grad^2 F(\tilde{\vecW}). \]
By \pref{eq:technicallemma2smoothness1}, \pref{eq:technicallemma2smoothness2}, and the above bounds on $\nrm*{\vecV_t}, \nrm*{\vecU_t-\tilde{\vecW}}$, we obtain for all $\theta \in [0,1]$ that
\[ \nrm*{\Delta'_t}_{\OPNORM} \le L_2(\vecW_0) \prn*{\theta\nrm*{\vecV_t}+\nrm*{\vecU_t - \tilde{\vecW}}} \le L_2(\vecW_0) \cL (450\hat{c}+1).\numberthis\label{eq:perturbedgdlem2deltatbound}\]
From here, exactly the same proof as that of Lemma 17, \citet{jin2017escape} lets us conclude. We detail it for completeness. Similar to the proof of Lemma 17, \citet{jin2017escape}, let $S$ be the subspace corresponding to eigenvectors of $\matH$ with eigenvalues larger or equal in absolute value to $\gamma$, and let $S^{\perp}$ be its orthogonal complement. Note $\vecE_1 \subseteq S$. Denote the norm of $\vecV_t$ projected onto $S$ by $\psi_t$, and the norm of $\vecV_t$ projected onto $S^{\perp}$ by $\phi_t$. 

Notice therefore from the assumptions of this Lemma that $\phi_0=0$ as $\vecV_0$ is a scalar multiple of $\vecE_1$. Similarly, note $\psi_0 = \nrm*{\vecV_0} \ge \frac{\delta}{2\sqrt{d}} \cdot r$ by \pref{eq:perturbedgdlem2v0bound}.

Let 
\[ B  := \eta L_2(\vecW_0) \cL (450\hat{c} + 1).\]
Observe $B \le 1$, as $\cL L_2(\vecW_0) \le 1$ and as $\eta \le \cmax \le \frac1{12100}$, $\hat{c}=11$.

Combining \pref{eq:perturbedgdlem2vtrecurse} with \pref{eq:perturbedgdlem2deltatbound} gives that
\[
\psi_{t+1} \geq (1 + \gamma\eta)\psi_t - B \sqrt{\psi_t^2 + \phi_t^2},
\phi_{t+1} \leq (1 + \gamma\eta)\phi_t + B \sqrt{\psi_t^2 + \phi_t^2}.\numberthis\label{eq:perturbedgdlem2recursions}
\]
\paragraph{The key induction:} Now we induct on $t$ to show that for all $t < T$,
\[
\phi_t \leq 4 B t \cdot \psi_t. 
\]
For the base case, recall by hypotheses of the Lemma that $\vecV_0$ is a scalar multiple of $\vecE_1$, thus $\phi_0=0$ and the base case holds. 

Now, for the inductive step, assume that the inductive hypothesis holds true for all $\tau \leq t$ for some $t$ such that $t+1 \le T$. Substituting the inequality \pref{eq:perturbedgdlem2recursions} for $\phi_{t+1}$ and applying the inductive hypothesis $\phi_t \leq 4 B t \cdot \psi_t$, we obtain
\[ \phi_{t+1} \le 4 B t(1+\gamma \eta) \psi_t + B\sqrt{\psi_t^2 + \phi_t^2}.\]
Also note \pref{eq:perturbedgdlem2recursions} gives
\begin{align*}
4B(t+1) \psi_{t+1} \ge 4 B (t+1)\prn*{(1+\gamma \eta)\psi_t - B \sqrt{\psi_t^2 + \phi_t^2}},
\end{align*}
which rearranges to 
\begin{align*}
4 B t(1+\gamma \eta) \psi_t \le 4 B(t+1)\psi_{t+1} + 4 B^2 (t+1)\sqrt{\psi_t^2+\phi_t^2} - 4 B (1+\gamma \eta)\psi_t.
\end{align*}
Therefore,
\[ \phi_{t+1} \le 4B(t+1)\psi_{t+1} + \prn*{4B^2(t+1)\sqrt{\psi_t^2+\phi_t^2} + B \sqrt{\psi_t^2 + \phi_t^2} - 4B (1+\gamma \eta)\psi_t}.\]
Thus, recalling $B\le 1$, to complete the induction it suffices to show the following:
\[
\prn*{1 + 4B^2(t + 1)}\sqrt{\psi_t^2 + \phi_t^2} \leq 4(1 + \gamma\eta)\psi_t.
\]
Choosing $\sqrt{c_{\max}} \leq \frac{1}{450\hat{c} + 1} \min\left\{\frac{1}{2\sqrt{2}}, \frac{1}{4\hat{c}}\right\}$ which is a universal constant, and choosing $\eta \leq \frac{c_{\max}}{L_1(\vecW_0)}$, we have:
\[
4B(t + 1) \leq 4B T \leq 4\eta L_2(\vecW_0) \cL (450\hat{c} + 1)\hat{c}\cF = 4\sqrt{\eta L_1(\vecW_0)} (450\hat{c} + 1)\hat{c} \leq 1.
\]
By the inductive hypothesis, this gives $\phi_t \le \psi_t$. In turn this implies that
\[
4(1 + \gamma\eta)\psi_t \geq 4\psi_t \geq 2\sqrt{2}\psi_t \geq (1 + 4B(t + 1))\sqrt{\psi_t^2 + \phi_t^2},
\]
finishing the induction.

\paragraph{Finishing the proof from here:} We thus obtain $\phi_t \le 4Bt \psi_t \le \psi_t$ for all $t$, where we use that $4BT \le 1$ as proven above, which just follows from our choice of parameters. Therefore,
\[ \psi_{t+1} \ge (1+\gamma \eta) \psi_t - B\sqrt{2}\psi_t > \prn*{1+\frac{\gamma \eta}2}\psi_t.\numberthis\label{eq:technicallemma2laststep}\]
The last step follows upon noting $B \le \eta L_2(\vecW_0) \cL(450\hat{c}+1) \le \sqrt{\cmax} (450\hat{c}+1) \gamma \eta \log^{-1}\prn*{\frac{d\kappa}{\delta}} < \frac{\gamma \eta}{2\sqrt{2}}$. The inequality is strict as $\gamma \eta >0$.

Finally, recalling that $\nrm*{\vecV_t} \le 300\cL\hat{c}$, $\psi_0 \ge \frac{\delta}{2\sqrt{d}} \cdot r$ and using \pref{eq:technicallemma2laststep}, we have for all $t < T$:
\begin{align*}
300(\cL \cdot \hat{c}) &\geq \nrm*{\vecV_t} \\
&\geq \psi_t \\
&> \prn*{1 + \frac{\gamma\eta}{2}}^t \psi_0 \\
&\ge \prn*{1 + \frac{\gamma\eta}{2}}^t \cdot \frac{\delta}{2\sqrt{d}} \cdot \frac{\cL}{\kappa \cdot \log\left(\frac{d \kappa}{\delta}\right)}.\numberthis\label{eq:perturbgdlem2laststepexponential}
\end{align*}
Note that $\delta \in \left(0, \frac{d \kappa}{e}\right]$ implies $\log\left(\frac{d \kappa}{\delta}\right) \geq 1$. Applying \pref{eq:perturbgdlem2laststepexponential} for $t=T-1$ we obtain:
\begin{align*}
T &< 1+ \log\left(600\kappa \sqrt{d}\delta^{-1} \cdot \hat{c} \log\left(\frac{d \kappa}{\delta}\right)\right) \cdot \log^{-1}\prn*{1 + \frac{\gamma\eta}{2}}\\
&\leq 1+2.01\log\left(600\kappa \sqrt{d}\delta^{-1} \cdot \hat{c} \log\left(\frac{d \kappa}{\delta}\right)\right) \cdot \frac{1}{\gamma\eta}\\
&\le 1+2.01\prn*{\log(600\hat{c})+1.01\log(d\kappa/\delta)} \cdot \frac1{\gamma \eta}\\
&\le \prn*{\frac1{40}+1+2.0301}\cF_2 \le \hat{c}\cF_2.
\end{align*}
These last steps follow by:
\begin{itemize}
    \item Taking $\cmax$ a small enough universal constant so that $\gamma \eta \le \frac{1}{60} \cdot \frac{\cmax}{L_1(\vecW_0)} \le \frac{\cmax}{60}$ satisfies $\frac{2.01}{x} > \log^{-1}(1+x/2)$, which is valid for all $0<x<0.02$. 
    \item \pref{rem:perturbedgdremarksimplify}, which states that we can assume WLOG $\log(d\kappa/\delta)$ is larger than a universal constant. In particular we can assume WLOG that $\log(d\kappa/\delta)$ solves $\log x < x^{0.01}$ (hence $\log(\kappa \sqrt{d} \delta^{-1} \log(d\kappa/\delta)) \le 1.01\log(d\kappa/\delta)$), that $2.01 \log(600\hat{c}) = 2.01 \log(6600) \le \log(d\kappa/\delta)$ (recall $\hat{c}=11$), and that $\cF_2=\frac{\log(d\kappa/\delta)}{\gamma \eta} \ge 40$.
\end{itemize}
This completes the proof.
\end{proof}

\section{Restarted SGD finding Second Order Stationary Points}\label{sec:restartedsgdproofs}
Here, we formally prove \pref{thm:sgdsecondorder}. We formally instantiate \pref{alg:restartedsgd} here. One may notice a slight difference in \pref{alg:restartedsgd} vs the algorithm of \citet{fang2019sharp}: we artificially inject bounded noise at a particular scale $\tilde{\sigma}$. This ensures we can escape saddle points that are in the $F(\vecW_0)$-sublevel set $\cL_{F,F(\vecW_0)}$. Note we may not be able to escape saddle points that are not in $\cL_{F,F(\vecW_0)}$, but that does not matter thanks to our framework \pref{thm:generalframework}, which effectively lets us consider only behavior within $\cL_{F, F(\vecW_0)}$. Also note a practitioner can find such a noise scaling $\tilde{\sigma}$ (depending on suboptimality at initialization $F(\vecW_0)$) via appropriate cross-validation.

The general proof strategy here is similar to the way we adapted the proof of \citet{jin2017escape} in \pref{sec:perturbedgdproofs}. Namely, we use the self-bounding regularity conditions to control the derivatives of $F$ in appropriate neighborhoods of the $F(\vecW_0)$-sublevel set $\cL_{F,F(\vecW_0)}$. 

\begin{algorithm}[tb]
   \caption{Restarted SGD, from \citet{fang2019sharp}}
   \label{alg:restartedsgd}
\begin{algorithmic}
   \STATE Initialize at $\vecW_0$, and consider $K_0 = \Tilde{\Theta}\prn*{\epsilon^{-2}}$, $\eta = \Tilde{\Theta}\prn*{\epsilon^{1.5}}$, $B = \Tilde{\Theta}\prn*{\epsilon^{0.5}}$, $\tilde{\sigma} = 2 \sigma'_1(\vecW_0)$, all explicitly defined in \pref{subsec:restartedsgdparamsnotation}.
   \STATE Let $t=0$ (the total number of iterates), $k=0$ (the restart counter), $\vecX^0 = \vecW_0$ (the point we consider the escape from).
   \WHILE{$k<K_0$}
        \STATE Let $\vecX^{t+1}=\vecX^t - \eta \prn*{\grad f(\vecX^t; \vecZeta_{t+1})+\tilde{\sigma} \vecLambda^{t+1}}$, where $\vecLambda^{t+1}$ is uniform from $\ball(\vecOrigin, 1)$ and independent of everything else, and $\vecZeta_{t+1}$ is an i.i.d. minibatch sample
        \STATE $t \leftarrow t+1$, $k \leftarrow k+1$
        \IF{$\nrm*{\vecX^k-\vecX^0} > B$}
            \STATE $\vecX^0 \leftarrow \vecX^k$, $k \leftarrow 0$
        \ENDIF
   \ENDWHILE
   \STATE Return $\frac1{K_0} \sum_{k=0}^{K_0-1} \vecX^k$
\end{algorithmic}
\end{algorithm}

\subsection{Notation and Parameters}\label{subsec:restartedsgdparamsnotation}
We set the parameters of the algorithm as follows. We will highlight the significance of these parameters in \pref{subsec:restartedsgdpreliminaries}.

\paragraph{Noise Parameters:}
Define
\begin{align*}
\sigma'(\vecW_0) &= \sigma\prn*{F(\vecW_0)+1}. \numberthis\label{eq:tildesigma1defnew} \\
\tilde{\sigma} &= 2\sigma'(\vecW_0).\numberthis\label{eq:tildesigmadefnew}\\
\sigma_1(\vecW_0) &= \max\crl*{\sigma'(\vecW_0)+\tilde{\sigma}, 1}. \numberthis\label{eq:sigma1defnew}
\end{align*}
Note this only depends on $\rho_0$ (and therefore only on $\rho_1)$ and $F(\vecW_0)$. Note $\tilde{\sigma} \in [\sigma'(\vecW_0), 2\sigma'(\vecW_0)]$.\footnote{In fact, this is the only condition we need on $\tilde{\sigma}$. In practice, such a $\tilde{\sigma}$ by fine-enough cross validation in terms of only $F(\vecW_0)$.} Also note $\sigma_1(\vecW_0) \le 3\sigma'(\vecW_0)$.

\paragraph{Update Rule:}
Define
\[ \grad \tilde{f}(\vecX^t;\vecZeta_{t+1}) := \grad f(\vecX^t;\vecZeta_{t+1})+\tilde{\sigma}\vecLambda^{t+1}.\label{eq:gradoraclenew}\]
Thus the SGD update rule in \pref{alg:restartedsgd} (without considering the restarts) is $\vecX^{t+1} = \vecX^t - \eta \grad \tilde{f}(\vecX^t;\vecZeta_{t+1})$. Note the slight abuse of notation; $\grad \tilde{f}(\vecX^t;\vecZeta_{t+1})$ is not necessarily an actual gradient.\footnote{This choice of notation is made to demonstrate the artificial noise injections $\tilde{\sigma}\vecLambda^{t+1}$ are not fundamentally needed. They are not necessary if the stochastic gradient $\grad f(\cdot;\cdot)$ enjoys suitable anticoncentration properties.} This will not cause issues or ambiguity for the rest of this section.

\paragraph{Effective Smoothness Parameters in $F(\vecW_0)$-sublevel set:} 
We define the `local smoothness parameters' as follows, slightly differently compared to the proof of \pref{thm:escapesecondordergd}. Define
\begin{align*}
L_1(\vecW_0) &:= \max\crl*{1, \rho_1\prn*{F(\vecW_0)+1}, \rho_3(\rho_0(F(\vecW_0)+1) + \sigma'(\vecW_0), F(\vecW_0)+1)}, \numberthis\label{eq:L1defnew}\\
L_2(\vecW_0) &:= \max\crl*{1, \rho_2\prn*{F(\vecW_0)+1}, \rho_0\prn*{F(\vecW_0)+1}^2 \max\crl*{4,\prn*{\sigma_1(\vecW_0)+\rho_0(F(\vecW_0)+1)}^2} }. \numberthis\label{eq:L2defnew}
\end{align*}
Note all of these parameters only depend on $F(\vecW_0)$, through $\rho_1(\cdot), \rho_2(\cdot), \rho_3(\cdot, \cdot)$ (recall $\rho_0(\cdot)$ can be defined in terms of $\rho_1(\cdot)$). 

\paragraph{Parameters of \pref{alg:restartedsgd}:} We define the remaining parameters of \pref{alg:restartedsgd} as follows. Consider any $\epsilon > 0$ and $p \in (0,1)$. We choose:
\begin{align*}
\tilde{C}_1 &= 2\lfloor \frac{\log (3/p)}{\log(0.8^{-1})}+1\rfloor \log \prn*{\frac{24\sqrt{d}}{\eta}},\\
\delta &= \sqrt{L_2(\vecW_0) \epsilon},\\
\delta_2 &= 16\delta, \\
B &= \frac{\delta}{L_2(\vecW_0) \tilde{C_1}},\\
K_0 &= \tilde{C}_1 \eta^{-1} \delta_2^{-1},\\
\eta &\le \frac{B^2 \delta}{512 \max(\sigma_1(\vecW_0)^2, 1) \tilde{C_1} \log(48 K_0/p)} \cdot \frac1{3(1+\log(K_0))}.\numberthis\label{eq:paramchoice}
\end{align*}
Also define 
\[ K_o = 2 \log \prn*{\frac{24\sqrt{d}}{\eta}}\eta^{-1} \delta_2^{-1}, \text{ thus } K_0 = \lfloor \frac{\log (3/p)}{\log(0.8^{-1})}+1\rfloor K_o.\]
\begin{remark}
To choose $\eta$ satisfying the above inequality, one can perform the same analysis as on footnote 4, page 7 of \citet{fang2019sharp}. We first choose $\tilde{\eta}$ appropriately by setting
\[ \tilde{\eta} = \frac{B^2 \delta}{4096 \max(\sigma_1(\vecW_0)^2, 1) \log(48/p) \log(p) \lfloor \frac{\log (3/p)}{\log(0.8^{-1})}+1\rfloor},\]
and then set $\eta = \tilde{\eta} \log^{-3}(1/\tilde{\eta})$. 
\end{remark}

\begin{remark}\label{rem:restartedsgdreduceparams}
Analogously to the proof of \pref{thm:escapesecondordergd}, note it suffices to show the result for $\epsilon \le \frac1{L_2(\vecW_0)}$; for $\epsilon > \frac1{L_2(\vecW_0)}$, we can just apply the result for $\epsilon=\frac1{L_2(\vecW_0)}$, and the result remains the same up to $F(\vecW_0)$-dependent parameters in the $O(\cdot)$. Thus we can suppose that $\delta_2$ (and $\delta$) are at most some universal constant. We also can take $L_1(\vecW_0), L_2(\vecW_0), \sigma_1(\vecW_0)$ to be the max between their currently definition and an appropriate universal constant. Thus due to the choice of parameters above, we may assume that
\begin{align*}
\tilde{C}_1, K_0 &\ge 1,\\
\log(K_0), \sigma_1(\vecW_0) &\ge 1,\\
B &\le \min\prn*{1, \frac{\sigma_1(\vecW_0)}{L_1(\vecW_0)}, \frac1{L_1(\vecW_0)}, \frac1{L_2(\vecW_0)}}, \\
\eta &\le \min\crl*{1, \frac1{\sigma_1(\vecW_0)^2}}.\\
\end{align*}
From here note we have $\eta L_1(\vecW_0) \le 1$. As \textit{these assumptions come with no loss of generality}, we make these assumptions for the rest of the proof.
\end{remark}

\paragraph{Notation:} Consider a sequence of iterates $\vecX^0, \vecX^1, \ldots$ beginning at $\vecX^0$ comprising an instance of the while loop in \pref{alg:restartedsgd}. For such a sequence, let $\mathfrak{F}^k$ be the $\sigma$-algebra defined by all the prior iterates and the noise up through $\vecX^k$, namely $\sigma\crl*{\vecX^0, \vecZeta_1, \vecLambda^1, \vecX^1, \ldots, \vecX^{k-1}, \vecZeta_{k}, \vecLambda^k}$. Let $\cK_0$ be a stopping time given by 
\[ \cK_0 = \inf_k \crl*{k \ge 0:\nrm*{\vecX^k - \vecX^0} \ge B}.\]
Note $\vecX^k$ and $1_{\cK_0 \ge k}, 1_{\cK_0 > k}$ are $\mathfrak{F}^k$-measurable. Thus, $1_{\cK_0 > k-1} \equiv 1_{\cK_0 \ge k}$ is $\mathfrak{F}^{k-1}$-measurable. 

\subsection{Result}
We now formally prove \pref{thm:sgdsecondorder}. The following \pref{thm:sgdsecondorderformal} can readily be seen to imply \pref{thm:sgdsecondorder}.
\begin{theorem}\label{thm:sgdsecondorderformal}
Suppose $F$ satisfies \pref{ass:thirdorderselfbounding} and the stochastic gradient oracle satisfies \pref{ass:noiseregularity} and \pref{ass:gradoraclecontrol}. Run \pref{alg:restartedsgd} initialized at $\vecW_0$, run with parameters chosen as per \pref{subsec:restartedsgdparamsnotation}. 

Consider any $p \in (0, 1)$. With probability at least $1-\frac74 p \cdot \frac{(F(\vecW_0)+1)7 \eta K_0}{B^2}$, upon making 
\[ K_0+\frac{7 \eta K_0^2 (F(\vecW_0)+1)}{B^2}\text{ oracle calls to }\grad f(\cdot;\cdot), \]
\pref{alg:restartedsgd} will output $\tilde{O}\prn*{\frac{7 \eta K_0^2 (F(\vecW_0)+1)}{B^2}}$ candidate vectors $\vecW$, one of which satisfies
\[ \nrm*{\grad F(\vecW)} \le 18 L_2(\vecW_0) B^2, \lambdamin(\grad^2 F(\vecW)) \ge -17 \delta.\]
\end{theorem}
\begin{remark}
Before proceeding, we justify why \pref{thm:sgdsecondorderformal} implies \pref{thm:sgdsecondorder}. Simply take $\epsilon\leftarrow\frac{\epsilon}{289L_2(\vecW_0)}$ in \pref{thm:sgdsecondorderformal}. 
Plugging this in, we obtain a result on finding a SOSP as per the definition in \pref{eq:sospproblem}.\footnote{Recall this definition refers to $\vecW$ such that $\nrm*{\grad F(\vecW)} \le \epsilon, \grad^2 F(\vecW) \succeq -\sqrt{\epsilon}\matI$.}
The oracle complexity has the desired dependence on $\epsilon$ and polylog dependence on $d,p$. The probability is at least $1-p\cdot \Tilde{\Theta}(\epsilon^{-1.5})$, where the $\Tilde{\Theta}$ are hiding polylog terms in $d,1/\epsilon,1/p$ and dependence on $F(\vecW_0)$ (through $\rho_1(\cdot), \rho_2(\cdot), \rho_3(\cdot), \sigma(\cdot)$).
This holds for any $p \in (0,1)$. 

Now consider the final desired success probability $1-\tilde{\delta}$ governed in terms of $\tilde{\delta} \in (0,1)$ in \pref{thm:sgdsecondorder}.
Let $p=\tilde{\delta}\epsilon^{1.5} \cdot \text{polylog}(d, 1/\epsilon)$ in the guarantee from the above paragraph. This gives \pref{thm:sgdsecondorder}, with the requested probability and oracle complexity.
\end{remark}

We now prove \pref{thm:sgdsecondorderformal} via our framework, \pref{thm:generalframework}. 
\begin{proof}[Proof of \pref{thm:sgdsecondorderformal} and thus \pref{thm:sgdsecondorder}]
We again use our framework \pref{thm:generalframework}. Consider any $p\in (0,1)$, and choose parameters as per \pref{subsec:restartedsgdparamsnotation}. 

Let 
\[ \cS=\{\vecW: \nrm*{\grad F(\vecW)} \le 18 L_2(\vecW_0) B^2, \lambdamin(\grad^2 F(\vecW)) \ge -17 \delta\}. \]
Define $\cA$ as follows, identically to how we defined them for Restarted SGD in \pref{subsec:exinframework}. Consider any given $\vecU_0 \in \mathbb{R}^d$. Let $\vecP_0=\vecU_0$. We define a sequence $(\vecP_i)_{0 \le i \le K_0}$ via $\vecP_i = \vecP_{i-1} - \eta (\grad f(\vecP_{i-1};\vecZeta_{i})+\tilde{\sigma} \Lambda^{i})$. Note this sequence can be equivalently defined by repeatedly composing the function $\vecU \rightarrow \vecU - \eta (\grad f(\vecU; \vecZeta)+\tilde{\sigma}\Lambda)$.

If it exists, let $i, 1 \le i \le K_0$ be the minimal index such that $\nrm*{\vecP_i - \vecP_0} > B$. Otherwise let $i=K_0$. In either case, we define
\[ \cA(\vecU_0) = \prn*{ \vecP_i, \frac1{i} \sum_{t=0}^{i-1} \vecP_t}, \text{ hence }\cA_1(\vecU_0)=\vecP_i, \cA_2(\vecU_0) = \frac1{i} \sum_{t=0}^{i-1} \vecP_t. \]
We now let 
\[ \Toracle(\vecU_0)=K_0,\text{ and }\Delta = \frac{B^2}{7\eta K_0}.\] 
Following the notation from \pref{alg:restartedsgd}, notice that $\cA(\vecU_0)$ corresponds to next vector set to $\vecX^0$ in the while loop of \pref{alg:restartedsgd}, when the while loop begins at $\vecX^0 = \vecU_0$.

Crucial to this proof are the following two Lemmas. While inspired from \citet{fang2019sharp}, a crucial difference is that \textit{they hold only in the $F(\vecW_0)$-sublevel set $\cL_{F,F(\vecW_0)}$}. 

\begin{lemma}[Equivalent of Proposition 10, \citet{fang2019sharp}]\label{lem:sgdfindingssp}
Consider $\vecX^0$ in the while loop of \pref{alg:restartedsgd}. Suppose $\vecX^0 \in \cL_{F,F(\vecW_0)}$. With probability at least $1-p$, if $\vecX^k$ does not move out of the ball $\ball(\vecX^0, B)$ within the first $K_0$ iterations in the while loop of \pref{alg:restartedsgd}, letting $\vecXbar=\frac1{K_0}\sum_{k=0}^{K_0-1} \vecX^k$, we have
\[ \nrm*{\grad F(\vecXbar)} \le 18 L_2(\vecW_0) B^2, \lambdamin(\grad^2 F(\vecXbar)) \ge -17 \delta.\]
\end{lemma}
\begin{lemma}[Equivalent of Proposition 9, \citet{fang2019sharp}]\label{lem:fasterdescentlemma}
Consider $\vecX^0$ in the while loop of \pref{alg:restartedsgd}. Suppose $\vecX^0\in \cL_{F,F(\vecW_0)}$. With probability at least $1-\frac34 p$, if $\vecX^k$ moves out of $\ball(\vecX^0, B)$ in $K_0$ iterations or fewer in the while loop of \pref{alg:restartedsgd}, we have 
\[ F(\vecX^{\cK_0}) < F(\vecX^0) - \frac{B^2}{7 \eta K_0}.\]
\end{lemma}

\paragraph{Finishing the proof:} The main point is to prove the following Claim.
\begin{claim}\label{claim:restartedsgddecreasesequence}
For any $\vecU_0 \in \cL_{F,F(\vecW_0)}$, $\cA$ is a $(\cS, K_0, \Delta, \frac74 p, \vecU_0)$-decrease procedure. 
\end{claim}
\begin{proof}[Proof of \pref{claim:restartedsgddecreasesequence}]
Apply \pref{lem:sgdfindingssp} and \pref{lem:fasterdescentlemma} to the sequence $(\vecP_i)_{0 \le i \le K_0}$, recalling that $\cA(\vecU_0)$ corresponds to next vector set to $\vecX^0$ in the while loop of \pref{alg:restartedsgd} when the while loop begins at $\vecX^0 = \vecP_0 = \vecU_0$. By a Union Bound over the events of \pref{lem:sgdfindingssp} and \pref{lem:fasterdescentlemma}, with probability at least $1-\frac74 p$, we have the following:
\begin{itemize}
\item Suppose there exists $t < K_0$ such that $\vecP_t \not\in \ball(\vecP_0, B)=\ball(\vecU_0, B)$. Let $t'$ be the minimal such $t$. By \pref{lem:fasterdescentlemma}, we have 
\[ F(\cA_1(\vecU_0)) = F(\vecP_{t'}) \le F(\vecP_0) - \frac{B^2}{7 \eta K_0}=F(\vecU_0)-\Delta.\]
\item Otherwise, we have $\cA_2(\vecU_0) =\overline{\vecP}$ where $\overline{\vecP} = \frac{1}{K_0} \sum_{k=0}^{K_0-1} \vecP_k$. In this case, by \pref{lem:sgdfindingssp}, we have 
\[ \cA_2(\vecU_0)  = \overline{\vecP} \in\cS.\]
\end{itemize}
Consequently, $\cA$ is a $(\cS, K_0, \Delta, \frac74 p, \vecU)$-decrease procedure. 
\end{proof}
Now with \pref{claim:restartedsgddecreasesequence}, directly applying \pref{thm:generalframework} and plugging in the relevant parameters, we obtain \pref{thm:sgdsecondorderformal}.
\end{proof}
\begin{remark}
To sanity check these results, note the rate from \pref{lem:fasterdescentlemma} will get worse as $\eta$ gets smaller because $K_0 \eta = 2\lfloor \frac{\log (3/p)}{\log(0.8^{-1})}+1\rfloor \log \prn*{\frac{24\sqrt{d}}{\eta}} \delta_2^{-1}$ will increase as $\eta$ gets smaller.
\end{remark}
The rest of \pref{sec:restartedsgdproofs} will now be devoted to the proofs of \pref{lem:sgdfindingssp} and \pref{lem:fasterdescentlemma}. For the rest of \pref{sec:restartedsgdproofs}, we suppose $F$ satisfies \pref{ass:thirdorderselfbounding} and the stochastic gradient oracle satisfies \pref{ass:noiseregularity} and \pref{ass:gradoraclecontrol}. These proofs are similar to that of \citet{fang2019sharp}, but hinges crucially on the fact that the analysis in \citet{fang2019sharp} is `local'.

\subsection{Preliminaries}\label{subsec:restartedsgdpreliminaries}
We now establish useful properties of the parameters of the algorithm defined in \pref{subsec:restartedsgdparamsnotation}, analogously to \pref{lem:formalizingperturbationopnrom}.

\paragraph{Locality of balls $\ball(\vecX^0, B)$:}
\begin{lemma}\label{lem:restartedsgdballcloseenough}
We have $B \le \frac1{2\rho_0(F(\vecW_0)+1)}$. In particular, for any $\vecU \in \ball(\vecW, B)$ for $\vecW \in \cL_{F,F(\vecW_0)}$, we have $\nrm*{\vecU-\vecW} \le \frac1{2\rho_0(F(\vecW_0)+1)} \le \frac1{2\rho_0(F(\vecW)+1)}$. 
\end{lemma}
\begin{proof}
As per \pref{rem:restartedsgdreduceparams}, we have $\epsilon \le 1$. Thus by the choice of parameters in \pref{eq:L2defnew},
\[ B \le \frac{\delta}{L_2(\vecW_0)} \le \frac1{\sqrt{L_2(\vecW_0)}} \le \frac1{2\rho_0(F(\vecW_0)+1)}.\]
This completes the proof.
\end{proof}

\paragraph{Control over the stochastic gradient oracle:}
\begin{lemma}\label{lem:noisechangeold}
For all $\vecU$ such that $\vecU \in \ball\prn*{\vecW, \frac1{\rho_0(F(\vecW_0)+1)}}$ for $\vecW \in \cL_{F,F(\vecW_0)}$, we have $\nrm*{\grad f(\vecU;\vecZeta) - \grad F(\vecU)} \le \sigma'(\vecW_0)$ for all $\vecZeta$.
\end{lemma}
\begin{proof}
By \pref{ass:noiseregularity}, we have 
\[ \nrm*{\grad f(\vecU;\vecZeta) - \grad F(\vecU)} \le \sigma(F(\vecU)).\]
Now as $\vecW \in \cL_{F,F(\vecW_0)}$, we have 
\[ \frac1{\rho_0(F(\vecW_0)+1)} \le \frac1{\rho_0(F(\vecW)+1)}.\]
Thus by \pref{lem:boundfuncvalueradius} and again as $\vecW \in \cL_{F,F(\vecW_0)}$, we have 
\[ F(\vecU) \le F(\vecW)+1 \le F(\vecW_0)+1.\]
Combining these gives \pref{lem:noisechangeold}.
\end{proof}
\begin{lemma}\label{lem:noisechange}
For all $\vecU$ such that $\vecU \in \ball\prn*{\vecW, \frac1{\rho_0(F(\vecW_0)+1)}}$ for $\vecW \in \cL_{F,F(\vecW_0)}$, $\nrm*{\grad \tilde{f}(\vecU;\vecZeta) - \grad F(\vecU)} \le \sigma_1(\vecW_0)$ for all $\vecZeta$.
\end{lemma}
\begin{proof}
This immediately follows from \pref{lem:noisechangeold} and the definition of $\grad \tilde{f}(\vecU;\vecZeta)$, as $\nrm*{\tilde{\sigma}\vecLambda^t}\le \tilde{\sigma}$.
\end{proof}

\paragraph{Locality after one step of SGD:}
\begin{lemma}\label{lem:restartedsgdonestepcloseenough}
Consider any $\vecU \in \ball(\vecW, B)$ for $\vecW \in \cL_{F,F(\vecW_0)}$. Then for all points $\vecP$ in the line segment between $\vecU$ and $\vecU - \eta \grad \tilde{f}(\vecU; \vecZeta)$ for any $\vecZeta$, we have $\vecP\in\ball\prn*{\vecW, \frac1{\rho_0(F(\vecW_0)+1)}}$.
\end{lemma}
\begin{proof}
It suffices to show $\vecU - \eta \grad \tilde{f}(\vecU; \vecZeta) \in \ball\prn*{\vecW, \frac1{2(\rho_0(F(\vecW_0)+1))}}$; after establishing this, the result then follows by Triangle Inequality and \pref{lem:restartedsgdballcloseenough}. To this end, by Triangle Inequality, it suffices to show that 
\[ \eta \nrm*{\grad \tilde{f}(\vecU; \vecZeta)} \le \frac1{2\rho_0(F(\vecW_0)+1)}. \]
Indeed, the same reasoning as in the proof of \pref{lem:restartedsgdballcloseenough} gives 
\[ F(\vecU) \le F(\vecW_0)+1.\]
Thus, \pref{ass:gradoraclecontrol} gives
\[ \nrm*{\grad F(\vecU)} \le \rho_0(F(\vecW_0)+1),\]
and so \pref{lem:noisechange} gives
\[ \nrm*{\grad \tilde{f}(\vecU;\vecZeta)} \le \sigma_1(\vecW_0)+\rho_0(F(\vecW_0)+1).\]
As per \pref{rem:restartedsgdreduceparams}, we have
\[ \eta \le \frac12 B^2 \delta \le \frac12 \cdot \frac{\delta^3}{L_2(\vecW_0)^2} \le \frac1{2L_2(\vecW_0)^{0.5}}.\]
Combining all the above gives
\begin{align*}
\eta \nrm*{\grad \tilde{f}(\vecU; \vecZeta)}  &\le  \frac1{2 L_2(\vecW_0)^{0.5}} \cdot \prn*{\sigma_1(\vecW_0)+\rho_0(F(\vecW_0)+1)} \\
&\le \frac1{2\rho_0(F(\vecW_0)+1)\prn*{\sigma_1(\vecW_0)+\rho_0(F(\vecW_0)+1)}} \cdot \prn*{\sigma_1(\vecW_0)+\rho_0(F(\vecW_0)+1)}\\
&\le \frac1{2\rho_0(F(\vecW_0)+1)},
\end{align*}
which by our earlier remarks completes the proof.
\end{proof}

\paragraph{Properties of the effective smoothness parameters:}
\begin{lemma}\label{lem:smoothnesschange}
Consider any $\vecX^0 \in \cL_{F,F(\vecW_0)}$. Then we have $\nrm*{\grad^2 F(\vecU)}_{\OPNORM} \le L_1(\vecW_0)$ for all $\vecU$ such that either:
\begin{itemize}
    \item $\vecU \in \ball(\vecX^0, B)$,
    \item Or $\vecU$ lies in the line segment between some $\vecU' \in \ball(\vecX^0, B)$ and $\vecU' - \eta \grad \tilde{f}(\vecU';\vecZeta)$, for any $\vecZeta$.
\end{itemize}
\end{lemma}
\begin{proof}
By \pref{lem:restartedsgdballcloseenough} and \pref{lem:restartedsgdonestepcloseenough}, irrespective of which case for $\vecU$ in the conditions of \pref{lem:smoothnesschange} holds, we have 
\[ \vecU \in \ball\prn*{\vecX^0, \frac1{\rho_0(F(\vecW_0)+1)}}.\]
As $\vecX^0 \in \cL_{F,F(\vecW_0)}$, this implies
\[ \nrm*{\vecU-\vecX^0} \le \frac1{\rho_0(F(\vecW_0)+1)} \le \frac1{\rho_0(F(\vecX^0)+1)}.\]
By \pref{lem:boundfuncvalueradius} and as $\vecX^0 \in \cL_{F,F(\vecW_0)}$, it follows that 
\[ F(\vecU) \le F(\vecX^0)+1 \le F(\vecW_0)+1.\]
The conclusion now follows by \pref{ass:selfbounding}.
\end{proof}

\begin{lemma}\label{lem:hessianlipschitzchange}
Consider any $\vecX^0 \in \cL_{F,F(\vecW_0)}$. Consider any $\vecU_1, \vecU_2$ such that each $\vecU_i$, $i=1,2$ is such that either:
\begin{itemize}
    \item $\vecU_i \in \ball(\vecX^0, B)$,
    \item Or $\vecU_i$ lies in the line segment between some $\vecU' \in \ball(\vecX^0, B)$ and $\vecU' - \eta \grad \tilde{f}(\vecU';\vecZeta)$, for any $\vecZeta$.
\end{itemize}
Then 
\[\nrm*{\grad^2 F(\vecU_1) - \grad^2 F(\vecU_2)}_{\OPNORM} \le L_2(\vecW_0)\nrm*{\vecU_1 - \vecU_2}.\]
\end{lemma}
\begin{proof}
Irrespective of which condition applies to $\vecU_i$, By \pref{lem:restartedsgdballcloseenough} and \pref{lem:restartedsgdonestepcloseenough}, we have 
\[ \vecU_i \in \ball\prn*{\vecX^0, \frac1{\rho_0(F(\vecW_0)+1)}}\]
for $i=1,2$. Thus the line segment $\overline{\vecU_1 \vecU_2}$ is contained in $\ball\prn*{\tilde{\vecW}, \frac1{\rho_0(F(\vecW_0)+1)}}$. As $\vecX^0 \in \cL_{F,F(\vecW_0)}$, the result now follows from applying \pref{lem:thirdorderlipschitzlem} and \pref{lem:boundfuncvalueradius}. 
\end{proof}
\begin{remark}\label{rem:restartedsgdpropertiesonestep}
The reason for the second case in the condition on $\vecU$ or $\vecU_i$ from \pref{lem:smoothnesschange}, \pref{lem:hessianlipschitzchange} will become clear in the proof of \pref{lem:fasterdescentlemma}. In particular, to prove \pref{lem:fasterdescentlemma}, we will consider $\vecU - \eta \grad \tilde{f}(\vecU; \vecZeta)$ for $\vecU \in \ball\prn*{\vecX^0, B}$ where $\vecX^0 \in\cL_{F,F(\vecW_0)}$.
\end{remark}

\begin{lemma}\label{lem:gradoraclelipschitzconstant}
Consider any $\vecX^0 \in \cL_{F,F(\vecW_0)}$. Then for any $\vecU \in \ball(\vecX^0, B)$ and any $\vecZeta$, 
\[ \nrm*{\grad^2 f(\vecU;\vecZeta)}_{\OPNORM} \le L_1(\vecW_0).\]
\end{lemma}
\begin{proof}
By \pref{lem:restartedsgdballcloseenough}, we have 
\[ \vecU \in \ball\prn*{\vecX^0, \frac1{\rho_0(F(\vecW_0)+1)}}.\]
By \pref{lem:boundfuncvalueradius}, because $\vecX^0 \in \cL_{F,F(\vecW_0)}$, we have
\[ F(\vecU) \le F(\vecW_0)+1.\]
Moreover, as $\vecX^0 \in \cL_{F,F(\vecW_0)}$ and by \pref{lem:noisechangeold} and \pref{corr:gradcontrol}, 
\[ \nrm*{\grad f(\vecU; \vecZeta)} \le \nrm*{\grad F(\vecU)} + \sigma'(\vecW_0) \le \rho_0(F(\vecW_0)+1) + \sigma'(\vecW_0).\]
Thus the result follows from \pref{ass:gradoraclecontrol}.
\end{proof}
\begin{remark}
While \pref{lem:gradoraclelipschitzconstant} is phrased as an upper bound on the operator norm of $\grad^2 f(\cdot;\vecZeta)$, it can be easily phrased in terms of the local Lipschitz constant of $\grad f(\cdot;\vecZeta)$, similar to one of the possibilities in \pref{ass:thirdorderselfbounding}.
\end{remark}

\paragraph{Enough noise to escape saddles:}
Now we verify that the noise scheme here gives us enough noise to escape saddle points in the $F(\vecW_0)$-sublevel set $\cL_{F,F(\vecW_0)}$.
\begin{definition}[$(q^*, \vecV)$-narrow property; Definition 2 in \citet{fang2019sharp}]\label{def:narrowproperty}
A Borel set $\cA \subset \mathbb{R}^d$ satisfies the $(q^{\star}, \vecV)$-narrow property if for any $\vecU \in \cA$, $q \ge q^{\star}$, $\vecU+q\vecV \in \cA^c$.
\end{definition}
Immediately, we obtain the following properties of this definition, as also noted in \citet{fang2019sharp}.
\begin{lemma}\label{lem:restartedsgdnarrowpropertyscaleshift}
If $\cA$ satisfies the $(q^{\star}, \vecV)$-narrow property, then for all $c_1 \in \mathbb{R}^d$, $c_2 \in \mathbb{R}$, $c_1+c_2 \cA$ satisfies the $(|c_2| q^{\star}, \vecV)$-narrow property.
\end{lemma}
We now introduce the following definition:
\begin{definition}[$\vecV$-dispersive Property; Equivalent of Definition 3 in \citet{fang2019sharp}]\label{def:dispersiveroperty}
We say that a random vector $\tilde{\vecXi}$ has the $\vecV$-dispersive property if for any $\cA$ satisfying the $\prn*{\frac{\sigma_1(\vecW_0)}{4\sqrt{d}}, \vecV}$-narrow property, we have 
\[ \mathbb{P}(\tilde{\vecXi}  \in \cA) \le \frac12.\]
\end{definition}
Note the slight change of the constant $\frac12$ rather than $\frac14$ in the above definition compared to that of \citet{fang2019sharp}; this subtle difference will appear in the following proofs, although this will not change too much conceptually. 

Now we prove the following Lemma, which shows that our update rule contains enough noise to escape saddle points:
\begin{lemma}[Dispersive Noise; see also Algorithm 3, \citet{fang2019sharp}]\label{lem:dispersivenoiseregpoints}
The update $\grad \tilde{f}(\vecX^t;\vecZeta_{t+1})$ admits the $\vecV$-dispersive property for all unit vectors $\vecV$, for any $\vecX^t$.
\end{lemma}
Note this does not necessarily hold for the stochastic gradient oracle itself under our assumptions, hence the artificial noise injection of $\tilde{\sigma}\vecLambda^t$.
\begin{proof}[Proof of \pref{lem:dispersivenoiseregpoints}]
First, we prove that the random vector $\tilde{\sigma}\vecLambda^{t+1}$ satisfies the Dispersive Noise property for all unit vectors $\vecV$. Consider any $\cA$ satisfying the $\prn*{\frac{\sigma_1(\vecW_0)}{4\sqrt{d}}, \vecV}$-narrow property. Note we have 
\begin{align*}
\mathbb{P}(\tilde{\sigma} \vecLambda^{t+1} \in \cA)  &= \mathbb{P}\prn*{\vecLambda^{t+1} \in \tilde{\sigma}^{-1}\cA}\\
&\le \frac{\sigma_1(\vecW_0)/4\sqrt{d}}{\tilde{\sigma}} \cdot \frac{\text{Vol}^{d-1}\ball(\vecOrigin, 1)}{\text{Vol}^{d}\ball(\vecOrigin, 1)} \\
&\le \frac{\sigma_1(\vecW_0)/4\sqrt{d}}{\tilde{\sigma}} \cdot \sqrt{d} = \frac{\sigma_1(\vecW_0)}{4\tilde{\sigma}}.    
\end{align*}
Here, the inequality follows from an elementary calculation with multivariate calculus, analogous to the calculation in the proof of \pref{lem:escapesaddlepoints}, which we detailed in full in this article. 
An analogous calculation can also be found in \citet{jin2017escape}, proof of Lemma 14, and in Appendix F, \citet{fang2019sharp}.

Now, note as $\tilde{\sigma} \ge \sigma'(\vecW_0)$, we have
\[ \frac{\sigma_1(\vecW_0)}{4\tilde{\sigma}} \le \frac{\sigma'(\vecW_0)+\tilde{\sigma}}{4\tilde{\sigma}} \le \frac12, \]
and so 
\[ \mathbb{P}(\tilde{\sigma} \vecLambda^{t+1} \in \cA) \le \frac12.\]
Consequently the random vector $\tilde{\sigma}\vecLambda^t$ satisfies the Dispersive Noise property for all unit vectors $\vecV$.

Now, we show that $\grad \tilde{f}(\vecX^t; \vecZeta_{t+1})$ satisfies the $\vecV$-dispersive property as wanted. The proof is analogous to part iii, Proposition 4 of \citet{fang2019sharp}. Consider any unit vector $\vecV$. Recall that $\Lambda^t$ and $\grad f(\vecX^t; \vecZeta_{t+1})$ are independent. Since the $(q^{\star}, \vecV)$-narrow property is evidently preserved with the same parameters by adding a fixed vector to $\cA$, we obtain the following bound on the following conditional probability:
\begin{align*}
\mathbb{P}\prn*{\grad \tilde{f}(\vecX^t; \vecZeta_{t+1}) \in \cA | \grad f(\vecX^t; \vecZeta_{t+1})} &= \mathbb{P}\prn*{\grad f(\vecX^t; \vecZeta_{t+1}) + \tilde{\sigma}\vecLambda^{t+1} 
 \in \cA|\grad f(\vecX^t; \vecZeta_{t+1})} \\
 &= \mathbb{P}\prn*{\tilde{\sigma}\vecLambda^{t+1} \in -\grad f(\vecX^t; \vecZeta_{t+1}) + \cA|\grad f(\vecX^t; \vecZeta_{t+1})} \le \frac12.
\end{align*}
This holds irrespective of conditioning, which implies that $\grad \tilde{f}(\vecX^t; \vecZeta_{t+1})$ satisfies the $\vecV$-dispersive property.
\end{proof}

\subsection{Escaping Saddles}\label{subsec:escapesaddles}
We first aim to prove that we can efficiently escape strict saddle points in the $F(\vecW_0)$-sublevel set, similarly to \citet{fang2019sharp}. In particular, we aim to prove the following \pref{lem:sgdescapesaddles}. The contrapositive of \pref{lem:sgdescapesaddles} will in turn be used to prove \pref{lem:sgdfindingssp}, which establishes that \pref{alg:restartedsgd} can find SOSPs.
\begin{lemma}[Equivalent of Proposition 7 in \citet{fang2019sharp}]\label{lem:sgdescapesaddles}
Consider a sequence of iterates $\vecX^0, \vecX^1, \ldots$ beginning at $\vecX^0$ comprising an instance of the while loop in \pref{alg:restartedsgd}. 
Suppose $\vecX^0 \in \cL_{F,F(\vecW_0)}$ and that $\lambdamin(\grad^2 F(\vecX^0)) \le -\delta_2$ for $\delta_2>0$. 
Then when the while loop of \pref{alg:restartedsgd} is initialized at $\vecX^0$, with probability at least $1-\frac{p}3$, we have
\[ \cK_0 \le K_0 = \lfloor \frac{\log (3/p)}{\log(0.8^{-1})}+1\rfloor K_o.\]
\end{lemma}
\begin{remark}
For $\delta_2$ very small, note the guarantee from \pref{lem:sgdescapesaddles} will deteriorate because $K_0$ scales with $\delta_2^{-1}$. 
\end{remark}

To prove \pref{lem:sgdescapesaddles}, we use the same strategy as in \citet{fang2019sharp}. However, as we do not have global Lipschitzness of the gradient and Hessian, we must be careful. We use that the strategy only requires control over points that are `local', i.e. near $\vecX^0$, since the proof strategy studies escape from the ball $\ball(\vecX^0, B)$. We then appeal to control over $F$ in $\ball(\vecX^0, B)$ that we have by \pref{subsec:restartedsgdpreliminaries}. 
\begin{remark}
In this section \pref{subsec:escapesaddles}, probability is over the samples $\vecZeta_{k}$ and the artificial noise injections $\Lambda^{k}$. 
\end{remark}

Now we go into the details. As in \citet{fang2019sharp}, let $\vecW^k(\vecU)$ be the iterates of SGD starting from a given $\vecU$ using the \textit{same} stochastic samples as $\vecX^k$ and the same noise additions $\tilde{\sigma}\vecLambda^k$. In particular
\[ \vecW^k(\vecU) = \vecW^{k-1}(\vecU) - \eta \grad \tilde{f}(\vecW^{k-1}(\vecU); \vecZeta_k).\]
Thus $\vecX^k = \vecW^k(\vecX^0)$. 

Also for all $\vecU$, let $\cKexit(\vecU)$ be the stopping time defined by
\[ \cKexit(\vecU) := \inf \{k \ge 0:\nrm*{\vecW^k(\vecU) - \vecX^0} > B\}.\]
Thus $\cK_0 = \cKexit(\vecX^0)$. 

The high-level idea from \citet{fang2019sharp}, similar to as in \citet{jin2017escape}, is to consider the `bad initialization region' around $\ball(\vecX^0, B)$ where iterates initialized in this bad region escape with low probability. We then prove that this bad initialization region is `narrow', and consequently we can escape the saddle point efficiently.

In particular, define
\[ \cS^B_{K_o}(\vecX^0) = \{ \vecU \in \mathbb{R}^d: \mathbb{P}(\cKexit(\vecU) < K_o) \le 0.4\}.\]
Note by definition that $\cS^B_{K_0}(\vecX^0) \subseteq \ball(\vecX^0, B)$.

First let $q_0 = \frac{\sigma_1(\vecW_0) \eta}{4\sqrt{d}}$. We establish the following Lemma, which verifies that $\cS^B_{K_o}(\vecX^0)$ is `narrow' in a suitable sense.
\begin{lemma}[Equivalent of Lemma 8 in \citet{fang2019sharp}; also similar to Lemma 15, \citet{jin2017escape}]\label{lem:sgdboundprob}
Suppose the assumptions of \pref{lem:sgdescapesaddles} hold. Let $\vecE_1$ be an arbitrary unit eigenvector of $\grad^2 F(\vecX^0)$ corresponding to its smallest eigenvalue $-\delta_m \le -\delta_2$. Then for any $q \ge q_0=\frac{\sigma_1(\vecW_0) \eta}{4\sqrt{d}}$ and any $\vecU, \vecU+q\vecE_1 \in \ball(\vecX^0, B)$, we have that
\[ \mathbb{P}\prn*{\prn*{\cKexit(\vecU) \ge K_o}\text{ and } \prn*{\cKexit(\vecU+q\vecE_1) \ge K_o}} \le 0.1.\]
Here probability is over the single sequence of samples used to compute stochastic gradients and the artificial noise injection. 
\end{lemma}
\begin{remark}
The proof of \pref{lem:sgdboundprob} crucially uses that $\grad^2 F(\vecX^0)$ has a negative eigenvector, as one would expect.
\end{remark}
Note we have, as in \citet{fang2019sharp}, that
\[ K_o = 2 \log \prn*{\frac{24\sqrt{d}}{\eta}}\eta^{-1} \delta_2^{-1} \ge \frac{\log(6/q_0)}{\log(1+\eta \delta_2)} \ge \frac{\log(6B/q_0)}{\log(1+\eta \delta_2)}.\numberthis\label{eq:paraminequality}\]
This follows evidently from the choice of parameters and definition of $q_0$, and \pref{rem:restartedsgdreduceparams} which states that it is enough to show the result for $\eta \delta_2$ at most a universal constant, namely one satisfying $\log(1+x) \ge \frac{x}2$. Now using \pref{lem:sgdboundprob}, we prove \pref{lem:sgdescapesaddles}:
\begin{proof}[Proof of \pref{lem:sgdescapesaddles} given \pref{lem:sgdboundprob}]
Given \pref{lem:sgdboundprob}, we first prove that the bad initialization region $\cS^B_{K_o}(\vecX^0)$ satisfies the $(q_0, \vecE_1)$-narrow property, i.e. that there are no points $\vecU, \vecU+q\vecE_1 \in \cS^B_{K_o}(\vecX^0)$ where $q \ge q_0=\frac{\sigma_1(\vecW_0) \eta}{4\sqrt{d}}$. This part of the proof is identical to Proposition 7, \citet{fang2019sharp}. If such points existed we would have 
\[ \mathbb{P}(\cKexit(\vecU) \ge K_o) \ge 0.6, \mathbb{P}(\cKexit(\vecU+q\vecE_1) \ge K_o) \ge 0.6.\]
This implies 
\begin{align*}
\mathbb{P}\prn*{\prn*{\cKexit(\vecU) \ge K_o}\text{ and } \prn*{\cKexit(\vecU+q\vecE_1) \ge K_o}} &\ge \mathbb{P}(\cKexit(\vecU) \ge K_o) + \mathbb{P}(\cKexit(\vecU+q\vecE_1))-1 \\
&\ge 0.2,
\end{align*}
which contradicts \pref{lem:sgdboundprob}.

From here, we prove \pref{lem:sgdescapesaddles}. For this rest of the proof of \pref{lem:sgdescapesaddles}, we only consider $\vecU$ and do not consider the iterates from $\vecU+q\vecE_1$. 
Recall $\cS^B_{K_o}$ satisfies the $(q_0, \vecE_1)$-narrow property with $q_0=\frac{\eta\sigma_1(\vecW_0)}{4\sqrt{d}}$ as shown above. 
Thus we have for any $\vecU\in \ball(\vecX^0, B)$, 
\begin{align*}
\mathbb{P}\prn*{\vecW^1(\vecU) \in \cS^B_{K_o}(\vecX^0)} &= \mathbb{P}\prn*{\vecU - \eta \grad \tilde{f}(\vecU;\vecZeta_1) \in \cS^B_{K_o}(\vecX^0)} \\
&= \mathbb{P}\prn*{\grad \tilde{f}(\vecU; \vecZeta_1) \in \eta^{-1}\prn*{-\cS^B_{K_o}(\vecX^0)+\vecU}} \le \frac12.\numberthis\label{eq:restartedsgdescapesaddlesstucksmall}
\end{align*}
The last step follows from the definition of the $\vecW^k(\vecU)$, the scale and translation properties of the $(q_0,\vecE_1)$-narrow property which implies that $\eta^{-1}\prn*{-S^B_{K_o}(\vecX^0)+\vecU}$ satisfies the $(\frac{\sigma_1(\vecW_0)}{4\sqrt{d}}, \vecE_1)$-narrow property, and that $\grad \tilde{f}\prn*{\vecU; \vecZeta_1}$ satisfies the $\vecE_1$-dispersive property by \pref{lem:dispersivenoiseregpoints}. 

Note as events we have $\{\cKexit(\vecW^1(\vecU)) < K_o\} \subseteq \{\cKexit(\vecU) \le K_o\}$. Thus by Law of Total Expectation, for all $\vecU\in \ball(\vecX^0, B)$, 
\begin{align*}
\mathbb{P}\prn*{\cKexit(\vecU) \le K_o} &\ge \mathbb{P}\prn*{\cKexit(\vecW^1(\vecU)) < K_o} \\
&\ge \mathbb{E}\brk*{\mathbb{P}\prn*{\cKexit(\vecW^1(\vecU)) < K_o | \mathfrak{F}^1} | \crl*{\vecW^1(\vecU) \in (\cS^B_{K_o}(\vecX^0))^c}}.\numberthis\label{eq:restartedsgdescapesaddlestotalexp}
\end{align*}
Conditioned on $\vecW^1(\vecU) \le (\cS^B_{K_o}(\vecX^0))^c$, we have by definition of $\cS^B_{K_o}(\vecX^0)$ that $\mathbb{P}\prn*{\cKexit(\vecW^1(\vecU)) < K_o | \mathfrak{F}^1} \ge 0.4$. By \pref{eq:restartedsgdescapesaddlesstucksmall}, for all $\vecU\in \ball(\vecX^0, B)$, we have
\[ \mathbb{P}\prn*{\vecW^1(\vecU) \in \cS^B_{K_o}(\vecX^0)^c} \ge \frac12. \]
Thus combining with \pref{eq:restartedsgdescapesaddlestotalexp} implies for all $\vecU\in \ball(\vecX^0, B)$,
\[ \mathbb{P}\prn*{\cKexit(\vecU) \le K_o} \ge0.4 \cdot \frac12=0.2.\numberthis\label{eq:lowerboundprobescapeinball}\]
Now consider any $N' \ge 1$. Notice as events,
\begin{align*}
\crl*{\cKexit(\vecU) > N'K_o} &= \crl*{\cKexit(\vecW^{(N'-1)K_o}(\vecU)) > K_o} \\
&= \crl*{\cKexit(\vecW^{(N'-1)K_o}(\vecU)) > K_o} \cap \crl*{\cKexit(\vecU)>(N'-1)K_o}.
\end{align*}
Therefore,
\[ \mathbb{P}\prn*{\cKexit(\vecU) > N'K_o} = \mathbb{E}\brk*{\mathbb{P}\prn*{\cKexit(\vecW^{(N'-1)K_o}(\vecU)) > K_o | \mathfrak{F}^{K_o}} | \crl*{\cKexit(\vecU)>(N'-1)K_o}}.\]
Note that conditioned on $\cKexit(\vecU)>(N'-1)K_o$, it follows that $\cKexit(\vecW^{(N'-1)K_o}(\vecU)) \in \ball(\vecX^0, B)$. Therefore $\mathbb{P}\prn*{\cKexit(\vecW^{(N'-1)K_o}(\vecU)) > K_o | \mathfrak{F}^{K_o}} \le  \sup_{\vecU' \in \ball(\vecX^0, B)}\mathbb{P}\prn*{\cKexit(\vecU') > K_o}$. Using \pref{eq:lowerboundprobescapeinball}, we can upper bound
\begin{align*}
\mathbb{P}\prn*{\cKexit(\vecU) > N'K_o} &\le \mathbb{P}\prn*{\cKexit(\vecU)>(N'-1)K_o} \cdot \sup_{\vecU' \in \ball(\vecX^0, B)}\mathbb{P}\prn*{\cKexit(\vecU') > K_o}\\
&\le 0.8 \mathbb{P}\prn*{\cKexit(\vecU)>(N'-1)K_o}.\numberthis\label{eq:restartedsgdescapesaddlesrecurseeq}
\end{align*}
Recall that $K_0 = \lfloor \frac{\log (3/p)}{\log(0.8^{-1})}+1\rfloor K_o$. Let $N=\lfloor \frac{\log (3/p)}{\log(0.8^{-1})}+1\rfloor$. We obtain by repeatedly applying \pref{eq:restartedsgdescapesaddlesrecurseeq} for $N'=N, N-1, \ldots$ that 
\[ \mathbb{P}\prn*{\cKexit(\vecU) > NK_o} \le 0.8^N \le p/3.\]
This gives the desired result.
\end{proof}
Now we prove \pref{lem:sgdboundprob}.
\begin{proof}[Proof of \pref{lem:sgdboundprob}]
Again, we proceed similarly as the proof of Lemma 8, \citet{fang2019sharp}. The main difference is we only have control over the relevant derivatives prior to the escape from $\ball(\vecX^0, B)$ (recall $\vecX^0 \in \cL_{F,F(\vecW_0)}$). However, it turns out that this is sufficient for the proof to go through.

\paragraph{Setup.} Recall that we have $\vecW^0(\vecU)=\vecU$, and
\begin{align*}
\vecW^k(\vecU) &= \vecW^{k-1}(\vecU) - \eta \grad \tilde{f}(\vecW^{k-1}(\vecU);\vecZeta_k),\\
\vecW^k(\vecU+q\vecE_1) &= \vecW^{k-1}(\vecU+q\vecE_1) - \eta \grad \tilde{f}(\vecW^{k-1}(\vecU+q\vecE_1);\vecZeta_k).
\end{align*}
Now define the following stopping time:
\[ \cK_1 = \cKexit(\vecU) \land \cKexit(\vecU+q\vecE_1).\]
For solely the purpose of analysis, consider the following sequence:
\[ \vecZ^k = \begin{cases} \vecW^k(\vecU+q\vecE_1) - \vecW^k(\vecU) &: k \le \cK_1 \\ \prn*{\matI - \eta \grad^2 F(\vecX^0)} \vecZ^{k-1}  &: k > \cK_1\end{cases}.\numberthis\label{eq:restartedsgdescapesaddleszk}\]
Clearly the $\vecZ^k$ are $\mathfrak{F}^k$-measurable, because the event $\{k \le \cK_1\}$ is $\mathfrak{F}^k$-measurable.
\begin{remark}
Note unlike \citet{fang2019sharp}, the first case holds when $k \le \cK_1$ rather than $k<\cK_1$. That being said we expect that if one uses the exact same definition as in \citet{fang2019sharp} for the $\vecZ^k$, the proof this generalized smooth setting will still work, with a slightly modified argument compared to the proof we present.
\end{remark}
Notice by definition of $\vecW^0(\vecU), \vecW^0(\vecU+q\vecE_1)$ and assumption of \pref{lem:sgdboundprob} that $\vecU, \vecU+q\vecE_1 \in \ball(\vecX^0, B)$, we have $\cK_1>0$. Thus,
\[ \vecZ^0 = q\vecE_1. \]

\paragraph{Controlling the $\vecZ^k$.}
Let $\matH = \grad^2 F(\vecX^0)$. We have the following lemma to control the $\vecZ^k$ from \pref{eq:restartedsgdescapesaddleszk}. 

For all $k$, define
\begin{align*}
\matD^{k} &:= \grad^2 F(\vecX^0) - \int_0^1 \grad^2 F\prn*{\vecW^{k}(\vecU) + \theta (\vecW^{k}(\vecU+q\vecE_1)-\vecW^{k}(\vecU))} \DERIV \theta,\numberthis\label{eq:restartedsgdescapesaddleDk}\\
\vecXi^k_d &:= \prn*{\grad F(\vecW^{k-1}(\vecU+q\vecE_1)) - \grad F(\vecW^{k-1}(\vecU))}
    - \prn*{\grad \tilde{f}(\vecW^{k-1}(\vecU+q\vecE_1); \vecZeta_k) - \grad \tilde{f}(\vecW^{k-1}(\vecU); \vecZeta_k)}.\numberthis\label{eq:restartedsgdescapesaddlexi1}
\end{align*}
Recall by definition of $\vecW^k(\vecU)$, we have 
\begin{align*}
\grad \tilde{f}(\vecW^{k-1}(\vecU+q\vecE_1); \vecZeta_k) &=  \grad f(\vecW^{k-1}(\vecU+q\vecE_1); \vecZeta_k) + \tilde{\sigma}\vecLambda^{k},\\
\grad \tilde{f}(\vecW^{k-1}(\vecU); \vecZeta_k) &=  \grad f(\vecW^{k-1}(\vecU); \vecZeta_k) + \tilde{\sigma}\vecLambda^{k},
\end{align*}
for the \textit{same} noise sequence $\vecLambda^{k}$. Thus we also have
\[ \vecXi^k_d = \prn*{\grad F(\vecW^{k-1}(\vecU+q\vecE_1)) - \grad F(\vecW^{k-1}(\vecU))}
    - \prn*{\grad f(\vecW^{k-1}(\vecU+q\vecE_1); \vecZeta_k) - \grad f(\vecW^{k-1}(\vecU); \vecZeta_k)}.\numberthis\label{eq:escapesaddlexidiff}\]
\begin{lemma}[Equivalent of Lemma 13, \citet{fang2019sharp}]\label{lem:controllingz_kescapesaddle}
We have that for all $k \le \cK_1$, 
\[
    \vecZ^k = (\matI - \eta \matH) \vecZ^{k-1} + \eta \matD^{k-1} \vecZ^{k-1} + \eta \vecXi^k_d.
\]
Furthermore, we have the following properties of the $\matD^k$ and $\vecXi^k_d$ defined in \pref{eq:restartedsgdescapesaddleDk}, \pref{eq:restartedsgdescapesaddlexi1}:
\begin{enumerate}
    \item For all such $k \le \cK_1$, we have
\[ \nrm*{\matD^{k-1}} \le L_2(\vecW_0) \max \prn*{\nrm*{\vecW^{k-1}(\vecU+q\vecE_1) - \vecX^0}, \nrm*{\vecW^{k-1}(\vecU) - \vecX^0}} \le L_2(\vecW_0) B.\]
\item For all $k$, we have 
\[ \mathbb{E}\brk*{\vecXi^k_d | \mathfrak{F}^{k-1}}=0.\]
\item For all $k \le \cK_1$, we have
\[ \nrm*{\vecXi^k_d} \le 2L_1(\vecW_0) \nrm*{\vecZ^{k-1}}.\]
\end{enumerate}
\end{lemma}
\begin{proof}
We prove each part one at a time:
\begin{enumerate}
    \item For $k \le \cK_1$, using the definition of $\vecZ^k$, it follows that
\begin{align*}
\vecZ^{k} &= \vecW^k(\vecU+q\vecE_1)- \vecW^k(\vecU) \\
&= \vecW^{k-1}(\vecU+q\vecE_1) - \vecW^{k-1}(\vecU) - \eta \prn*{\grad \tilde{f}(\vecW^{k-1}(\vecU+q\vecE_1); \vecZeta_k) - \grad \tilde{f}(\vecW^{k-1}(\vecU); \vecZeta_k)}\\
&= \vecZ^{k-1} - \eta \prn*{\grad F(\vecW^{k-1}(\vecU+q\vecE_1)) - \grad F(\vecW^{k-1}(\vecU))}
    \\
    &+ \eta \brk*{\prn*{\grad F(\vecW^{k-1}(\vecU+q\vecE_1)) - \grad F(\vecW^{k-1}(\vecU))}
    - \prn*{\grad \tilde{f}(\vecW^{k-1}(\vecU+q\vecE_1); \vecZeta_k) - \grad \tilde{f}(\vecW^{k-1}(\vecU); \vecZeta_k)}} \\
&= \vecZ^{k-1} - \eta \brk*{\int_0^1 \grad^2 F\prn*{\vecW^{k-1}(\vecU)+\theta (\vecW^{k-1}(\vecU+q\vecE_1)-\vecW^{k-1}(\vecU))} \DERIV \theta}\vecZ^{k-1}+\eta \vecXi^k_d \\
&= \vecZ^{k-1} - \eta \prn*{\matH - \matD^{k-1}} \vecZ^{k-1}+\eta \vecXi^k_d.
\end{align*}
This proves the desired property of the $\vecZ^k$.

\item For the required properties of the $\matD^{k-1}$, consider any $k \le \cK_1$. First, notice $\vecW^{k-1}(\vecU) + \theta (\vecW^{k-1}(\vecU+q\vecE_1)-\vecW^{k-1}(\vecU)) = \theta \vecW^{k-1}(\vecU+q\vecE_1) + (1-\theta)\vecW^{k-1}(\vecU)$ for any $\theta \in [0,1]$. For $k \le \cK_1$, both $\vecW^{k-1}(\vecU+q\vecE_1), \vecW^{k-1}(\vecU) \in \ball(\vecX^0, B)$. Note this still remains true for $k = \cK_1$ because for $k-1=\cK_1-1 < \cK_1$, the definition of $\cK_1$ implies that the iterates $\vecW^{k-1}(\vecU+q\vecE_1), \vecW^{k-1}(\vecU) \in \ball(\vecX^0, B)$.

Thus for any $\theta \in [0,1]$, $\vecW^{k-1}(\vecU) + \theta (\vecW^{k-1}(\vecU+q\vecE_1)-\vecW^{k-1}(\vecU)) \in \ball(\vecX^0, B)$, and so all points $\vecP$ on the line segment between $\vecX^0$ and $\vecW^{k-1}(\vecU) + \theta (\vecW^{k-1}(\vecU+q\vecE_1)-\vecW^{k-1}(\vecU))$ lie in $\ball(\vecX^0, B)$. Thus by \pref{lem:hessianlipschitzchange}, 
\begin{align*}
\nrm*{\matD^{k-1}} &= \nrm*{\grad^2 F(\vecX^0) - \int_0^1 \grad^2 F(\vecW^{k-1}(\vecU) + \theta (\vecW^{k-1}(\vecU+q\vecE_1)-\vecW^{k-1}(\vecU))) \DERIV \theta} \\
&\le \int_0^1 \nrm*{\grad^2 F(\vecX^0) - \grad^2 F\prn*{\vecW^{k-1}(\vecU) + \theta (\vecW^{k-1}(\vecU+q\vecE_1)-\vecW^{k-1}(\vecU))}} \DERIV \theta \\
&\le L_2(\vecW_0) \int_0^1 \nrm*{\theta \prn*{\vecW^{k-1}(\vecU+q\vecE_1) - \vecX^0} + (1-\theta)\prn*{\vecW^{k-1}(\vecU) - \vecX^0}} \DERIV \theta \\
&\le L_2(\vecW_0) \max \crl*{\nrm*{\vecW^{k-1}(\vecU+q\vecE_1) - \vecX^0}, \nrm*{\vecW^{k-1}(\vecU) - \vecX^0}} \\
&\le L_2(\vecW_0) B.
\end{align*}
The last line follows since $k \le \cK_1$, hence $k-1<\cK_1$, thus $\vecW^{k-1}(\vecU+q\vecE_1), \vecW^{k-1}(\vecU) \in \ball(\vecX^0, B)$. 

\item Next as the stochastic gradient oracle $\grad f(\cdot;\vecZeta)$ is unbiased, applying Linearity of Expectation on \pref{eq:escapesaddlexidiff}, it follows that $\mathbb{E}\brk*{\vecXi^k_d | \mathfrak{F}^{k-1}}=0$ for all $k$.

For the bound on the magnitude of $\vecXi^k_d$, again recall by the above that for $k\le \cK_1$, we have 
\[ \vecW^{k-1}(\vecU+q\vecE_1), \vecW^{k-1}(\vecU) \in \ball(\vecX^0, B).\]
Thus for all $\vecP$ on the line segment between $\vecW^{k-1}(\vecU+q\vecE_1), \vecW^{k-1}(\vecU)$, we have $\vecP \in \ball(\vecX^0, B)$. Thus by \pref{lem:smoothnesschange}, $\nrm*{\grad^2 F(\vecP)} \le L_1(\vecW_0)$. By \pref{lem:gradoraclelipschitzconstant}, for any $\vecZeta$, $\nrm*{\grad^2 f(\vecP;\vecZeta)} \le L_1(\vecW_0)$. Recalling \pref{eq:escapesaddlexidiff} gives
\begin{align*}
&\nrm*{\vecXi^k_d} \\
&\le \nrm*{\grad F(\vecW^{k-1}(\vecU+q\vecE_1)) - \grad F(\vecW^{k-1}(\vecU))}
+\nrm*{\grad f(\vecW^{k-1}(\vecU+q\vecE_1); \vecZeta_k) - \grad f(\vecW^{k-1}(\vecU); \vecZeta_k)} \\
&\le 2L_1(\vecW_0)\nrm*{\vecW^{k-1}(\vecU+q\vecE_1)-\vecW^{k-1}(\vecU)} \\
&= 2L_1(\vecW_0)\nrm*{\vecZ^{k-1}}.
\end{align*}
In the last step, we used the definition of $\vecZ^{k}$ for $k \le \cK_1$. 
\end{enumerate}
This proves all the desired parts of \pref{lem:controllingz_kescapesaddle}.
\end{proof}

\paragraph{Controlling iterates under a high probability event.} 
We now consider a rescaled iteration as considered in \citet{fang2019sharp}. Recall the definition of $\delta_m \geq \delta_2$ in the statement of \pref{lem:sgdboundprob}. For each $k = 0, 1, \ldots$, we define:
\[
\vecPsi_k := q^{-1}(1 + \eta \delta_m)^{-k} \vecZ_k. 
\]
\begin{lemma}[Equivalent of the first part of Lemma 14, \citet{fang2019sharp}]\label{lem:restartedsgdescapesaddleshighprobcontrol1}
Define $\hatmatD_k := (1 + \eta \delta_m)^{-1} \matD_k$, and slightly overloading notation, define 
\[ \vecZeta^k_d := q^{-1}(1 + \eta \delta_m)^{-k} \vecXi^k_d. \]
Then for $k \le \cK_1$, we have $\vecPsi^0 = \vecE_1$ and
\[
\vecPsi^k =
\frac{\matI - \eta \matH}{1 + \eta \delta_m} \vecPsi^{k-1} + \eta \hatmatD^{k-1} \vecPsi^{k-1} + \eta \vecZeta^k_d,
\]
as well as the properties
\begin{align*}
\nrm*{\hatmatD^{k}} &\le L_2(\vecW_0) B \text{ for all } 0  \le k < \cK_1, \\
\nrm*{\vecZeta^k_d} &\leq 2 L_1(\vecW_0) \nrm*{\vecPsi^{k-1}} \text{ for all } 1  \le k \le \cK_1.
\end{align*}
\end{lemma}
\begin{proof}
We prove all the desired parts of \pref{lem:restartedsgdescapesaddleshighprobcontrol1}.
\begin{itemize}
    \item The fact that $\vecPsi^0=\vecE_1$ follows immediately, because $\vecZ^0 = q\vecE_1$. For the general recursion for $\vecPsi^k$, consider any $k \le \cK_1$. First note that by the recursion for the $\vecZ^k$ for $k \le \cK_1$ in \pref{lem:controllingz_kescapesaddle}, we have
\begin{align*}
\vecPsi^k &= q^{-1}(1 + \eta \delta_m)^{-k} \vecZ^k \\
&= \frac{\matI - \eta \matH}{1 + \eta \delta_m} q^{-1}(1 + \eta \delta_m)^{-(k-1)} \vecZ^{k-1} \\
&\hspace{1in}+ \eta \frac{\matD^{k-1}}{1 + \eta \delta_m} q^{-1}(1 + \eta \delta_m)^{-(k-1)} \vecZ^{k-1} + \eta q^{-1}(1 + \eta \delta_m)^{-k} \vecXi^k_d \\
&= \frac{\matI - \eta \matH}{1 + \eta \delta_m} \vecPsi^{k-1} + \eta \hatmatD^{k-1} \vecPsi^{k-1} + \eta \vecZeta^k_d.
\end{align*}

\item Consider any $k \le \cK_1$. For the requisite properties of $\hatmatD^{k}$ for $k<\cK_1$, the upper bound on the norm of $\hatmatD^{k}$ follows immediately from \pref{lem:controllingz_kescapesaddle}. 

Next from the definition of $\vecZeta^k_d$ and \pref{lem:controllingz_kescapesaddle}, for $k \le \cK_1$ we have that
\begin{align*}
\nrm*{\vecZeta^k_d} &\le q^{-1}(1 + \eta \delta_m)^{-k} \nrm*{\vecXi^k_d} \\
&\leq 2 L_1(\vecW_0) q^{-1} \frac{(1 + \eta \delta_m)^{-(k-1)}}{1 + \eta \delta_m} \nrm*{\vecZ^{k-1}} \\
&\le 2L_1(\vecW_0) \nrm*{\vecPsi^{k-1}}.    
\end{align*}
\end{itemize}
This proves \pref{lem:restartedsgdescapesaddleshighprobcontrol1}.
\end{proof}

\begin{lemma}[Equivalent of the rest of Lemma 14, \citet{fang2019sharp}]\label{lem:restartedsgdescapesaddleshighprobcontrol2}
With the step size $\eta$ from \pref{eq:paramchoice}, there exists an event $\cH_{o}$ (namely, from \pref{eq:highprobevent_escapesaddle}) with probability at least 0.9, such that for all $k \le \min(\cK_1-1, K_0)$ we have
\[
\nrm*{\vecPsi^k}^2 \le 4,\numberthis\label{eq:restartedsgdescapesaddleshighprobcontrol2event1}\]
and 
\[ \vecE_1^\top \vecPsi^k > \frac{1}{2}.\numberthis\label{eq:restartedsgdescapesaddleshighprobcontrol2event2} \]
\end{lemma}

\begin{proof}
Define
\[
\vechatPsi^{k-1} = \frac{\matI - \eta \matH}{1 + \eta \delta_m} \vecPsi^{k-1}. 
\]
Recall that $\matH = \grad^2 F(\vecX^0)$ and $\vecX^0$ is in the $F(\vecW_0)$-sublevel set $\cL_{F,F(\vecW_0)}$. Therefore, from \pref{ass:selfbounding}, $\nrm*{\matH} \le L_1(\vecW_0)$. By definition of $\delta_m$, it follows that 
\[ -\delta_m \matI \preceq \matH \preceq L_1(\vecW_0) \matI.\]
Since $\eta L_1(\vecW_0) \leq 1$, it follows that the matrix $\matI - \eta \matH$ is symmetric and has all eigenvalues in $[0, 1 + \eta \delta_m]$. This implies 
\[ \nrm*{\vechatPsi^{k-1}} \le \nrm*{\vecPsi^{k-1}}.\numberthis\label{eq:restartedsgdescapesadleshatpsibound}\]
Note that $\vechatPsi^{k-1}$ and $\vecPsi^{k-1}$ are measurable on $\mathfrak{F}^{k-1}$. This combined with \pref{lem:controllingz_kescapesaddle} and \pref{lem:restartedsgdescapesaddleshighprobcontrol1} implies that for all $1 \le k \le \cK_1$,
$$
\mathbb{E}\brk*{(\vechatPsi^{k-1})^{\top} \vecZeta^k_d \cdot 1_{\nrm*{\vecPsi^{k-1}} \leq 2} | \mathfrak{F}^{k-1}}
= 1_{\nrm*{\vecPsi^{k-1}} \leq 2} \cdot \mathbb{E}[(\vechatPsi^{k-1})^{\top} \vecZeta^k_d | \mathfrak{F}^{k-1}] = 0, 
$$
and
\[
|(\vechatPsi^{k-1})^\top \vecZeta^k_d \cdot 1_{\nrm*{\vecPsi^{k-1}} \leq 2}|^2 
\le 1_{\nrm*{\vecPsi^{k-1}} \leq 2} \cdot 4L^2_1(\vecW_0) \nrm*{\vecPsi^{k-1}}^4 \le (8L_1(\vecW_0))^2.
\]
Now define the following real-valued stochastic process:
\[ Y_k = (\vechatPsi^{k-1})^\top \vecZeta^k_d  1_{\nrm*{\vecPsi^{k-1}} \le 2} 1_{k-1 < \cK_1}= \begin{cases} (\vechatPsi^{k-1})^\top \vecZeta^k_d \cdot 1_{\nrm*{\vecPsi^{k-1}} \leq 2}&:k\le \cK_1\\ 0&:k>\cK_1.\end{cases}\]
Note $Y_k$ is $\mathfrak{F}_k$-measurable, and that $(\vechatPsi^{k-1})^\top, 1_{\nrm*{\vecPsi^{k-1}} \le 2}, 1_{k-1 < \cK_1} \equiv 1_{k \le \cK_1}$ are all $\mathfrak{F}_{k-1}$-measurable. Thus, by \pref{lem:controllingz_kescapesaddle} and the definition of $\vecZeta^k_d$ from \pref{lem:restartedsgdescapesaddleshighprobcontrol1},
\[\mathbb{E}\brk*{Y_k|\mathfrak{F}_{k-1}}=0 .\]
Furthermore combining the above justification with the trivial case $k >\cK_1$, we obtain
\[ |Y_k|  \le 8L_1(\vecW_0).\]
By the (standard) Azuma's Inequality, with probability $1 - 0.1 / (2K_0)$, for any given $l, 1 \le l \le K_0$:
\[
\abs*{\sum_{k=1}^l Y_k}
\leq 8L_1(\vecW_0) \sqrt{2l \log(40K_0)} \leq 8L_1(\vecW_0) \sqrt{2 K_0 \log(40K_0)} 
\le \frac{1}{\eta},
\]
where the last inequality follows from the given choice of parameters.

Analogously, by \pref{lem:controllingz_kescapesaddle} and \pref{lem:restartedsgdescapesaddleshighprobcontrol1}, we also have for $1 \le k \le \cK_1$:
\[
\mathbb{E}\brk*{\vecE_1^\top \vecZeta^k_d \cdot 1_{\nrm*{\vecPsi^{k-1}} \le 2} | \mathfrak{F}^{k-1}} = 0, |\vecE_1^\top \vecZeta^k_d \cdot 1_{\nrm*{\vecPsi^{k-1}} \le 2}| \leq 4L_1(\vecW_0).
\]
Define
\[ Y'_k := \vecE_1^\top \vecZeta^k_d \cdot 1_{\nrm*{\vecPsi^{k-1}} \le 2} 1_{k \le \cK_1}.\]
The (standard) Azuma's Inequality now implies that with probability at least $1 - 0.1 / (2K_0)$, for any given $l, 1 \le l \le K_0$:
\[
\abs*{\sum_{k=1}^l Y'_k} 
\leq 4L_1(\vecW_0)\sqrt{2 l \log(40K_0)} \le \frac{1}{4\eta}.
\]
By the Union Bound, there exists an event $\cH_{o}$ happening with probability at least $0.9$ such that the following inequalities hold for each $l = 1, 2, \ldots, K_0$:
\[ \abs*{\sum_{k=1}^l Y_k} \le \frac{1}{\eta}, \abs*{\sum_{k=1}^l Y'_k}\leq \frac{1}{4\eta}.\numberthis\label{eq:highprobevent_escapesaddle}\]
In particular under the event $\cH_o$, for any $l \le \min(\cK_1-1, K_0)$, using the definitions of $Y_k, Y'_k$ we obtain
\[
\abs*{\sum_{k=1}^l \vechatPsi_{k-1}^\top \vecZeta^k_d \cdot 1_{\nrm*{\vecPsi^{k-1}} \leq 2}}\leq \frac{1}{\eta},
\abs*{\sum_{k=1}^l \vecE_1^\top \vecZeta^k_d \cdot 1_{\nrm*{\vecPsi_{k-1}} \leq 2}} \leq \frac{1}{4\eta}.\numberthis\label{eq:restartedsgdescapesaddleshighprobeventimplication}
\]
Now from \pref{lem:restartedsgdescapesaddleshighprobcontrol1}, it follows for all $k \le \cK_1$ that
\begin{align*}
\nrm*{\vecPsi^k}^2 &= 
\nrm*{ \frac{\matI - \eta \matH}{1 + \eta \delta_m} \vecPsi^{k-1} + \eta \hatmatD^{k-1} \vecPsi^{k-1} + \eta \vecZeta^k_d }^2 \\
&= \nrm*{\vechatPsi^{k-1}}^2 + 2\eta (\vechatPsi_{k-1})^\top \hatmatD_{k-1} \vecPsi_{k-1} 
+ \eta^2 \nrm*{ \hatmatD_{k-1} \vecPsi^{k-1} + \vecZeta^k_d }^2 + 2\eta (\vechatPsi^{k-1})^\top \vecZeta^k_d \\
&= \nrm*{\vecPsi^{k-1}}^2 + Q_{1,k} + Q_{2,k} + Q_{3,k}\numberthis\label{eq:restartedsgdescapesaddlespsikrecursebound}
\end{align*}
where we define
\[ Q_{1,k} := 2\eta (\vechatPsi^{k-1})^\top \hatmatD^{k-1} \vecPsi^{k-1}, Q_{2,k} := \eta^2 \nrm*{ \hatmatD_{k-1} \vecPsi_{k-1} + \vecZeta^{k}_d }^2, Q_{3,k} := 2\eta (\vechatPsi^{k-1})^\top \vecZeta_{d}^{k}.\]
For $k \le\cK_1$, we have $k-1<\cK_1$. Thus by \pref{lem:restartedsgdescapesaddleshighprobcontrol1} and \pref{eq:restartedsgdescapesadleshatpsibound}, we have
\[
Q_{1,k} \leq 2\eta L_2(\vecW_0) B \nrm*{\vecPsi^{k-1}}^2,\numberthis\label{eq:restartedsgdescapesaddlesQ1kbound}
\]
and
\begin{align*}
Q_{2,k} &\le 2\eta^2 \nrm*{\hatmatD^{k-1} \vecPsi^{k-1}}^2 + 2\eta^2 \nrm*{\vecZeta_{d}^{k}}^2 \\
&\le 2\eta^2 \cdot L_2(\vecW_0)^2 B^2 \nrm*{\vecPsi^{k-1}}^2 + 8\eta^2 L_1(\vecW_0)^2 \nrm*{\vecPsi^{k-1}}^2 \\
&\le 16\eta^2 L_1(\vecW_0)^2 \nrm*{\vecPsi^{k-1}}^2.\numberthis\label{eq:restartedsgdescapesaddlesQ2kbound}
\end{align*}
The last inequality above follows as per \pref{rem:restartedsgdreduceparams}.

Now we complete the proof. Under the event $\cH_{o}$ from \pref{eq:highprobevent_escapesaddle}, we prove \pref{eq:restartedsgdescapesaddleshighprobcontrol2event1} by induction on $k$ (recall our condition for $k$ for \pref{lem:restartedsgdescapesaddleshighprobcontrol2} is that $0 \le k \le\min(\cK_1-1, K_0)$). 

When $k = 0$, by \pref{lem:restartedsgdescapesaddleshighprobcontrol1}, $\vecPsi^0 = \vecE_1$, so $\nrm*{\vecPsi^0} = \nrm*{\vecE_1} =1 \le 2$ and $\vecE_1^{\top} \vecPsi^0 = \nrm*{\vecE_1}^2=1$ (recall $\vecE_1$ is a \textit{unit} eigenvector), proving the base case. 

Now for the inductive step, consider some $k \le\min(\cK_1-1, K_0)$. Suppose $\nrm*{\vecPsi^l} \le 2$ holds for all $l, 0 \le l \le k-1$. Then because $k<\cK_1$, upon applying the above bounds \pref{eq:restartedsgdescapesaddlespsikrecursebound}, \pref{eq:restartedsgdescapesaddlesQ1kbound}, \pref{eq:restartedsgdescapesaddlesQ2kbound} we have:
\begin{align*}
\nrm*{\vecPsi^k}^2 &\le \nrm*{\vecPsi^0}^2 + \sum_{s=1}^k Q_{1,s} + \sum_{s=1}^k Q_{2,s} + \sum_{s=1}^k Q_{3,s} \\
&\le 1 + 2\eta \sum_{s=1}^k L_2(\vecW_0) B \nrm*{\vecPsi^{s-1}}^2 + 16\eta^2 L_1(\vecW_0)^2 \sum_{s=1}^k \nrm*{\vecPsi^s}^2 + 2\eta \sum_{s=1}^k (\vechatPsi^{s-1})^\top \vecZeta_d^s \\
&\le 1 + 2 L_2(\vecW_0) B \cdot 4 \cdot \eta k + 16\eta^2 \cdot L_1(\vecW_0)^2 \cdot 4 \cdot k + 2\eta \sum_{s=1}^k (\vechatPsi^{s-1})^\top \vecZeta_d^s \cdot 1_{\nrm*{\vecPsi^{s-1}} \le 2} \\
&\le 1 + 16L_2(\vecW_0) B \cdot \eta K_0 + 2\eta \sum_{s=1}^k {\vechatPsi_{s-1}}^\top \vecZeta_d^s \cdot 1_{\nrm*{\vecPsi_{s-1}} \le 2} \le 1 + 1 + 2\eta \cdot \frac{1}{\eta} = 4.
\end{align*}
To upper bound the above, we used our choice of step size $\eta \le \frac{L_2(\vecW_0) B}{8L_1(\vecW_0)^2}$ and $B \le \frac{1}{L_1(\vecW_0)}$ as per \pref{rem:restartedsgdreduceparams}, our above upper bounds on $Q_{1,s}, Q_{2,s}$, and that the event $\cH_o$ implies \pref{eq:restartedsgdescapesaddleshighprobeventimplication}. 

This completes the induction and proves \pref{eq:restartedsgdescapesaddleshighprobcontrol2event1}.

With \pref{eq:restartedsgdescapesaddleshighprobcontrol2event1}, we prove \pref{eq:restartedsgdescapesaddleshighprobcontrol2event2}. Namely note for $k \le \min(\cK_1-1, K_0)$, summing and telescoping the recursion for $\vecPsi^k$ from \pref{lem:restartedsgdescapesaddleshighprobcontrol1}, we have:
\begin{align*}
\vecE_1^\top \vecPsi_k &= \vecE_1^\top \vecPsi_0 + \sum_{s=0}^{k-1} \eta \vecE_1^\top \hatmatD_s \vecPsi^s + \sum_{s=0}^{k-1} \eta \vecE_1^\top \vecZeta_d^s\\
&\ge 1 - \eta \sum_{s=0}^{k-1} 2L_2(\vecW_0) B \nrm*{\vecPsi^{s}} + \eta \sum_{s=0}^{k-1} \vecE_1^\top \vecZeta_d^s \cdot 1_{\nrm*{\vecPsi^{s-1}} \leq 2} \\
&\ge 1 - \eta \cdot K_0 \cdot 2L_2(\vecW_0) B \cdot 2 + \eta \sum_{s=0}^{k-1} \vecE_1^\top \vecZeta_d^s \cdot 1_{\|\vecPsi^{s-1}\| \leq 2} \geq 1 - \frac{1}{8} - \frac{2}{8} \geq \frac{1}{2}.
\end{align*}
Here to lower bound the final sum, we used that $\vecPsi_0=\vecE_1$ and the upper bound on $\nrm*{\hatmatD_s}$ from \pref{lem:restartedsgdescapesaddleshighprobcontrol1}, the fact that we have already established $\nrm*{\vecPsi^s} \le 2$ for all $s < k$ as we showed \pref{eq:restartedsgdescapesaddleshighprobcontrol2event1}, and that the event $\cH_o$ implies \pref{eq:restartedsgdescapesaddleshighprobeventimplication}.

This proves all parts of \pref{lem:restartedsgdescapesaddleshighprobcontrol2}.
\end{proof}

\paragraph{Finish.} 
Now we prove \pref{lem:sgdboundprob} via the same high-level strategy as the proof of Lemma 8, \citet{fang2019sharp}. Note on the event $\{ \cK_1 > K_o\}$, we have 
\[ \vecZ^{K_o} = \vecW^{K_o}(\vecU+q\vecE_1) - \vecW^{K_o}(\vecU) = (\vecW^{K_o}(\vecU+q\vecE_1) - \vecX^0)) - (\vecW^{K_o}(\vecU) - \vecX^0).\]
Thus by definition of $\cK_1$, the event $\{ \cK_1 > K_o\}$ implies that
\[ \nrm*{\vecZ^{K_o}} \le \nrm*{\vecW^{K_o}(\vecU+q\vecE_1) - \vecX^0} + \nrm*{\vecW^{K_o}(\vecU)- \vecX^0} \le 2B.\]
That is, 
\[ \{ \cK_1 > K_o\} \subseteq \{ \nrm*{\vecZ^{K_o}} \le 2B\}.\]
However, consider the event $\cH_o$, \pref{eq:highprobevent_escapesaddle} from \pref{lem:restartedsgdescapesaddleshighprobcontrol2}. On the event $\{ \cK_1 > K_o\} \cap \cH_o$, we have $K_o \le \min(\cK_1-1, K_0)$, and so by \pref{lem:restartedsgdescapesaddleshighprobcontrol2}, we have
\[ \vecE_1^\top \vecPsi^{K_o} > \frac12.\]
Thus by definition of $\vecPsi^k$ and recalling $\delta_m \ge \delta_2 >0$, on the event $\{ \cK_1 > K_o\} \cap \cH_o$ we have
\[ \nrm*{\vecZ^{K_o}} = q(1+\eta \delta_m)^{K_o} \nrm*{\vecPsi^{K_o}} \ge q_0 (1+\eta \delta_2)^{K_o} \abs*{\vecE_1^\top \vecPsi^{K_o}} > q_0 \cdot \frac{6B}{q_0} \cdot \frac12 = 3B,\]
where the last inequality uses \pref{eq:paraminequality}. This means that
\[ \{ \cK_1 > K_o\} \cap \cH_o \subseteq \{ \nrm*{\vecZ^{K_o}} \ge 3B\}.\]
Putting our work together, we see that
\[ \{\cK_1 > K_o\} \cap \cH_o \subseteq \{ \nrm*{\vecZ^{K_o}} \ge 3B\} \cap \{ \nrm*{\vecZ^{K_o}} \le 2B\}=\emptyset.\]
Therefore
\[ \{\cK_1 > K_o\} \subseteq \cH_o^c \implies \mathbb{P}\prn*{\cK_1 > K_o} \le \mathbb{P}(\cH^c_0)\le 0.1.\]
Recalling the definition of $\cK_1$, we conclude \pref{lem:sgdboundprob}.
\end{proof}
\begin{remark}
Note we only have $\vecE_1^T \vecPsi^{k} > \frac12$ for $k < \cK_1$ due to the lack of global Lipschitz bounds on the graedient and Hessian of $F$, unlike in the proof of Lemma 8, \citet{fang2019sharp}.
\end{remark}

\subsection{Faster Descent}\label{subsec:fasterdescent}
\paragraph{Setup:} As in \pref{subsec:escapesaddles}, let $\cK_0$ denote the escape time of $\ball(\vecX^0, B)$ for while loop of \pref{alg:restartedsgd} when the while loop begins at $\vecX^0$. In this section, we aim to prove \pref{lem:fasterdescentlemma}.

As in \pref{subsec:escapesaddles}, the difference between \pref{lem:fasterdescentlemma} and Proposition 9 of \citet{fang2019sharp} is that \textit{this result only holds at points in the $F(\vecW_0)$-sublevel set $\cL_{F,F(\vecW_0)}$}. 
For the rest of this section, we work under the assumptions of \pref{lem:fasterdescentlemma}; thus for the rest of this section, $\vecX^0$ is in the $F(\vecW_0)$-sublevel set $\cL_{F,F(\vecW_0)}$.

The idea here is similar to that of \pref{subsec:escapesaddles}. At a high level, we have the requisite control over the gradient and Hessian since the iterates we consider are in a neighborhood of a point $\vecX^0 \in \cL_{F,F(\vecW_0)}$. 
As in the previous part and as in \citet{fang2019sharp}, we let
\[ \matH := \grad^2 F(\vecX^0),\]
and let
\[
\vecXi^{k+1} := \grad \tilde{f}(\vecX^k; \vecZeta_{k+1}) - \grad F(\vecX^k), \quad k \geq 0.\numberthis\label{eq:vecxikdef}
\]
Note as $\vecLambda^{k+1}$ has mean 0 and as the stochastic gradient oracle is unbiased, we have that for all $k \ge 0$,
\[ \mathbb{E}\brk*{\vecXi^{k+1} | \mathfrak{F}^k} = 0.\]
Let $\cS$ be the subspace spanned by all eigenvectors of $\grad^2 F(\vecX^0)$ whose eigenvalue is greater than $0$, and $\cS^\perp$ denotes the complement space. Also, let $\projS \in \mathbb{R}^{d \times d}$ and $\projSperp \in \mathbb{R}^{d \times d}$ denote the projection matrices onto the spaces $\cS$ and $\cS^\perp$, respectively. Let $\vecU^k = \projS(\vecX^k - \vecX^0)$, and $\vecV^k = \projSperp(\vecX^k - \vecX^0)$. We can decompose the update equation of SGD as:
\[
\vecU^{k+1} = \vecU^k - \eta \projS \grad F(\vecX^k) - \eta \projS \vecXi^{k+1}, 
\]
\[
\vecV^{k+1} = \vecV^k - \eta \projSperp \grad F(\vecX^k) - \eta \projSperp \vecXi^{k+1}, 
\]
for $k \geq 0$. Clearly $\vecU^0 = \pmb{0}$, $\vecV^0 = \pmb{0}$. 

Now decompose $\matH = \matU \matLambda \matU^T$ by the Spectral Theorem where $\matU \in \mathbb{R}^{d\times d}$ is unitary and $\matLambda \in \mathbb{R}^{d\times d}$ is diagonal. Let $\matLambda_{>0}$ denote the diagonal matrix with diagonal entries equal to the positive (diagonal) entries of $\matLambda$. Let $\matLambda_{\le 0}$ denote the diagonal matrix with diagonal entries equal to the zero or negative (diagonal) entries of $\matLambda$. Now define
\[ \matHS := \matU \matLambda_{>0} \matU^T, \matHSperp := \matU \matLambda_{\le 0} \matU^T.\]
Thus $\matHS$ has range in $\cS$, and $\matHSperp$ has range in $\cS^\perp$. Note $\matHS, \matHSperp$ are both symmetric.

From here, define the following quadratic approximations:
\[ G_{\cS}(\vecU) := \brk*{\projS \grad F(\vecX^0)}^\top\vecU + \frac12 \vecU^\top \matHS \vecU, G_{\cS^{\perp}}(\vecV) := \brk*{\projSperp \grad F(\vecX^0)}^\top\vecV + \frac12 \vecV^\top \matHSperp \vecV. \]
Now define the quadratic approximation 
\[
G(\vecX) = G_{\cS}(\vecU) + G_{\cS^\perp}(\vecV) \text{ where }\vecU = \projS(\vecX - \vecX^0), \vecV = \projSperp(\vecX - \vecX^0). 
\]
It is easy to see that
\[ G(\vecX) =\brk*{\grad F(\vecX^0)}^\top  (\vecX - \vecX^0) + \frac12 (\vecX - \vecX^0)^\top \matH (\vecX - \vecX^0).\]
For convenience, let
\[ \projSgradF(\vecX^k) = \projS \grad F(\vecX^k), \projSperpgradF(\vecX^k) = \projSperp \grad F(\vecX^k).\]
Similarly, let 
\[ \vecXiS^k = \projS \vecXi^k, \vecXiSperp^k = \projSperp \vecXi^k. \]
Also denote the stopping time 
\[ \cK = \cK_0 \land K_0 .\]
Due to its `local' nature around the $\vecX^0$ in the $F(\vecW_0)$-sublevel set, we still have the following result from \citet{fang2019sharp}:
\begin{lemma}[Equivalent of Lemma 15, \citet{fang2019sharp}]\label{lem:localsmoothineq}
Consider any $\vecU \in \cL_{F,F(\vecW_0)}$, and consider any $\vecX\in \ball(\vecU, B)$. Then we have
\[
\nrm*{\grad F(\vecX) - \grad G(\vecX)} \le \frac{L_2(\vecW_0) B^2}{2}. 
\]
Furthermore, for any symmetric matrix $\matA$, with $0 < a \leq \frac{1}{\nrm*{\matA}_2}$, for any $i = 0, 1, \ldots$, and $j = 0, 1, \ldots$, we have
\[
\nrm*{(\matI - a\matA)^i \matA (\matI - a\matA)^j}_2 \leq \frac{1}{a(i + j + 1)}. 
\]
\end{lemma}
\begin{proof}    
Notice that for all $0 \le \theta \le 1$, $\theta \vecX + (1-\theta)\vecU \in \ball(\vecU, B)$. Thus as $\vecU \in \cL_{F,F(\vecW_0)}$, by \pref{lem:hessianlipschitzchange}, we have 
\[ \nrm*{\grad^2 F(\theta \vecX + (1-\theta)\vecU) - \grad^2 F(\vecU)} \le L_2(\vecW_0) \cdot \theta\nrm*{\vecX-\vecU}\text{ for all }0\le \theta \le 1.\]
Thus we have
\begin{align*}
\nrm*{\grad F(\vecX) - \grad G(\vecX)} 
&= \nrm*{\grad F(\vecX) - \grad F(\vecX^0) - \grad^2 F(\vecU)(\vecX - \vecU)} \\
&= \nrm*{\crl*{\int_0^1 \prn*{\nabla^2 F(\vecX^0 + \theta(\vecX - \vecU))-\grad^2 F(\vecU)} \DERIV\theta}\prn*{\vecX - \vecU}}\\
&\le \nrm*{\int_0^1 \crl*{L_2(\vecW_0) \cdot \theta \nrm*{\vecX - \vecU}} \DERIV\theta} \cdot \nrm*{\vecX - \vecU}\\
&\le \frac{L_2(\vecW_0) B^2}{2}.
\end{align*}
The second part of the Lemma follows from the exact same proof of \pref{lem:technicallemma1} in \pref{sec:perturbedgdproofs}. It is also proved in the proofs of Lemma 15, \citet{fang2019sharp}, and in the proof of Lemma 16 of \citet{jin2017escape}. For more detail, let the eigenvalues of $\matA$ be $\{\lambda_k\}$. Thus for any $i, j \geq 0$, the eigenvalues of $(\matI - a \matA)^i  \matA(\matI - a \matA)^j$ are $\{\lambda_k (1 - a \lambda_k)^{i+j}\}$. We now detail a calculation from \citet{jin2017escape}. Letting $g_t(\lambda) := \lambda(1 - a \lambda)^t$ and setting its derivative to zero yields
\[
\nabla g_t(\lambda) = (1 - a \lambda)^t - t a \lambda (1 - a \lambda)^{t-1} = 0.
\]
It is easy to check that $\lambda_t^\star = \frac{1}{(1 + t)a}$ is the unique maximizer, and $g_t(\lambda)$ is monotonically increasing in $(-\infty, \lambda_t^\star]$.

This gives:
\[
\nrm*{ (\matI - a \matA)^i \matA (\matI - a \matA)^{j} } = \max_{k} \lambda_i (1 - a \lambda_k)^{i+j} \leq \hat{\lambda}(1 - a \hat{\lambda})^{i+j} \leq \frac{1}{(1 + i+j)a},
\]
where $\hat{\lambda} = \min\{\ell, \lambda_{i+j}^\star\}$.
\end{proof} 
\begin{lemma}\label{lem:restartedsgdfasterdescentxikcontrol}
For any $k \le \cK_0$, we have 
\[ \nrm*{\vecXi^{k}} \le \sigma_1(\vecW_0).\]
\end{lemma}
\begin{proof}
Note for $k \le \cK_0$, we have $k-1<\cK_0$ and so $\vecX^{k-1} \in \ball(\vecX^0, B)$. Recall furthermore that $\vecX^0 \in \cL_{F,F(\vecW_0)}$. Thus, by \pref{lem:noisechange} and \pref{lem:restartedsgdballcloseenough},
\[ \nrm*{\vecXi^{k}} = \nrm*{\grad \tilde{f}(\vecX^{k-1}; \vecZeta_{k}) - \grad F(\vecX^{k-1})} \le \sigma_1(\vecW_0), \]
as desired.
\end{proof}

\paragraph{Analyzing the Quadratic Approximation:} We now analyze the quadratic approximation $G(\vecX)$ as done in \citet{fang2019sharp}. First we analyze the part in $\cS$:
\begin{lemma}[Equivalent of Lemma 16, \citet{fang2019sharp}]\label{lem:quaddescentanalysis}
Set hyperparameters from (8). With probability at least $1 - p/4$, we have
\begin{align*}
&G_{\cS}(\vecU^{\cK}) - G_{\cS}(\vecU^0) \\
&\le -\frac{25\eta}{32} \sum_{k=0}^{\cK-1} \nrm*{\grad G_{\cS}(\vecY^k)}^2 
+ 4\eta\sigma_1(\vecW_0)^2 \prn*{\log(K_0) + 3} \log\prn*{\frac{48K_0}{p}} 
+ \eta L_2(\vecW_0)^2B^4K_0
\\
&= -\frac{25\eta}{32} \sum_{k=0}^{\cK-1} \nrm*{\grad G_{\cS}(\vecY^k)}^2 + \tilde{O}(\epsilon^{1.5}).
\end{align*}
\end{lemma}
\begin{proof}    
We follow a similar strategy as before of combining the proof of \citet{fang2019sharp} with our self-bounding framework. To analyze $G_{\cS}(\cdot)$  we first consider an auxiliary Gradient Descent trajectory, which performs the update:
\[
\vecY^{k+1} = \vecY^k - \eta \grad G_{\cS}(\vecY^k), \quad k \geq 0, 
\]
and $\vecY^0 = \vecU^0$. $\vecY^k$ performs Gradient Descent on $G_{\cS}(\cdot)$, which is deterministic given $\vecX^0$. 

Noting $G_{\cS}$ has Hessian $\matH_{\cS}$, and that $\matH$ is the Hessian of $F$ at the point $\vecX^0 \in \cL_{F,F(\vecW_0)}$, we obtain from \pref{ass:selfbounding} that
\[ \nrm*{\matH_{\cS}} \le \nrm*{\matH} \le L_1(\vecW_0). \]
Since the following only concern $G_{\cS}$, then identically to the proof of Lemma 16, \citet{fang2019sharp}, we obtain the following:
\begin{itemize}
    \item By $L_1(\vecW_0)$-smoothness of $G_{\cS}$ (recall $G_{\cS}$ has Hessian $\matH_S$), we obtain the so-called `Descent Lemma':
    \begin{align*}
    G_{\cS}(\vecY^{k+1}) &\le G_{\cS}(\vecY^k) + \langle \grad G_{\cS}(\vecY^k), \vecY^{k+1} - \vecY^k \rangle + \frac{L_1(\vecW_0)}{2}\nrm*{\vecY^{k+1} - \vecY^k}^2. \\
    &= G_{\cS}(\vecY^{k}) - \eta\prn*{1 - \frac{L_1(\vecW_0)\eta}{2}} \nrm*{\grad G_{\cS}(\vecY^{k})}^2. 
    \end{align*}
    \item Telescoping the above for $0 \le k \le \cK-1$, and by our choice of $\eta$ which satisfies $\eta L_1(\vecW_0) \le \frac1{16}$ as per \pref{rem:restartedsgdreduceparams}, we obtain
    \[
    G_{\cS}(\vecY^{\cK})\leq G_{\cS}(\vecY^{0}) - \frac{31\eta}{32} \sum_{k=0}^{\cK-1} \nrm*{\grad G_{\cS}(\vecY^{k})}^2.\numberthis\label{eq:preliminarystep} 
    \]
\end{itemize}
To obtain \pref{lem:quaddescentanalysis}, we upper bound the difference between $\vecU^{\cK}$ and $\vecY^{\cK}$. For all $k \ge 0$, define
\[
\vecZ^{k} := \vecU^{k} - \vecY^{k}.
\]
We aim to upper bound $\vecZ^{\cK}$ (in an appropriate sense) using the concentration argument of \citet{fang2019sharp}:
\begin{lemma}[Equivalent of Lemma 17, \citet{fang2019sharp}]\label{lem:concentrationdifferencedescentsgd}  
With probability at least $1 - p/6$, we have
\[
\nrm*{\vecZ^{k}} \leq \frac{3B}{32}\approx \tilde{\Theta}\prn*{\epsilon^{0.5}}, \numberthis\label{eq:restartedsgdfasterdescentconcentrationdiff1} \]
and 
\[ \vecZ^{k^\top} \matH_{\cS} \vecZ^{k} \le 8\sigma_1(\vecW_0)^2 \eta \prn*{\log(K_0) + 1} \log\prn*{\frac{48K_0}{p}} + \eta L_2(\vecW_0)^2 B^4 K_0 \approx \tilde{\Theta}\prn*{\epsilon^{0.5}}.\numberthis\label{eq:restartedsgdfasterdescentconcentrationdiff2}
\]
Here $\tilde{\Theta}(\cdot)$ hides $F(\vecW_0)$-dependence.
\end{lemma}
\begin{proof}[Proof of \pref{lem:concentrationdifferencedescentsgd}]
Clearly $\vecZ^0=\pmb{0}$. From the definitions of $\vecU^k, \vecY^k$, we have
\begin{align*}
\vecZ^{k+1} &= \vecZ^{k} - \eta(\grad G_{\cS}(\vecU^{k}) - \grad G_{\cS}(\vecY^{k})) 
- \eta(\projSgradF(\vecX^{k}) - \grad G_{\cS}(\vecU^{k})) - \eta \vecXiS^{k+1} \\
&= (\matI - \eta \matH_{\cS})\vecZ^k - \eta(\projSgradF(\vecX^k) - \grad G_{\cS}(\vecU^k)) - \eta \vecXiS^{k+1}, \quad k \geq 0.\numberthis\label{eq:zrecursionsgd}
\end{align*}
Unraveling the above recursion gives:
\[
\vecZ^k = -\sum_{j=1}^k \eta (\matI- \eta \matH_{\cS})^{k-j} \vecXiS^j 
- \eta \sum_{j=0}^{k-1} (\matI- \eta \matH_{\cS})^{k-1-j}(\projSgradF(\vecX^j) - \grad G_{\cS}(\vecU^j)), \quad k \geq 0. \numberthis\label{eq:zrecursionsgdexpanded}
\]
Setting $k = \cK$, Triangle Inequality gives
\[
\nrm*{\vecZ^{\cK}} \leq \nrm*{\sum_{j=1}^{\cK} \eta (\matI - \eta \matH_{\cS})^{\cK-j} \vecXiS^j} 
+ \nrm*{\eta \sum_{j=0}^{\cK-1} (\matI - \eta \matH_{\cS})^{\cK-1-j}(\projSgradF(\vecX^j) - \grad G_{\cS}(\vecU^j))}.
\]
We separately bound these two terms:
\begin{itemize}
    \item For the first term, for any fixed $l$ from $1$ to $K_0$, and any $j$ from $1$ to $\min(l, \cK_0)$, we have
\[
\mathbb{E}\brk*{\eta(\matI - \eta \matH_{\cS})^{l-j} \vecXiS^j | \mathfrak{F}^{j-1}} = 0,
\nrm*{\eta(\matI - \eta \matH_{\cS})^{l-j} \vecXiS^j} \leq \eta\sigma_1(\vecW_0).
\]
The first equality uses $\nrm*{\vecXiS^j} = \nrm*{\projS \vecXi^j}$ and that the stochastic gradient oracle is unbiased. The inequality uses that $\pmb{\cP}$ is a projection matrix, $\nrm*{\vecXiS^j} = \nrm*{\projS \vecXi^j} \le \sigma_1(\vecW_0)$ which follows as $j \le \cK_0$ and \pref{lem:restartedsgdfasterdescentxikcontrol}, and $\nrm*{(\matI - \eta \matH_{\cS})^{l-j}} \le 1$ which follows as $l \ge j$ and $\matH_{\cS} \succeq 0$. (Note the importance that $j \le \cK_0$, which gives us enough control over the noise term $\vecXiS^j$.)

Now to deal with the fact that the above control only applies for certain $j$, we define a stochastic process as follows, analogously to our proof of \pref{lem:sgdboundprob}. For all fixed $1 \le l \le K_0$, define a stochastic process $Y_{l,j}$ over all $1 \le j \le l$ by:
\[ Y_{l,j} = \eta(\matI - \eta \matH_{\cS})^{l-j} \vecXiS^j 1_{j-1 < \cK} = \begin{cases} \eta(\matI - \eta \matH_{\cS})^{l-j} \vecXiS^j &: j \le \cK  \\ 0 &: j > \cK. \end{cases}\]
Recalling $\cK = \cK_0 \wedge K_0$, it's easy to check that for any \textit{fixed} $l$, $Y_{l,j}$ is $\mathfrak{F}^{j}$-measurable. Furthermore, $\eta(\matI - \eta \matH_{\cS})^{l-j}, 1_{j-1 < \cK}$ are both $\mathfrak{F}^{j-1}$-measurable. Thus combining with the earlier observations, we obtain that
\[ \mathbb{E}\brk*{Y_{l,j} | \mathfrak{F}^{j-1}} = 0, \nrm*{Y_{l,j}} \le \eta \sigma_1(\vecW_0). \]
Thus, by the Vector-Martingale Concentration Inequality \pref{thm:subgaussianazumahoeffding}, we have with probability $1 - p/(12K_0)$,
\[
\nrm*{\sum_{j=1}^l Y_{l,j}} \leq 2\eta\sigma_1(\vecW_0) \sqrt{l \log\prn*{\frac{48K_0}{p}}} \leq 2\eta\sigma_1(\vecW_0) \sqrt{K_0 \log\prn*{\frac{48K_0}{p}}} \le \frac{B}{16}.\numberthis\label{eq:restartedsgdfasterdescentconcentrationdiff1event1}
\]
The last inequality uses our choice of parameters.

By a Union Bound, with probability at least $1 - p/12$, \pref{eq:restartedsgdfasterdescentconcentrationdiff1event1} holds for all $l$ from $1$ to $K_0$. In particular, with probability at least $1-p/12$ we have for $\cK$ (recall $\cK \le K_0$) that
\[ \nrm*{\sum_{j=1}^{\cK} \eta(\matI - \eta \matH_{\cS})^{\cK-j} \vecXiS^j} = \nrm*{\sum_{j=1}^{\cK} Y_{\cK, j}} \le \frac{B}{16},\]
where we define $Y_{\cK, j}$ the obvious way. This holds because with probability at least $1-p/12$, we have the bound \pref{eq:restartedsgdfasterdescentconcentrationdiff1event1} on $\nrm*{\sum_{j=1}^l Y_{l,j}}$ irrespective of which value of $1 \le l \le K_0$ that $\cK$ takes on. 
The first equality holds by our definition of $Y_{l,j}$ for $j \le l = \cK$.
\item For the second term, we have
\begin{align*}
\nrm*{\eta \sum_{j=0}^{\cK-1} (\matI - \eta \matH_{\cS})^{\cK-1-j}(\projSgradF(\vecX^j) - \grad G_{\cS}(\vecU^j))} &\le \eta \sum_{j=0}^{\cK-1} \nrm*{\projSgradF(\vecX^j) - \grad G_{\cS}(\vecU^j)} \\
&\le \eta \sum_{j=0}^{\cK-1} \nrm*{\grad F(\vecX^j) - \grad G(\vecX^j)} \\
&\le \frac{\eta L_2(\vecW_0) B^2 K_0}2 \le \frac{B}{32}.
\end{align*}
The first inequality uses the Triangle Inequality and that $\nrm*{(\matI - \eta \matH_{\cS})^{\cK-1-j}}_2 \le 1$ for $j$ from $0$ to $\cK-1$; this follows because $\nrm*{\matH_{\cS}} \le L_1(\vecW_0)$ and as $\eta \le \frac1{L_1(\vecW_0)}$. The second inequality uses $\nrm*{\projS(\grad F(\vecX) - \grad G(\vecX))} \leq \nrm*{\grad F(\vecX) - \grad G(\vecX)}$ because $\projS$ is a projection matrix. The third inequality follows from \pref{lem:localsmoothineq}, and the fact that for all $j \le \cK-1$, $\vecX^j \in \ball(\vecX^0, B)$. The last inequality uses the choice of parameters.
\end{itemize}
Combining the above gives \pref{eq:restartedsgdfasterdescentconcentrationdiff1}, the first part of \pref{lem:concentrationdifferencedescentsgd}.

Now prove the second part of \pref{lem:concentrationdifferencedescentsgd}, namely \pref{eq:restartedsgdfasterdescentconcentrationdiff2}. Using the fact that $(\vecA+\vecb)^\top \matA(\vecA+\vecb) \leq 2\vecA^\top \matA \vecA + 2\vecb^\top \matA \vecb$ for any symmetric positive definite matrix $\matA$ and the recursion \pref{eq:zrecursionsgdexpanded} for $\vecZ^k$, we have
\begin{align*}
&\prn*{\vecZ^{\cK}}^\top \matH_{\cS} \vecZ^{\cK} \\
&\leq 2\eta^2 \prn*{\sum_{j=1}^{\cK} (\matI - \eta \matH_{\cS})^{\cK-j-1}}^\top \matH_{\cS} \prn*{\sum_{j=1}^{\cK} (\matI - \eta \matH_{\cS})^{K-j} \vecXi_u^j} \\
&+ 2\eta^2 \left(\sum_{j=0}^{\cK-1} (\matI - \eta \matH_{\cS})^{\cK-1-j} \left(\projSgradF(\vecX^j) - \grad G_{\cS}(\vecU^j)\right)\right)^\top \matH_{\cS} 
\left(\sum_{j=0}^{\cK-1} (\matI - \eta \matH_{\cS})^{\cK-1-j} \left(\projSgradF(\vecX^j) - \grad G_{\cS}(\vecU^j)\right)\right)\\
&=2\nrm*{\eta \sum_{j=1}^{\cK} \matH_{\cS}^{1/2} (\matI - \eta \matH_{\cS})^{\cK-j} \vecXiS^j}^2 \\
&+2\eta^2 \sum_{j=0}^{\cK-1} \sum_{l=0}^{\cK-1} \left(\projSgradF(\vecX^j) - \grad G_{\cS}(\vecU^j)\right)^\top (\matI - \eta \matH_{\cS})^{\cK-1-j} \matH_{\cS} (\matI- \eta \matH_{\cS})^{\cK-1-l} 
\left(\projSgradF(\vecX^l) - \grad G_{\cS}(\vecU^l)\right)  \\
&\leq 2\nrm*{\eta \sum_{j=1}^{\cK} \matH_{\cS}^{1/2} (\matI - \eta \matH_{\cS})^{\cK-j} \vecXiS^j}^2+ 2\eta^2 \frac{L_2(\vecW_0)^2 B^4}{4} \sum_{j=0}^{\cK-1} \sum_{l=0}^{\cK-1} 
\nrm*{(\matI - \eta \matH_{\cS})^{\cK-1-j} \matH_{\cS} (\matI- \eta \matH_{\cS})^{\cK-1-l}}.
\end{align*}
The last inequality follows by properties of projection matrices and by \pref{lem:localsmoothineq}, recalling that for $j \le \cK-1$, $\vecX^j \in \ball(\vecX^0, B)$.

Now we bound each of these two terms separately:
\begin{itemize}
\item For the first term, for any fixed $l, 1 \le l \le K_0$, again we define a stochastic process for any $j, 1 \le j \le l$ by:
\[ Y_{l,j} = \eta \prn*{\matH_{\cS}^{1/2} (\matI - \eta \matH_{\cS})^{l-j} \vecXiS^j} 1_{j-1 < \cK} = \begin{cases} \eta \prn*{\matH_{\cS}^{1/2} (\matI - \eta \matH_{\cS})^{l-j} \vecXiS^j} &: j \le \cK \\ 0&: j > \cK. \end{cases}\]
Analogously to earlier, recalling $\cK \le \cK_0$, for \textit{fixed} $l$, it is evident that $Y_{l,j}$ is $\mathfrak{F}^{j}$-measurable, $\eta \matH_{\cS}^{1/2} (\matI - \eta \matH_{\cS})^{l-j} 1_{j-1 < \cK}$ is $\mathfrak{F}^{j-1}$-measurable, and thus
\[ \mathbb{E}\brk*{Y_{l,j} | \mathfrak{F}^{j-1}} = 0. \]
We furthermore have
\[ \nrm*{Y_{l,j}}^2 \le \frac{\eta \sigma_1(\vecW_0)^2}{1+2(l-j)},\]
which follows by noting for any $1 \le l \le K_0$ and $j \le \cK \le \cK_0$,
\begin{align*}
\nrm*{\eta (\matH_{\cS}^{1/2} (\matI - \eta \matH_{\cS})^{l-j} \vecXiS^j)}^2 
&\le \eta^2 \nrm*{\vecXiS^j}^2 \nrm*{\matH_{\cS}^{1/2} (\matI - \eta \matH_{\cS})^{l-j} \matH_{\cS} (\matI - \eta \matH_{\cS})^{l-j}} \nrm*{\vecXiS^j}^2 \\
&\le \frac{\eta \sigma_1(\vecW_0)^2}{1+2(l-j)}.
\end{align*}
This uses the second part of \pref{lem:localsmoothineq}, that $\nrm*{\matH_{\cS}} \le L_1(\vecW_0)$, that $j \le \cK_0$ which gives $\nrm*{\vecXiS^j} \le \sigma_1(\vecW_0)$ by \pref{lem:restartedsgdfasterdescentxikcontrol}, and our choice of $\eta$ (which cancels one of the $\sigma_1(\vecW_0)^2$ factors).

For a given $l$, by the Vector-Martingale Concentration Inequality \pref{thm:subgaussianazumahoeffding}, we have with probability $1 - p/(12K_0)$ that
\begin{align*}
\nrm*{\sum_{j=1}^l Y_{l,j} }^2 &\leq 4\eta \sigma_1(\vecW_0)^2 \log\left(\frac{48K_0}{p}\right) \sum_{j=1}^l \frac{1}{1 + 2(l-j)} \\
&\le 4\eta \sigma_1(\vecW_0)^2 (\log(K_0) + 1) \log\left(\frac{48K_0}{p}\right).\numberthis\label{eq:restartedsgdfasterdescentconcentrationdiff1event2}
\end{align*}
The last step above uses $l \le K_0$, $\sum_{j=1}^l \frac{1}{1 + j} \leq \log(K_0) + 1$.

By the Union Bound, with probability at least $1 - \frac{p}{12}$, \pref{eq:restartedsgdfasterdescentconcentrationdiff1event2} holds for all $l$ from $1$ to $K_0$. Because $1 \leq \cK \leq K_0$, using the definition of $Y_{l,j}$ for $l \le \cK$, we obtain with probability at least $1 - \frac{p}{12}$ that
\[
\eta \nrm*{\sum_{j=1}^{\cK} \matHS^{1/2} (\matI - \eta \matH_{\cS})^{K-j} \vecXiS^j}^2 = \nrm*{\sum_{j=1}^{\cK} Y_{\cK,j}}^2 
\le 4\eta \sigma_1(\vecW_0)^2 (\log(K_0) + 1) \log\left(\frac{48K_0}{p}\right).
\]
\item For the second term, using the second part of \pref{lem:localsmoothineq} and that $\cK \le K_0$, and then rearranging order of the sum and performing explicit calculation yields
\begin{align*}
&\eta^2 \frac{L_2(\vecW_0)^2 B^4}{4} \sum_{j=0}^{\cK-1} \sum_{l=0}^{\cK-1} 
\nrm*{(\matI - \eta \matH_{\cS})^{\cK-1-j} \matH_{\cS} (\matI- \eta \matH_{\cS})^{\cK-1-l}}\\
&\le \eta \frac{L_2(\vecW_0)^2 B^4}{4}  \sum_{j=0}^{K_0-1} \sum_{l=0}^{K_0-1} \frac{1}{1 + j + l} \\
&\le \eta \frac{L_2(\vecW_0)^2 B^4}{4} \sum_{l=0}^{2(K_0-1)} \frac{\min(1+j, 2K_0-1-j)}{1+j} \\
&\le \frac{\eta L_2(\vecW_0)^2 B^4 K_0}{2}.
\end{align*}
\end{itemize}
Combining the above two bounds proves \pref{eq:restartedsgdfasterdescentconcentrationdiff2}, the second part of \pref{lem:concentrationdifferencedescentsgd}.
\end{proof}

We introduce one more Lemma, an intermediate step in the proof of \citet{fang2019sharp}. 
\begin{lemma}\label{lem:concentrationdifferenceinnerprodsgd}
We have with probability at least $1-p/12$ that
\[\tri*{\grad G_{\cS}(\vecY^{\cK}), \vecU^{\cK}-\vecY^{\cK}} \le \frac{3 \eta \sum_{k=0}^{\cK} \nrm*{\grad G_{\cS}(\vecY^k)}^2}{16} + 8\eta \sigma_1(\vecW_0)^2 \log(48K_0 / p) + \eta L_2(\vecW_0)^2 B^4 K_0 / 2. \]
\end{lemma}
\begin{proof}[Proof of \pref{lem:concentrationdifferenceinnerprodsgd}]
Let $\vecY^* = \arg\min_{\vecY} G_{\cS}(\vecY)$; this exists as $G$ is convex in the subspace $\cS$, by the definition of $\cS$. By the optimality condition of $\vecY^*$, we have:
\[
\projSgradF(\vecX^0) = -\matH_{\cS} \vecY^*. \numberthis\label{eq:restartedsgdconcentrationdifferenceinnerprodsgdoptcondition}
\]
Let $\tilde{\vecY}^k = \vecY^k - \vecY^*$. From the update rule of $\vecY^k$ and the optimality condition \pref{eq:restartedsgdconcentrationdifferenceinnerprodsgdoptcondition}, we obtain:
\[
\matH_{\cS} \tilde{\vecY}^k = \grad G_{\cS}(\vecY^k), \tilde{\vecY}^{k+1} = \tilde{\vecY}^k - \eta \matH_{\cS} \tilde{\vecY}^k.\numberthis\label{eq:tildeyrecursion} 
\]
Consequently, using \pref{eq:tildeyrecursion} and \pref{eq:zrecursionsgdexpanded}, we have:
\begin{align*}
&\tri*{\grad G_{\cS}(\vecY^{\cK}), \vecU^{\cK}-\vecY^{\cK}} \\
&= \tri*{ \vecYtilde^{\cK}, \vecZ^{\cK} }_{\matH_{\cS}} \\
&= \eta \sum_{k=1}^{\cK} \tri*{ \tilde{\vecY}^{k-1}, \vecXi_u^k }_{\matH_{\cS} (\matI - \eta \matH_{\cS})^{\cK-k+1}} - \eta \sum_{k=0}^{\cK-1} \tri*{ \tilde{\vecY}^k, \projSgradF(\vecX^k) - \grad G_{\cS}(\vecU^k) }_{\matH_{\cS} (\matI - \eta \matH_{\cS})^{\cK-k}}.
\end{align*}
Now we bound both of these sums in a manner similar to the proof of \pref{lem:concentrationdifferencedescentsgd}:
\begin{itemize}
\item For the first term: 
For any fixed $l$, $1 \le l \le K_0$, define a real-valued stochastic process for any $k$, $1 \le k \le \min(l, \cK_0)$ by:
\[ Y_{l,k} = \langle \vecYtilde^{k-1}, \vecXiS^k \rangle_{\matH_{\cS}(\matI - \eta \matH_{\cS})^{l-k+1}}1_{k-1<\cK}=\begin{cases} \langle \vecYtilde^{k-1}, \vecXiS^k \rangle_{\matH_{\cS}(\matI - \eta \matH_{\cS})^{l-k+1}}&: k \le \cK\\ 0 &: k > \cK.\end{cases}\]
Analogously to earlier, recalling $\cK \le \cK_0$, it's easy to check that for any \textit{fixed} $l$, $Y_{l,k}$ is $\mathfrak{F}^{k}$ measurable, and that all terms defining $Y_{l,k}$ are $\mathfrak{F}^{k-1}$ measurable except $\vecXiS^k$. 
Thus,
\[ \mathbb{E}\brk*{Y_{l,k} | \mathfrak{F}^{k-1}} = 0. \]
We furthermore have for any fixed $l, 1 \le l \le K_0$ and $k, 1 \le k \le l$,
\[ \nrm*{Y_{l,k}}^2 \le \sigma_1(\vecW_0)^2 \nrm*{\grad G_{\cS}(\vecY^{k-1})}^2.\]
To justify why the above holds, clearly this is evident for $k>\cK$. For $k \le \cK \le \cK_0$, note that
\begin{align*}
\abs*{Y_{l,k}}^2 = \nrm*{\langle \tilde{\vecY}^{k-1}, \vecXiS^k \rangle_{\matH_{\cS}(\matI - \eta \matH_{\cS})^{l-k+1}}}^2 &= \abs*{\tri*{\matHS \vecYtilde^{k-1}, \vecXiS^k}}_{(\matI - \eta \matH_{\cS})^{l-k+1}}^2 \\
&= \abs*{\tri*{\grad G_{\cS}(\vecY^{k-1}), \vecXiS^k}}_{(\matI - \eta \matH_{\cS})^{l-k+1}}^2 \\
&\le \sigma_1(\vecW_0)^2 \nrm*{\grad G_{\cS}(\vecY^{k-1})}^2 \nrm*{(\matI - \eta \matH_{\cS})^{l-k+1}}^2\\
&\le \sigma_1(\vecW_0)^2 \nrm*{\grad G_{\cS}(\vecY^{k-1})}^2.
\end{align*}
Here we used that $\matH_S$ is symmetric, that \pref{eq:tildeyrecursion}, that $\nrm*{\matI - \eta \matH_{\cS}^{l-k+1}} \le 1$ which we have argued earlier in the proof of \pref{lem:concentrationdifferencedescentsgd}, and that $\nrm*{\vecXiS^k} \le \sigma_1(\vecW_0)$ as $k \le l \le \cK_0$ by \pref{lem:restartedsgdfasterdescentxikcontrol} and properties of projection matrices.

Now for any $l, 1 \le l\le K_0$, by the Azuma–Hoeffding inequality, we have with probability at least $1 - p/(12K_0)$ that
\[ \abs*{\eta \sum_{k=1}^l Y_{l,k}} \le \sqrt{2\eta^2 \sigma_1(\vecW_0)^2 \log(24K_0/p) \sum_{k=0}^{l-1} \nrm*{\grad G_{\cS}(\vecY^k)}^2}. \]
Taking a Union Bound, it follows that with probability at least $1-p/12$, the above holds for all $l$ with $1 \le l \le K_0$. 

Because $1 \leq \cK \leq K_0$ always holds, using the definition of $Y_{l,k}$ for $k \le \cK$, we obtain with probability at least $1 - \frac{p}{12}$ that
\begin{align*}
\abs*{\eta \sum_{k=1}^{\cK} \langle \tilde{\vecY}_{k-1}, \vecXiS^k \rangle_{\matH_{\cS}(\matI - \eta \matH_{\cS})^{\cK-k+1}}}&= \abs*{\eta \sum_{k=1}^{\cK} Y_{\cK, k}} \\
&\le \sqrt{2\eta^2 \sigma_1(\vecW_0)^2 \log(24K_0/p) \sum_{k=0}^{\cK-1} \nrm*{\grad G_{\cS}(\vecY^k)}^2} \\
&\le \frac{\eta }{16} + 8 \eta \sigma_1(\vecW_0)^2 \log(48 K_0/p)
\end{align*}
where we used AM-GM in the last step. This holds because we have this upper bound on $\abs*{\sum_{k=1}^l Y_{l,k}}$ irrespective of which value of $l, 1 \le l \le K_0$ that $\cK$ takes on. The first equality holds by our definition of $Y_{l,k}$ for $k \le \cK$.

\item For the second term: note
\begin{align*}
&\eta \sum_{k=0}^{\cK-1} \tri*{ \tilde{\vecY}^k, \projSgradF(\vecX^k) - \grad G_{\cS}(\vecU^k) }_{\matH_{\cS} (\matI - \eta \matH_{\cS})^{K-k}} \\
&= \eta \sum_{k=0}^{\cK-1} \tri*{\grad G_{\cS}(\vecY^{\cK}), \projSgradF(\vecX^k) - \grad G_{\cS}(\vecU^k) }_{(\matI - \eta \matH_{\cS})^{\cK-k}} \\
&\le \eta \sum_{k=0}^{\cK-1} \nrm*{\grad G_{\cS}(\vecY^{\cK})} \nrm*{\projSgradF(\vecX^k) - \grad G_{\cS}(\vecU^k)} \\
&\le \frac{\eta \sum_{k=0}^{\cK-1}\nrm*{\grad G_{\cS}(\vecY^{\cK})}^2}{8} + 2\eta \sum_{k=0}^{\cK-1}\nrm*{\projSgradF(\vecX^k) - \grad G_{\cS}(\vecU^k)}^2 \\
&\le \frac{\eta \sum_{k=0}^{\cK-1}\nrm*{\grad G_{\cS}(\vecY^{\cK})}^2}{8} + \frac12 \eta L_2(\vecW_0)^2 B^4 K_0.
\end{align*}
The first step above uses that $\matHS$ is symmetric and \pref{eq:tildeyrecursion}. The second step uses that $k \le \cK$ and that $\nrm*{\matI- \eta \matH_{\cS}} \le 1$, as argued in the proof of \pref{lem:concentrationdifferencedescentsgd}. The third step uses AM-GM. The last step uses that $\cK \le K_0$ and \pref{lem:localsmoothineq}; for $k < \cK$, we have $\vecX^k \in \ball(\vecX^0, B)$.
\end{itemize}
Combining these above two bounds proves \pref{lem:concentrationdifferenceinnerprodsgd}.
\end{proof}

Now we finish the proof of \pref{lem:quaddescentanalysis}. As done in \citet{fang2019sharp}, we combine \pref{lem:concentrationdifferencedescentsgd}, \pref{lem:concentrationdifferenceinnerprodsgd} with \pref{eq:preliminarystep} to prove \pref{lem:quaddescentanalysis} as follows. In particular, taking a Union Bound over the events from \pref{lem:concentrationdifferencedescentsgd} and \pref{lem:concentrationdifferenceinnerprodsgd}, we obtain with probability at least $1-p/4$ that
\begin{align*}
G_{\cS}(\vecU^{\cK}) &= G_{\cS}(\vecY^{\cK}) + \langle \grad G_{\cS}(\vecY^{\cK}), \vecU^{\cK} - \vecY^{\cK} \rangle + \frac{1}{2}(\vecU^{\cK} - \vecY^{\cK})^\top \matH (\vecU^{\cK} - \vecY^{\cK}) \\
&\le G_{\cS}(\vecY^{\cK}) + \langle \grad G_{\cS}(\vecY^{\cK}), \vecU^{\cK} - \vecY^{\cK} \rangle + \frac{1}{2}(\vecU^{\cK} - \vecY^{\cK})^\top \matHS (\vecU^{\cK} - \vecY^{\cK}) \\
&\le G_{\cS}(\vecY^{\cK}) + \frac{3\eta}{16} \sum_{k=0}^{{\cK}-1} \nrm*{\nabla G_{\cS}(\vecY^k)}^2 \\
&\hspace{1in}+ 4\eta \sigma_1(\vecW_0)^2 (\log(K_0) + 3) \log(48K_0/p) + L_2(\vecW_0)^2 \eta B^4 K_0.
\end{align*}
Here the first two lines used the definition of $G_{\cS}$ and $\cS$. The last line above applied \pref{lem:concentrationdifferenceinnerprodsgd} together with the second part of \pref{lem:concentrationdifferencedescentsgd}. 

Now combining the above with \pref{eq:preliminarystep}, we obtain
\begin{align*}
G_{\cS}(\vecU^{\cK}) &\le G_{\cS}(\vecY^{\cK}) + \frac{3\eta}{16} \sum_{k=0}^{{\cK}-1} \nrm*{\nabla G_{\cS}(\vecY^k)}^2 \\
&\hspace{1in}+ 4\eta \sigma_1(\vecW_0)^2 (\log(K_0) + 3) \log(48K_0/p) + L_2(\vecW_0)^2 \eta B^4 K_0 \\
&\le G_{\cS}(\vecU^0) - \frac{25}{32} \sum_{k=0}^{{\cK}-1} \nrm*{\nabla G_{\cS}(\vecY^k)}^2 \\
&\hspace{1in}+ 4\eta \sigma_1(\vecW_0)^2 (\log(K_0) + 3) \log(48K_0/p) + \eta L_2(\vecW_0)^2 B^4 K_0,
\end{align*}
where we also used $\vecY^0=\vecU^0$. This proves \pref{lem:quaddescentanalysis}.
\end{proof}
We now analyze the orthogonal complement of $\cS$, $\cS^{\perp}$ as in \citet{fang2019sharp}, where the analysis again goes through since the iterates are `local', being prior to the escape time $\cK$:
\begin{lemma}[Equivalent of Lemma 18, \citet{fang2019sharp}]\label{lem:perpdescentsgd}
Deterministically, we have:
\[
G_{\cS^\perp}(\vecV^{\cK}) \leq G_{\cS^\perp}(\vecV^0) - \sum_{k=1}^{\cK} \eta \tri*{\grad G_{\cS^\perp}(\vecV_{\cK-1}), \vecXiSperp^k}- \frac{7\eta}{8} \sum_{k=0}^{\cK-1} \nrm*{\grad G_{\cS^\perp}(\vecX^k)}^2 + L_2(\vecW_0)^2 B^4 \eta K_0^2.
\]
Note by choice of parameters that $L_2(\vecW_0)^2 B^4 \eta K_0^2 = \tilde{O}(\epsilon^{1.5})$, where again the $\tilde{O}(\cdot)$ hides $F(\vecW_0)$-dependence.
\end{lemma}
\begin{proof}
By definition of $G_{\cS^{\perp}}$, and using definition of $\cS^{\perp}$ which implies $\matHSperp \preceq 0$, we obtain
\begin{align*}
G_{\cS^{\perp}}\left(\vecV^{k+1}\right) &=  G_{\cS^{\perp}}\left(\vecV^{k}\right) + \tri*{ \grad G_{\cS^{\perp}}\left(\vecV^{k}\right), \vecV^{k+1} - \vecV^k  } +  \frac12\prn*{\vecV^{k+1} - \vecV^{k}}^{\top}\matHSperp\prn*{\vecV^{k+1}-\vecV^{k}}\\
&\le G_{\cS^{\perp}}\left(\vecV^{k}\right) + \tri*{ \grad G_{\cS^{\perp}}\left(\vecV^{k}\right), \vecV^{k+1} - \vecV^k  }\\
&= G_{\cS^{\perp}}\left(\vecV^{k}\right) -\eta  \tri*{ \grad G_{\cS^{\perp}}\left(\vecV^{k}\right), \projSperpgradF(\vecX^k) + \vecXiSperp^{k+1} } \\
&= G_{\cS^{\perp}}\left(\vecV^{k}\right) - \eta \nrm*{\grad G_{\cS^{\perp}}(\vecV^k)}^2 - \tri*{\eta \grad G_{\cS^{\perp}}\left(\vecV^{k}\right), \projSperpgradF(\vecX^k)-\grad G_{\cS^{\perp}}\left(\vecV^{k}\right)} \\
&\hspace{1in}- \eta \tri*{\grad G_{\cS^{\perp}}\left(\vecV^{k}\right),\vecXiSperp^{k+1}} \\
&\le G_{\cS^{\perp}}\left(\vecV^{k}\right)- \eta \tri*{\grad G_{\cS^{\perp}}\left(\vecV^{k}\right),\vecXiSperp^{k+1}}-\frac{7\eta}{8}  \nrm*{\grad G_{\cS^{\perp}}(\vecV^k)}^2 + 2\eta \nrm*{\projSperpgradF(\vecX^k)-\grad G_{\cS^{\perp}}\left(\vecV^{k}\right)}^2.
\end{align*}
The last step uses AM-GM. 

Substituting and telescoping the above for $k$ from $0$ to $\cK-1$, we have:
\begin{align*}
&G_{\cS^{\perp}}(\vecV^{\cK})  \\
&\le G_{\cS^{\perp}}(\vecV^{0}) - \sum_{k=1}^{\cK} \eta \langle \grad G_{\cS^{\perp}}(\vecV^{k-1}), \vecXi_{\vecV}^k \rangle 
- \frac{7\eta}{8} \sum_{k=0}^{\cK-1} \nrm*{\grad G_{\cS^{\perp}}(\vecX^{k})}^2 
+ 2\eta \sum_{k=0}^{\cK-1} \nrm*{\grad_{\vecV} F(\vecX^{k}) - \grad G_{\cS^{\perp}}(\vecV^{k})}^2 \notag\\
&\le G_{\cS^{\perp}}(\vecV^{0}) - \sum_{k=1}^{K} \eta \langle \grad G_{\cS^{\perp}}(\vecV^{k-1}), \vecXi_{\vecV}^k \rangle 
- \frac{7\eta}{8} \sum_{k=0}^{K-1} \nrm*{\grad G_{\cS^{\perp}}(\vecX^{k})}^2 + \frac{L_2(\vecW_0)^2 B^4 \eta K_0}{2}. 
\end{align*}
Here, the second inequality uses that by \pref{lem:localsmoothineq}, for all $k \le \cK-1$, we have $\vecX^k \in \ball(\vecX^0, B)$ and so
\[
\nrm*{\projSperpgradF(\vecX^k) - G_{\cS^{\perp}}(\vecV^{k})} = \nrm*{\cP_{\cS_{\perp}}(\grad F(\vecX^k)- \grad G(\vecX^{k}))}
\leq \nrm*{\grad F(\vecX^k)- \grad G(\vecX^{k})} \le \frac{L_2(\vecW_0) B^2}{2}.
\]
This completes the proof.
\end{proof}

\paragraph{Completing the Proof:} Now we have all the ingredients in hand to prove \pref{lem:fasterdescentlemma}.
\begin{proof}[Proof of \pref{lem:fasterdescentlemma}]
Again, we follow the strategy of \citet{fang2019sharp} and adapt it to our setting here where we do not have global bounds on the Lipschitz constants of the gradient and Hessian. With \pref{lem:quaddescentanalysis} and \pref{lem:perpdescentsgd} in hand, the idea will be to show
\[
\sum_{k=0}^{\cK_0 - 1} \nrm*{\grad G_{\cS_{\perp}}(\vecV^k)}^2 + \sum_{k=0}^{\cK_0 - 1} \nrm*{\grad G_{\cS}(\vecY^k)}^2 = \tilde{\Omega}(1),
\]
and to bound the noise term
\[
-\sum_{k=1}^K \eta \langle \grad G_{\cS^\perp}(\vecV^{k-1}), \vecXiSperp^k \rangle.
\]
We break the proof of \pref{lem:fasterdescentlemma} into two cases:
\begin{enumerate}
    \item $\nrm*{\grad F(\vecX^0)} > 5\sigma_1(\vecW_0)$. 
    \item $\nrm*{\grad F(\vecX^0)} \le 5\sigma_1(\vecW_0)$.
\end{enumerate}

\textbf{Case 1:}
This case is more straightforward as the gradient is large, and will not use the quadratic approximation we developed earlier.

Consider any $k, 0 \le k \le \cK-1$. Thus $\vecX^k \in \ball(\vecX^0, B)$, and so $\vecU \in \ball(\vecX^0, B)$ for all $\vecU \in \overline{\vecX^0 \vecX^k}$. By \pref{lem:smoothnesschange}, as $\vecX^0 \in \cL_{F,F(\vecW_0)}$, we have $\nrm*{\grad^2 F(\vecU)} \le L_1(\vecW_0)$ for all such $\vecU$. Thus as $\nrm*{\grad F(\vecX^0)} > 5\sigma_1(\vecW_0)$ and by our choice of parameters, 
\[
\nrm*{\grad F(\vecX^k)} \geq \nrm*{\grad F(\vecX^0)} - \nrm*{\grad F(\vecX^k) - \grad F(\vecX^0)} \ge 5\sigma_1(\vecW_0) - L_1(\vecW_0) B \ge \frac{9}2 \sigma_1(\vecW_0).\numberthis\label{eq:noisevsgradbound}
\]
Similarly, as $\vecX^{k+1}=\vecX^k - \eta \grad \tilde{f}(\vecX^k;\vecZeta_{k+1})$ and again as $\vecX^0 \in \cL_{F,F(\vecW_0)}$, we have $\nrm*{\grad^2 F(\vecU)} \le L_1(\vecW_0)$ for all $\vecU \in \overline{\vecX^k \vecX^{k+1}}$ by \pref{lem:smoothnesschange}. Applying \pref{lem:secondordersmoothnessineq}, for all $0 \le k \le \cK-1$, we obtain:
\begin{align*}
F(\vecX^{k+1}) - F(\vecX^k) &\le \tri*{\grad F(\vecX^k), \vecX^{k+1} - \vecX^k }+ \frac{L_1(\vecW_0)}{2} \nrm*{\vecX^{k+1} - \vecX^k}^2 \\
&= -\eta \nrm*{\grad F(\vecX^k)}^2 - \eta \tri*{\grad F(\vecX^k), \vecXi^{k+1} } + \frac{L_1(\vecW_0)\eta^2}{2} \nrm*{\grad F(\vecX^k) + \vecXi^{k+1}}^2. \\
&\le -\eta \nrm*{\grad F(\vecX^k)}^2 - \eta \tri*{ \grad F(\vecX^k), \vecXi^{k+1} } + L_1(\vecW_0) \eta^2 \nrm*{\grad F(\vecX^k)}^2 + L_1(\vecW_0) \eta^2 \nrm*{\vecXi^{k+1}}^2. \\
&\le \eta\prn*{-\frac{15 }{16} +\frac{5 }{32}} \nrm*{\grad F(\vecX^k)}^2 + \frac{8}{5} \eta \sigma_1(\vecW_0)^2 + L_1(\vecW_0) \eta^2 \sigma_1(\vecW_0)^2 \\
&\le -\frac{25 \eta}{32} \nrm*{\grad F(\vecX^k)}^2 + 2 \eta \sigma^2. \\
&\le -\eta \left(\frac{25}{32} - \frac{8}{81}\right)  \nrm*{\grad F(\vecX^k)}^2. 
\end{align*}
Note here that we need to consider a bound on the Lipschitz constant of the gradient between $\vecX^{\cK-1}$ and $\vecX^{\cK}$; see \pref{rem:restartedsgdpropertiesonestep}. Here, we used the update rule of SGD, AM-GM and Young's Inequality, that $L_1(\vecW_0) \eta \le \frac1{16}$ by our choice of hyperparameters, \pref{lem:restartedsgdfasterdescentxikcontrol}, and finally \pref{eq:noisevsgradbound} in the last step.

Telescoping the above inequality from $k = 0$ to $\cK-1$, we get:
\[ F(\vecX^{\cK}) - F(\vecX^0) \leq -\eta \left(\frac{25}{32} - \frac{8}{81}\right) \sum_{k=0}^{\cK-1} \nrm*{\grad F(\vecX^k)}^2.\numberthis\label{eq:easycaselowerboundgradsumpt0} \]
To upper bound the right hand side above, note by Triangle Inequality that
\begin{align*}
\nrm*{\eta \sum_{k=0}^{\cK-1} \grad F(\vecX^k)} &= \nrm*{-\eta \sum_{k=0}^{\cK-1} \grad F(\vecX^k)} \\
&= \nrm*{\vecX^{\cK} - \vecX^0 + \eta\sum_{k=1}^{\cK} \vecXi^k} \\
&\ge \nrm*{\vecX^{\cK} - \vecX^0} - \nrm*{\eta\sum_{k=1}^{\cK} \vecXi^k} .\numberthis\label{eq:easycaselowerboundgradsumpt1}
\end{align*}
By the Vector-Martingale Concentration Inequality \pref{thm:subgaussianazumahoeffding} and the bound $\nrm*{\vecXi^k} \le \sigma_1(\vecW_0)$ for all $k \le \cK$ by \pref{lem:noisechange}, we obtain with probability at least $1-p/12$:
\[
\nrm*{\eta\sum_{k=1}^{\cK} \vecXi^k} = \nrm*{\eta\sum_{k=1}^{K_0} \vecXi^k 1_{k \le \cK}} \le 2\eta\sigma_1(\vecW_0)\sqrt{K_0 \log(48/p)} \le \frac{B}{16}. \numberthis\label{eq:easycaselowerboundgradsumpt2}
\]
Here, we used the fact that $1_{k \le \cK} \equiv 1_{k-1 < \cK}$ and consequently $1_{k \le \cK}$ is $\mathfrak{F}^{k-1}$-measurable, and that $\mathbb{E}\brk*{\vecXi^k|\mathfrak{F}^{k-1}}=0$, $\nrm*{\vecXi^k} \le \sigma_1(\vecW_0)$ for all $k \le \cK$.

Suppose the above event implying \pref{eq:easycaselowerboundgradsumpt2} occurs, which has probability at least $1-\frac{p}{12}$. Under this event, suppose that $\vecX^k$ is able to leave the ball $\ball(\vecX^0, B)$ in $K_0$ iterations or less. If this is the case, then we have $\cK=\cK_0 \le K_0$, and so $\nrm*{\vecX^{\cK}-\vecX^0} \ge B$. Thus conditioned on the aforementioned event implying \pref{eq:easycaselowerboundgradsumpt2}, if $\vecX^k$ is able to leave the ball $\ball(\vecX^0, B)$ in $K_0$ iterations or less, we obtain
\begin{align*}
\eta \sum_{k=0}^{\cK-1} \nrm*{\grad F(\vecX^k)}^2 \ge \frac{1}{\eta \cK} \nrm*{\sum_{k=0}^{\cK-1} \eta \grad F(\vecX^k)}^2 \ge \frac{1}{\eta \cK} \prn*{B-\frac{1}{16}B}^2 \ge \frac{15^2 B^2}{16^2 \eta \cK} \ge \frac{15^2 B^2}{16^2 \eta K_0},
\end{align*}
where we combined \pref{eq:easycaselowerboundgradsumpt1}, \pref{eq:easycaselowerboundgradsumpt2} to lower bound $\nrm*{\sum_{k=0}^{\cK-1} \eta \grad F(\vecX^k)}$. 
Here the first step holds by the elementary inequality $\nrm*{\sum_{i=0}^{l} \vecA_i}^2 \le l\sum_{i=0}^{l} \nrm*{\vecA_i}^2$, and the last step uses $K_0 \ge \cK$.

Consequently by combining with \pref{eq:easycaselowerboundgradsumpt0}, with probability at least $1-\frac{p}{12}$, if $\vecX^k$ is able to leave the ball $\ball(\vecX^0, B)$ in $K_0$ iterations or less, we have
\[
F(\vecX^{\cK}) \leq F(\vecX^0) - \left(\frac{25}{32} - \frac{8}{81}\right) \cdot \frac{15^2 B^2}{16^2 \eta K_0} < F(\vecX^0) - \frac{B^2}{7 \eta K_0}.
\]

\textbf{Case 2:} Suppose $\nrm*{\grad F(\vecX^0)} \le 5\sigma_1(\vecW_0)$. To obtain the desired result, we first define and prove the following Lemmas. Proving these Lemmas in turn utilizes the Lemmas on quadratic approximation we have established earlier.

\begin{lemma}\label{lem:upperboundgradperpsmallnoise}
For all $0 \le k \le \cK-1$, we have
\[ \nrm*{\grad G_{\cS^\perp}(\vecV^k)} \le \frac{11}2\sigma_1(\vecW_0).\]
\end{lemma}
\begin{proof}
By the condition in this case, properties of projection matrices, and as $\vecV^0=0$,
\[
\nrm*{\grad G_{\cS^\perp}(\vecV^0)} = \nrm*{\projSperpgradF(\vecX^0)} \leq \nrm*{\grad F(\vecX^0)} \leq 5\sigma_1(\vecW_0).
\]
Note for $k \le \cK-1$, we have 
\[
\nrm*{\vecV^k - \vecV^0} = \nrm*{\projSperp(\vecX^k - \vecX^0)} \leq B.
\]
Thus
\begin{align*}
\nrm*{\grad G_{\cS^\perp}(\vecV^k)} &\leq \nrm*{\grad G_{\cS^\perp}(\vecV^0)} + \nrm*{\grad G_{\cS^\perp}(\vecV^k) - \grad G_{\cS^\perp}(\vecV^0)} \\
&\leq 5\sigma_1(\vecW_0) + L_1(\vecW_0) B \\
&\leq \frac{11}{2}\sigma.    
\end{align*}
The above uses our choice of hyperparameters, and that 
\[ \nrm*{\grad G_{\cS^\perp}(\vecV^k) - \grad G_{\cS^\perp}(\vecV^0)} = \nrm*{\matHSperp (\vecV^k - \vecV^0)} \le \nrm*{\matH} \nrm*{\vecV^k-\vecV^0} \le L_1(\vecW_0)\nrm*{\vecV^k-\vecV^0},\]
which in turn follows because $\vecX^0 \in \cL_{F,F(\vecW_0)}$ and by \pref{ass:selfbounding}.
\end{proof}
The next Lemma is obtained by combining \pref{lem:quaddescentanalysis} and \pref{lem:perpdescentsgd}, and it gives us a way to upper bound $F(\vecX^k)-F(\vecX^0)$.
\begin{lemma}[Equivalent of Lemma 19 in \citet{fang2019sharp}]\label{lem:sgdsecondorderdescentboundhardcase}
If $\nrm*{\grad F(\vecX^0)} \leq 5\sigma_1(\vecW_0)$, with probability $1 - \frac{p}4$, we have
\begin{align*}
F(\vecX^{\cK}) \leq &F(\vecX^0) - \eta \sum_{k=1}^{\cK} \tri*{\grad G_{\cS^\perp}(\vecV^{k-1}), \vecXiSperp^k } 
+ \left(\frac{3}{256} + \frac{1}{80}\right) \frac{B^2}{\eta K_0} \\
&- \frac{7\eta}{8} \sum_{k=0}^{\cK-1} \nrm*{\grad G_{\cS^\perp}(\vecV^k)}^2 
- \frac{25\eta}{32} \sum_{k=0}^{\cK-1} \nrm*{\grad G_{\cS}(\vecY^k)}^2. 
\end{align*}
\end{lemma}
\begin{proof} 
For $k \le \cK-1$, we have $\vecX^k \in \ball(\vecX^0, B)$. Consequently the entire line segment $\overline{\vecX^0 \vecX^k}$ lies in $\ball(\vecX^0, B)$. As $\vecX^0 \in \cL_{F,F(\vecW_0)}$, by \pref{lem:smoothnesschange}, we have
\[ \nrm*{\grad F(\vecX^k) - \grad F(\vecX^0)} \le L_1(\vecW_0) \nrm*{\vecX^k - \vecX^0} \le L_1(\vecW_0) B.\]
Thus by our choice of parameters, as per \pref{rem:restartedsgdreduceparams},
\[ \nrm*{\grad F(\vecX^k)} \leq \nrm*{\grad F(\vecX^0)} + \nrm*{\grad F(\vecX^k) - \grad F(\vecX^0)} 
\leq 5\sigma_1(\vecW_0) + L_1(\vecW_0) B \leq \frac{11}{2}\sigma_1(\vecW_0).\]
Recalling $\nrm*{\vecXi^{\cK}} \le \sigma_1(\vecW_0)$ by \pref{lem:restartedsgdfasterdescentxikcontrol}, we obtain from our choice of parameters as per \pref{rem:restartedsgdreduceparams} that
\[
\nrm*{\vecX^{\cK} - \vecX^0} \leq \nrm*{\vecX^0 - \vecX^{\cK-1}} + \eta\nrm*{\grad F(\vecX^{\cK-1}) + \vecXi^{\cK}} \leq B + \frac{13}{2}\eta\sigma_1(\vecW_0) 
\leq B + \frac{B}{100}. \numberthis\label{eq:hardcasesgdsecondorderbounddistance}
\]
Using this, we then bound the difference between $F(\vecX^{\cK})$ and $G(\vecX^{\cK})$. As $\vecX^{\cK}=\vecX^{\cK-1} - \eta \grad \tilde{f}(\vecX^{\cK-1};\vecZeta_{\cK})$, as $\vecX^{\cK-1}\in\ball(\vecX^0, B)$, and as $\vecX^0 \in \cL_{F,F(\vecW_0)}$, we have $\nrm*{\grad^2 F(\vecU) - \grad^2 F(\vecX^0)} \le L_2(\vecW_0)\nrm*{\vecU - \vecX^0}$ for all $\vecU \in \overline{\vecX^{\cK-1} \vecX^{\cK}}$ by \pref{lem:hessianlipschitzchange}. Applying \pref{lem:thirdordersmoothnessineq} and recalling that $G_{\cS}(\vecU^{\cK}) + G_{\cS^\perp}(\vecV^{\cK}) = G(\vecX^{\cK} - \vecX^0)$, we obtain
\[
F(\vecX^{\cK}) - F(\vecX^0) - G_{\cS}(\vecU^{\cK}) - G_{\cS^\perp}(\vecV^{\cK}) \leq \frac{L_2(\vecW_0)}{6}\nrm*{\vecX^{\cK} - \vecX^0}^3 
\leq \frac{L_2(\vecW_0) B^3}{5}.\numberthis\label{eq:restartedsgdfasterdescenthardcasethirdorderbound1}
\]
Here, we used \pref{eq:hardcasesgdsecondorderbounddistance} in the last step. Note here that we need to consider a bound on the Lipschitz constant of the Hessian between $\vecX^{\cK-1}$ and $\vecX^{\cK}$; see \pref{rem:restartedsgdpropertiesonestep}.

Now, take a Union Bound over \pref{lem:quaddescentanalysis} and \pref{lem:perpdescentsgd}. We now add the bounds from \pref{lem:quaddescentanalysis} and \pref{lem:perpdescentsgd} to upper bound $G_{\cS}(\vecU^{\cK}) + G_{\cS^\perp}(\vecV^{\cK})$ and use that $G_{\cS}(\vecU^0) + G_{\cS^\perp}(\vecV^0) = 0$. Combining with \pref{eq:restartedsgdfasterdescenthardcasethirdorderbound1}, we obtain with probability at least $1 - p/4$ that
\begin{align*}
F(\vecX^{\cK}) &\leq F(\vecX^0) - \eta \sum_{k=1}^{\cK} \tri*{ \grad G_{\cS^\perp}(v_{k-1}), \vecXi_v^k } 
+ 4\eta\sigma_1(\vecW_0)^2 (1 + 3\log(K_0)) \log\left(\frac{48}{p}\right) \\
&- \frac{7\eta}{8} \sum_{k=0}^{\cK-1} \nrm*{\grad G_{\cS^\perp}(\vecV^k)}^2 - \frac{25\eta}{32} \sum_{k=0}^{\cK-1} \nrm*{\grad G_{\cS}(\vecY^k)}^2 
+ \frac{3L_2(\vecW_0) B^4\eta K_0}{2} + \frac{L_2(\vecW_0) B^3}{5}.\numberthis\label{eq:restartedsgdfasterdescenthardcasethirdorderbound2}
\end{align*}
Note by our choice of hyperparameters (analogous to the choice of hyperparameters from \citet{fang2019sharp}), we have the following bounds:
$
4\eta\sigma_1(\vecW_0)^2 (1 + 3\log(K_0)) \log\left(\frac{48}{p}\right) \leq \frac{B^2}{256\eta K_0}, 
\frac{3L_2(\vecW_0) B^4\eta K_0}{2} \leq \frac{B^2}{128\eta K_0}, 
\frac{L_2(\vecW_0) B^3}{5} \leq \frac{B^2}{80\eta K_0}. 
$

Combining these above inequalities with \pref{eq:restartedsgdfasterdescenthardcasethirdorderbound2}, with probability at least $1 - p/4$, we obtain
\begin{align*}
F(\vecX^{\cK}) &\leq F(\vecX^0) - \eta \sum_{k=1}^{\cK} \tri*{ \grad G_{\cS^\perp}(\vecV_{k-1}), \vecXiSperp^k }
+ \left(\frac{3}{256} + \frac{1}{80}\right) \frac{B^2}{\eta K_0} \\
&- \frac{7\eta}{8} \sum_{k=0}^{{\cK}-1} \nrm*{\grad G_{\cS^\perp}(\vecV^k)}^2 
- \frac{25\eta}{32} \sum_{k=0}^{{\cK}-1} \nrm*{\grad G_{\cS}(\vecY^k)}^2.
\end{align*}
This implies \pref{lem:sgdsecondorderdescentboundhardcase}.
\end{proof}
By \pref{lem:sgdsecondorderdescentboundhardcase}, we want to lower bound the gradient norm of $G_{\cS^{\perp}}, G_{\cS}$. We do this in the following Lemma, assuming $\vecX^k$ leaves the ball $\ball(\vecX^0, B)$ in $K_0$ iterations.
\begin{lemma}[Equivalent of Lemma 20 in \citet{fang2019sharp}]\label{lem:sgdsecondorderdescentvalboundhardcase}
With probability $1 - \frac{p}{6}$, if $\vecX^k$ exits $\ball(\vecX^0, B)$ in $K_0$ iterations (i.e. $\cK=\cK_0 \le K_0$), we have
\[
\eta \sum_{k=0}^{\cK-1} \nrm*{\grad G_{\cS^\perp}(\vecV^k)}^2 
+ \eta \sum_{k=0}^{\cK-1} \nrm*{\grad G_{\cS}(\vecY^k)}^2 
\geq \frac{169B^2}{512\eta K_0}. 
\]
\end{lemma}
\begin{proof} 
At a high level, the proof idea is similar to the proof of Case 1 earlier. Telescoping the recursions $\vecV^k = \vecV^{k-1} - \eta \vecXiSperp^{k} - \eta \projSperpgradF(\vecX^k)$ and $\vecY^k = \vecY^{k-1} - \eta \grad G_{\cS}(\vecY^k)$, we obtain
\begin{align*}
\nrm*{\eta \sum_{k=0}^{\cK-1} \prn*{ \grad G_{\cS^\perp}(\vecV^k) + \grad G_{\cS}(\vecY^k) }}  &= \nrm*{-\eta \sum_{k=0}^{\cK-1} \prn*{\grad G_{\cS^\perp}(\vecV^k) + \grad G_{\cS}(\vecY^k) }} \\
&= \nrm*{\vecV^{\cK} - \vecV^0 + \eta \sum_{k=0}^{\cK-1} \prn*{ \vecXiSperp^{k+1} - \grad G_{\cS^\perp}(\vecV^k) + \projSperpgradF(\vecX^k) } + \vecY^{\cK} - \vecY^0} \\
&\ge \nrm*{\vecV^{\cK} - \vecV^0 + \eta \sum_{k=0}^{\cK-1} \vecXiSperp^{k+1} + \prn*{ \vecU^{\cK} - \vecU^0 } - \prn*{ \vecZ^{\cK}- \vecZ^0 }} \\
&\hspace{1in}- \nrm*{\eta \sum_{k=0}^{\cK-1} \prn*{ \grad G_{\cS^\perp}(\vecV^k) - \projSperpgradF(\vecX^k) }}. 
\end{align*}
Here, we used that $\vecZ^k = \vecU^k - \vecY^k$ and the Triangle Inequality.

Next, recall $\vecX^k - \vecX^0 = \vecU^k + \vecV^k$ for all $k \ge 0$, and $\vecU^0 = \vecV^0=0$. Thus $\vecX^k - \vecX^0 = \vecV^k - \vecV^0 + \vecU^k - \vecU^0$. Furthermore notice 
\[ \grad G_{\cS^\perp}(\vecV^k) - \projSperpgradF(\vecX^k) = \matHSperp \prn*{\grad G(\vecX^k) -\grad F(\vecX^k)}.\]
For all $k \le \cK-1$ we have $\vecX^k \in \ball(\vecX^0, B)$, so as $\vecX^0 \in \cL_{F,F(\vecW_0)}$, \pref{lem:localsmoothineq} gives
\[ \nrm*{\eta \sum_{k=0}^{\cK-1} \prn*{ \grad G_{\cS^\perp}(\vecV^k) - \projSperpgradF(\vecX^k) }} \le \eta K_0 \cdot \frac{L_2(\vecW_0) B^2}2.\]
Applying these observations and Triangle Inequality again, we obtain
\begin{align*}
\nrm*{\eta \sum_{k=0}^{\cK-1} \prn*{ \grad G_{\cS^\perp}(\vecV^k) + \grad G_{\cS}(\vecY^k) }} &\ge \nrm*{\vecX^{\cK} - \vecX^0} - \nrm*{\vecZ^{\cK} - \vecZ^0} - \eta \nrm*{\sum_{k=1}^{\cK} \vecXiSperp^{k}} 
- \frac{\eta K_0 L_2(\vecW_0) B^2}{2} \\
&\ge \nrm*{\vecX^{\cK} - \vecX^0} - \nrm*{\vecZ^{\cK} - \vecZ^0} - \frac{B}{32} - \eta \nrm*{\sum_{k=1}^{\cK} \vecXiSperp^{k}}.\numberthis\label{eq:restartedsgdfasterdescentlastlem1}
\end{align*}
 and \pref{lem:localsmoothineq} combined with the fact that projection matrices do not increase norm and that $\vecX^k \in \ball(\vecX^0, B)$ for $k < \cK$, and the final statement is by the choice of hyperparameters.

Using \pref{lem:concentrationdifferencedescentsgd} and that $\vecZ^0=0$, we obtain with probability at least $1 - \frac{p}{12}$ that 
\[ \nrm*{\vecZ^{\cK} - \vecZ^0} \le \frac{3B}{32}.\numberthis\label{eq:restartedsgdfasterdescentlastlem2} \]
Now recall that $1_{k \leq \cK} \equiv 1_{k-1 < \cK}$ is $\mathfrak{F}^{k-1}$-measurable, which implies 
\[ \mathbb{E}\brk*{\vecXiSperp^k 1_{\{k \le \cK\}} | \mathfrak{F}^{k-1}}=\pmb{0}, \]
as the stochastic gradient oracle is unbiased. Furthermore, recall $\nrm*{\vecXi^k} \le \sigma_1(\vecW_0)$ for $k\le\cK$, and projection matrices do not increase norm.
Thus by the Vector-Martingale Concentration Inequality \pref{thm:subgaussianazumahoeffding}, with probability at least $1 - \frac{p}{12}$, we have 
\[
\nrm*{\eta \sum_{k=1}^{\cK} \vecXiSperp^{k}} 
= \nrm*{\eta \sum_{k=1}^{K_0} \vecXiSperp^k 1_{\{k \le \cK\}}} \leq 2\eta \sigma_1(\vecW_0) \sqrt{K_0 \log\left(\frac{48}{p}\right)}  \leq \frac{B}{16}.\numberthis\label{eq:restartedsgdfasterdescentlastlem3}
\]
Thus taking a Union Bound over the events implying \pref{eq:restartedsgdfasterdescentlastlem2}, \pref{eq:restartedsgdfasterdescentlastlem2} and combining with the earlier display \pref{eq:restartedsgdfasterdescentlastlem1}, with probability at least $1 - \frac{p}{6}$, we have
\[
\nrm*{\eta \sum_{k=0}^{\cK-1} \grad G_{\cS^\perp}(\vecV^k) + \grad G_{\cS}(\vecY^k)}
\geq \nrm*{\vecX^{\cK} - \vecX^0} - \frac{3B}{16}. 
\]
Thus with probability at least $1 - \frac{p}{6}$, if $\vecX^k$ exits $\ball(\vecX^0, B)$ in $K_0$ iterations (that is, if we have $K_0 \ge \cK$), we have
\[\nrm*{\eta \sum_{k=0}^{\cK-1} \grad G_{\cS^\perp}(\vecV^k) + \grad G_{\cS}(\vecY^k)}
\geq \nrm*{\vecX^{\cK} - \vecX^0} - \frac{3B}{16} \ge B-\frac{3B}{16}, \]
and so
\begin{align*}
\eta \sum_{k=0}^{\cK-1} \nrm*{\grad G_{\cS^\perp}(\vecV^k)}^2 + \eta \sum_{k=0}^{\cK-1} \nrm*{\grad G_{\cS}(\vecY^k)}^2 &\ge\frac{1}{2\eta \cK} \nrm*{\eta \sum_{k=0}^{\cK-1} \prn*{\grad G_{\cS^\perp}(\vecV^k) + \grad G_{\cS}(\vecY^k)}}^2 \\
&\ge \frac{1}{2\eta \cK} \prn*{B - \frac{3B}{16}}^2 = \frac{169B^2}{512\eta \cK} \ge\frac{169B^2}{512\eta K_0}.
\end{align*}
In the first step above we used the elementary inequality $\nrm*{\sum_{i=1}^l \vecA_i}^2 \leq l \sum_{i=1}^l \nrm*{\vecA_i}^2$ and Young's Inequality. This proves \pref{lem:sgdsecondorderdescentvalboundhardcase}.
\end{proof}

We now combine \pref{lem:sgdsecondorderdescentboundhardcase}, \pref{lem:sgdsecondorderdescentvalboundhardcase} to prove \pref{lem:fasterdescentlemma}. First recall by \pref{lem:sgdsecondorderdescentboundhardcase}, with probability $1 - p/4$, we have
\begin{align*}
F(\vecX^{\cK}) \leq &F(\vecX^0) - \eta \sum_{k=1}^{\cK} \tri*{\grad G_{\cS^\perp}(\vecV^{k-1}), \vecXiSperp^k } 
+ \left(\frac{3}{256} + \frac{1}{80}\right) \frac{B^2}{\eta K_0} \\
&- \frac{7\eta}{8} \sum_{k=0}^{\cK-1} \nrm*{\grad G_{\cS^\perp}(\vecV^k)}^2 
- \frac{25\eta}{32} \sum_{k=0}^{\cK-1} \nrm*{\grad G_{\cS}(\vecY^k)}^2. \numberthis\label{eq:finalupperboundsecondordersgd}
\end{align*}
We first control $\sum_{k=1}^{\cK} \tri*{ \grad G_{\cS^\perp}(\vecV^{k-1}), \vecXiSperp^{k}}$ by concentration. For all $k$ from $1$ to $K_0$, note
\[
\mathbb{E}\brk*{\eta \tri*{ \grad G_{\cS^\perp}(\vecV^{k-1}), \vecXiSperp^{k} } 1_{k \leq \cK} | \mathfrak{F}_{k-1}} = 0,
\]
because $1_{k \leq \cK} \equiv 1_{k-1 \leq \cK}$, so all terms in $\eta \tri*{ \grad G_{\cS^\perp}(\vecV^{k-1}), \vecXiSperp^{k} } 1_{k \leq \cK}$ except $\vecXiSperp^k$ are $\mathfrak{F}^{k-1}$-measurable.

Furthermore, by \pref{lem:upperboundgradperpsmallnoise} and \pref{lem:restartedsgdfasterdescentxikcontrol}, for all $k \le \cK$, we have
\[
\nrm*{\eta \tri*{ \grad G_{\cS^\perp}(\vecV^{k-1}), \vecXiSperp^{k} } 1_{k \leq \cK}} \le \frac{11 \eta \sigma_1(\vecW_0)^2}2,
\]
and
\[
\mathbb{E}\brk*{\crl*{\eta \tri*{ \grad G_{\cS^\perp}(\vecV^{k-1}), \vecXiSperp^{k} } 1_{k \leq \cK}}^2 | \mathfrak{F}^{k-1}}
\leq \eta^2 \sigma_1(\vecW_0)^2 1_{k \leq K} \nrm*{\grad G_{\cS^\perp}(\vecV^k)}^2.
\]
Taking $\delta = \frac{p}{3 \log(K_0)}$ in the Data-Dependent Bernstein Inequality \pref{thm:subgaussiandatadependentconcentration}, we obtain with probability at least $1 - \frac{p}{3}$,
\begin{align*}
&\sum_{k=1}^{\cK} -\eta \tri*{ \grad G_{\cS^\perp}(\vecV^{k-1}), \vecXiSperp^{k} } \\
&=\sum_{k=1}^{K_0} -\eta \tri*{ \grad G_{\cS^\perp}(\vecV^{k-1}), \vecXiSperp^{k} } 1_{k \leq \cK} \\
&\leq \max \left\{
11\eta \sigma_1(\vecW_0)^2 \log\left(\frac{3 \log(K_0)}{p}\right),
4 \sqrt{\eta^2 \sigma_1(\vecW_0)^2 \sum_{k=0}^{\cK-1} \nrm*{\grad G_{\cS^\perp}(\vecV^k)}^2 \log\left(\frac{3 \log(K_0)}{p}\right)}
\right\}.\numberthis\label{eq:secondordersgdlastconcentration}
\end{align*}
We upper bound each of these terms in the maximum. With our choice of parameters and one application of AM-GM, we have
\[
11\eta \sigma_1(\vecW_0)^2 \log\left(\frac{3 \log(K_0)}{p}\right) \leq \frac{B^2}{100\eta K_0},
\]
and
\begin{align*}
4\sqrt{\eta^2\sigma^2 \sum_{k=0}^{K-1} \nrm*{\grad G_{\cS^\perp}(\vecV^{k})}^2 \log\left(\frac{3 \log(K_0)}{p}\right)} 
&\le 32 \log\left(\frac{3 \log(K_0)}{p}\right)\eta \sigma_1(\vecW_0)^2 + \frac{\eta}{8} \sum_{k=0}^{\cK-1} \nrm*{\grad G_{\cS^\perp}(\vecV^k)}^2 \\
&\le \frac{B^2}{32\eta K_0} + \frac{\eta}{8} \sum_{k=0}^{\cK-1} \nrm*{\grad G_{\cS^\perp}(\vecV^k)}^2.
\end{align*}
Consequently the second upper bound dominates the maximum from \pref{eq:secondordersgdlastconcentration}. Substituting the above into \pref{eq:secondordersgdlastconcentration}, with probability at least $1 - \frac{p}{3}$, we obtain
\[
\sum_{k=1}^{\cK} -\eta \tri*{ \grad G_{\cS^\perp}(\vecV^{k-1}), \vecXiSperp^{k} }
\leq \frac{B^2}{32\eta K_0} + \frac{\eta}{8} \sum_{k=0}^{\cK-1} \nrm*{\grad G_{\cS^\perp}(\vecV^k)}^2. 
\]
Combining with \pref{eq:finalupperboundsecondordersgd}, we obtain with probability at least $1-\frac{7p}{12}$ that
\begin{align*}
F(\vecX^{\cK}) - F(\vecX^0) &\le 
\left(\frac{3}{256} + \frac{1}{80} + \frac{1}{32}\right)\frac{B^2}{\eta K_0} 
- \frac{3\eta}{4} \sum_{k=0}^{\cK-1} \nrm*{\grad G_{\cS^\perp}(\vecV^k)}^2 
- \frac{3\eta}{4} \sum_{k=0}^{\cK-1} \nrm*{\grad G_{\cS}(\vecY^k)}^2 
\end{align*}
Taking a Union Bound with the event from \pref{lem:sgdsecondorderdescentvalboundhardcase}, we obtain with probability at least $1-\frac34 p$, if $\vecX^k$ moves out of the ball $\vecB(\vecX^0, B)$ within $K_0$ iterations (i.e. $\cK=\cK_0 \le K_0$), then
\begin{align*}
F(\vecX^{\cK_0}) - F(\vecX^0) &=F(\vecX^{\cK}) - F(\vecX^0) \le -\left(\frac{3}{4} \cdot \frac{169}{512} - \frac{3}{256} - \frac{1}{80} - \frac{1}{32}\right)\frac{B^2}{\eta K_0} < -\frac{B^2}{7\eta K_0}.
\end{align*}
This proves \pref{lem:fasterdescentlemma} in Case 2. 

Combining Case 1 and Case 2, we obtain \pref{lem:fasterdescentlemma}.
\end{proof}

\subsection{Finding Second Order Stationary Points}
Here, we finish the proof by showing with high probability, if the algorithm does not escape $\ball(\vecX^0, B)$ in $K_0$ iterates, then the average of the $K_0$ iterates is a SOSP. In particular, we aim to prove \pref{lem:sgdfindingssp}. Here is where \pref{lem:sgdescapesaddles} is used. In the following, we define $\vecXi^k$ as in \pref{eq:vecxikdef}.
Furthermore, note the proofs of \pref{lem:localsmoothineq} and \pref{lem:restartedsgdfasterdescentxikcontrol} still go through under the conditions of \pref{lem:sgdfindingssp}, so we may apply those Lemmas in our proof here.

\begin{proof}
We adopt the proof strategy of \cite{fang2019sharp} in a similar way as we have thus far. 
\begin{itemize}
\item By \pref{lem:sgdescapesaddles}, with probability $1-\frac{p}3$ (namely if the event \pref{eq:highprobevent_escapesaddle} from \pref{lem:sgdescapesaddles} occurs), then if $\lambdamin(\grad^2 F(\vecXbar)) \le -\delta_2$, $\vecX^k$ will move out of the ball $\ball(\vecX^0, B)$ within $K_0$ iterations. By taking the contrapositive, we see that with probability $1-\frac{p}3$, if $\vecX^k$ does not move out of the ball $\ball(\vecX^0, B)$ in $K_0$ iterations, then $\lambdamin(\grad^2 F(\vecX^0)) \ge -\delta_2$. In this case, we have $\vecX^k \in \ball(\vecX^0, B)$ for all $1 \le k \le K_0$, so $\vecXbar \in \ball(\vecX^0, B)$. Thus by \pref{lem:hessianlipschitzchange} and as $\vecX^0 \in \cL_{F,F(\vecW_0)}$,
\[ \lambdamin(\grad^2 F(\vecXbar)) \ge \lambdamin(\grad^2 F(\vecX^0)) - L_2(\vecW_0) \nrm*{\vecXbar - \vecX^0} \ge -\delta_2 - L_2(\vecW_0) B \ge -17\delta,\]
where the final inequality follows from our choice of parameters. That is, with probability $1-\frac{p}3$, if $\vecX^k$ does not move out of the ball $\ball(\vecX^0, B)$ in $K_0$ iterations, then $\lambdamin(\grad^2 F(\vecXbar)) \ge -17\delta$.

\item To complete the proof and show $\vecXbar$ is a SOSP, we will show that $\nrm*{\grad F(\vecXbar)}$ is small. To this end, we upper bound $\frac1{K_0} \nrm*{ \sum_{k=1}^{K_0} \vecXi^k}$ using concentration. In deriving this bound we do \textit{not} yet suppose that $\vecX^k$ does not move out of $\ball(\vecX^0, B)$ in its first $K_0$ iterations. Consider
\[ \nrm*{ \sum_{k=1}^{K_0} \vecXi^k 1_{k \le \cK_0}}=\nrm*{ \sum_{k=1}^{K_0} \vecXi^k 1_{k-1 < \cK_0}}.\]
As $1_{k-1<\cK_0}$ is $\mathfrak{F}^{k-1}$-measurable, 
\[ \mathbb{E}\brk*{\vecXi^k 1_{k \le \cK_0}|\mathfrak{F}^{k-1}}=\pmb{0}.\]
Furthermore by \pref{lem:restartedsgdfasterdescentxikcontrol}, for $k \le \cK_0$ we have
\[ \nrm*{\vecXi^k 1_{k \le \cK_0}} \le \sigma_1(\vecW_0). \]
Thus the Vector-Martingale Concentration Inequality \pref{thm:subgaussianazumahoeffding} gives with probability at least $1-2p/3$ that
\[ \frac1{K_0}\nrm*{ \sum_{k=1}^{K_0} \vecXi^k 1_{k \le \cK_0}} \le \frac{2\sigma_1(\vecW_0)\sqrt{K_0 \log(6/p)}}{K_0} \le L_2(\vecW_0) B^2.\numberthis\label{eq:restartedsgdfindingsospeq1}\]
The last inequality follows from our choice of parameters. 

Now conditioning on the above event implying \pref{eq:restartedsgdfindingsospeq1} which occurs with probability at least $1-2p/3$, suppose $\vecX^k$ does not move out of the ball $\ball(\vecX^0, B)$ in $K_0$ iterations. Then we have $\cK_0 > K_0$, and so from \pref{eq:restartedsgdfindingsospeq1}, we have
\[ \frac1{K_0} \nrm*{ \sum_{k=1}^{K_0} \vecXi^k}=\frac1{K_0}\nrm*{ \sum_{k=1}^{K_0} \vecXi^k 1_{k \le \cK_0}} \le L_2(\vecW_0) B^2.\]
Furthermore, if $\vecX^k$ does not move out of the ball $\ball(\vecX^0, B)$ in $K_0$ iterations, then we have $\vecXbar \in \ball(\vecX^0, B)$. We find an upper bound $\nrm*{\grad F(\vecXbar)}^2$. We again consider the quadratic approximation $G(\vecX)$ at $\vecX^0$ defined in \pref{subsec:fasterdescent}, and follow the notation from there. Noting $G(\cdot)$ is a quadratic and so its gradient is a linear map, we obtain
\begin{align*}
\nrm*{G(\vecXbar)} &= \nrm*{\frac1{K_0} \sum_{k=0}^{K_0-1} \grad G(\vecX^k)} \\
&\le \nrm*{\frac1{K_0} \sum_{k=0}^{K_0-1} \grad F(\vecX^k)}+\nrm*{\frac1{K_0} \sum_{k=0}^{K_0-1} \grad G(\vecX^k)-\grad F(\vecX^k)} \\
&= \frac1{K_0 \eta} \nrm*{\vecX^{K_0-1} - \vecX^0 - \eta \sum_{k=1}^{K_0} \vecXi^k} + \nrm*{\frac1{K_0} \sum_{k=0}^{K_0-1} \grad G(\vecX^k)-\grad F(\vecX^k)} \\
&\le \frac{B}{K_0 \eta} + \frac1{K_0} \nrm*{ \sum_{k=1}^{K_0} \vecXi^k} + \frac1{K_0} \cdot K_0 \cdot \frac{L_2(\vecW_0) B^2}{2} \\
&\le \prn*{\frac{16}{\tilde{C_1}}+\frac12} L_2(\vecW_0) B^2 + \frac1{K_0} \nrm*{ \sum_{k=1}^{K_0} \vecXi^k}.
\end{align*}
Here we used the choice of parameters, that $\vecX^k \in \ball(\vecX^0, B)$ for all $0 \le k \le K_0$ combined with \pref{lem:localsmoothineq} and that $\vecX^0 \in \cL_{F,F(\vecW_0)}$, and Triangle Inequality repeatedly.

Note because $\vecX^0 \in \cL_{F,F(\vecW_0)}$ and as $\vecXbar \in \ball(\vecX^0, B)$, by \pref{lem:localsmoothineq}, the above implies
\[ \nrm*{\grad F(\vecXbar)} \le \nrm*{\grad G(\vecXbar)} + \frac{L_2(\vecW_0)B^2}2 \le 17 L_2(\vecW_0) B^2 + \frac1{K_0}\nrm*{ \sum_{k=1}^{K_0} \vecXi^k} \le 18 L_2(\vecW_0) B^2.\]
Consequently, with probability at least $1-2p/3$, if $\vecX^k$ does not move out of the ball $\ball(\vecX^0, B)$ within $K_0$ iterations, then 
\[ \nrm*{\grad F(\vecXbar)} \le18 L_2(\vecW_0) B^2.\]
\end{itemize}
Taking a Union Bound, it follows that with probability at least $1-p$, if $\vecX^k$ does not escape $\ball(\vecX^0, B)$ within the first $K_0$ iterations, we have both
\[ \nrm*{\grad F(\vecXbar)} \le 18 L_2(\vecW_0) B^2, \lambdamin(\grad^2 F(\vecXbar)) \ge -17 \delta.\]
This proves \pref{lem:sgdfindingssp}.
\end{proof}

\section{Examples}\label{sec:examplesproofs}
\subsection{Phase Retrieval}\label{subsec:examplesappendixphaseretrieval}
By \pref{thm:escapesecondordergd} and \pref{thm:sgdsecondorder}, it suffices to show that 1) $F_{\text{pr}}$ satisfies \pref{ass:thirdorderselfbounding} and 2) $F_{\text{pr}}$ is a strict saddle problem (that is, all SOSPs are near-optima in a suitable sense). In the rest of this subsection, denote $F_{\text{pr}}$ by $F$ for short. As shown in \citet{candes2015phase, priorpaper}, Section 2.3 and Lemma 16 part a respectively, direct calculation shows $F(\vecW)$ takes the form
\[ F(\vecW) = \vecW^{\top}(\matI - (\vecW^{\star})(\vecW^{\star})^{\top})\vecW + \frac34(\nrm*{\vecW}^2-1)^2.\numberthis\label{eq:phaseretrievalloss}\]
As $\nrm*{\vecW^{\star}}=1$, we have $F(\vecW) \ge 0$. Furthermore, we have $\inf_{\vecW\in\mathbb{R}^d}F(\vecW)=0$, attained for example at $\vecW=\pm \vecW^{\star}$. Also note for any fixed $\vecW$, $F$ is absolutely continuous on a compact neighborhood of $\vecW$.

\paragraph{$F$ satisfies \pref{ass:thirdorderselfbounding}:} By \citet{priorpaper}, Lemma 20, we have that
\[ \nrm*{\grad^2 F(\vecW)} \le \rho_1(F(\vecW))\]
for $\rho_1(x)=9\sqrt{x}+10$. It remains to show that 
\[ \nrm*{\grad^3 F(\vecW)} \le \rho_2(F(\vecW))\]
for some increasing, non-negative $\rho_2$, where $\nrm*{\grad^3 F(\vecW)}$ refers to operator norm of the third order tensor. Equivalently, we will show that for any $\vecW$ and any unit vector $\vecU$, we have
\[ \lim_{\delta \rightarrow 0}\frac{\nrm*{\grad^2 F(\vecW+\delta\vecU) - \grad^2 F(\vecW)}_{\OPNORM}}{\delta \nrm*{\vecU}} \le  \rho_2(F(\vecW)). \]
As shown in the proof of Lemma 20, \citet{priorpaper}, we obtain from direct calculation that
\[ \grad^2 F(\vecW)=2\matI - 2 (\vecW^{\star})(\vecW^{\star})^{\top}+ 3(\nrm*{\vecW}^2-1)\matI+6\vecW \vecW^{\top}.\numberthis\label{eq:phaseretrievalhessian}\]
Thus, by repeatedly applying Triangle Inequality and \pref{lem:upperboundouterproductopnorm} and as $\nrm*{\vecU}=1$, 
\begin{align*}
&\nrm*{\grad^2 F(\vecW+\delta\vecU) - \grad^2 F(\vecW)}_{\OPNORM} \\
&= \nrm*{3(\nrm*{\vecW+\delta\vecU}^2-\nrm*{\vecW}^2)\matI + 6(\vecW+\delta\vecU)(\vecW+\delta\vecU)^{\top} - 6 \vecW \vecW^{\top}}_{\OPNORM} \\
&\le 3\abs*{\nrm*{\vecW+\delta \vecU}-\nrm*{\vecW}}\cdot \prn*{\nrm*{\vecW+\delta \vecU}+\nrm*{\vecW}} \\
&\hspace{1in}+ 6\nrm*{(\vecW+\delta\vecU)(\vecW+\delta\vecU)^{\top}  - \vecW (\vecW+\delta\vecU)^{\top}+\vecW (\vecW+\delta\vecU)^{\top} -  \vecW \vecW^{\top}}_{\OPNORM} \\
&\le 3 \delta \nrm*{\vecU} \prn*{2\nrm*{\vecW}+\delta} + 6\prn*{\nrm*{\delta \vecU (\vecW+\delta \vecU)^{\top}}_{\OPNORM}+\nrm*{\vecW (\delta \vecU)^{\top}}_{\OPNORM}} \\
&\le \delta\nrm*{\vecU} \prn*{3\prn*{2\nrm*{\vecW}+\delta}+6\nrm*{\vecW+\delta\vecU}+6\nrm*{\vecW}} \\
&\le \delta\nrm*{\vecU} \prn*{18\nrm*{\vecW}+9\delta}.
\end{align*}
Here, we used the inequality $\abs*{\nrm*{\vecX+\vecY} - \nrm*{\vecX}} \le \nrm*{\vecY}$. 

Consequently, 
\[ \lim_{\delta \rightarrow 0}\frac{\nrm*{\grad^2 F(\vecW+\delta\vecU) - \grad^2 F(\vecW)}_{\OPNORM}}{\delta \nrm*{\vecU}} \le \lim_{\delta \rightarrow 0} 18\nrm*{\vecW}+9\delta \le 18\nrm*{\vecW}+1.\]
By Lemma 16 part d, \citet{priorpaper}, using Jensen's Inequality we have
\[ F(\vecW) \ge (\nrm*{\vecW}^2-1)^2.\]
Note for $\nrm*{\vecW} \ge 2$, this implies
\[ 18\nrm*{\vecW}+1 \le 18 \prn*{\nrm*{\vecW}+1}^2 \prn*{\nrm*{\vecW}-1}^2 \le 18 F(\vecW).\]
Combining with the case $\nrm*{\vecW}<2$, we obtain
\[\lim_{\delta \rightarrow 0}\frac{\nrm*{\grad^2 F(\vecW+\delta\vecU) - \grad^2 F(\vecW)}_{\OPNORM}}{\delta \nrm*{\vecU}} \le 18\nrm*{\vecW}+1 \le 18F(\vecW)+37,\]
so we can just take $\rho_2(x)=18x+37$.

\paragraph{Next, we check that $F$ is a strict saddle problem:} We check this here. Similar results, in slightly different of a setting where we solve phase retrieval from samples from data, are shown in \citet{sun2018geometric}.

Suppose $\nrm*{\grad F(\vecW)} \le \delta$ for $\delta \le (\frac1{20})^4$. Note by Lemma 16 part b, \citet{priorpaper}, $\tri*{\vecW^{\star}, \grad F(\vecW)} = 3(\nrm*{\vecW}^2-1)\tri*{\vecW, \vecW^{\star}}$. By Cauchy-Schwartz and recalling $\vecW^{\star}$ is a unit vector, this gives
\[ \delta \ge \nrm*{\vecW^{\star}}\nrm*{\grad F(\vecW)} \ge \abs*{\tri*{\vecW^{\star}, \grad F(\vecW)}} = 3\abs*{\nrm*{\vecW}^2-1} \cdot \abs*{\tri*{\vecW, \vecW^{\star}}}.\numberthis\label{eq:phaseretrievalprovesseq}\]
\begin{itemize}
\item Suppose $\abs*{\tri*{\vecW, \vecW^{\star}}} \ge \sqrt{\delta}$. Combining this with \pref{eq:phaseretrievalprovesseq} gives
\[ \abs*{\nrm*{\vecW}^2-1} \le \frac{\sqrt{\delta}}{3}.\]
By Lemma 16 part c, \citet{priorpaper},
\begin{align*}
\nrm*{\grad F(\vecW)}^2 &= 12\nrm*{\vecW}^2 F(\vecW)-8(\nrm*{\vecW}^2-\tri*{\vecW, \vecW^{\star}}^2)  \\
&=(12\nrm*{\vecW}^2-8)F(\vecW)+6(\nrm*{\vecW}^2-1)^2,
\end{align*}
where the last equality follows from the explicit form $F(\vecW)$ from \pref{eq:phaseretrievalloss}. Thus using $\abs*{\nrm*{\vecW}^2-1} \le \frac{\sqrt{\delta}}3$, we obtain
\[ \delta^2 \ge \nrm*{\grad F(\vecW)}^2 = (12\nrm*{\vecW}^2-8)F(\vecW)+6(\nrm*{\vecW}^2-1)^2 \ge (4-4\sqrt{\delta})F(\vecW).\]
For $\delta \le \frac14$, this gives
\[ F(\vecW) \le \frac{\delta^2}{4-4\sqrt{\delta}} \le \frac{\delta^2}2. \]
\item Otherwise, suppose $\abs*{\tri*{\vecW, \vecW^{\star}}} \le \sqrt{\delta}$. Note by differentiating \pref{eq:phaseretrievalloss}, as shown in the proof of Lemma 16 part b, \citet{priorpaper},
\[ \grad F(\vecW) = 2\vecW - 2\tri*{\vecW, \vecW^{\star}}\vecW^{\star}+3(\nrm*{\vecW}^2-1) \vecW = - 2\tri*{\vecW, \vecW^{\star}}\vecW^{\star}+(3\nrm*{\vecW}^2-1) \vecW.\]
Thus by Triangle Inequality,
\[ \abs*{3\nrm*{\vecW}^2-1} \cdot \nrm*{\vecW} \le \nrm*{\grad F(\vecW)} + 2\abs*{\tri*{\vecW, \vecW^{\star}}} \nrm*{\vecW^{\star}} \le \delta + 2\sqrt{\delta} \le 4\sqrt{\delta}.\]
Consequently either $\nrm*{\vecW} \le 2\delta^{1/4}$ or $\abs*{3\nrm*{\vecW}^2-1} \le 2\delta^{1/4}$. 

In the first case, by Cauchy Schwartz and \pref{eq:phaseretrievalhessian}, notice for any unit vector $\vecU$ that
\begin{align*}
\vecU^{\top} \grad^2 F(\vecW) \vecU &= \vecU^{\top}\prn*{2\matI - 2(\vecW^{\star})(\vecW^{\star})^{\top}+ 3(\nrm*{\vecW}^2-1)\matI+6 \vecW \vecW^T}]\vecU \\
&\le -\nrm*{\vecU}^2 + 3\nrm*{\vecU}^2 \cdot (2\delta^{1/4})^2 + 6 \nrm*{\vecU}^2\cdot (2\delta^{1/4})^2 \\
&\le -1+36\delta^{1/2} \le -\frac{9}{10}, 
\end{align*}
since $\delta \le (\frac1{20})^4$. 

In the second case, using \pref{eq:phaseretrievalhessian}, notice as $\nrm*{\vecW^{\star}}=1$, we have
\begin{align*}
\vecW^{\star^\top}\grad^2 F(\vecW) \vecW^{\star} &= \vecW^{\star^\top}(3\nrm*{\vecW}^2-1) \vecW^{\star} - 2 \nrm*{\vecW^{\star}}^2 + 6\abs*{\tri*{\vecW, \vecW^{\star}}}^2 \\
&\le 2\delta^{1/4} - 2 + 6\delta \le -\frac{9}{5}.
\end{align*}
Consequently in either case, $\grad^2 F(\vecW)$ has at least one negative eigenvalue with value at most $-\frac9{10}$. 
\end{itemize}
Consider $\epsilon$ smaller than a universal constant, and take $\delta=\sqrt{\epsilon}$ in the above result. It follows from the analysis here that if we find an SOSP to tolerance $\epsilon$ as per the definition \pref{eq:sospproblem}, we obtain $\vecW$ with $F(\vecW) \le \frac{\epsilon}2$.

Thus, it follows that running Perturbed GD or Restarted SGD as described in \pref{thm:escapesecondordergd} or \pref{thm:sgdsecondorder} respectively, we will obtain $\vecW$ with suboptimality $F(\vecW) \le \epsilon$, where the number of oracle calls depends on $1/\epsilon, d, F(\vecW_0)$ in the same way as in \pref{thm:escapesecondordergd} or \pref{thm:sgdsecondorder} respectively. 

\subsection{Matrix PCA}\label{subsec:examplesappendixpca}
Again by \pref{thm:escapesecondordergd}, \pref{thm:sgdsecondorder}, it suffices to show that 1) $F_{\text{pca}}$ satisfies \pref{ass:thirdorderselfbounding} and 2) is a strict saddle problem (that is, all SOSPs are near-optima in a suitable sense). We will show this, with the parameters governing the strict saddle property depending on the spectral gap 
$\lambda_1(\matM)-\lambda_2(\matM)$.\footnote{Thus our result will be vacuous when the spectral gap is 0.} In the rest of this subsection, denote $F_{\text{pca}}$ by $F$ for short. Recall the loss function for PCA takes the form
\[ F(\vecW) = \frac12 \nrm*{\vecW \vecW^{\top} - \matM}_F^2,\]
where $\matM$ is a symmetric PD matrix. Note for any fixed $\vecW$, $F$ is absolutely continuous on a compact neighborhood of $\vecW$. Note $F(\vecW) \ge 0$ always holds. While it is not true that $\inf_{\vecW\in\mathbb{R}^d}F(\vecW) = 0$, to enforce this, we can consider the shifted function $G := F - \inf_{\vecW\in\mathbb{R}^d}F(\vecW)$. The derivatives of $G$ are identical to those of $F$, and furthermore $G(\vecX) - G(\vecY) = F(\vecX) - F(\vecY)$ for all $\vecX,\vecY$. Thus to apply \pref{thm:escapesecondordergd}, \pref{thm:sgdsecondorder} and show that Perturbed GD or Restarted SGD can globally optimize $G$ and therefore $F$ by finding SOSPs, it remains to show $F$ satisfies \pref{ass:thirdorderselfbounding} and is strict saddle.

\paragraph{$F$ satisfies \pref{ass:thirdorderselfbounding}:} Direct calculation, also in \citet{jin2021nonconvex}, yields
\[ \grad F(\vecW) = (\vecW \vecW^{\top} - \matM) \vecW, \grad^2 F(\vecW) = \nrm*{\vecW}^2 \matI + 2\vecW \vecW^{\top} - \matM.\numberthis\label{eq:gradhessianpcacalc}\]
We now check self-bounding regularity for the Hessian and third order derivative tensor. First observe
\[ \vecW^{\top} (\vecW \vecW^{\top}) \vecW = \nrm*{\vecW}^4.\]
Combining with \pref{lem:upperboundouterproductopnorm}, we obtain
\begin{align*}
\nrm*{\vecW} &= \nrm*{\vecW \vecW^{\top}}_{\OPNORM}^{1/2} \\
&\le \prn*{\nrm*{\vecW \vecW^{\top} - \matM}_{\OPNORM} + \nrm*{\matM}_{\OPNORM}}^{1/2} \\
&\le \nrm*{\vecW \vecW^{\top} - \matM}_F^{1/2} + \nrm*{\matM}_{\OPNORM}^{1/2} \\
&\le 2F(\vecW)^{1/4}+ \nrm*{\matM}_{\OPNORM}^{1/2}.\numberthis\label{eq:pcaupperboundeuclideannorm}
\end{align*}
Now we check the self bounding conditions. 
For the Hessian, note from \pref{eq:gradhessianpcacalc} and \pref{eq:pcaupperboundeuclideannorm} and using \pref{lem:upperboundouterproductopnorm}, 
\[ \nrm*{\grad^2 F(\vecW)}_{\OPNORM} \le 3\nrm*{\vecW}^2 + \nrm*{\matM}_{\OPNORM} \le 3(2F(\vecW)^{1/4}+ \nrm*{\matM}_{\OPNORM}^{1/2})^2+ \nrm*{\matM}_{\OPNORM}.\]
Thus we can take $\rho_1(x) = 3(2x^{1/4}+\nrm*{\matM}_{\OPNORM}^{1/2})^2+ \nrm*{\matM}_{\OPNORM}$.

For the third order derivative tensor, following the strategy in \pref{subsec:examplesappendixphaseretrieval}, we will show that for any $\vecW$ and any unit vector $\vecU$, we have
\[ \lim_{\delta \rightarrow 0}\frac{\nrm*{\grad^2 F(\vecW+\delta\vecU) - \grad^2 F(\vecW)}_{\OPNORM}}{\delta \nrm*{\vecU}} \le  \rho_3(F(\vecW)). \]
Applying \pref{eq:gradhessianpcacalc} and \pref{lem:upperboundouterproductopnorm} and note 
\begin{align*}
(\vecW + \delta \vecU)(\vecW + \delta \vecU)^{\top} - \vecW \vecW^{\top} &= (\vecW + \delta \vecU)(\vecW + \delta \vecU)^{\top} - (\vecW + \delta \vecU)\vecW^{\top} + (\vecW + \delta \vecU)\vecW^{\top}- \vecW \vecW^{\top} \\
&= (\vecW + \delta \vecU) (\delta\vecU)^{\top} + \delta \vecU \vecW^{\top}.
\end{align*}
This gives
\begin{align*}
&\lim_{\delta \rightarrow 0}\frac{\nrm*{\grad^2 F(\vecW+\delta\vecU) - \grad^2 F(\vecW)}_{\OPNORM}}{\delta \nrm*{\vecU}} \\
&= \lim_{\delta \rightarrow 0} \frac{(\nrm*{\vecW + \delta \vecU}^2 - \nrm*{\vecW}^2) + 2\nrm*{(\vecW + \delta \vecU)(\vecW + \delta \vecU)^{\top} - \vecW \vecW^{\top}}_{\OPNORM}}{\delta \nrm*{\vecU}} \\
&\le \lim_{\delta \rightarrow 0} \frac{\abs*{\nrm*{\vecW + \delta \vecU} - \nrm*{\vecW}} \cdot (2\nrm*{\vecW} + \delta\nrm*{\vecU})  + \delta\nrm*{\vecU} (2\nrm*{\vecW}+\delta\nrm*{\vecU})}{\delta \nrm*{\vecU}}  \\
&\le \lim_{\delta \rightarrow 0} \frac{\delta\nrm*{\vecU}(2\nrm*{\vecW} + \delta\nrm*{\vecU})  + \delta\nrm*{\vecU} (2\nrm*{\vecW}+\delta\nrm*{\vecU})}{\delta \nrm*{\vecU}} \\
&= \lim_{\delta \rightarrow 0}4\nrm*{\vecW} + 2\delta\nrm*{\vecU} \\
&= 4\nrm*{\vecW} \\
&\le 8F(\vecW)^{1/4}+4\nrm*{\matM}^{1/2}_{\OPNORM}.
\end{align*}
Here we used the inequality $\abs*{\nrm*{\vecX+\vecY} - \nrm*{\vecX}} \le \nrm*{\vecY}$. The last step used \pref{eq:pcaupperboundeuclideannorm}. Thus we can take $\rho_2(x) = 8{x}^{1/4}+4\nrm*{\matM}^{1/2}_{\OPNORM}$.

\paragraph{Next, we check $F$ is a strict saddle problem:} We check this here. A similar verification is done in \citet{ge2017no}.

Let $\vecV_1, \ldots, \vecV_d$ be the (unit) eigenvectors of $\matM$ corresponding to $\lambda_1(\matM) \ge \lambda_2(\matM) \ge\cdots \ge \lambda_d(\matM)>0$ respectively (recall $\matM$ is assumed to be PD). Thus the $\vecV_i$ form an orthonormal basis of $\mathbb{R}^d$. Furthermore for convenience let $\lambda_i := \lambda_i(\matM)$ for all $1 \le i \le d$. As $\matM$ is symmetric and PD, by the Spectral Theorem, we can write 
\[ \matM = \sum_{i=1}^d \lambda_i \vecV_i \vecV_i^\top.\]
Suppose $\vecW$ is a SOSP to tolerance $\epsilon$ for $\epsilon < \min\crl*{1, \frac{(\lambda_1-\lambda_2)^2}{16}, \frac38 (\lambda_1-\lambda_2)^{5/2}}$. Note the minimizers of $F$ are $\vecW=\pm \sqrt{\lambda_1} \vecV_1$. We will show that $\vecW$ is close to these minimizers: in particular, that $\min\crl*{\nrm*{\vecW-\sqrt{\lambda_1} \vecV_1}^2, \nrm*{\vecW+\sqrt{\lambda_1} \vecV_1}^2} \le \epsilon$.

Write $\vecW = c_1 \vecV_1 + \cdots + c_d \vecV_d$. Thus, our goal is to show that $\abs*{(c_1^2+\cdots+c_d^2)-\lambda_1} < \sqrt{\epsilon}$. By \pref{eq:gradhessianpcacalc}, we have 
\begin{align*}
\epsilon &\ge \nrm*{\grad F(\vecW)} = \nrm*{\matM \vecW - \nrm*{\vecW}^2 \vecW} = \nrm*{\sum_{i=1}^d \prn*{(c_1^2+\cdots+c_d^2) - \lambda_i }c_i\vecV_i}.
\end{align*}
That is, we have 
\[ \sum_{i=1}^d c_i^2 \prn*{(c_1^2+\cdots+c_d^2) - \lambda_i }^2 \le \epsilon^2.\numberthis\label{eq:expcagradimply1}\]
Furthermore by \pref{eq:gradhessianpcacalc}, we have 
\[ \grad^2 F(\vecW) = (c_1^2+\cdots+c_d^2) \matI + 2\sum_{i, j} c_i c_j \vecV_i \vecV_j^\top - \sum_{i=1}^d \lambda_i \vecV_i \vecV_i^\top.\]
Since $\vecW$ is a SOSP, for all $\vecV_k, 1 \le k \le d$, we have 
\[ -\sqrt{\epsilon} \le \vecV_k^\top \grad^2 F(\vecW) \vecV_k = (c_1^2+\cdots+c_d^2) + 2c_k^2 - \lambda_k.\numberthis\label{eq:expcahessimply1}\]
We now break into cases:
\begin{itemize}
\item Suppose for all $i$, we have $\abs*{(c_1^2+\cdots+c_d^2)-\lambda_i} \ge \sqrt{\epsilon}$. From \pref{eq:expcagradimply1}, this gives $\sum_{i=1}^d c_i^2 \le \epsilon$. Taking $k=1$ in \pref{eq:expcahessimply1}, we obtain
\[ -\sqrt{\epsilon} \le 3\sum_{i=1}^d c_i^2 - \lambda_1 \le 3\epsilon-\lambda_1 \implies \lambda_1 \le \sqrt{\epsilon}+3\epsilon,\]
contradicting that $\epsilon < \min\crl*{1, \frac{(\lambda_1-\lambda_2)^2}{16}}$.

\item Else, suppose there exists $i$ such that $\abs*{(c_1^2+\cdots+c_d^2)-\lambda_i} < \sqrt{\epsilon}$. Suppose that $i \ge 2$. Then taking $k=1$ in \pref{eq:expcahessimply1}, we obtain
\[ -\sqrt{\epsilon} \le \lambda_i+\sqrt{\epsilon} + 2c_1^2 - \lambda_1 \implies c_1^2 \ge \frac{\lambda_1-\lambda_i}2 - \sqrt{\epsilon} \ge \frac{\lambda_1-\lambda_2}4,\]
where the last inequality uses $\lambda_i \le \lambda_2$ and $\epsilon < \prn*{\frac{\lambda_1-\lambda_2}4}^2$.

Note furthermore that as $\epsilon \le \prn*{\frac{\lambda_1-\lambda_2}4}^2$, as $\abs*{(c_1^2+\cdots+c_d^2)-\lambda_i} < \sqrt{\epsilon}$, and as $\lambda_i\le \lambda_2 < \lambda_1$, we have $\abs*{(c_1^2+\cdots+c_d^2)-\lambda_1} > \frac{3(\lambda_1-\lambda_2)}4$. Thus \pref{eq:expcagradimply1} implies
\[ \epsilon^2 > 0 + \frac{\lambda_1-\lambda_2}4 \cdot \frac{9}{16} (\lambda_1-\lambda_2)^2,\]
contradicting that $\epsilon < \frac38 (\lambda_1-\lambda_2)^{5/2}$.
\end{itemize}
Therefore, we must have $i=1$ in the second case above. That is, $\abs*{(c_1^2+\cdots+c_d^2)-\lambda_1} < \sqrt{\epsilon}$, as desired.

Thus, it follows that running Perturbed GD or Restarted SGD as described in \pref{thm:escapesecondordergd} or \pref{thm:sgdsecondorder} respectively, we will obtain $\vecW$ that is distance at most $\sqrt{\epsilon}$ from a global minimizer of $F$ for $\epsilon < \min\crl*{1, \frac{(\lambda_1-\lambda_2)^2}{16}, \frac38 (\lambda_1-\lambda_2)^{5/2}}$. Here the number of oracle calls depends on $1/\epsilon, d, F(\vecW_0)$ the same way as in \pref{thm:escapesecondordergd} or \pref{thm:sgdsecondorder} respectively. For $\epsilon\ge \min\crl*{1, \frac{(\lambda_1-\lambda_2)^2}{16}, \frac38 (\lambda_1-\lambda_2)^{5/2}}$, we can replace $\epsilon$ by any real strictly smaller than $\min\crl*{1, \frac{(\lambda_1-\lambda_2)^2}{16}, \frac38 (\lambda_1-\lambda_2)^{5/2}}$ in the guarantees from \pref{thm:escapesecondordergd} or \pref{thm:sgdsecondorder}. 

\section{Simulations}\label{sec:experiments}
\graphicspath{ {./neurips25_submission/Graphs/} }

Our algorithmic results \pref{thm:gdfirstorder}, \pref{thm:adaptivegdguarantee}, \pref{thm:sgdfirstorder}, \pref{thm:escapesecondordergd}, and \pref{thm:sgdsecondorder} have strong practical implications. They directly suggest that under generalized smoothness, the step sizes $\eta$ that lead to convergence/successful optimization become smaller for larger initialization $F(\vecW_0)$ and larger self-bounding functions $\rho_1(\cdot), \rho_2(\cdot)$. For example in \pref{thm:gdfirstorder}, we set $\eta=\frac1{L_1(\vecW_0)}$ where $L_1(\vecW_0) = \max\crl*{1, \rho_0(F(\vecW_0)+1), \rho_0(F(\vecW_0))\rho_0(F(\vecW_0)+1), \rho_1(F(\vecW_0)+1)}$ was defined in \pref{eq:L1def}.

That is, our work suggests that larger suboptimality at initialization and larger self-bounding functions shrink the `window' for choosing a working $\eta$ in practice, when the loss function satisfies generalized smoothness. This has strong practical implications: it implies that for losses with non-Lipschitz gradient/Hessian, one should tune $\eta$ based on suboptimality at initialization. This contrasts sharply with the Lipschitz gradient/Hessian case, see e.g. \citep{bubeck2015convex, jin2017escape, fang2019sharp}, where the range of working $\eta$ is fixed in terms of the Lipschitz constant of the gradient and/or Hessian, and does not depend on the initialization. 

In this section, we empirically validate this implication of our work.

\subsection{Synthetic Simulations with GD}\label{subsec:experimentstoygd}
\paragraph{Simulation Details:} We consider $F(\vecW) = \nrm*{\matA \vecW}^p$ for $p=2,3,4,5,6$, where $\matA = \text{diag}\prn*{\frac1{20}, \frac1{19}, \ldots, \frac12, 1}$. When $p=2$, $F(\vecW)$ is smooth. When $p \ge 3$, $F(\vecW)$ is not smooth, but it is straightforward to verify that it satisfies \pref{ass:selfbounding}, similar to our verifications in \pref{subsec:detailsforcomparewithlit}. One can furthermore verify that as $p$ increases, the corresponding self-bounding function $\rho_1(\cdot)$ from \pref{ass:selfbounding} increase. This choice of generalized smooth function was motivated by \citet{gaash2025convergence}, who used $\nrm*{\matA \vecW}^4$ with the exact same $\matA$ in their experiments to study optimization with first-order methods under generalized smoothness. 

For each $p=2,3,4,5,6$, we consider the following settings for GD:
\begin{itemize}
    \item Step sizes: We consider 30 step sizes $\{\eta_i\}_{i=1}^{30}, \eta_1 < \cdots < \eta_{30}$ evenly spaced on a log scale between $10^{-8}$ and $10^1$, inclusive.
    \item Initialization: For each step size $\eta_i$, we initialize GD at 4 distributions $\pi_j=\cN(\vecOrigin, c_j\matI_{20})$ for $c_j\in\{2.5, 5, 7.5, 10\}$. For each of these 4 distributions $\pi_j$, we draw 100 points $\vecW_0 \sim \pi_j$ to use as our initialization.
    \item Number of steps: For each $\eta_i$ and each $\vecW_0 \sim \pi_j$, we run GD initialized at $\vecW_0$ with step size $\eta_i$ for $T=1000$ iterations. Here as $F$ is known, we analytically compute the gradient. 
\end{itemize}
For each $p$ and initialization $\pi_j$, we consider all 30 possible $\eta_i$, which we plot on the $x$-axis. For each $\eta_i$, we consider all 100 initializations $\vecW_0 \sim \pi_j$. For each initialization $\vecW_0$, letting $\{\vecW_t\}$ be the resulting sequence of iterates of GD, we compute $\tfrac{\nrm*{\grad F(\vecW_T)}}{F(\vecW_0)}$ for $T=1000$. For $\eta_i$ that led to faithful convergence of GD, on the $y$-axis, we then plot the mean of $\tfrac{\nrm*{\grad F(\vecW_T)}}{F(\vecW_0)}$ over those 100 initializations as a blue dot, with blue vertical error bars indicating $\pm 2$ standard deviations. We considered the ratio $\tfrac{\nrm*{\grad F(\vecW_T)}}{F(\vecW_0)}$ because for $L$-smooth functions, established optimization theory predicts that this converges at a rate independent of $F(\vecW_0)$ and only depending on $T$ and $L$ \citep{bubeck2015convex}.

The simulations for \pref{subsec:experimentstoygd} were run on a Jupyter notebook in Python in Google Colab Pro, connected to a single NVIDIA T4 GPU. Our code can be found in the attached files.

\paragraph{Divergence of GD and working step sizes:} We observe that for some $\eta_i$ larger than some threshold depending on $p$ and $\pi_j$, the iterates of GD diverge. In particular, the resulting ratio $\tfrac{\nrm*{\grad F(\vecW_T)}}{F(\vecW_0)}$ becomes massive, often on the order of $10^5$ or more, indicating that $\eta_i$ was too large for GD to converge. To identify the smallest $\eta_i$ where this first occurs, or equivalently find the largest working step size among $\{\eta_i\}_{i=1}^{30}$, for a given $\pi_j$ and $\eta_i$, we computed the average $\tfrac{\nrm*{\grad F(\vecW_T)}}{F(\vecW_0)}$ over the 100 initializations. If this average was 100 or more times larger than this average for $\eta_{i-1}$, we took this as an indication that the iterates of GD with this step size $\eta_i$ or larger step sizes diverge, and for this $p$ and $\pi_j$, we stopped considering any larger $\eta_{i'}$, $i'>i$. We then save this $\eta_i$ to indicate the smallest $\eta_i$ for which divergence occurred. This $\eta_i$ is indicated with a red line in the following plots.

This smallest $\eta_i$ for which divergence occurred plays a crucial role in validating our theoretical claims. Established optimization theory predicts that for smooth functions (here, when $p=2$), this $\eta_i$ is identical across different initializations \citep{bubeck2015convex}. Meanwhile for generalized smooth functions, as per our remarks earlier and from \pref{subsec:examplesmainbody}, we predict that as $F(\vecW_0)$ increases, the range of working step sizes, and consequently also the smallest $\eta_i$ for which divergence occurs, will \textit{decrease}. Note as $c_j$ increases (recall $\pi_j \sim \cN(\vecOrigin, c_j \matI_{20})$ and $c_j \in \{2.5, 5, 7.5, 10\}$), we expect $F(\vecW_0)$ to increase, at least on average or with high probability over the 100 initializations $\vecW_0 \sim \pi_j$.

\begin{figure}
    \centering
    \begin{subfigure}[b]{0.46\textwidth}
     \includegraphics[width=\textwidth]{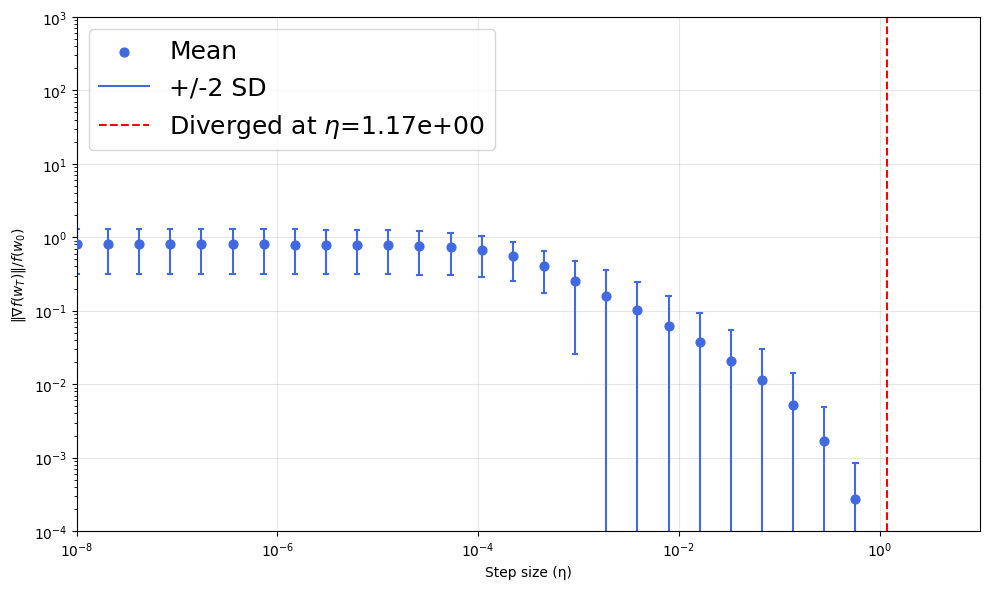}
     \caption{$\pi_j = \cN(\vecOrigin, 2.5\matI_{20})$. The first divergence is at $\eta_i\approx 1.17$.}
     \label{fig:gd_p2_cov2.5}
 \end{subfigure}
 \hfill
 \begin{subfigure}{0.46\textwidth}
     \includegraphics[width=\textwidth]{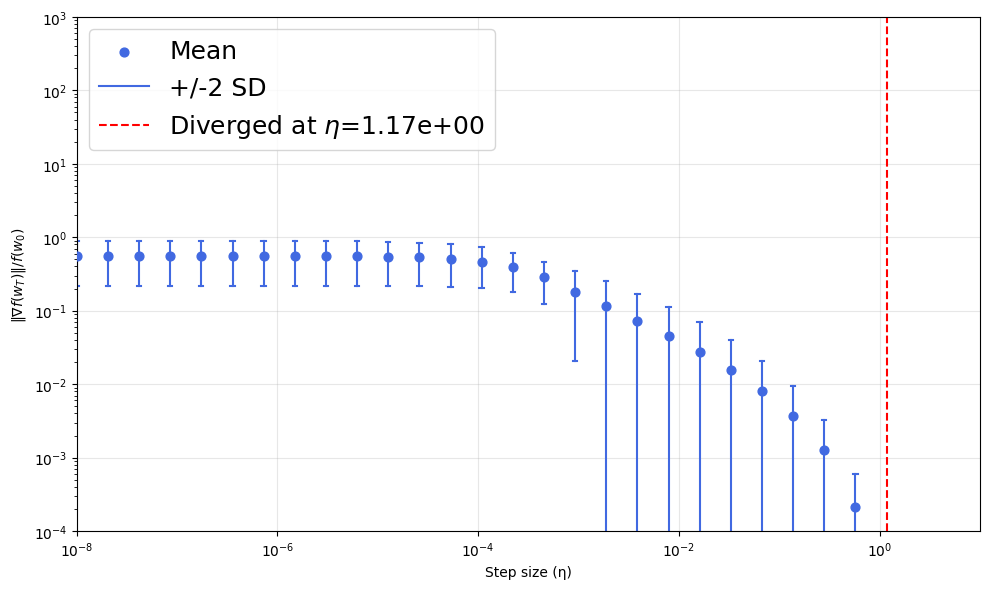}
     \caption{$\pi_j = \cN(\vecOrigin, 5.0\matI_{20})$. The first divergence is at $\eta_i\approx 1.17$.}
     \label{fig:gd_p2_cov5}
 \end{subfigure}
 \hfill
 \begin{subfigure}{0.46\textwidth}
     \includegraphics[width=\textwidth]{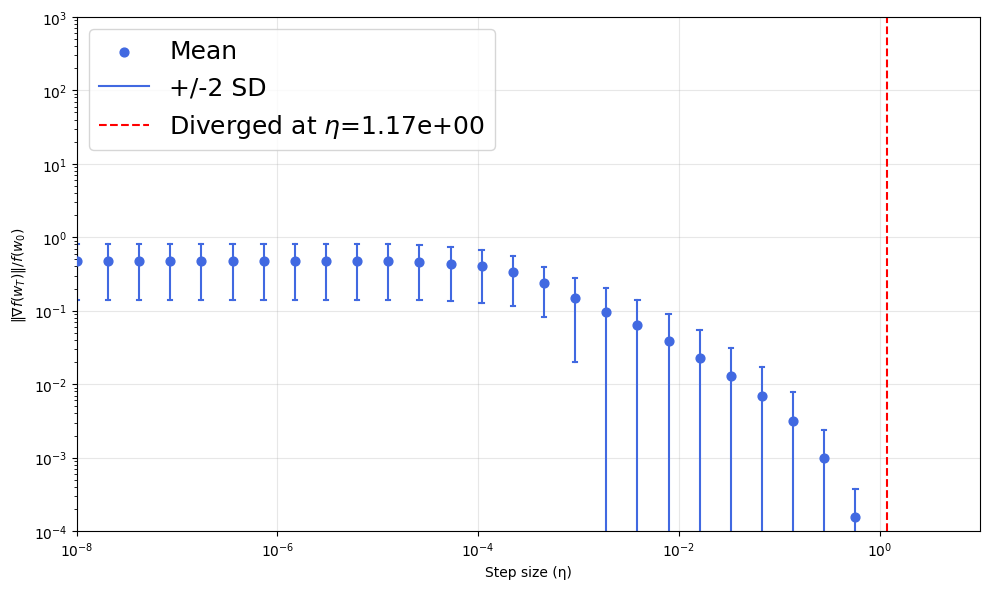}
     \caption{$\pi_j = \cN(\vecOrigin, 7.5\matI_{20})$. The first divergence is at $\eta_i\approx 1.17$.}
     \label{fig:gd_p2_cov7.5}
 \end{subfigure}
 \hfill
 \begin{subfigure}{0.46\textwidth}
     \includegraphics[width=\textwidth]{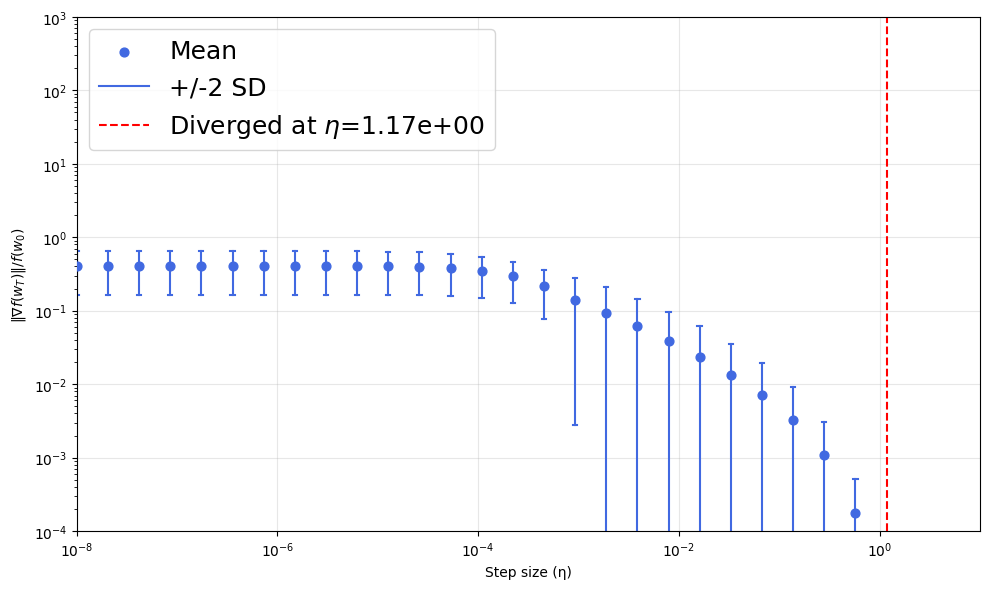}
     \caption{$\pi_j = \cN(\vecOrigin, 10\matI_{20})$. The first divergence is at $\eta_i\approx 1.17$.}
     \label{fig:gd_p2_cov10}
 \end{subfigure} 
 \caption{GD simulation results for $p=2$. For all $\pi_j$, the smallest $\eta_i$ leading to divergence is $\approx 1.17$.}
 \label{fig:gd_p2}
\end{figure}

\begin{figure}
    \centering
    \begin{subfigure}{0.46\textwidth}
     \includegraphics[width=\textwidth]{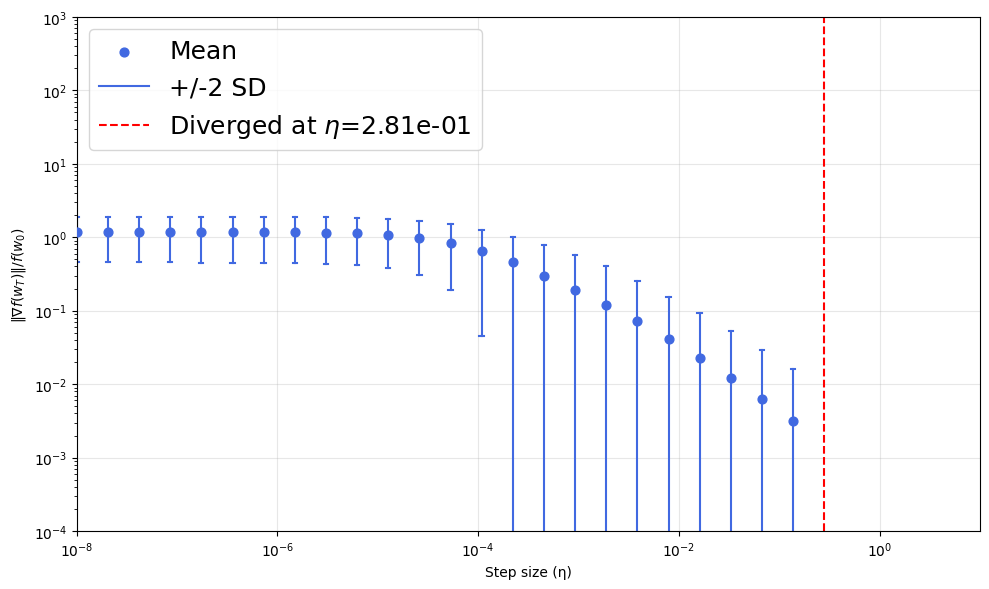}
     \caption{$\pi_j = \cN(\vecOrigin, 2.5\matI_{20})$. The first divergence is at $\eta_i\approx 0.281$.}
     \label{fig:gd_p3_cov2.5}
 \end{subfigure}
 \hfill
 \begin{subfigure}{0.46\textwidth}
     \includegraphics[width=\textwidth]{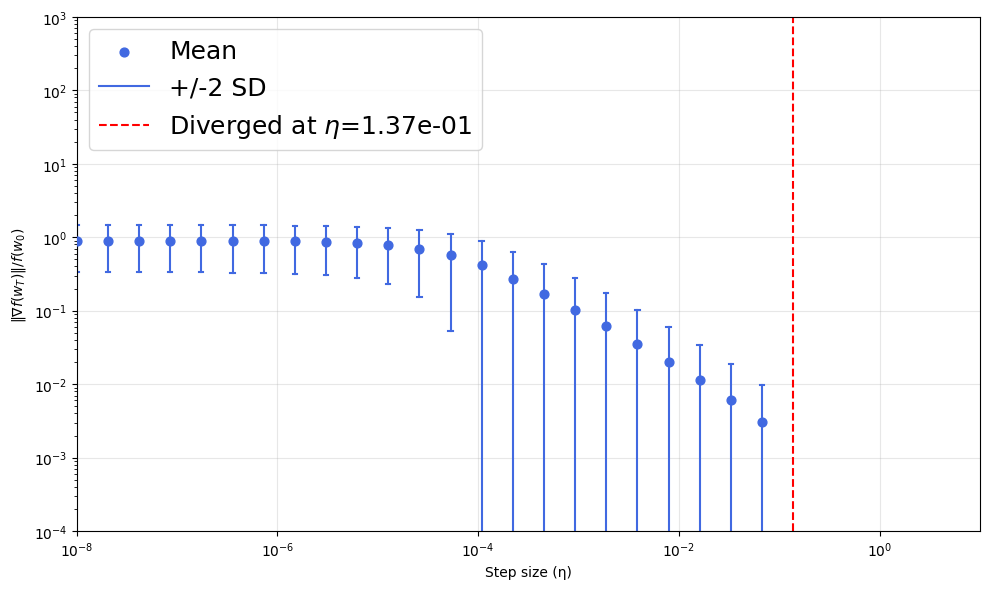}
     \caption{$\pi_j = \cN(\vecOrigin, 5.0\matI_{20})$. The first divergence is at $\eta_i\approx 0.137$.}
     \label{fig:gd_p3_cov5}
 \end{subfigure}
 \hfill
 \begin{subfigure}{0.46\textwidth}
     \includegraphics[width=\textwidth]{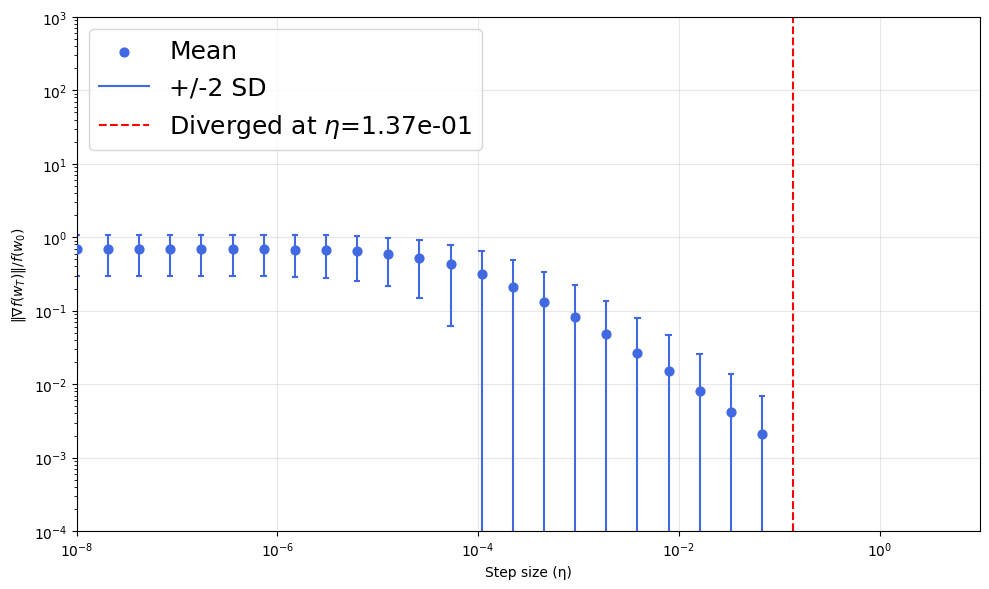}
     \caption{$\pi_j = \cN(\vecOrigin, 7.5\matI_{20})$. The first divergence is at $\eta_i\approx 0.137$.}
     \label{fig:gd_p3_cov7.5}
 \end{subfigure}
 \hfill
 \begin{subfigure}{0.46\textwidth}
     \includegraphics[width=\textwidth]{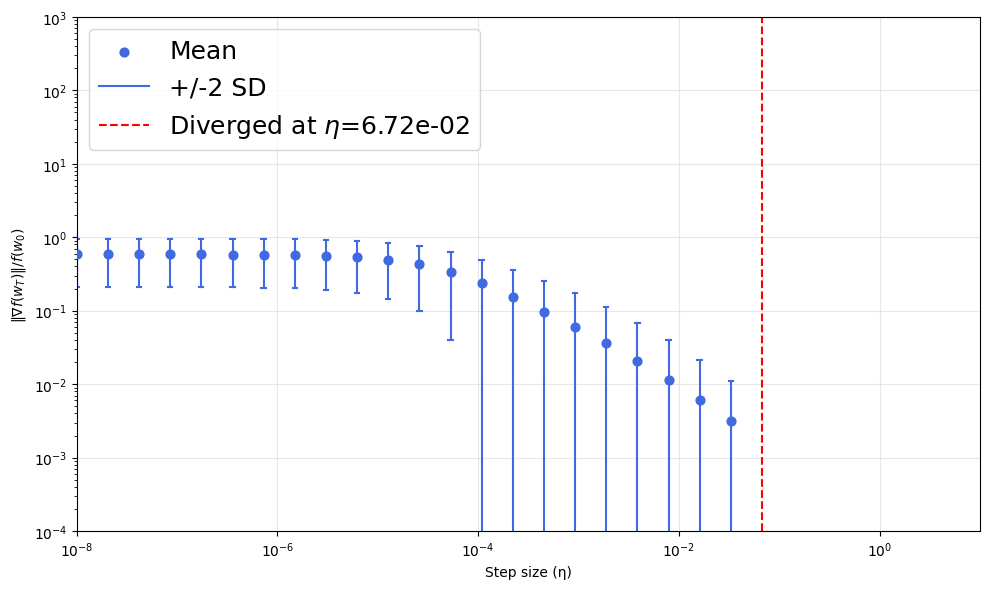}
     \caption{$\pi_j = \cN(\vecOrigin, 10\matI_{20})$. The first divergence is at $\eta_i\approx 0.0672$.}
     \label{fig:gd_p3_cov10}
 \end{subfigure} 
 \caption{GD simulation results for $p=3$. For $\pi_j = \cN(\vecOrigin, 2.5 \matI_{20})$, the first divergence is at $\eta_i\approx 0.281$. For $\pi_j = \cN(\vecOrigin, 5 \matI_{20}), \cN(\vecOrigin, 7.5 \matI_{20})$, the first divergence is at $\eta_i\approx 0.137$. 
 For $\pi_j = \cN(\vecOrigin, 10\matI_{20})$, the first divergence is at $\eta_i\approx 0.0672$.}
 \label{fig:gd_p3}
\end{figure}

\begin{figure}
    \centering
    \begin{subfigure}{0.46\textwidth}
     \includegraphics[width=\textwidth]{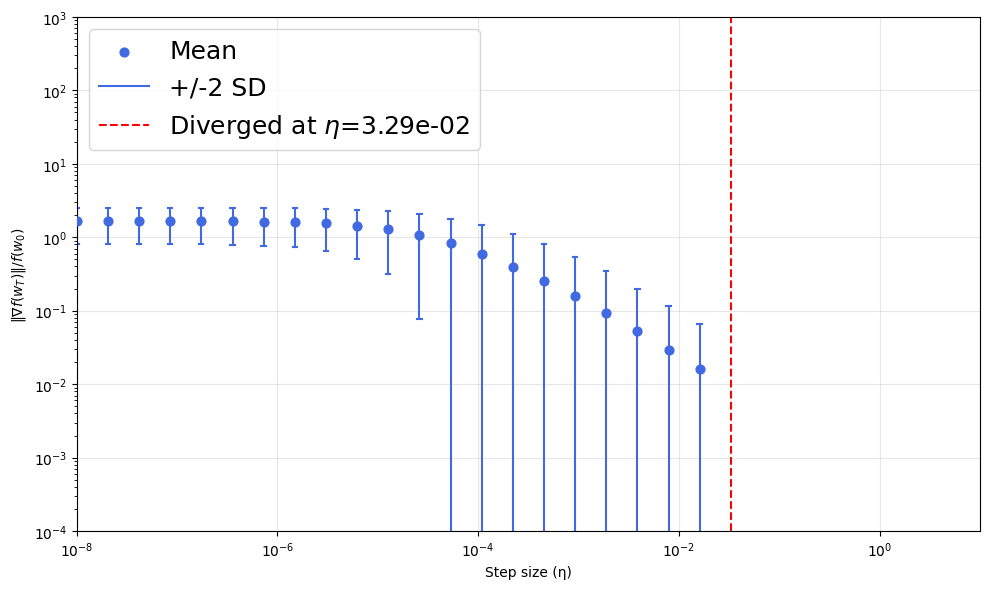}
     \caption{$\pi_j = \cN(\vecOrigin, 2.5\matI_{20})$. The first divergence is at $\eta_i\approx 3.29 \cdot 10^{-2}$.}
     \label{fig:gd_p4_cov2.5}
 \end{subfigure}
 \hfill
 \begin{subfigure}{0.46\textwidth}
     \includegraphics[width=\textwidth]{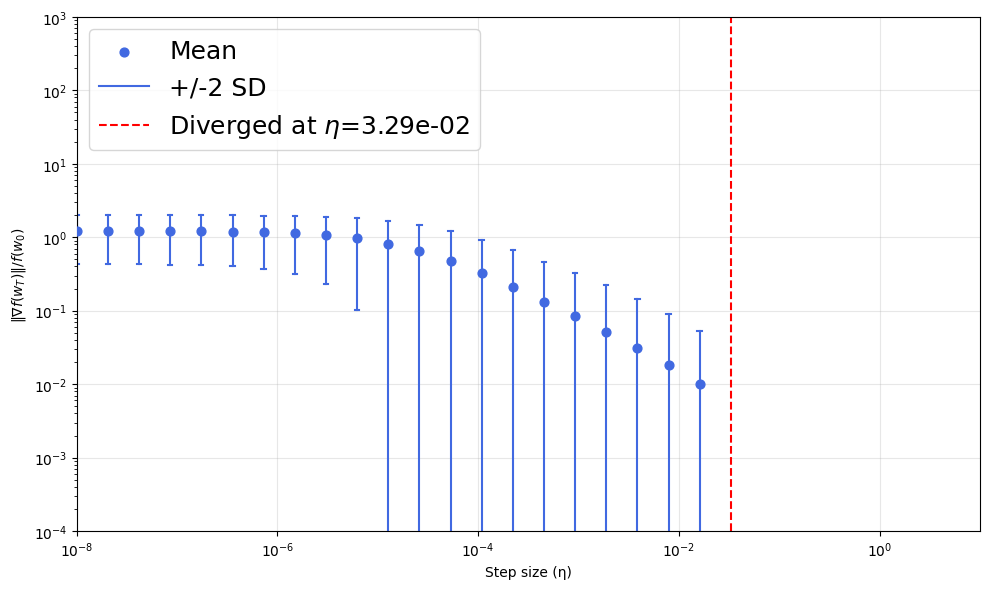}
     \caption{$\pi_j = \cN(\vecOrigin, 5.0\matI_{20})$. The first divergence is at $\eta_i\approx 3.29 \cdot 10^{-2}$.}
     \label{fig:gd_p4_cov5}
 \end{subfigure}
 \hfill
 \begin{subfigure}{0.46\textwidth}
     \includegraphics[width=\textwidth]{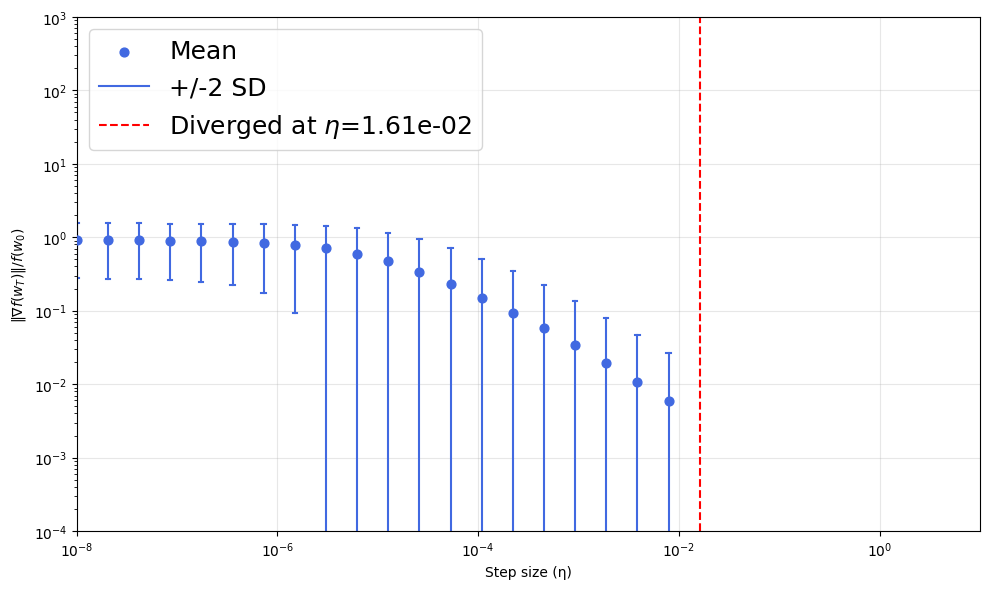}
     \caption{$\pi_j = \cN(\vecOrigin, 7.5\matI_{20})$. The first divergence is at $\eta_i\approx 1.61 \cdot 10^{-2}$.}
     \label{fig:gd_p4_cov7.5}
 \end{subfigure}
 \hfill
 \begin{subfigure}{0.46\textwidth}
     \includegraphics[width=\textwidth]{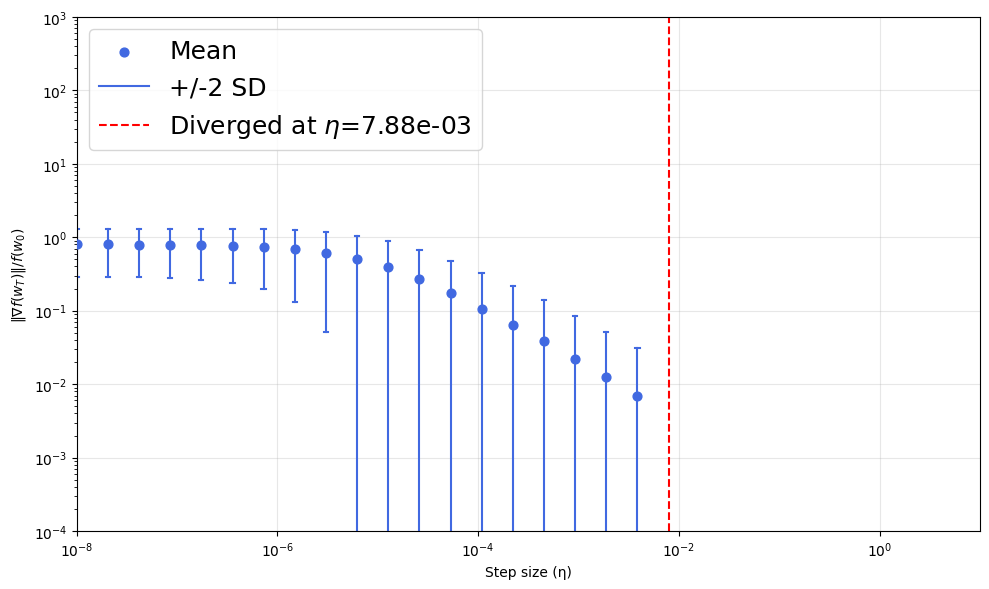}
     \caption{$\pi_j = \cN(\vecOrigin, 10\matI_{20})$. The first divergence is at $\eta_i\approx 7.88 \cdot 10^{-3}$.}
     \label{fig:gd_p4_cov10}
 \end{subfigure} 
 \caption{GD simulation results for $p=4$. For $\pi_j = \cN(\vecOrigin, 2.5 \matI_{20}), \cN(\vecOrigin, 5 \matI_{20})$, the first divergence is at $\eta_i\approx 3.29 \cdot 10^{-2}$. For $\pi_j = \cN(\vecOrigin, 7.5\matI_{20})$, the first divergence is at $\eta_i \approx 1.61 \cdot 10^{-2}$. For $\pi_j = \cN(\vecOrigin, 10\matI_{20})$, the first divergence is at $\eta_i\approx 7.88 \cdot 10^{-3}$.}
 \label{fig:gd_p4}
\end{figure}

\begin{figure}
    \centering
    \begin{subfigure}{0.46\textwidth}
     \includegraphics[width=\textwidth]{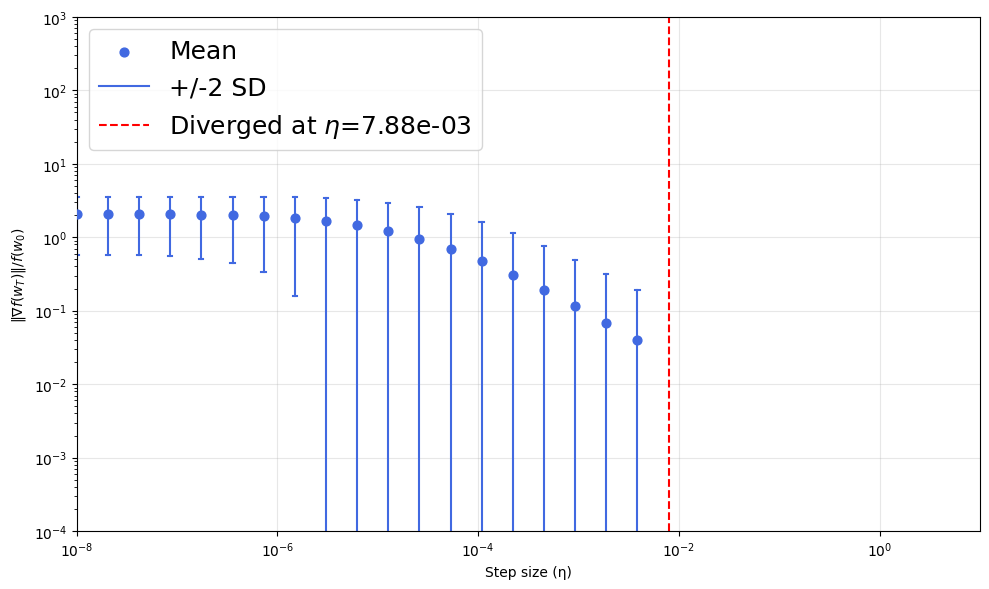}
     \caption{$\pi_j = \cN(\vecOrigin, 2.5\matI_{20})$. The first divergence is at $\eta_i\approx 7.88 \cdot 10^{-3}$.}
     \label{fig:gd_p5_cov2.5}
 \end{subfigure}
 \hfill
 \begin{subfigure}{0.46\textwidth}
     \includegraphics[width=\textwidth]{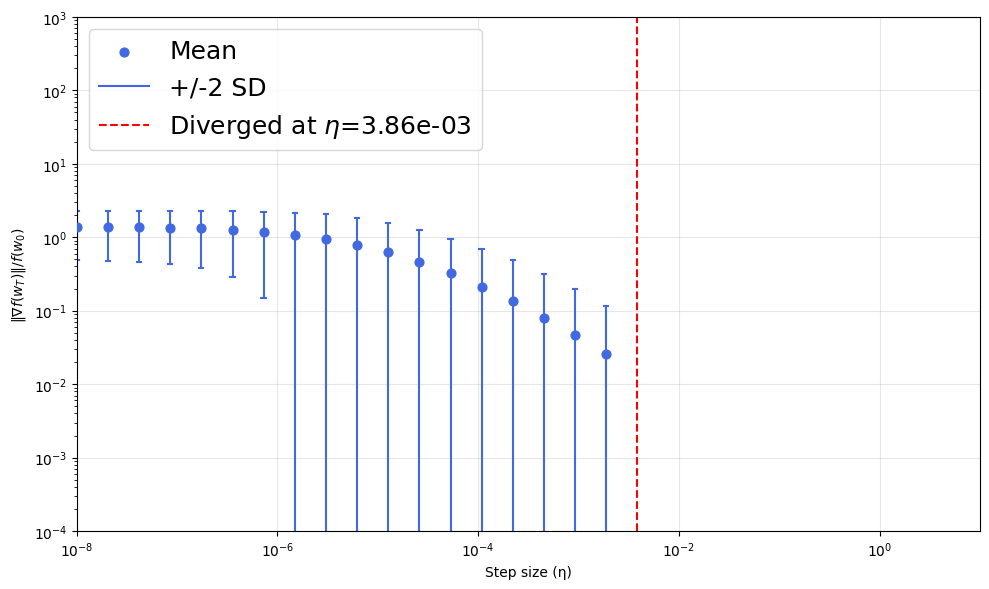}
     \caption{$\pi_j = \cN(\vecOrigin, 5.0\matI_{20})$. The first divergence is at $\eta_i\approx 3.86 \cdot 10^{-3}$.}
     \label{fig:gd_p5_cov5}
 \end{subfigure}
 \hfill
 \begin{subfigure}{0.46\textwidth}
     \includegraphics[width=\textwidth]{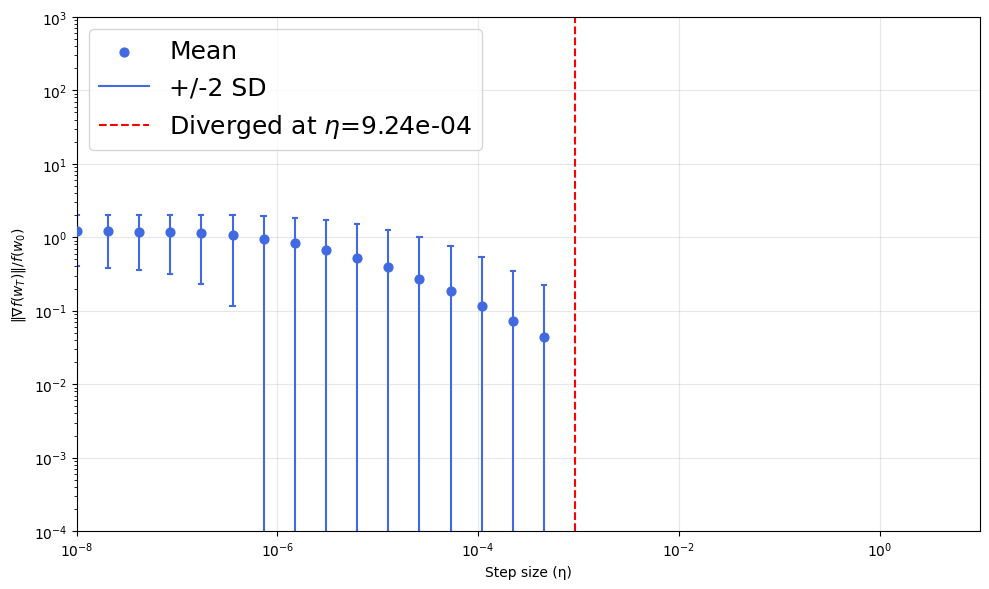}
     \caption{$\pi_j = \cN(\vecOrigin, 7.5\matI_{20})$. The first divergence is at $\eta_i\approx 9.24 \cdot 10^{-4}$.}
     \label{fig:gd_p5_cov7.5}
 \end{subfigure}
 \hfill
 \begin{subfigure}{0.46\textwidth}
     \includegraphics[width=\textwidth]{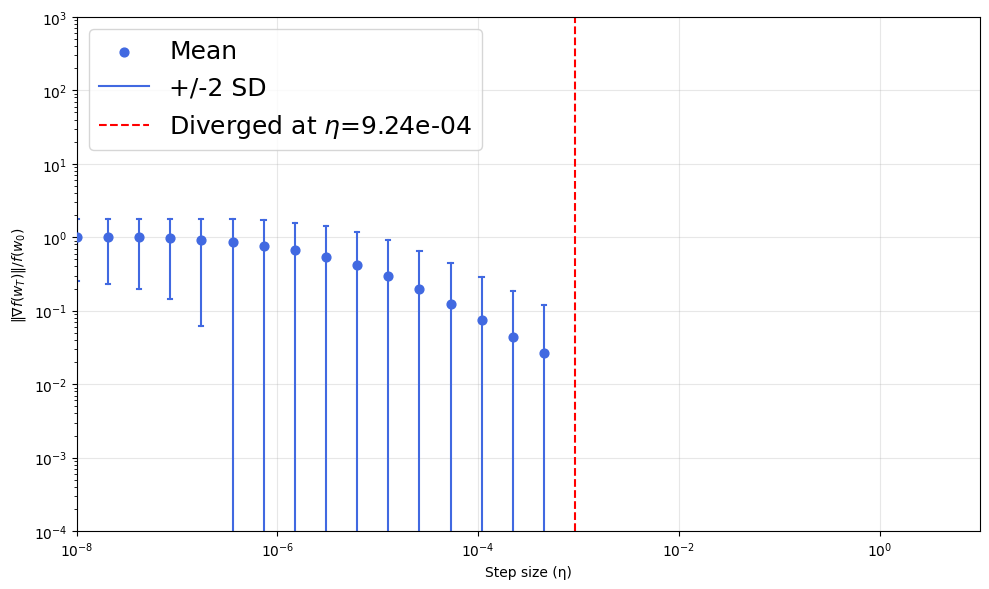}
     \caption{$\pi_j = \cN(\vecOrigin, 10\matI_{20})$. The first divergence is at $\eta_i\approx 9.24 \cdot 10^{-4}$.}
     \label{fig:gd_p5_cov10}
 \end{subfigure} 
 \caption{GD simulation results for $p=5$. For $\pi_j = \cN(\vecOrigin, 2.5 \matI_{20})$, the first divergence is at $\eta_i\approx 7.88 \cdot 10^{-3}$. For $\pi_j = \cN(\vecOrigin, 5 \matI_{20})$, the first divergence is at $\eta_i\approx 3.86 \cdot 10^{-3}$. For $\pi_j = \cN(\vecOrigin, 7.5 \matI_{20}), \cN(\vecOrigin, 10\matI_{20})$, the first divergence is at $\eta_i\approx 9.24 \cdot 10^{-4}$.}
 \label{fig:gd_p5}
\end{figure}

\begin{figure}
    \centering
    \begin{subfigure}{0.46\textwidth}
     \includegraphics[width=\textwidth]{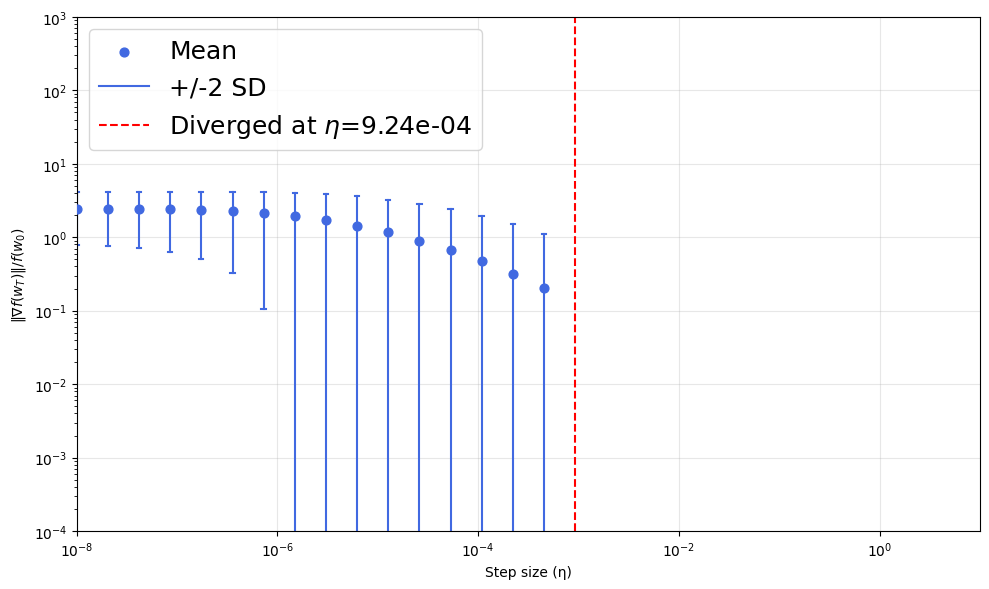}
     \caption{$\pi_j = \cN(\vecOrigin, 2.5\matI_{20})$. The first divergence is at $\eta_i\approx 9.24 \cdot 10^{-4}$.}
     \label{fig:gd_p6_cov2.5}
 \end{subfigure}
 \hfill
 \begin{subfigure}{0.46\textwidth}
     \includegraphics[width=\textwidth]{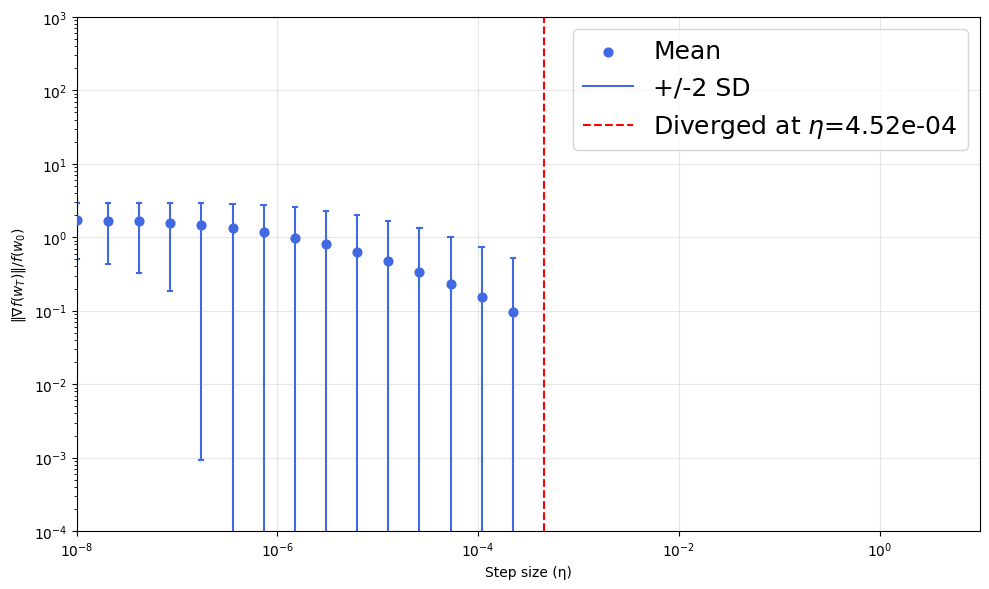}
     \caption{$\pi_j = \cN(\vecOrigin, 5.0\matI_{20})$. The first divergence is at $\eta_i\approx 4.52 \cdot 10^{-4}$.}
     \label{fig:gd_p6_cov5}
 \end{subfigure}
 \hfill
 \begin{subfigure}{0.46\textwidth}
     \includegraphics[width=\textwidth]{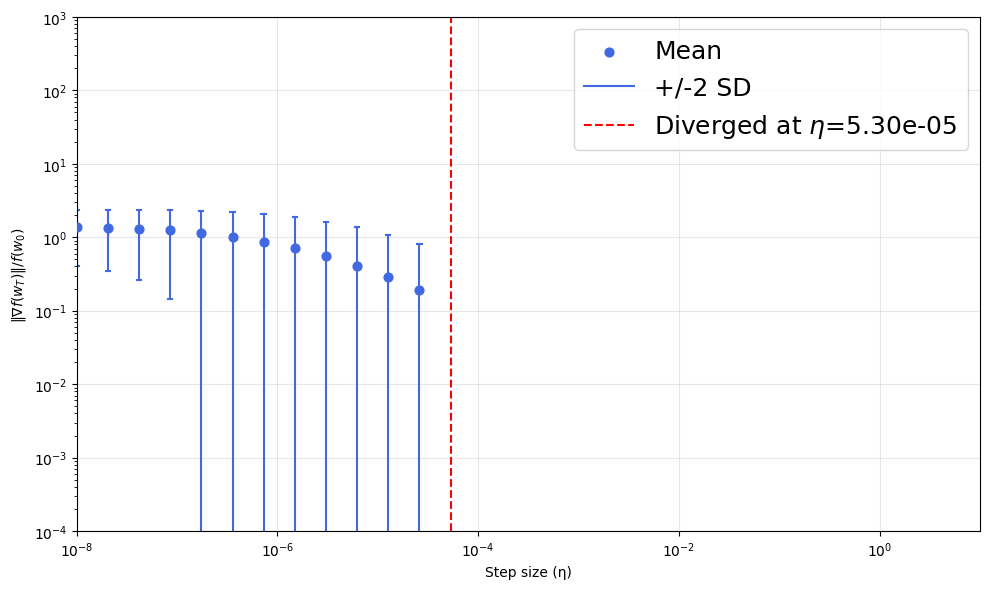}
     \caption{$\pi_j = \cN(\vecOrigin, 7.5\matI_{20})$. The first divergence is at $\eta_i\approx 5.30\cdot 10^{-5}$.}
     \label{fig:gd_p6_cov7.5}
 \end{subfigure}
 \hfill
 \begin{subfigure}{0.46\textwidth}
     \includegraphics[width=\textwidth]{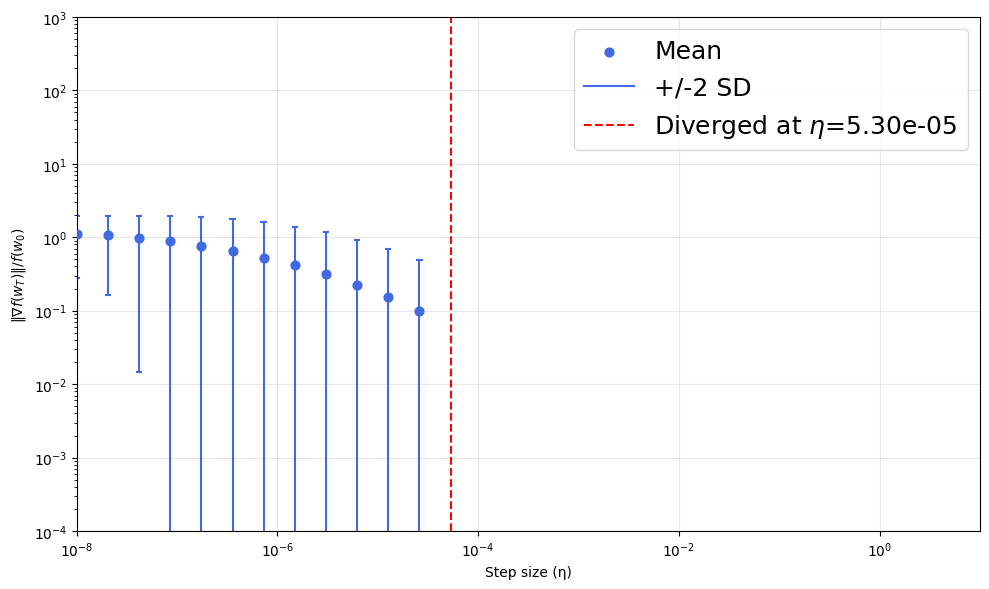}
     \caption{$\pi_j = \cN(\vecOrigin, 10\matI_{20})$. The first divergence is at $\eta_i\approx 5.30 \cdot 10^{-5}$.}
     \label{fig:gd_p6_cov10}
 \end{subfigure} 
 \caption{GD simulation results for $p=6$. For $\pi_j = \cN(\vecOrigin, 2.5 \matI_{20})$, the first divergence is at $\eta_i\approx 9.24 \cdot 10^{-4}$. For $\pi_j = \cN(\vecOrigin, 5 \matI_{20})$, the first divergence is at $\eta_i\approx 4.52 \cdot 10^{-4}$. For $\pi_j = \cN(\vecOrigin, 7.5\matI_{20}), \cN(\vecOrigin, 10\matI_{20})$, the first divergence is at $\eta_i\approx 5.30 \cdot 10^{-5}$.}
 \label{fig:gd_p6}
\end{figure}

\paragraph{Results:} Our simulations validate this theory very accurately. Note in the following figures that the $y$-axis is normalized, as we plot $\tfrac{\nrm*{\grad F(\vecW_T)}}{F(\vecW_0)}$ where $T=1000$. Thus larger $c_j$ lead to comparable values on the $y$-axis.
\begin{itemize}
\item When $p=2$: In \pref{fig:gd_p2}, we plot the results in the manner described above for all 4 initializations $\pi_j$. As is predicted by established optimization theory for smooth functions \citep{bubeck2015convex}, the first step size leading to divergence $\eta_i$ is identical across all the $\pi_j$.

\item When $p=3, 4, 5, 6$: We plot the results in the manner described above for all 4 initializations $\pi_j$ in \pref{fig:gd_p3}, \pref{fig:gd_p4}, \pref{fig:gd_p5}, \pref{fig:gd_p6} respectively. Unlike the $p=2$ case, in all of these cases, the first step size leading to divergence $\eta_i$ generally decreases as the covariance $c_j \matI_{20}$ of $\pi_j$ increases from $2.5$ to 10. 
\end{itemize}
We also notice the following, both in line with our theoretical claims:
\begin{itemize}
    \item For a given $p$, consider how this first step size $\eta_i$ leading to divergence decreases as the covariance $c_j \matI_{20}$ of $\pi_j$ increases from $2.5$ to 10. We find that the rate of this decrease increases as $p$ increases. The ratio of the first $\eta_i$ leading to divergence for $\pi_1$ vs $\pi_4$ is approximately $4.18, 4.18, 8.53, 17.43$ for $p=3,4,5,6$ respectively.

    As remarked earlier, for larger $p$, the corresponding self-bounding function $\rho_1(\cdot)$ is larger for $F(\vecW)=\nrm*{\matA\vecW}^p$ (see \pref{subsec:detailsforcomparewithlit} for a similar verification). Thus this behavior is consistent with our results, as the step size from all of our results depends on $F(\vecW_0)$ through $\rho_1(\cdot)$.

    \item Fixing $\pi_j$ and comparing across $p$, we see that the first step size leading to divergence $\eta_i$ decreases as $p$ increases. Again this is not a surprise considering our theoretical results, as for larger $p$, both $F(\vecW_0)$ for $\vecW_0 \sim \pi_j$ and the self-bounding function $\rho_1(\cdot)$ become larger.
\end{itemize}
For each $p \in \{2,3,4,5,6\}$ and $\pi_j$, we also record the smallest $\eta_i$ for which divergence occurred in \pref{tab:gdstepsizes}, which highlights the aforementioned trends. 

\begin{table}[h!]
\centering
\begin{tabular}{lcccc}
\toprule
 & $\pi_j=\cN(\vecOrigin, 2.5\matI_{20})$ & $\pi_j=\cN(\vecOrigin, 5.0\matI_{20})$ & $\pi_j=\cN(\vecOrigin, 7.5\matI_{20})$ & $\pi_j=\cN(\vecOrigin, 10\matI_{20})$ \\
\midrule
$p=2$ & $1.17 \cdot 10^0$ & $1.17 \cdot 10^0$  & $1.17 \cdot 10^0$  & $1.17  \cdot 10^0$  \\
$p=3$ & $2.81 \cdot 10^{-1}$ & $1.37 \cdot 10^{-1}$ & $1.37 \cdot 10^{-1}$ & $6.72 \cdot 10^{-2}$ \\
$p=4$ & $3.29 \cdot 10^{-2}$ & $3.29 \cdot 10^{-2}$ & $1.61 \cdot 10^{-2}$ &  $7.88 \cdot 10^{-3}$\\
$p=5$ & $7.88 \cdot 10^{-3}$ & $3.86 \cdot 10^{-3}$ & $9.24 \cdot 10^{-4}$ & $9.24 \cdot 10^{-4}$ \\
$p=6$ &  $9.24 \cdot 10^{-4}$ & $4.52 \cdot 10^{-4}$ & $5.30 \cdot 10^{-5}$ & $5.30 \cdot 10^{-5}$
\end{tabular}
\caption{The smallest $\eta_i$ leading to divergence for a given $p$ and initialization $\pi_j$.}
\label{tab:gdstepsizes} 
\end{table}

\subsection{Synthetic Simulations with SGD}\label{subsec:experimentstoysgd}
\paragraph{Simulation Details:} We adopt the exact same settings as in \pref{subsec:experimentstoygd}. The only difference is that we study SGD rather than GD, and hence we simulate stochastic gradients. We do so similarly to \citet{gaash2025convergence}: we artificially add $\cN\prn*{\vecOrigin, 0.01 \matI_{20}}$ to $\grad F$ at each iteration of SGD.\footnote{Note our result for convergence of SGD to FOSPs, \pref{thm:sgdfirstorder}, applies for Gaussian noise as per \pref{rem:sgdfirstordersubgaussian}.} The simulations for \pref{subsec:experimentstoysgd} were again run on a Jupyter notebook in Python in Google Colab Pro, connected to a single NVIDIA T4 GPU. Our code is in the attached files.

\begin{figure}
    \centering
    \begin{subfigure}{0.46\textwidth}
     \includegraphics[width=\textwidth]{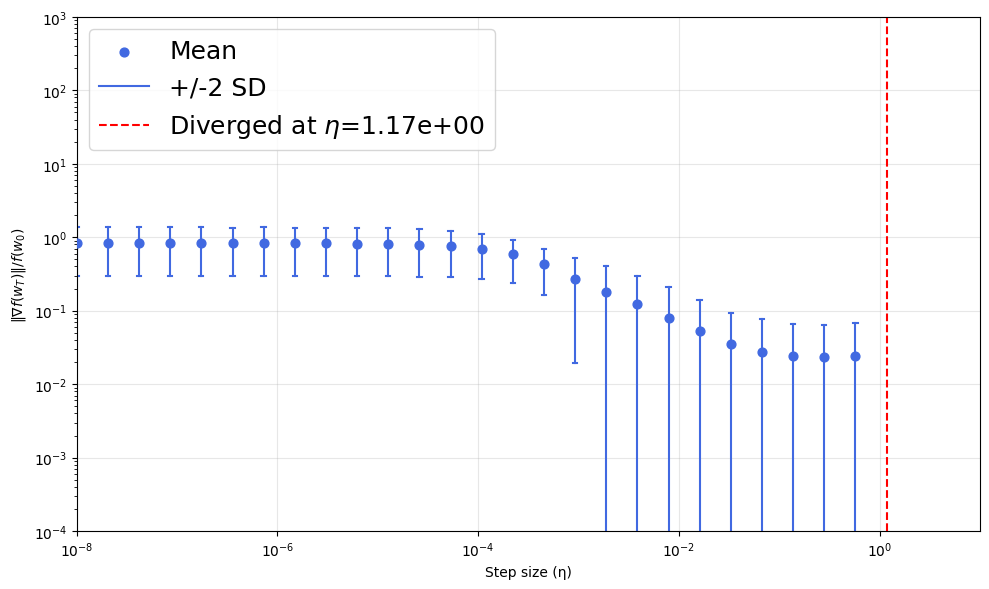}
     \caption{$\pi_j = \cN(\vecOrigin, 2.5\matI_{20})$. The first divergence is at $\eta_i\approx 1.17$.}
     \label{fig:sgd_p2_cov2.5_naive}
 \end{subfigure}
 \hfill
 \begin{subfigure}{0.46\textwidth}
     \includegraphics[width=\textwidth]{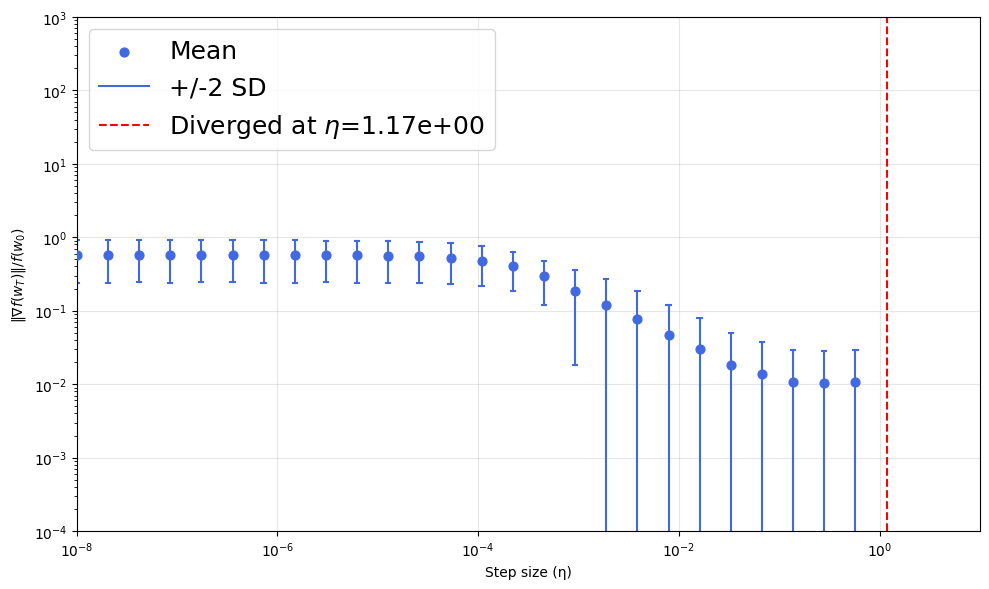}
     \caption{$\pi_j = \cN(\vecOrigin, 5.0\matI_{20})$. The first divergence is at $\eta_i\approx 1.17$.}
     \label{fig:sgd_p2_cov5_naive}
 \end{subfigure}
 \hfill
 \begin{subfigure}{0.46\textwidth}
     \includegraphics[width=\textwidth]{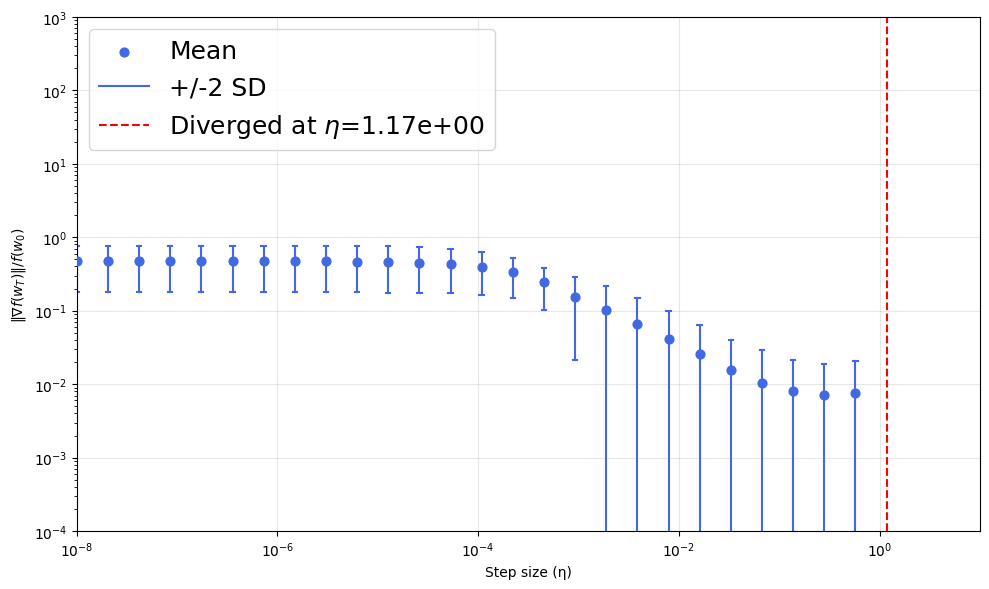}
     \caption{$\pi_j = \cN(\vecOrigin, 7.5\matI_{20})$. The first divergence is at $\eta_i\approx 1.17$.}
     \label{fig:sgd_p2_cov7.5_naive}
 \end{subfigure}
 \hfill
 \begin{subfigure}{0.46\textwidth}
     \includegraphics[width=\textwidth]{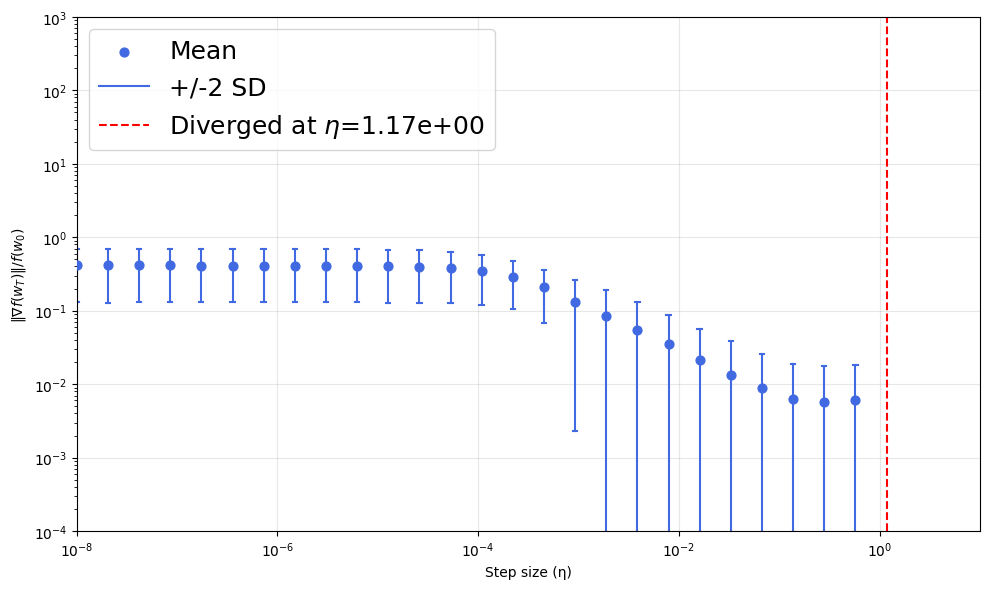}
     \caption{$\pi_j = \cN(\vecOrigin, 10\matI_{20})$. The first divergence is at $\eta_i\approx 1.17$.}
     \label{fig:sgd_p2_cov10_naive}
 \end{subfigure} 
 \caption{SGD simulation results for $p=2$. For all $\pi_j$, the smallest $\eta_i$ leading to divergence is $\approx 1.17$.}
 \label{fig:sgd_p2_naive}
\end{figure}

\begin{figure}
    \centering
    \begin{subfigure}{0.46\textwidth}
     \includegraphics[width=\textwidth]{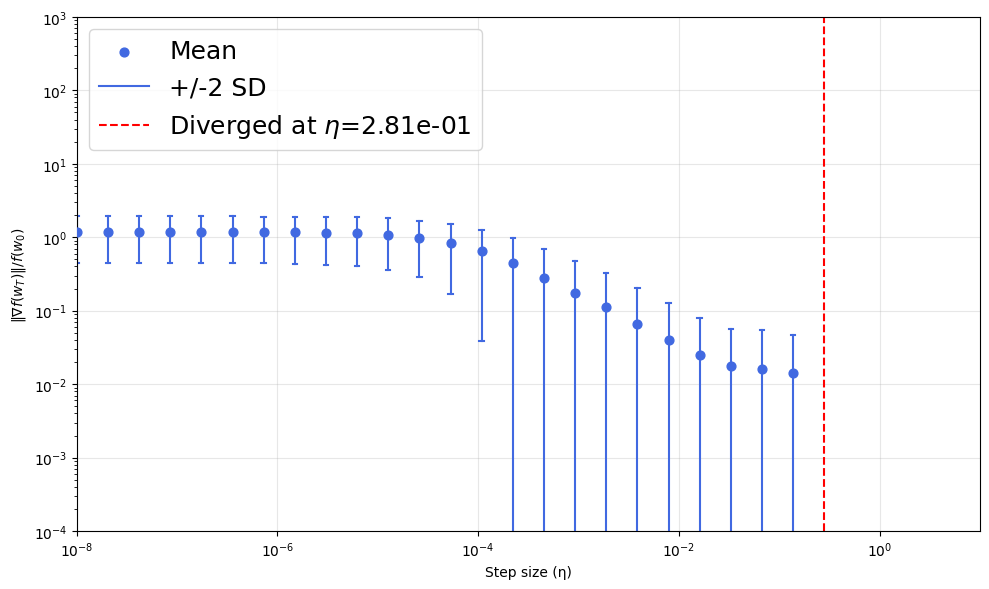}
     \caption{$\pi_j = \cN(\vecOrigin, 2.5\matI_{20})$. The first divergence is at $\eta_i\approx 0.281$.}
     \label{fig:sgd_p3_cov2.5_naive}
 \end{subfigure}
 \hfill
 \begin{subfigure}{0.46\textwidth}
     \includegraphics[width=\textwidth]{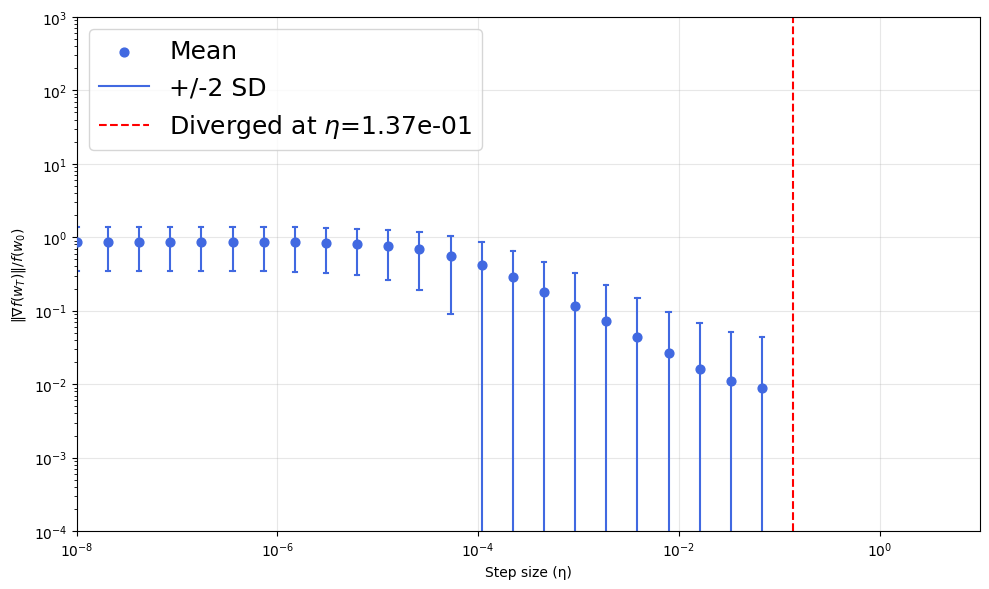}
     \caption{$\pi_j = \cN(\vecOrigin, 5.0\matI_{20})$. The first divergence is at $\eta_i\approx 0.137$.}
     \label{fig:sgd_p3_cov5_naive}
 \end{subfigure}
 \hfill
 \begin{subfigure}{0.46\textwidth}
     \includegraphics[width=\textwidth]{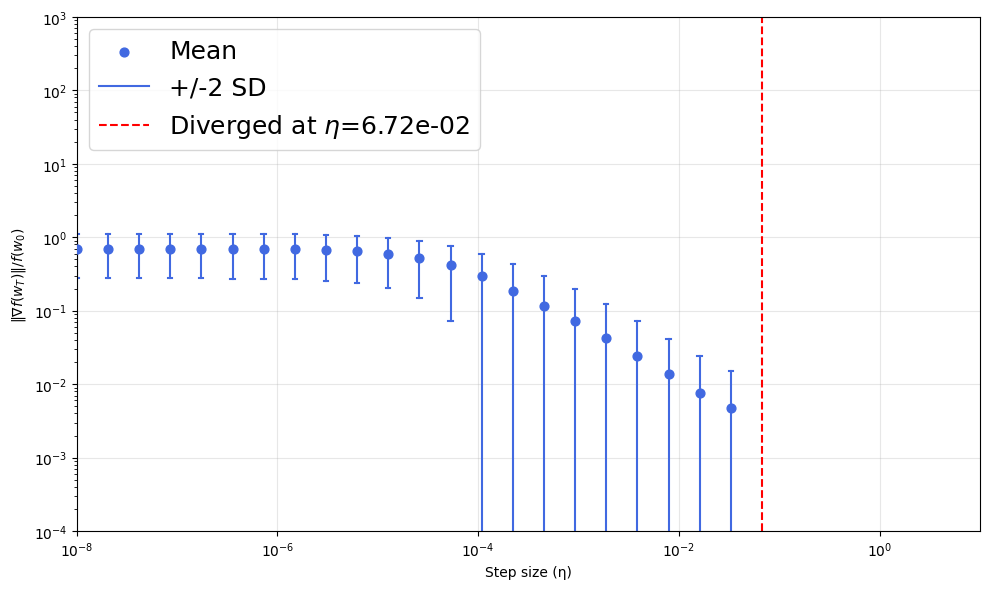}
     \caption{$\pi_j = \cN(\vecOrigin, 7.5\matI_{20})$. The first divergence is at $\eta_i\approx 6.72 \cdot 10^{-2}$.}
     \label{fig:sgd_p3_cov7.5_naive}
 \end{subfigure}
 \hfill
 \begin{subfigure}{0.46\textwidth}
     \includegraphics[width=\textwidth]{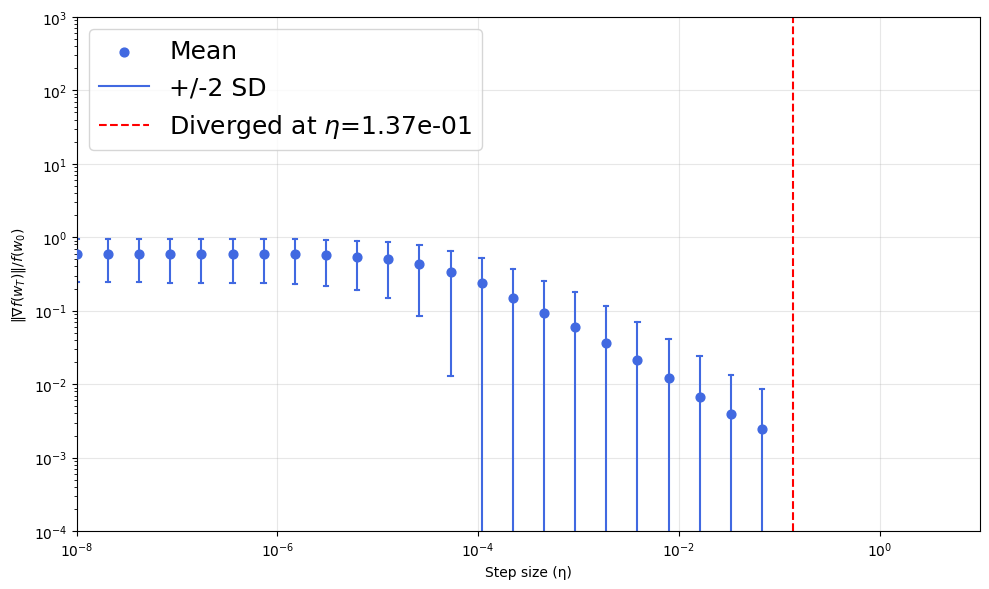}
     \caption{$\pi_j = \cN(\vecOrigin, 10\matI_{20})$. The first divergence is at $\eta_i\approx 0.137$.}
     \label{fig:sgd_p3_cov10_naive}
 \end{subfigure} 
 \caption{SGD simulation results for $p=3$. For $\pi_j = \cN(\vecOrigin, 2.5 \matI_{20})$, the first divergence is at $\eta_i\approx 0.281$. For $\pi_j = \cN(\vecOrigin, 7.5 \matI_{20})$, the first divergence is at $\eta_i\approx 6.72 \cdot 10^{-2}$. 
 For the other $\pi_j$, the first divergence is at $\eta_i\approx 0.137$.}
 \label{fig:sgd_p3_naive}
\end{figure}

\begin{figure}
    \centering
    \begin{subfigure}{0.46\textwidth}
     \includegraphics[width=\textwidth]{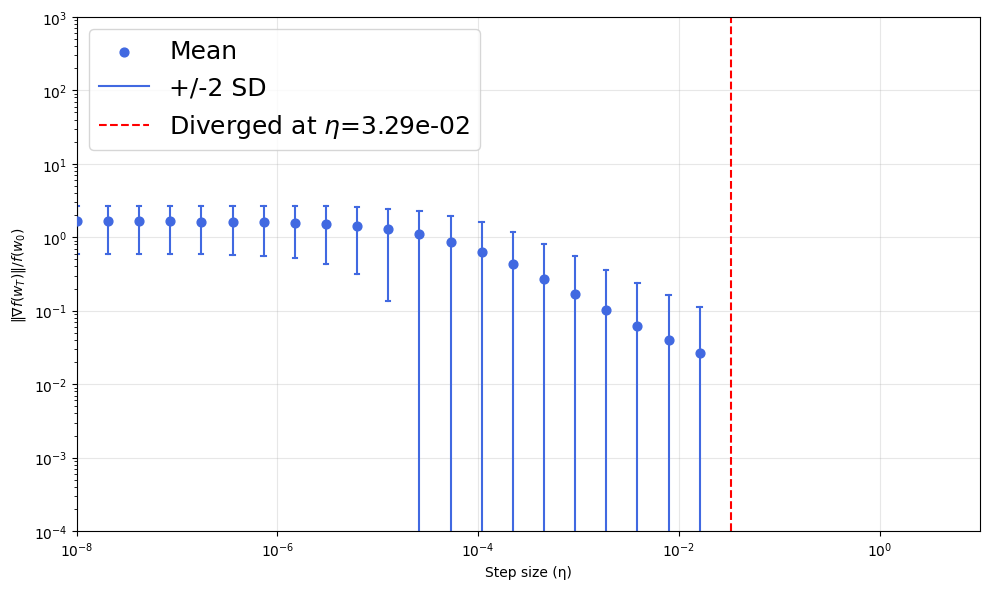}
     \caption{$\pi_j = \cN(\vecOrigin, 2.5\matI_{20})$. The first divergence is at $\eta_i\approx 3.29 \cdot 10^{-2}$.}
     \label{fig:sgd_p4_cov2.5_naive}
 \end{subfigure}
 \hfill
 \begin{subfigure}{0.46\textwidth}
     \includegraphics[width=\textwidth]{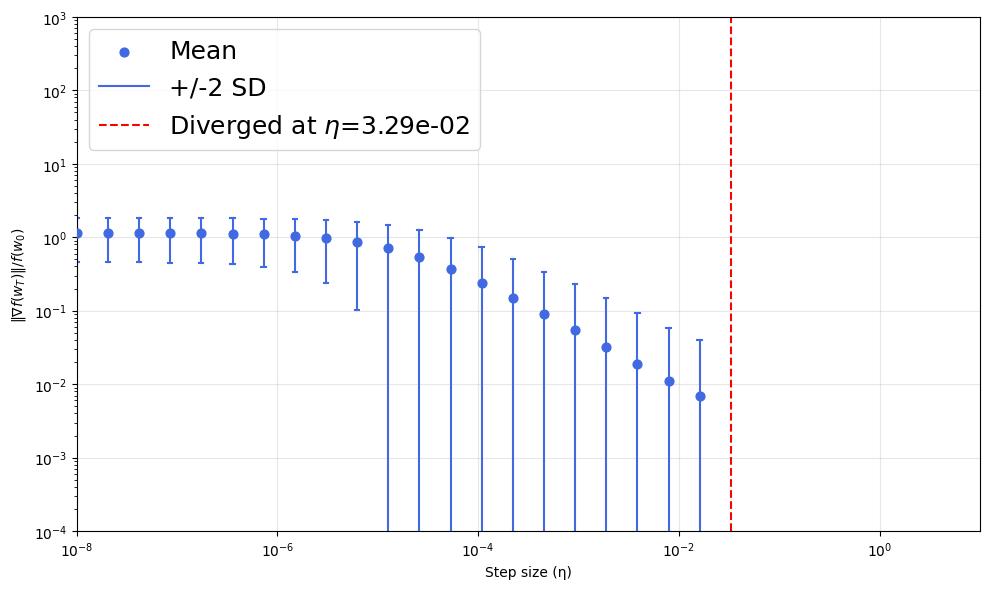}
     \caption{$\pi_j = \cN(\vecOrigin, 5.0\matI_{20})$. The first divergence is at $\eta_i\approx 3.29 \cdot 10^{-2}$.}
     \label{fig:sgd_p4_cov5_naive}
 \end{subfigure}
 \hfill
 \begin{subfigure}{0.46\textwidth}
     \includegraphics[width=\textwidth]{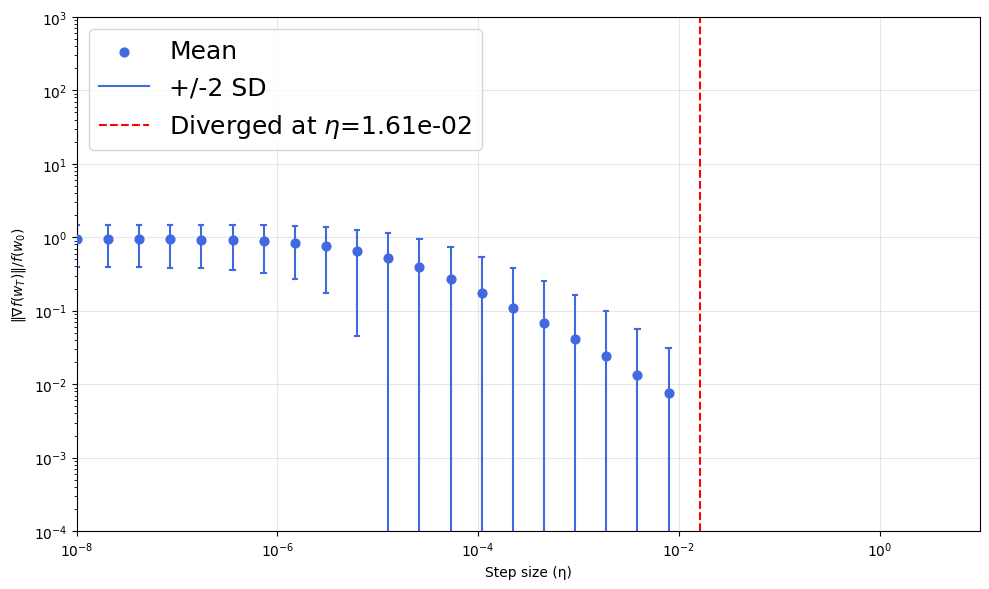}
     \caption{$\pi_j = \cN(\vecOrigin, 7.5\matI_{20})$. The first divergence is at $\eta_i\approx 1.61 \cdot 10^{-2}$.}
     \label{fig:sgd_p4_cov7.5_naive}
 \end{subfigure}
 \hfill
 \begin{subfigure}{0.46\textwidth}
     \includegraphics[width=\textwidth]{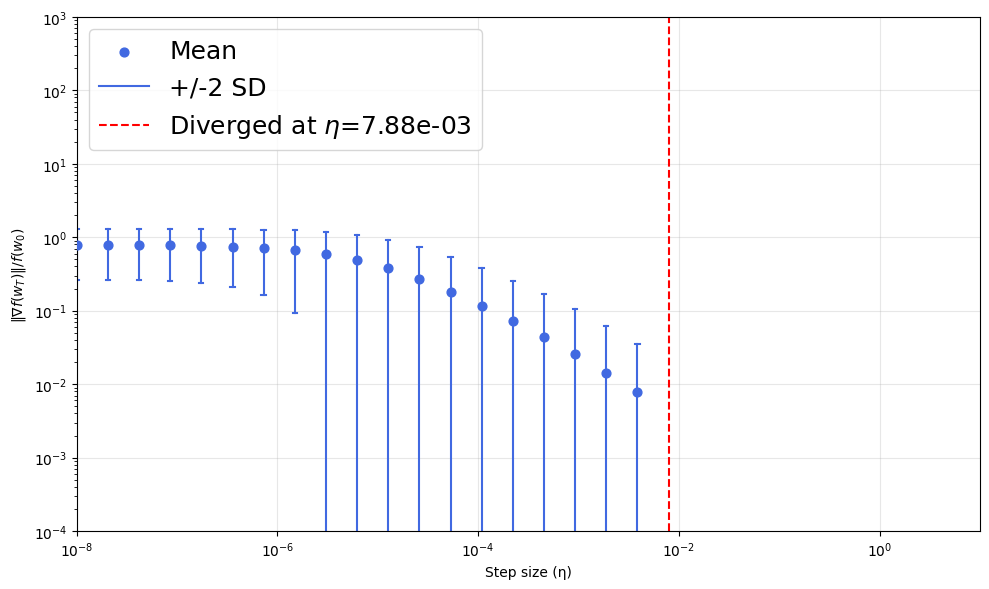}
     \caption{$\pi_j = \cN(\vecOrigin, 10\matI_{20})$. The first divergence is at $\eta_i\approx 7.88 \cdot 10^{-3}$.}
     \label{fig:sgd_p4_cov10_naive}
 \end{subfigure} 
 \caption{SGD simulation results for $p=4$. For $\pi_j = \cN(\vecOrigin, 2.5 \matI_{20}), \cN(\vecOrigin, 5.0 \matI_{20})$, the first divergence is at $\eta_i\approx 3.29 \cdot 10^{-2}$. For $\pi_j = \cN(\vecOrigin, 7.5 \matI_{20})$, the first divergence is at $\eta_i \approx 1.61 \cdot 10^{-2}$. For $\pi_j\sim \cN(\vecOrigin, 10 \matI_{20})$, the first divergence is at $\eta_i\approx 7.88 \cdot 10^{-3}$.}
 \label{fig:sgd_p4_naive}
\end{figure}

\begin{figure}
    \centering
    \begin{subfigure}{0.46\textwidth}
     \includegraphics[width=\textwidth]{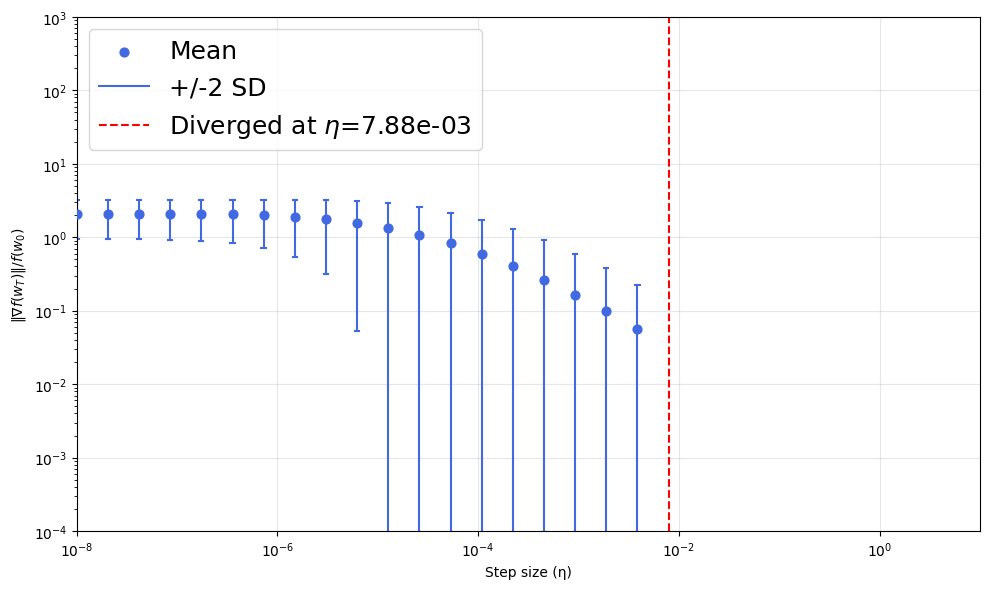}
     \caption{$\pi_j = \cN(\vecOrigin, 2.5\matI_{20})$. The first divergence is at $\eta_i\approx 7.88 \cdot 10^{-3}$.}
     \label{fig:sgd_p5_cov2.5_naive}
 \end{subfigure}
 \hfill
 \begin{subfigure}{0.46\textwidth}
     \includegraphics[width=\textwidth]{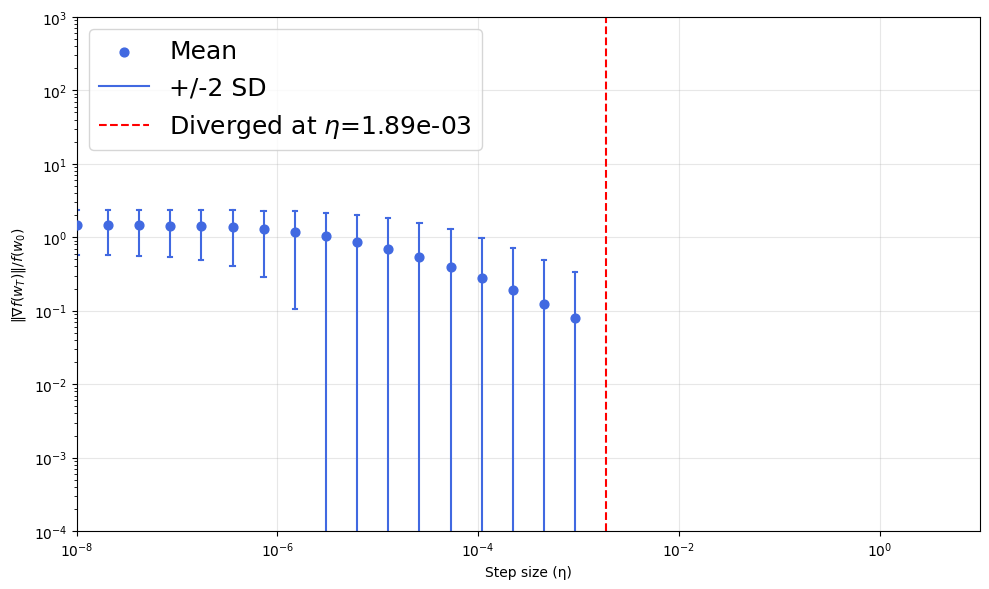}
     \caption{$\pi_j = \cN(\vecOrigin, 5.0\matI_{20})$. The first divergence is at $\eta_i\approx 1.89 \cdot 10^{-3}$.}
     \label{fig:sgd_p5_cov5_naive}
 \end{subfigure}
 \hfill
 \begin{subfigure}{0.46\textwidth}
     \includegraphics[width=\textwidth]{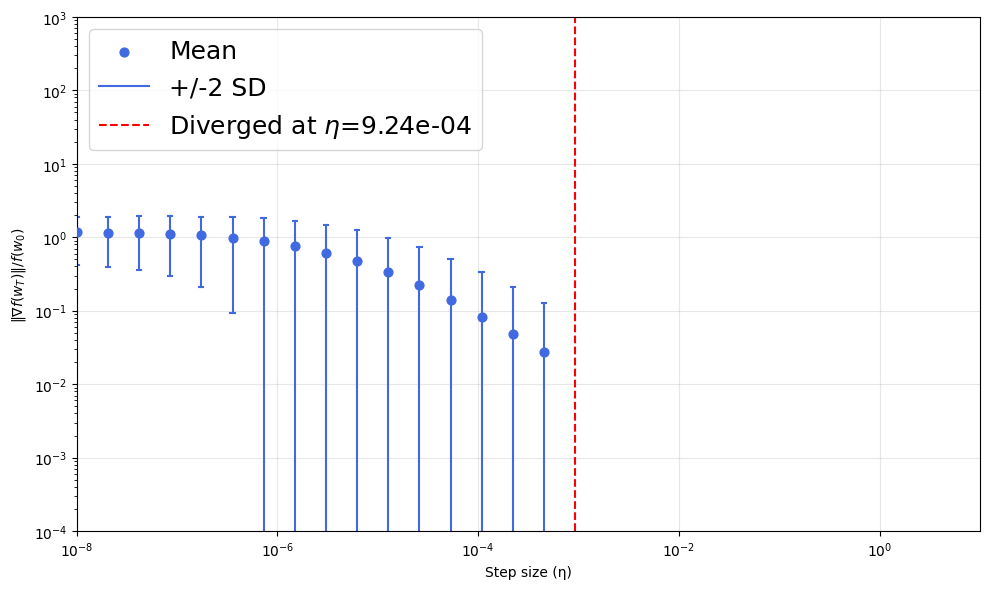}
     \caption{$\pi_j = \cN(\vecOrigin, 7.5\matI_{20})$. The first divergence is at $\eta_i\approx 9.24 \cdot 10^{-4}$.}
     \label{fig:sgd_p5_cov7.5_naive}
 \end{subfigure}
 \hfill
 \begin{subfigure}{0.46\textwidth}
     \includegraphics[width=\textwidth]{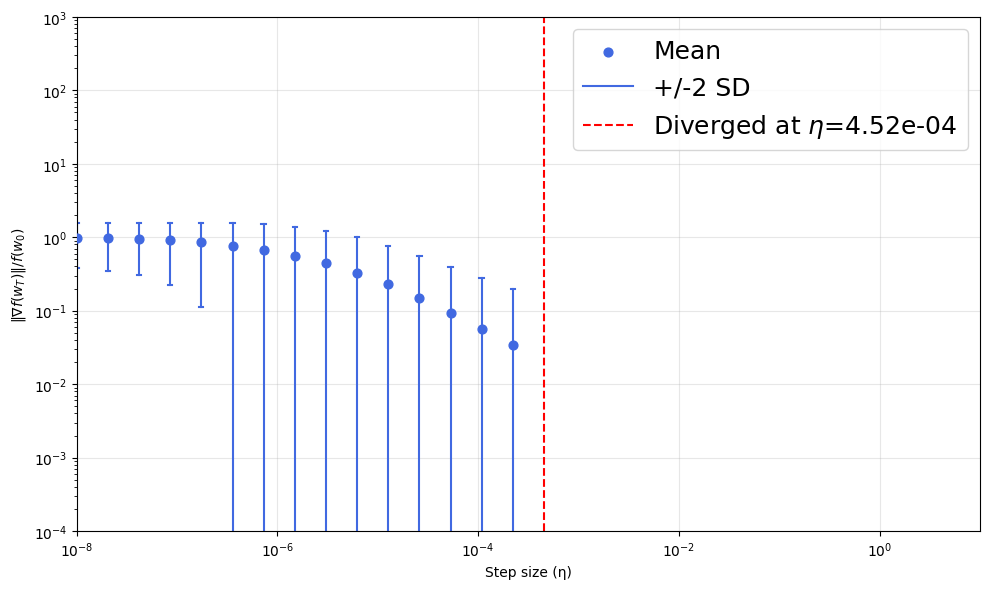}
     \caption{$\pi_j = \cN(\vecOrigin, 10\matI_{20})$. The first divergence is at $\eta_i\approx 4.52 \cdot 10^{-4}$.}
     \label{fig:sgd_p5_cov10_naive}
 \end{subfigure} 
 \caption{SGD simulation results for $p=5$. 
 For $\pi_j = \cN(\vecOrigin, 2.5 \matI_{20}), \cN(\vecOrigin, 5.0 \matI_{20}), \cN(\vecOrigin, 7.5 \matI_{20}), \cN(\vecOrigin, 10 \matI_{20})$, the first divergence is at $\eta_i\approx 7.88 \cdot 10^{-3}, 1.89 \cdot 10^{-3}, 9.24 \cdot 10^{-4}, 4.52 \cdot 10^{-4}$ respectively.}
 \label{fig:sgd_p5_naive}
\end{figure}

\begin{figure}
    \centering
    \begin{subfigure}{0.46\textwidth}
     \includegraphics[width=\textwidth]{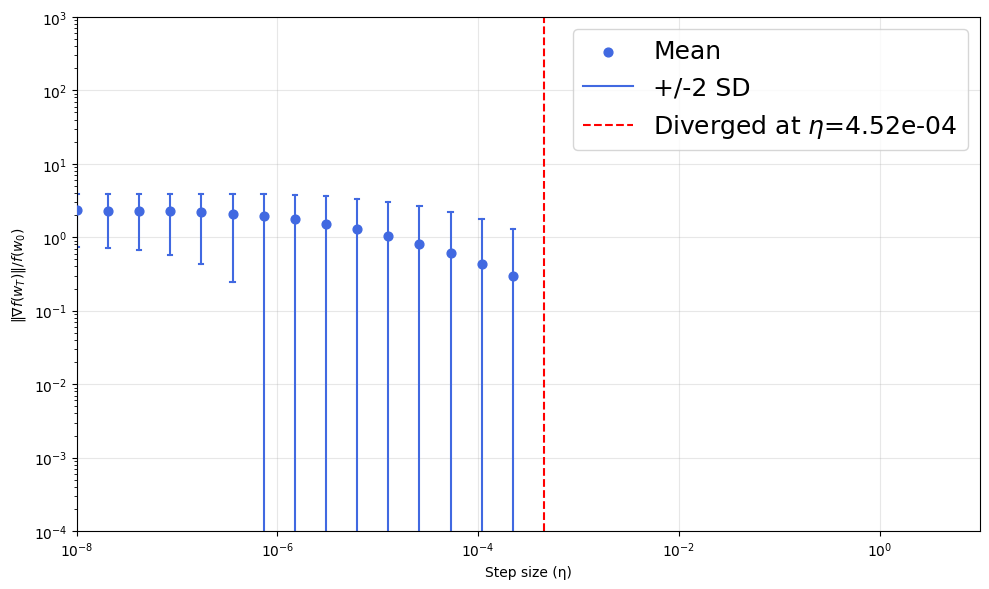}
     \caption{$\pi_j = \cN(\vecOrigin, 2.5\matI_{20})$. The first divergence is at $\eta_i\approx 4.52 \cdot 10^{-4}$.}
     \label{fig:sgd_p6_cov2.5_naive}
 \end{subfigure}
 \hfill
 \begin{subfigure}{0.46\textwidth}
     \includegraphics[width=\textwidth]{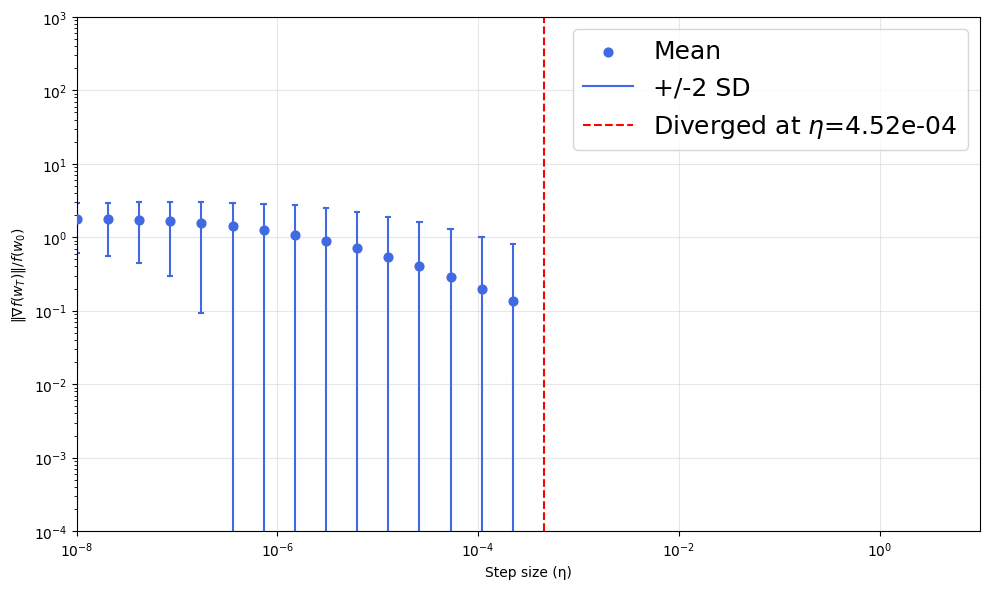}
     \caption{$\pi_j = \cN(\vecOrigin, 5.0\matI_{20})$. The first divergence is at $\eta_i\approx 4.52 \cdot 10^{-4}$.}
     \label{fig:sgd_p6_cov5_naive}
 \end{subfigure}
 \hfill
 \begin{subfigure}{0.46\textwidth}
     \includegraphics[width=\textwidth]{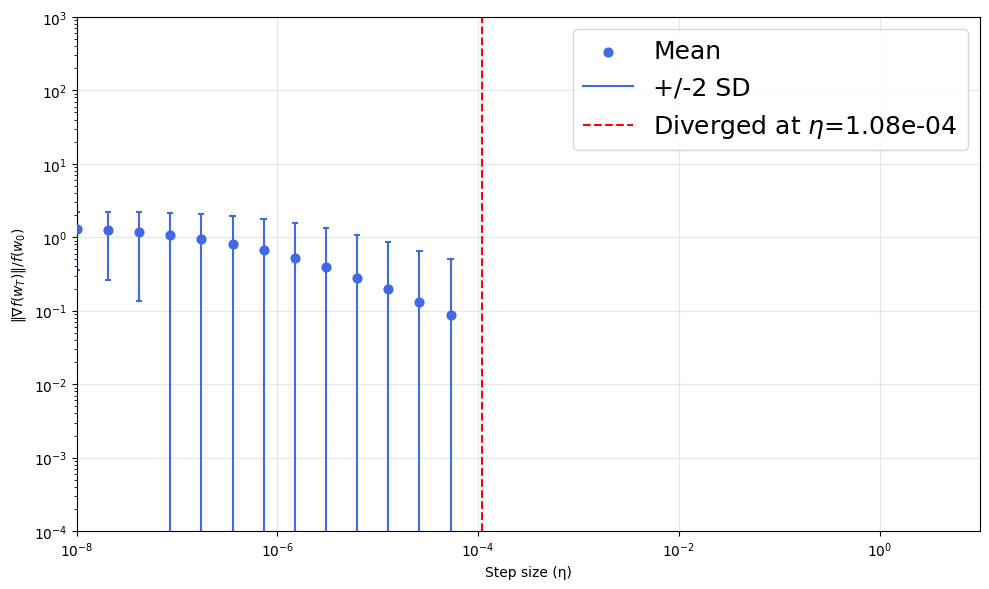}
     \caption{$\pi_j = \cN(\vecOrigin, 7.5\matI_{20})$. The first divergence is at $\eta_i\approx 1.08 \cdot 10^{-4}$.}
     \label{fig:sgd_p6_cov7.5_naive}
 \end{subfigure}
 \hfill
 \begin{subfigure}{0.46\textwidth}
     \includegraphics[width=\textwidth]{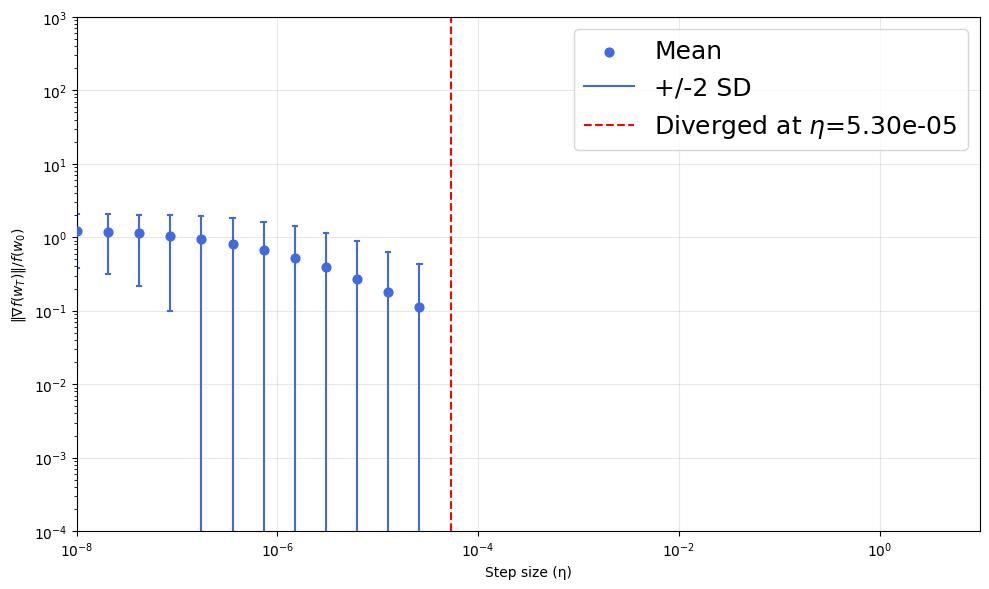}
     \caption{$\pi_j = \cN(\vecOrigin, 10\matI_{20})$. The first divergence is at $\eta_i\approx 5.30 \cdot 10^{-5}$.}
     \label{fig:sgd_p6_cov10_naive}
 \end{subfigure} 
 \caption{SGD simulation results for $p=6$. 
 For $\pi_j = \cN(\vecOrigin, 2.5 \matI_{20}), \cN(\vecOrigin, 5.0 \matI_{20}), \cN(\vecOrigin, 7.5 \matI_{20}), \cN(\vecOrigin, 10 \matI_{20})$, the first divergence are at $\eta_i\approx 4.52 \cdot 10^{-4}, 4.52 \cdot 10^{-4}, 1.08 \cdot 10^{-4}, 5.30 \cdot 10^{-5}$ respectively.}
 \label{fig:sgd_p6_naive}
\end{figure}

\paragraph{Results:} Our conclusions are similar to those from \pref{subsec:experimentstoygd}. When $p=2$, as predicted by established optimization theory for smooth functions, the first step size leading to divergence $\eta_i$ is identical across the $\pi_j$ (see \pref{fig:sgd_p2_naive}). In contrast for $p=3,4,5,6$, this $\eta_i$ generally decreases as the covariance $c_j \matI_{20}$ of $\pi_j$ increases from $2.5$ to 10 (see \pref{fig:sgd_p3_naive}, \pref{fig:sgd_p4_naive}, \pref{fig:sgd_p5_naive}, \pref{fig:sgd_p6_naive}). We note that while the general trends are similar to those from \pref{subsec:experimentstoygd}, we can clearly see the presence of the stochastic gradients in these plots. In many of these plots, $\tfrac{\nrm*{\grad F(\vecW_T)}}{F(\vecW_0)}$ becomes roughly constant for $\eta$ large enough such that $T=1000$ yields reasonable convergence; for such $\eta$, by $T=1000$, the true gradients are small enough and the noise from the stochastic gradients takes over.

Once more, consider how the first step size leading to divergence $\eta_i$ decreases as the covariance $c_j \matI_{20}$ of $\pi_j$ increases from $2.5$ to 10. We find that the rate of this decrease generally increases as $p$ increases. We also again see that fixing $\pi_j$ and comparing across $p$, the first step size leading to divergence $\eta_i$ decreases as $p$ increases. As discussed in \pref{subsec:experimentstoygd}, both of these phenomena are consistent with our theoretical results.
For each $p \in \{2,3,4,5,6\}$ and $\pi_j$, we again record the smallest $\eta_i$ for which divergence occurred in \pref{tab:sgdnaivestepsizes}. 

\begin{table}[h!]
\centering
\begin{tabular}{lcccc}
\toprule
 & $\pi_j=\cN(\vecOrigin, 2.5\matI_{20})$ & $\pi_j=\cN(\vecOrigin, 5.0\matI_{20})$ & $\pi_j=\cN(\vecOrigin, 7.5\matI_{20})$ & $\pi_j=\cN(\vecOrigin, 10\matI_{20})$ \\
\midrule
$p=2$ & $1.17 \cdot 10^0$ & $1.17 \cdot 10^0$  & $1.17 \cdot 10^0$  & $1.17  \cdot 10^0$  \\
$p=3$ & $2.81 \cdot 10^{-1}$ & $1.37 \cdot 10^{-1}$ & $6.72 \cdot 10^{-2}$ & $1.37 \cdot 10^{-1}$ \\
$p=4$ & $3.29 \cdot 10^{-2}$ & $3.29 \cdot 10^{-2}$ & $1.61 \cdot 10^{-2}$ &  $7.88 \cdot 10^{-3}$\\
$p=5$ & $7.88 \cdot 10^{-3}$ & $1.89 \cdot 10^{-3}$ & $9.24 \cdot 10^{-4}$ & $4.52 \cdot 10^{-4}$ \\
$p=6$ &  $4.52 \cdot 10^{-4}$ & $4.52 \cdot 10^{-4}$ & $1.08 \cdot 10^{-4}$ & $5.30 \cdot 10^{-5}$
\end{tabular}
\caption{Smallest $\eta_i$ leading to divergence for a given $p$ and initialization $\pi_j$.}
\label{tab:sgdnaivestepsizes} 
\end{table}

\end{document}